\documentclass[a4paper]{article}

\usepackage{MyLaTeX,MyArticle}
\usepackage{lscape}

\RequirePackage{pgf,pgffor}
\usepackage{tikz}
\newcommand{\Proj}[1]{\mathcal{P}({#1})}
\numberwithin{equation}{section}
\numberwithin{table}{section}
\numberwithin{figure}{section}
\newcommand{\CSp}{\operatorname{CSp}}
\newcommand{\CSO}{\operatorname{CSO}}
\newcommand{\Irr}{\mathrm{Irr}}
\newcommand{\MOD}[1]{#1\mathrm{-mod}}

\newcommand{\Deg}{\mathrm{Deg}}
\newcommand{\Arg}{\mathrm{Arg}}
\renewcommand{\b}[1]{\textbf{#1}}
\renewcommand{\bo}[1]{\text{\mathversion{bold}$#1$\mathversion{normal}}}
\renewcommand{\L}{\b L}
\newcommand{\aA}{\mathrm{aA}}

\def\withapp{1}

\title{Perverse Equivalences and Brou\'e's Conjecture II: The Cyclic Case}
\author{David A.\ Craven}
\date{June 2012}
\tikzstyle{every node}=[circle, fill=black!0,
                        inner sep=0pt, minimum width=4pt]

\begin{document}
\maketitle

\begin{abstract} We study Brou\'e's abelian defect group conjecture for groups of Lie type using the recent theory of perverse equivalences and Deligne--Lusztig varieties. Our approach is to analyze the perverse equivalence induced by certain Deligne--Lusztig varieties (the geometric form of Brou\'e's conjecture) directly; this uses the cohomology of these varieties, together with information from the cyclotomic Hecke algebra. We start with a conjecture on the cohomology of these Deligne--Lusztig varieties, prove various desirable properties about it, and then use this to prove the existence of the perverse equivalences predicted by the geometric form of Brou\'e's conjecture whenever the defect group is cyclic (except possibly for two blocks whose Brauer tree is unknown). This is a necessary first step to proving Brou\'e's conjecture in general, as perverse equivalences are built up inductively from various Levi subgroups.

This article is the latest in a series by Rapha\"el Rouquier and the author with the eventual aim of proving Brou\'e's conjecture for unipotent blocks of groups of Lie type.
\end{abstract}

\section{Introduction}

Brou\'e's abelian defect group conjecture is one of the deepest conjectures in modular representation theory of finite groups, positing the existence of a derived equivalence between a block $B$ of a finite group $G$ and its Brauer correspondent, whenever the block has abelian defect groups. If $G$ is a group of Lie type and $B$ is a unipotent block (e.g., the principal block) then there is a special form of Brou\'e's conjecture, the \emph{geometric form}, in which the derived equivalence is given by the complex of cohomology of a particular variety associated with $G$, a \emph{Deligne--Lusztig variety}. Various properties of this derived equivalence arise from properties of this cohomology, and this offers another avenue in which these varieties have become important, beyond their original application in classifying unipotent characters of groups of Lie type, and their intrinsic interest.

The first objective of this article is to provide a conjecture giving the precise cohomology of these Deligne--Lusztig varieties over an algebraically closed field of characteristic $0$. This is the information required for the derived equivalence and so, equipped with this information, we can search directly for the derived equivalence without analyzing the geometry of Deligne--Lusztig varieties. Previously, only the cases where the prime $\ell$ divides $q\pm 1$ were conjectured \cite{dmr2007}, and the case where $\ell$ divides $\Phi_d(q)$ with $d$ the Coxeter number was solved by Lusztig in \cite{lusztig1976}, so this conjecture is a considerable extension of this work. We give the precise conjecture later in this introduction, and then give the theorems that we prove about it afterwards.

We then turn our attention to the applications to Brou\'e's conjecture. The majority of the article is spent proving the following theorem.

\begin{thm}\label{thm:combbroueconj} Let $B$ be a unipotent block of a finite group of Lie type, not of type $E_8$. If $B$ has cyclic defect groups, then the combinatorial form of Brou\'e's conjecture holds for $B$.
\end{thm}

The `combinatorial version' of Brou\'e's conjecture, at least for blocks with cyclic defect group, will be given in Section \ref{sec:combbroueconj}, with its rather more delicate extension to all groups to appear in a later paper in this series. In fact, the restriction on the type of the group in this theorem is largely not necessary, as there are only two unipotent blocks of $E_8$ for which the Brauer tree, or equivalently the combinatorial form of Brou\'e's conjecture, is not known. Along the way, we give a complete description of all perverse equivalences between a block with cyclic defect group and its Brauer correspondent in Theorem \ref{thm:allperversecyclic}.

\medskip

We now describe in more detail the results given in this paper. We start with the conjecture on the cohomology of Deligne--Lusztig varieties. Let $\ell\neq p$ be primes, $q$ a power of $p$, write $d$ for the multiplicative order of $q$ modulo $\ell$, and let $G=G(q)$ be a finite group of Lie type. (We are more precise about our setup in Section \ref{sec:gensetup}.) We assume that $\ell$ is large enough that the Sylow $\ell$-subgroup of $G$ is abelian. The exact varieties that we consider are given in Section \ref{sec:DLvars}; if $\kappa\geq 1$ is prime to $d$ then to the fraction $\kappa/d$ we attach a variety $Y_{\kappa/d}$ in a natural way; it is this variety whose cohomology over $\bar\Q_\ell$ that we wish to describe.

Let $\mc F$ denote the set of all polynomials in $\R[q]$ whose zeroes are either roots of unity or $0$. Notice that the generic degree of any unipotent character of a group of Lie type, including the Ree and Suzuki groups, which are polynomials in $q$, lie in the set $\mc F$. (It also includes the `unipotent degrees' of the real reflection groups $H_3$, $H_4$ and $I_2(p)$, see \cite{lusztig1993}.) If $\xi$ is a non-zero complex number, write $\Arg_{\kappa/d}(\xi)$ for the set of all positive numbers $\lambda$ such that $\lambda$ is an argument for $\xi$ and $\lambda\leq 2\pi \kappa/d$. If $f$ is a polynomial, write $\Arg_{\kappa/d}(f)$ for the multiset that is the union of $\Arg_{\kappa/d}(\xi)$ for $\xi$ all non-zero roots of $f$, with multiplicity.

\begin{defn} For coprime integers $d,\kappa\geq 1$ and $f\in \mc F$, write $a(f)$ for the multiplicity of $0$ as a zero of $f$, $A(f)=\deg(f)$, and $\phi_{\kappa/d}(f)$ for the sum of $|\Arg_{\kappa/d}(f)|$ and half the multiplicity of $1$ as a root of $f$. Set $\pi_{\kappa/d}(f)=(a(f)+A(f))\kappa/d+\phi_{\kappa/d}(f)$.
\end{defn}

If $\chi$ is a unipotent character lying in a block with $d$-cuspidal pair $(\b L,\bo\lambda)$ (see \cite{bmm1993} for a definition), and $\Deg(\chi)$ denotes the generic degree of $\chi$, then we write $\pi_{\kappa/d}(\chi)$ for the difference $\pi_{\kappa/d}(\Deg(\chi))-\pi_{\kappa/d}(\Deg(\bo\lambda))$. (For those unfamiliar with $d$-cuspidal pairs, as an example, for the principal block $\bo\lambda$ is the trivial character, and so $\pi_{\kappa/d}(\Deg(\bo\lambda))=1$ and $\pi_{\kappa/d}(\chi)=\pi_{\kappa/d}(\Deg(\chi))$.) We are now able to state the conjecture on cohomology for unipotent characters of $G$.
\begin{conj}\label{conj:DLcohom} If $\chi$ is a unipotent character of $\bar{\Q}_\ell G$ then $\pi_{\kappa/d}(\chi)$ is the unique degree of the cohomology of the Deligne--Lusztig variety $H^\bullet(Y_{\kappa/d},\bar{\Q}_\ell)$ in which $\chi$ appears.
\end{conj}

As we have mentioned before, one reason for interest in the cohomology of Deligne--Lusztig varieties is Brou\'e's conjecture: for unipotent blocks of groups of Lie type, it provides a more explicit version -- the \emph{geometric version} of Brou\'e's conjecture -- of a derived equivalence between the block and its Brauer correspondent. We will describe this in more detail in Section \ref{sec:DLvars}. In particular, this derived equivalence should be \emph{perverse} (see \cite{chuangrouquierun} and Section \ref{sec:pervequiv} below). The cohomology of the varieties $Y_{\kappa/d}$ should provide perverse equivalences for Brou\'e's conjecture, and the geometric version of Brou\'e's conjecture implies the following.

\begin{conj}\label{conj:perverse} If $\chi_1,\dots,\chi_s$ are the unipotent ordinary characters in the unipotent $\ell$-block $B$ of $kG$ with abelian defect group, then there is a perverse equivalence from $B$ to $B'$ with perversity function given by $\pi_{\kappa/d}(\chi_i)$, where $B'$ is the Brauer correspondent of $B$.
\end{conj}

Again, we are more specific about when this conjecture should hold in Section \ref{sec:DLvars}. The firming up of this conjecture, into the full \emph{combinatorial form} of Brou\'e's conjecture, where all aspects of the perverse equivalence are given, is the subject of a later paper, but in the case of cyclic defect groups it is completed here. The precise description is complicated, and will be given in Section \ref{sec:combbroueconj}.

The first test that Conjectures \ref{conj:DLcohom} and \ref{conj:perverse} might hold is to prove that $\pi_{\kappa/d}(\chi)$ is always an integer, which is the content of our first theorem. This result also holds for the unipotent degrees of the Coxeter groups that are not Weyl groups, by a case-by-case check.

\begin{thm}\label{thm:pifnisinteger} Let $d\geq 1$ be such that $\Phi_d(q)$ divides $|G(q)|$, and let $\kappa\geq 1$ be prime to $d$. If $\chi$ is a unipotent character of $G$ then $\pi_{\kappa/d}(\chi)$ is an integer.
\end{thm}

The next theorem checks that in a bijection with signs arising from a perfect isometry between a unipotent block and its Brauer correspondent, the sign attached to $\chi$ is $(-1)^{\pi_{\kappa/d}(\chi)}$.

\begin{thm}\label{thm:bijectionwithsigns} Let $B$ be a unipotent $\ell$-block of $kG$, with Brauer correspondent $B'$. In a bijection with signs $\Irr_K(B)\to \Irr_K(B')$ arising from a perfect isometry, the sign attached to a unipotent character $\chi$ is $(-1)^{\pi_{\kappa/d}(\chi)}$.
\end{thm}

We prove Theorems \ref{thm:pifnisinteger} and \ref{thm:bijectionwithsigns} simultaneously in Section \ref{sec:integrality}; the proof is not case-by-case, and is remarkably short, needing no facts about groups of Lie type beyond the statement that $\Deg(\chi)/\Deg(\bo\lambda)$ is a constant modulo $\Phi_d(q)$, which is known \cite[\S5]{bmm1993}. In particular, we get a geometric interpretation of $\pi_{\kappa/d}(f)$; the quantity $\pi_{\kappa/d}(f)$ is (modulo $2$) the argument of the complex number $f(\e^{2\kappa\pi\I/d})$ divided by $\pi$. This proof gives some meaning behind the somewhat obscure function $\pi_{\kappa/d}$.

\medskip

We move on to perverse equivalences: we firstly prove that the structure of a perverse equivalence is in some sense independent of $\ell$ when the defect group is cyclic, a fact closely related to the statement that the Brauer tree of a unipotent $\ell$-block only depends on the $d$ such that $\ell\mid\Phi_d(q)$, but not $\ell$ or $q$. The general statement that perverse equivalences should in some sense be independent of the characteristic $\ell$ of the field is still ongoing research of Rapha\"el Rouquier and the author. The next stage is to classify all possible perverse equivalences between a block $B$ with cyclic defect groups and its Brauer correspondent $B'$, which we do in Section \ref{ssec:perverseBrauer2}. It turns out that two obvious conditions -- one being that the perversity function satisfies the conclusion of Theorem \ref{thm:bijectionwithsigns} on the parity of the perversity function, the other that the perversity function, which is defined on simple modules of the block, increase towards the exceptional node -- are sufficient, and so there is a nice parametrization of all perverse equivalences in this situation.

\if\withapp1
This is enough to prove Conjecture \ref{conj:perverse} for blocks with cyclic defect group whenever the Brauer tree is known, but for applying to derived equivalences for higher-rank groups, which will be done inductively, we need more complete information about the derived equivalence, and prove the complete combinatorial Brou\'e's conjecture; this task takes the remainder of the article. For exceptional groups we only perform a few representative cases here in full detail, but in the appendix we list all unipotent blocks of weight $1$ for all exceptional groups, together with the parameters of the cyclotomic Hecke algebra and the Brauer tree, including the conjectures for the two currently unknown trees.
\else
This is enough to prove Conjecture \ref{conj:perverse} for blocks with cyclic defect group whenever the Brauer tree is known, but for applying to derived equivalences for higher-rank groups, which will be done inductively, we need more complete information about the derived equivalence, and prove the complete combinatorial Brou\'e's conjecture; this task takes the remainder of the article. For exceptional groups we only perform a few representative cases here, but full details (which is too extensive to publish here at 100 pages) are available on the author's website.
\fi

\medskip

The structure of this article is as follows: Section \ref{sec:gensetup} introduces the general setup and the following section introduces the Deligne--Lusztig varieties under study. We prove Theorems \ref{thm:pifnisinteger} and \ref{thm:bijectionwithsigns} in Section \ref{sec:integrality}, and look at some evidence in favour of the conjecture on Deligne--Lusztig varieties in Section \ref{sec:previouswork}.

A long section on perverse equivalences in next, in which we determine all perverse equivalences between a block with cyclic defect group and its Brauer correspondent, among other results. Section \ref{sec:combbroueconj} gives the final form of the combinatorial Brou\'e conjecture for blocks with cyclic defect group, which we will prove in the remaining sections. Section \ref{sec:evalpi} gives some formulae regarding calculating the $\pi_{\kappa/d}$-function, and the section afterward introduces cyclotomic Hecke algebras for the cyclic group $Z_e$, as well as proving the important Proposition \ref{prop:increasingperv}, which enables us to compute with a different function to the $\pi_{\kappa/d}$-function in classical groups.

\if\withapp1
We then have two sections that give the standard combinatorial devices of partitions and symbols and the unipotent character degrees, then studies the character degrees of blocks with cyclic defect group to prove one part of the combinatorial Brou\'e conjecture; the succeeding two sections wrap up the proof. The final section gives three example computations with the unipotent blocks of exceptional groups, with the rest being summarized in the appendix.
\else
We then have two sections that give the standard combinatorial devices of partitions and symbols and the unipotent character degrees, then studies the character degrees of blocks with cyclic defect group to prove one part of the combinatorial Brou\'e conjecture; the succeeding two sections wrap up the proof. The final section gives three example computations with the unipotent blocks of exceptional groups, with the rest being summarized on the author's website.
\fi

\section{General Setup and Preliminaries}
\label{sec:gensetup}

Let $q$ be a power of a prime $p$, and let $\b G$ be a connected, reductive algebraic group over the field $\bar\F_p$. Let $F$ be an endomorphism of $\b G$, with $F^\delta$ a Frobenius map for some $\delta\geq 1$ relative to an $\F_{q^\delta}$-structure on $\b G$, and write $G=\b G^F$ for the $F$-fixed points. (We may normally take $\delta=1$ unless $G$ is a Ree or Suzuki group, in which case $q$ is an odd power of $\sqrt2$ or $\sqrt3$ and $\delta=2$.) Let $W$ denote the Weyl group of $\b G$, $B^+$ the braid monoid of $W$, and let $\phi$ denote the automorphism of $W$ (and hence $B^+$) induced by $F$. We let $\ell\neq p$ be a good prime, and write $d$ for the multiplicative order of $q$ modulo $\ell$, so that $\ell\mid\Phi_d(q)$. Suppose that $\ell$ does not divide any other $\Phi_{d'}(q)$ for $d'\neq d$, so that a Sylow $\ell$-subgroup of $G$ is abelian; in particular, $\ell$ is odd. Finally, we let $\mc O$, $K$ and $k$ be, as usual, a complete discrete valuation ring, its field of fractions, and its residue field; we assume that $\mc O$ is an extension of the $\ell$-adic integers $\Z_\ell$, so that $K$ is an extension of $\Q_\ell$ and $k$ is an extension of $\F_\ell$; we assume, again as usual, that these extensions are sufficiently large, for example the algebraic closures. (The assumption that $\Q_\ell\subset K$ makes it easier for the theory of Deligne--Lusztig varieties.)

We make a few remarks about the particular groups of Lie type we are studying: since we are interested in unipotent blocks only, we may be quite flexible about the precise form of the group involved; the centre of a group always lies in the kernel of any unipotent character, and the set of unipotent characters is independent of taking or removing \emph{diagonal} automorphisms, although the defect group of a unipotent block might change.  For example, as long as $\ell$ does not divide $q-1$, the restriction map from $\GL_n(q)$ to $\SL_n(q)$ induces Morita equivalences of unipotent blocks; therefore, if we term the blocks of $\PSL_n(q)$ whose inflation to $\SL_n(q)$ to be unipotent, the unipotent blocks of $\PSL_n(q)$, $\SL_n(q)$, $\PGL_n(q)$ and $\GL_n(q)$ are all Morita equivalent, with simple modules with isomorphic Green correspondents, so for Brou\'e's conjecture it is irrelevant which one is considered.

For definiteness, when $G$ is classical we take it to be one of the groups $\GL_n(q)$ (which is important if $\ell\mid(q-1)$), $\GU_n(q)$ (which is important if $\ell\mid(q+1)$), $\SO_{2n+1}(q)$ (where $q$ is odd), $\Sp_{2n}(q)$, and $(\CSO_{2n}^{\pm})^0(q)$, where this last group is the subgroup of $\CSO_{2n}^\pm(q)$ of index $2$, where the outer automorphisms induced on the simple group are diagonal. (For $q$ odd, we could take $\SO_{2n}^\pm(q)$ as well, but for $q$ even the $\SO$-action induces the graph automorphism on the simple group, so we cannot take this group.)

Let $\kappa$ be a non-negative integer prime to $d$ and write $\zeta=\e^{2\kappa\pi\I/d}$, so that $\zeta$ is a primitive $d$th root of unity. (In previous work in this area it has sometimes been assumed that $0\leq\kappa\leq d-1$, but in this and subsequent papers we will need to also consider the case $\kappa\geq d$.) Let $B$ be a unipotent $\ell$-block of $G$ with defect group $D$, and let $T$ be a $\Phi_d$-torus containing $D$ with $D$ and $T$ of the same rank. Write $e$ for the number of unipotent characters of $d$: in almost all cases where the defect group $D$ is cyclic, $e=d$, $e=2d$ or $e=d/2$. To $B$ we associate a $d$-cuspidal pair $(\L,\bo\lambda)$, and for any unipotent character $\chi$ in $B$ we write $\Deg(\chi)$, or simply $\chi(1)$, for the generic degree of the unipotent character $\chi$, a polynomial in $q$. Write $E$ for the $\ell'$-group $\Norm_G(D)/\Cent_G(D)$, which is a complex reflection group, and its natural action on the $\Phi_d$-torus $T$ is as complex reflections.

As usual, if $f$ is a polynomial, $A(f)$ and $a(f)$ denote $\deg(f)$ and the multiplicity of $0$ as a zero of $f$ respectively: these are usually called Lusztig's $A$- and $a$-functions, or often simply the $A$- and $a$-functions. For a unipotent character $\chi$ in $B$, we introduce the notation
\[ \aA(\chi)=(a(\Deg(\chi))+A(\Deg(\chi)))-(a(\Deg(\bo\lambda))+A(\Deg(\bo\lambda))).\]
The quantity $\aA(\chi)$ is closely related to the parameters of the cyclotomic Hecke algebra of $B$ (see Section \ref{sec:cyclohecke}), and the eigenvalues of the Frobenius map.

Write $B'$ for the Brauer correspondent of $B$, a block of $H=\Norm_G(D)$. The simple $B$-modules will usually be denoted by $S_i$ and the simple $B'$-modules will be denoted by $T_i$. If $D$ is cyclic then the Brauer tree of $B'$ is a star, which we can envisage as being embedded in $\C$, with the exceptional node positioned at $0$ and the $e$ non-exceptional characters being equally spaced around the exceptional node on the circle $|z|=1$. We choose our orientation of the Brauer tree to be anti-clockwise, so that in the following example the projective cover of the trivial module has second radical layer $T_2$.
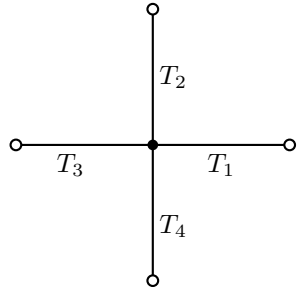
\begin{center}\begin{tikzpicture}[thick,scale=1.8]
\draw (0,-1) -- (0,1);
\draw (-1,0) -- (1,0);

\draw (0.5,-0.15) node{$T_1$};
\draw (0.15,0.5) node{$T_2$};
\draw (-0.6,-0.15) node{$T_3$};
\draw (0.15,-0.6) node{$T_4$};

\draw (1,0) node [draw] (l0) {};
\draw (-1,0) node [draw] (l0) {};
\draw (0,1) node [draw] (l0) {};
\draw (0,-1) node [draw] (l0) {};
\draw (0,0) node [fill=black!100] (ld) {};
\end{tikzpicture}\end{center}
In order to save space, we use the `$/$' character to delineate radical layers in a module, so that for example we  write $T_1/T_2/T_3/T_4/T_1$ for the radical layers of the projective cover of the trivial module above (assuming exceptionality $1$).

Write $\Proj{M}$ for the projective cover of the module $M$, and $\Omega(M)$ for the kernel of the natural map $\Proj{M}\to M$. In the opposite direction, write $\Omega^{-1}$ for the cokernel of the morphism mapping a module into its injective hull. Notice that $\Omega^2(T_i)=T_{i+1}$ (with indices taken modulo $e$) so that $\Omega^2$ acts like a rotation by $2\pi/e$ on the Brauer tree, and hence on the complex plane. It makes sense therefore to place $\Omega(T_i)$ on the circle of unit radius halfway between $T_i$ and $T_{i+1}$, so that $\Omega$ acts like a rotation of $\pi/e$ on the \emph{doubled Brauer tree} (this terminology, and concept, is not standard).

\section{Deligne--Lusztig Varieties}
\label{sec:DLvars}

In this section we give information on the varieties that we deal with in Conjecture \ref{conj:DLcohom}. In the geometric form of Brou\'e's conjecture, for each unipotent block $B$ of $kG$, where $\ell\mid\Phi_d(q)$, and each $\kappa\geq 1$ prime to $d$, there is a variety $Y_{\kappa/d}$, which has an action of $G$ on the one side and an action of the torus $T$ on the other: its complex of cohomology inherits this action, and the action of $T$ may be extended to an action of $\Norm_G(D)$, so that this complex provides a derived equivalence between $B$ and its Brauer correspondent $B'$. We now describe the variety $Y_{\kappa/d}$.

We first define the Deligne--Lusztig variety $Y(b)$, for $b\in B^+$. Let $w\mapsto \b w$ be the length-preserving lift $W\to B^+$ of the canonical map $B^+\to W$. Let $\b B$, $\b T$ and $\b U$ be, as usual, a fixed $F$-stable Borel subgroup, an $F$-stable torus $\b T$ contained in $\b B$ and the unipotent radical $\b U$ of $\b B$. Fix an $F$-equivariant morphism $\tau:B^+\to \Norm_{\b G}(\b T)$ that lifts the canonical map $\Norm_{\b G}(\b T)\to W$. For $w\in W$, write $\dot w=\tau(\b w)$, and 
given $w_1,\dots,w_m\in W$, we set
\begin{align*} Y(w_1,\dots,&w_m)=\{\,(g_1\b U,\dots,g_m\b U)\in (\b G/\b U)^m\mid
\\ &g_1^{-1}g_2^{\vphantom{-1}}\in\b U\dot w_1\b U,\dots,g_{r-1}^{-1}g_r^{\vphantom{-1}}\in\b U\dot w_{r-1}\b U,\dots,g_r^{-1}F(g_1)\in \b U\dot w_r\b U\,\}.\end{align*}
Up to isomorphism this variety depends only on the product $b=\b w_1\ldots\b w_m\in B^+$, and we write $Y(b)$ for this variety.

\medskip

We now describe the cases in which $Y_{\kappa/d}$ has been identified. Recall that $(\b L,\bo\lambda)$ is a $d$-cuspidal pair associated to the block $B$. If $\b L$ is a torus then the variety $Y_{\kappa/d}$ was identified in \cite{brouemichel1997}, and we briefly describe this case (see also \cite[\S3.4]{cravenrouquier2010un}). Let $\b w_0$ be the lift of the longest element of $W$ in $B^+$. Choose $b_d\in B^+$ such that $(b_d\phi)^d=(\b w_0)^2\phi^d$; the variety $Y_{\kappa/d}$ should be the Deligne--Lusztig variety $Y((b_d)^k)$.

Recently \cite{dignemichel2011un} a generalization of this construction of $Y(b_d)$ was given, producing so-called \emph{parabolic Deligne--Lusztig varieties}. The construction of these is more technical, and we do not give it here. In \cite{dignemichel2011un} a candidate variety $Y_{\kappa/d}$ is identified in the case where $\b L$ is minimal (i.e., the trivial character of $\b L$ is $d$-cuspidal). Thus in these cases the variety $Y_{\kappa/d}$ has been found, but in general the identification has not been explicitly worked out, although it seems as though it can be from the information contained in \cite{dignemichel2011un}.

\section{Integrality of $\pi_{\kappa/d}$ and a Bijection with Signs}
\label{sec:integrality}
In this section we prove that the $\pi_{\kappa/d}$-function, evaluated at a unipotent character, is always an integer, and demonstrate that, in a bijection with signs $\Irr_K(B)\to \Irr_K(B')$ that arises from a perfect isometry, that the sign attached to $\chi$ is $(-1)^{\pi_{\kappa/d}(\chi)}$. In this section we may assume that $d\geq 2$, and write $\zeta=\e^{2\kappa\pi\I/d}$, a primitive $d$th root of unity. (The case where $d=1$ is easy, since it is clear that it is an integer and we will see in Section \ref{sec:previouswork} that the integer is always even, tallying with \cite{bmm1993}.) Let $\mc F$, as before, denote the set of all polynomials in $q$ with real coefficients that have as zeroes either roots of unity of finite order or $0$. 

A preliminary result is needed to simplify some of the proofs that follow, and it will be useful in its own right; it describes the relationship between different fractions $\kappa/d$ that describe the same root of unity $\zeta=\e^{2\kappa\pi\I/d}$.

\begin{lem}\label{lem:add2pi} Let $f$ be a polynomial in $\mc F$. If $\kappa$ and $d$ are coprime positive integers, then
\[ \pi_{(\kappa+d)/d}(f)=\pi_{\kappa/d}(f)+2A(f).\]
\end{lem}
\begin{pf} The difference in the sets $\Arg_{(k+d)/d}(f)$ and $\Arg_{\kappa/d}(f)$ is one copy of an argument for each non-zero root of $f$, so that the difference in cardinalities is $A(f)-a(f)$. Obviously the remaining contribution to the $\pi_{\kappa/d}$-function -- $(A(f)+a(f))\kappa/d$ -- yields a difference of $A(f)+a(f)$, and the sum of these two is $2A(f)$, as claimed.
\end{pf}

Since $2A(f)$ is always an integer, we see that $\pi_{\kappa/d}(f)\equiv \pi_{(\kappa+d)/d}(f)$ modulo $2$ for any polynomial $f$, so in proving integrality and correct parity, we may assume that $\kappa$ is less than $d$.

We can of course extend the domain of $\pi_{\kappa/d}(-)$ to include all polynomials, and for the proof of the next result we extend the domain to include all polynomials with \emph{complex} coefficients that have as zeroes either roots of unity of finite order or $0$.

\begin{thm} Let $1\leq \kappa<d$ be coprime integers. Let $f$ be a polynomial in $\mc F$ such that $f(\zeta)\neq 0$. Writing $\arg(z)$ for the argument of the complex number $z$, taken in $[0,2\pi)$, modulo $2$ we have that $\arg(f(\zeta))/\pi\equiv\pi_{\kappa/d}(f)$.
\end{thm}
\begin{pf} Since $f$ is a polynomial with real coefficients, if $\omega$ is a complex zero of $f$ then so is $\bar\omega$. Since $\arg(zw)=\arg(z)+\arg(w)$ and $\pi_{\kappa/d}(fg)=\pi_{\kappa/d}(f)+\pi_{\kappa/d}(g)$ it suffices to prove the result for $f=q$, $f=(q\pm 1)$ and $f=(q-\omega)(q-\bar\omega)$, where $\omega$ is a root of unity of finite order. If $f=q$ then the result is obvious, since $\arg(f(\zeta))\equiv 2\pi \kappa/d\bmod 2\pi$ and $\pi_{\kappa/d}(f)=2\kappa/d$.

If $f=(q-\omega)$ for some $\omega\neq \zeta$ then 
\[ \pi_{\kappa/d}(f)=\begin{cases} \kappa/d+1/2& \omega=1,\\ \kappa/d & \arg(\zeta)\leq \arg(\omega),\\ \kappa/d+1 &\arg(\zeta)\geq \arg(\omega).\end{cases}\]
It is easy to see that if $z$ has norm $1$ then $z-1$ has argument $\arg(z)/2+\pi/2$ for $\arg(z)\in [0,2\pi)$, proving the result for $\omega=1$. As $q+1=-(-q-1)$, this does $\omega=-1$ as well. In fact, since, $\zeta-\omega=\omega(\zeta/\omega-1)$, we have
\[ \arg(\zeta-\omega)\equiv\arg(\omega)+\arg(\zeta/\omega-1)\equiv\arg(\omega)+\arg(\zeta/\omega)/2+\pi/2\mod 2\pi.\]
As we have declared that $\arg(-)$ lies in $[0,2\pi)$, $\arg(\zeta/\omega)$ is equal to $\arg(\zeta)-\arg(\omega)$ if $\arg(\zeta)\geq \arg(\omega)$, and $\arg(\zeta)-\arg(\omega)+2\pi$ if $\arg(\zeta)\leq \arg(\omega)$. Hence $(\zeta-\omega)(\zeta-\bar\omega)$ has argument $\arg(\zeta)$ if $\arg(\zeta)>\arg(\omega),\arg(\bar\omega)$ or $\arg(\zeta)<\arg(\omega),\arg(\bar\omega)$, and $\arg(\zeta)\pm\pi$ (to stay in $[0,2\pi)$) otherwise. Hence this argument divided by $\pi$ is either $2\kappa/d$ or $2\kappa/d+1$ (modulo $2$), as needed.
\end{pf}

Let $\chi$ be a unipotent ordinary character in a block $B$, with associated $d$-cuspidal pair $(\b L,\bo\lambda)$. It is known \cite[\S5]{bmm1993} that (as polynomials in $q$) $\Deg(\chi)\equiv (-1)^\ep\alpha\cdot\Deg(\bo\lambda)\bmod \Phi_d(q)$, for some positive $\alpha\in \Q$ and $\ep\in\Z$. Hence $\Deg(\chi)/\Deg(\bo\lambda)$ is a rational function which, when $q$ is evaluated at a primitive $d$th root of unity, becomes $\pm\alpha$, a real number. Thus $\pi_{\kappa/d}(\chi)$, which modulo $2$ is the argument of $\pm\alpha$ divided by $\pi$, must be $\ep$ modulo $2$; in particular, $\pi_{\kappa/d}(\chi)$ is always an integer, proving Theorem \ref{thm:pifnisinteger}.

If $\ell$ is large then it is also proved in \cite[\S5]{bmm1993} that $(-1)^\ep=(-1)^{\pi_{\kappa/d}(\chi)}$ is the sign in a perfect isometry between $B$ and $B'$, so this proves Theorem \ref{thm:bijectionwithsigns} as well.

\section{Previous Work and Known Cases}

\label{sec:previouswork}

In this section we will summarize some of the previous work on this problem, and how it interacts with Conjecture \ref{conj:DLcohom}.

In the cases of $d=1$ and $d=2$ with $\kappa=1$ there is already a conjecture from \cite{dmr2007}, which states that the degree should be $2\deg(\Deg(\chi))/d$, i.e., $2A(\Deg(\chi))/d$ for the principal block. Note that in these cases, there is only one option for $\zeta$.

\begin{prop}\label{prop:conjd12} If $d=1$ or $d=2$ then for $\chi$ in the principal $\ell$-block, $\pi_{\kappa/d}(\chi)=2\kappa A(\Deg(\chi))/d$.
\end{prop}
\begin{pf} Let $f=\Deg(\chi)$. If $\kappa=d=1$ then $f(1)\neq 0$, and hence $a(f)+\phi_{\kappa/d}(f)=\deg(f)$, since any zero of $f$ must either be $0$, so is counted in $a(f)$, or non-zero and not $1$, so counted in $\phi_{\kappa/d}(f)$. Hence $\pi_{\kappa/d}(f)=2A(f)$, and so $\pi_{\kappa/d}(\chi)=2A(f)$. The case of $\kappa$ arbitrary follows now from Lemma \ref{lem:add2pi}.

Now let $d=2$ and $\kappa=1$; if $\omega$ is a complex zero of $f$ then so is $\bar\omega$, since $f\in\R[q]$; hence exactly one of $\omega$ and $\bar\omega$ contributes to $\phi_{\kappa/d}(f)$. Finally, we count half of each zero that is $+1$, and since $d=2$ we cannot have that $\Phi_2(q)$ divides $f$, so that $\phi_{\kappa/d}(f)$ counts half of the number of zeroes of $f$ not equal to $0$. Thus $\pi_{\kappa/d}(f)=(a(f)+A(f))/2+\phi_{\kappa/d}(f)=A(f)$, as needed. The case of $\kappa$ an arbitrary odd integer follows again from Lemma \ref{lem:add2pi}.
\end{pf}

Notice that the case $d=1$ can also have $\kappa=0$, in which case $\pi_{\kappa/d}(f)=0$ as long as $f(1)\neq 0$. This suggests that there is a Deligne--Lusztig variety that is a collection of points, and indeed this is the case: this is the case proved by Puig in \cite{puig1990}, that establishes a Puig equivalence (in particular Morita equivalence) between the blocks $B$ and $B'$.

The other case where much is known about the structure of the Deligne--Lusztig variety is when $d$ is the Coxeter number, which for the groups other than the Ree and Suzuki groups is simply the largest integer $d$ such that $\Phi_d(q)\mid |G(q)|$ (in the Ree and Suzuki groups case it is $d''$ for $d$ the largest integer such that $\Phi_d(q)\mid|G(q)|$, so that the associated polynomial $\Phi_d''(q)$ has as a zero the root of unity with smallest argument). In this case, both the structure of the cohomology of the Deligne--Lusztig variety and the geometric version of Brou\'e's conjecture are known.

\begin{thm}[Lusztig \cite{lusztig1976}] Conjecture \ref{conj:DLcohom} on the cohomology of Deligne--Lusztig varieties holds whenever $d$ is the Coxeter number and $\kappa=1$.
\end{thm}

If $d$ is the Coxeter number then the Sylow $\Phi_d$-subgroups (or $\Phi_d''$-subgroups) are cyclic, so Rickard's theorem holds and there is a perverse equivalence (see Section \ref{ssec:perverseBrauer1}). In this case, it is actually seen that the perversity function for $d$ the Coxeter number and $\kappa=1$ is the canonical perversity function in Secction \ref{ssec:perverseBrauer1}. It is easy to see that, in this case, $\phi_{\kappa/d}(\chi)$ is half the multiplicity of $(q-1)$ in $\chi(1)$, where $\chi$ is a unipotent $B$-character, so it is straightforward to show that $\pi_{\kappa/d}(\chi)=(A(\chi)+a(\chi))/d+\phi_{\kappa/d}(\chi(1))$ is equal to $\pi_0(\chi)$, the canonical perversity function. (To see the last step we need the structure of the Brauer tree for $d$ the Coxeter number, but this is now known for all groups \cite{dudasrouquier2012un}; it satisfies the conjecture of Hiss, L\"ubeck and Malle from \cite{hlm1995}.) Hence we get the following result.

\begin{thm} Conjecture \ref{conj:perverse} holds whenever $d$ is the Coxeter number and $\kappa=1$.
\end{thm}

By work of Olivier Dudas \cite[Theorem B]{dudas2010un} and Dudas and Rouquier \cite{dudasrouquier2012un}, it is known that the complex of cohomology of the Deligne--Lusztig variety, over $\mc O$, does indeed induce a perverse equivalence, and so even the geometric version of Brou\'e's conjecture holds in this case.

In other work, Dudas has proved the geometric version of Brou\'e's conjecture for the principal $\ell$-block of $\GL_n(q)$ whenever $d>n/2$ and $\kappa=1$, and proved for all $\GL_n(q)$ that Conjecture \ref{conj:DLcohom} holds for all $d$ if and only if it holds for $d=1$ \cite{dudas2011un}, where we require $\kappa=1$ in these cases.

In addition to these results, Dudas and Jean Michel have calculated the cohomology of various Deligne--Lusztig varieties, and the results are consistent with the conjecture here. A non-exhaustive list, with $\kappa=1$ (except for the first), is the following:
\begin{enumerate}
\item $\GL_3(q)$, all $d$ and $\kappa$;
\item $G=\GU_4(q)$, $d=4$;
\item $G=\GU_6(q)$, $d=6$;
\item $G=F_4(q)$, $d=8$;
\item $G=E_6(q)$, $d=9$;
\item $G={}^2\!E_6(q)$, $d=12$;
\item $G=E_7(q)$, $d=14$ (principal series and cuspidal modules only);
\item $G=E_8(q)$, $d=24$ (principal series and cuspidal modules only).
\end{enumerate}

Finally, since Conjecture \ref{conj:perverse} can be thought of as a shadow of the geometric form of Brou\'e's conjecture and Conjecture \ref{conj:DLcohom}, so it is also of interest to know when this conjecture is known, particularly for non-cyclic defect groups. In \cite{cravenrouquier2010un}, Rapha\"el Rouquier and the author proved this version for the principal blocks of all groups of Lie type, $\ell=3$ (so in particular, dividing the order of the Weyl group), and $d=1,2$ whenever the Sylow $\ell$-subgroup is elementary abelian of order $9$. In addition, in as-yet unpublished work, we have verified it for $\ell=5$ and $d=4$ for the principal blocks of $(\CSO^+_8)^0(2)$ and $\Sp_8(2)$, $\ell=5$ and $d=8'$ for ${}^2\!F_4(2)$, and for $\ell=7$ and $d=3$ for the principal block of $G={}^3D_4(2)$. (To extend this to all appropriate $q$ we need to know that the Green correspondents of the simple $B$-modules in the principal blocks do not depend on $q$, a widely believed, but unproven, statement.)

\section{Perverse Equivalences}
\label{sec:pervequiv}

In this section we will give some of the theory of perverse equivalences, as developed in \cite{chuangrouquierun} originally, and \cite{cravenrouquier2010un}. We begin with a definition of (a special case of) perverse equivalences relying on cohomology (rather than equivalences of Serre subcategories as per the original definition) and describe an algorithm that computes all perverse equivalences. We then prove a theorem that the output of the algorithm is `generic' in $\ell$ in a suitable sense whenever the defect group is cyclic (the statement is general is ongoing research of Rapha\"el Rouquier and the author), before constructing all perverse equivalences between any Brauer tree algebra (for example, a block with cyclic defect group for a finite group) and that of the star with exceptional vertex at the centre (for example, the Brauer correspondent of a block with cyclic defect group). This infinite family will then contain all of the perverse equivalences that should arise from Deligne--Lusztig varieties $Y_{\kappa/d}$.

\subsection{Definition and Algorithm}
\label{ssec:defnalg}

We begin with the definition of a special type of perverse equivalence, which includes those equivalences expected for groups of Lie type.

\begin{defn} Let $R$ be one of $\mc O$ and $k$, and let $A$ and $A'$ be $R$-algebras. A derived equivalence $f:D^b(\MOD A)\to D^b(\MOD{A'})$ is \emph{perverse} if there exists a bijection between the simple $A$-modules $S_1,\dots,S_e$ and simple $A'$-modules $T_1,\dots,T_e$ (relabelled so that the bijection sends $S_i$ to $T_i$), and a function $\pi:\{1,\dots,r\}\to \Z_{\geq 0}$ such that, in the cohomology of $f(S_i)$, the only composition factors of $H^{-j}(f(S_i))$ are $T_\alpha$ for those $\alpha$ such that $\pi(\alpha)<j\leq \pi(i)$, and a single copy of $T_i$ in $H^{-\pi(i)}(f(S_i))$.
\end{defn}

Another way of viewing this is, if we construct a table with the cohomology in the $j$th place of the module $f(S_i)$ from left to right, then there is a single copy of $T_i$ on the $i$th row, at the position $-\pi(i)$, and if a module $T_\alpha$ appears for $\alpha\neq i$ in the table then the column it appears in is strictly to the left of $-\pi(\alpha)$. Hence, if we order the $S_i$ so that $\pi(i)$ weakly increases with $i$, the table is triangular in shape.

\medskip

Perverse equivalences are of interest because there is an algorithm to compute them. We will describe this algorithm only in the case of a block of a finite group, since this is the case that concerns us, and refer to \cite{chuangrouquierun} for the general case. This algorithm takes as inputs the following:
\begin{enumerate}
\item the Brauer correspondent $B'$ of the block $B$ of $kG$;
\item a function $\pi$ from the simple $B'$-modules to non-negative integers, or, labelling the simple $B'$-modules $T_1,\dots,T_e$, a function $\pi$ from $\{1,\dots,e\}\to\Z_{\geq 0}$;
\item a collection $\mc R$ of sequences of relatively $Q$-projective $B'$-modules for various $1<Q<D$ -- one sequence $\mc R_i$ for each simple $B'$-module $T_i$;
\end{enumerate}
it returns a collection of complexes in the derived category of $B'$, which are meant to represent the images of simple $B$-modules in a derived equivalence from $B$ to $B'$. Notice that the output -- a series of complexes -- makes no reference to $B$, and so what we actually get is a derived self-equivalence on $B'$. We say `we apply the algorithm to the triple $(B',\pi,\mc R)$' when we perform the algorithm on this triple, or if $\mc R$ is empty as it in the case where $D$ is cyclic, to the pair $(B',\pi)$.

The collection $\mc R$ should come from a stable equivalence between $B$ and $B'$, but can be an arbitrary collection for the statement of the algorithm, except the algorithm will fail for many such sequences $\mc R$. We assume that the number of terms in $\mc R_i$ is less than $\pi(i)$. The output is a set of complexes $X_i$ of $B'$-modules, with the set of degree $0$ terms (hopefully) being the Green correspondents of the simple $B$-modules.

The first term of the complex $X_i$ is the injective hull $\proj(T_i)$ of $T_i$, in degree $-\pi(i)$. The cohomology $H^{-\pi(i)}(X_i)$ consists of $T_i$ in the socle, and the largest submodule of $\Proj{T_i}/T_i$ consisting of those $T_\alpha$ such that $\pi(\alpha)<\pi(i)$. This module $M_{\pi(i)}$ will be the kernel of the map from degree $-\pi(i)$ to $-\pi(i)+1$; let $L_{\pi(i)}=\Omega^{-1}(M_{\pi(i)})$, i.e., $\Proj{T_i}/M_{\pi(i)}$.

Now let $0<j<\pi(i)$. Write $\mc R_i=(\mc R_{i,1},\dots,\mc R_{i,n_i})$. The $-j$th term of $X_i$ is the module $P_j$: this is the direct sum of $\mc R_{i,j}$ (which is $0$ if $j>n_i$) and the smallest injective module $P_j'$ such that the socle of $\mc R_{i,j}\oplus P_j'$ contains that of $L_{j+1}$ (so that, if $j>n_i$ then $P_j'$ is simply the injective hull of $L_{j+1}$). \textbf{At this stage, it is not generally true that $L_{j+1}$ is isomorphic to a submodule of $P_j$, or even if it is, that the quotient $P_j/L_{j+1}$ is independent of the choice of injective map, but we assume that $\mc R_{i,j}$ is chosen so that these conditions are satisfied; in particular, these hold if $\mc R_{i,j}=0$. (If these statements do not hold, we say that the algorithm fails, noting that the algorithm cannot fail for $\mc R=\emptyset$.)} The submodule $L_{j+1}$ is the image of the previous map, and define $M_j$ to be the largest submodule of $P_j$, containing $L_{j+1}$, such that $M_j/L_{j+1}$ has composition factors only those $T_\alpha$ such that $\pi(\alpha)<j$. The module $M_j/L_{j+1}$ is $H^{-j}(X_i)$, and $M_j$ is the kernel of the map from degree $-j$ to degree $-j+1$. Again, write $L_j=P_j/M_j$.

Finally, the $0$th term of $X_j$ is the module $L_1$, which should be the Green correspondent of a simple $B$-module, which we denote by $S_i$.

\medskip

An important remark is that, if the injective module $P_j$ in degree $-j$ has a simple module $T_\alpha$ in its socle, then $\pi(\alpha)>j$, since otherwise in degree $-j-1$ the module $T_\alpha$, which lies in the socle of $L_{j+1}$, would have been subsumed into $M_{j+1}$.

\medskip

We now discuss the cohomology of the complexes $X_i$, and how this may be used to reconstruct the decomposition matrix of the block $B$. Let $\pi$ and the $S_i$ and $T_i$ be as above, and let $X_i$ be the complex in $D^b(\MOD {B'})$ obtained by running the algorithm. The \emph{alternating sum of cohomology} $H(T_i)$ of $X_i$ is the virtual $B'$-module
\[ \bigoplus_{j=0}^{\pi(T_i)} \;\mathop{\bigoplus_{T\in \mathrm{cf}(H^{-j}(X_i))}} (-1)^{j-\pi(T)}T,\]
where $\mathrm{cf}(M)$ is the set of composition factors of $M$. These virtual $B'$-modules determine $r$ rows of the decomposition matrix in an easy way, and can determine the rest of the decomposition matrix if the corresponding rows of $B'$ are known (it is an easy task to determine these rows for $B'$).

We will explain this description via an example.

\begin{example} Let $G=G_2(3)$ and $\ell=13\mid\Phi_3(3)$, so that $d=3$, and let $\kappa=1$. Let $P$ denote a (cyclic) Sylow $\ell$-subgroup of $G$, $H=\Norm_G(P)\cong\Z_\ell\rtimes \Z_6$, and order the simple $kH$-modules so that the $i$th radical layer of $\Proj{k}$ is $T_i$ for $1\leq i\leq 6$, where $k$ denotes the trivial module as well as the field. Using the notation of \cite{carterfinite} for the unipotent characters of $G_2(q)$, the ordering on the simples for the principal block $B$ of $kG$ is $\phi_{1,0}$, $G_2[\theta^2]$, $\phi_{2,2}$, $G_2[\theta]$, $\phi_{1,6}$, $G_2[1]$, yielding a particular bijection between the unipotent characters of $B$ and the simple $B'$-modules. For the reason why we chose this particular ordering, see Section \ref{sec:combbroueconj}. This allows us to transfer the $\pi_{\kappa/d}$-function to the simple $B'$-modules, meaning we can run the algorithm.

Applying the formula for the $\pi_{1/3}$-function, with this ordering, we get $0,3,3,3,4,4$. (It is a coincidence that, in this case, the ordering on the simple $kH$-modules makes the $\pi$-function weakly increasing, and in general this does not happen.) By Theorem \ref{thm:allperversecyclic} below this is the perversity function for a perverse equivalence between $B$ and its Brauer correspondent, and the particular bijection needed is that above.

The Green correspondents of the simple $B$-modules have dimensions $1$, $12$, $11$, $12$, $5$ and $1$, and have radical layers (writing $i$ for $T_i$)
\[ C_1=1,\;\; C_2=6/\cdots/5,\;\; C_3=2/\cdots/6,\;\; C_4=3/\cdots/2,\;\; C_5=5/\cdots/3,\;\; C_6=4.\]
(We can delete the inner radical layers since a $kH$-module is determined by its dimension and socle (or top).)

Applying the algorithm to the pair $(B',\pi)$, where $B'=kH$ and the $\pi$-function given by $\pi_{1/3}(-)$ on the simple $kH$-modules, we get six complexes, of the form:
\[\begin{array}{lr}
    X_1:&\;\; C_1.
\\ X_2:&\;\; \Proj2\to\Proj6\to\Proj6\to C_2.
\\ X_3:&\;\; \Proj3\to\Proj2\to\Proj2\to C_3.
\\ X_4:&\;\; \Proj4\to\Proj3\to\Proj3\to C_4.
\\ X_5:&\;\; \Proj5\to\Proj6\to\Proj4\to\Proj5\to C_5.
\\ X_6:&\;\; \Proj6\to\Proj5\to\Proj5\to\Proj4\to C_6.
\end{array}\]
The cohomology of the complexes above is displayed in the following table.
\smallskip
\begin{center}\begin{tabular}{cccccc}
\hline $X_i$ & $H^{-4}$ & $H^{-3}$ & $H^{-2}$ & $H^{-1}$ & Total
\\ \hline $2$&& $1/2$& $1$& & $2$
\\ $3$& & $3$& & $1$& $3-1$
\\ $4$& & $4$ & & & $4$
\\ $5$& $1/2/3/4/5$ & & & & $5-4-3-2+1$
\\ $6$& $6$ & & & & $6$
\\ \hline \end{tabular}\end{center}

The column `Total' gives the alternating sum of cohomology. To construct the first six rows of the decomposition matrix for $B$, we stipulate that the vector consisting of $0$ everywhere except a $1$ in the $i$th position should be the sum of the rows (with signs) given in the Total column. Hence the third row, minus the first row, should be $(0,0,1,0,0,0)$, and hence the third row is $(1,0,1,0,0,0)$. Continuing this, we get the matrix below.
\[ \begin{tabular}{ccccccccc}
\hline Name & Degree & $\pi_{\kappa/d}$ & $S_1$ & $S_2$ & $S_3$ & $S_4$ & $S_5$ & $S_6$
\\ \hline $\phi_{1,0}$ & $1$ & $0$ & $1$ &&&&&
\\ $G_2[\theta^2]$ & $q\Phi_1^2\Phi_2^2/3$ & $3$&&1&&&&
\\ $\phi_{2,2}$ & $q\Phi_2^2\Phi_6/2$  & $3$&1&&1&&&
\\ $G_2[\theta]$ & $q\Phi_1^2\Phi_2^2/3$  & $3$&&&&1&&
\\ $\phi_{1,6}$ & $q^6$ & $4$ & &1&1&1&1&
\\ $G_2[1]$ & $q\Phi_1^2\Phi_6/6$ & $4$ &&&&&&1
\\ \hline\end{tabular}\]
To construct the rest of the rows, we take a `non-unipotent' row of $B'$ -- in this case it is $(1,1,1,1,1,1)$ -- multiply it by the row $(-1)^{\pi_{\kappa/d}(S_i)}$ -- yielding in this case $\b v=(1,-1,-1,-1,1,1)$ -- and take the sum of the $i$th row of the matrix multiplied by the $i$th entry of $\b v$ -- yielding $(0,0,0,0,1,1)$. In the cyclic case, the non-unipotent rows are those of the exceptional characters, and for $B'$ these are always $(1,1,\dots,1)$.
\end{example}

In the definition of the perverse equivalence there is a bijection between the simple $B$- and $B'$-modules, and this was the assignment $S_i\mapsto T_i$ above given by identifying the Green correspondent in the degree $0$ term.

\medskip

We end this section with a remark about the $\pi_{\kappa/d}$-function. We have defined the $\pi_{\kappa/d}$-function on the unipotent $B$-characters, and we need it on the simple $B'$-modules. There are many potential bijections, and finding the correct one is non-trivial; we state the correct bijection in this article for the cyclic case, but in general we need technical information provided by the cyclotomic Hecke algebra to find this bijection. This topic will be explored in a later paper in this series. Similarly, in the cyclic case the collection $\mc R$ of relative projective modules is empty, and the description of this in the general is the subject of a later paper in this series.

\subsection{Genericity}
\label{ssec:genericity}

Let $R$ be the ring of integers of some algebraic number field, and let $E$ be finite subgroup of $\GL_n(R)$. The fundamental example of this for our purposes is $E$ a complex reflection group, for example, $R=\Z$ and $E$ a Weyl group, or the group  $R=\Z[\I]$ and $E=G_8\leq \GL_2(R)$.  Let $\bar\ell$ be an integer (not necessarily prime, nor even a prime power) with $(|E|,\bar\ell)=1$, and suppose that $\bar\ell$ is chosen so that there is a surjective homomorphism $R\to \Z_{\bar\ell}$ (the ring $\Z/\bar\ell\Z$), inducing the map $\alpha:\GL_n(R)\to \GL_n(\Z_{\bar\ell})$, whose restriction to $E$ is injective: such $\bar\ell$ are \emph{admissible} integers for $E$. This yields a map $E\to\Aut(\Z_{\bar\ell}^n)$ (where here $\Z_{\bar\ell}$ is considered simply as a group), so we may form the group $G_{\bar\ell}=(\Z_{\bar\ell})^n\rtimes E$; this group is in some sense generic in the integer $\bar\ell$. These groups can be found as the normalizers of $\Phi_d$-tori in groups of Lie type, where $|\Phi_d|=\bar\ell$.

Now specify $\bar\ell$ to be a power of a prime $\ell$, and let $k$ be a field of characteristic $\ell$ (as is our convention). In the situation of Brou\'e's conjecture, there is an isomorphism between $B'$ and $kG_{\bar\ell}$ (since if $\bar\ell$ is a prime power, this group algebra has only one block), so the simple $B'$-modules are in one-to-one correspondence with the simple $kG_{\bar\ell}$-modules: one of the key difficulties is to define a canonical such bijection, which is a fundamental part of the problem discussed in the remark at the end of the previous section.

The simple $kG_{\bar\ell}$-modules are `independent' of $\bar\ell$, in the sense that there is a natural identification of the simple $kG_{\bar\ell}$-modules with the ordinary $E$-characters,
and hence and identification between the simple $kG_{\bar\ell}$- and simple $k'G_{\bar\ell'}$-modules, where $\bar\ell'$ is a power of a different prime $\ell'$, $k'$ is a field of characteristic $\ell'$, and $\bar\ell'$ is chosen to have the same above properties as $\bar\ell$. We say that the simple $kG_{\bar\ell}$- and $k'G_{\bar\ell'}$-modules are \emph{identified}. With this identification, it is clear that we can also identify the projective $kG_{\bar\ell}$- and $k'G_{\bar\ell'}$-modules, and we do so. An obvious remark, but worth making, is that the projective modules have dimension $\bar\ell$, and also the defect group $D$ of the block $kG_{\bar\ell}$ is cyclic of order $\bar\ell$.

Let $T_1,\dots,T_e$ be an ordering of the simple $kG_{\bar\ell}$-modules, and pass this ordering onto the $k'G_{\bar\ell'}$-modules through the identification. The main philosophy of genericity is the following: given a fixed $\pi$-function $\pi:\{1,\dots,e\}\to\Z_{\geq 0}$, the outputs when applying the algorithm to the pairs $(kG_{\bar\ell},\pi)$ -- yielding complexes $X_i$ -- and $(k'G_{\bar\ell'},\pi)$ -- yielding complexes $X_i'$ -- should be `the same' (note we are assuming that the collections $\mc R$ are empty, although a version should exist with these included). By `the same', we mean
\begin{enumerate}
\item in the complexes $X_i$ and $X_i'$, the projective modules appearing in each degree are identified;
\item the multisets of composition factors of the cohomologies $H^{-j}(X_i)$ and $H^{-j}(X_i')$ are equal up to identification.
\end{enumerate}
If these two conditions hold for all $\bar\ell$ and $\bar\ell'$ at least $m$, then we say that the algorithm \emph{is generic} for $(E,\pi)$ with lower bound $m$.

In general, for any $(E,\pi)$ there should exist $m\in \N$ such that the algorithm is generic with lower bound $m$, although this is ongoing research of Rapha\"el Rouquier and the author. In the cyclic case however, i.e., $n=1$, it can fairly easily be proved with no restriction on $\bar\ell$ and $\bar\ell'$ (except that they are admissible for $(E,\rho)$ of course), i.e., $m=1$, and we give this now.

\medskip

Let $R=\Z[\e^{2\pi\I/e}]$ and $E$ be the cyclic subgroup of $R^*$ of order $e$. An admissible prime power $\bar\ell$ is one where the prime $\ell$ satisfies $e\mid(\ell-1)$. Before we start, we want to extend our definition of identified modules: consider the indecomposable $kG_{\bar\ell}$-modules. The group algebra $kG_{\bar\ell}$ has a single block, with cyclic defect group, and the Brauer tree of $kG_{\bar\ell}$ is a star, with $e$ vertices on the boundary. The projective cover of any simple module is uniserial: label the simple $kG_{\bar\ell}$-modules so that $T_1$ is the trivial module, and the first $e$ radical layers of $\Proj{T_1}$ are the simple modules $T_1$, $T_2$, \dots, $T_e$. For any $1\leq i,j\leq e$, there exists a unique uniserial module with $j$ layers and socle $T_i$: write $U_{i,j}$ for this indecomposable module. If $\bar\ell'$ is a power of another prime $\ell'$ with $e\mid(\ell'-1)$, then we can perform the same construction, and produce uniserial modules $U_{i,j}'$; we identify $U_{i,j}$ and $U_{i,j}'$.

\begin{prop}\label{prop:genericity} Let $E$ be as in the previous paragraph, and let $\pi:\{1,\dots,e\}\to\Z_{\geq0}$ be arbitrary. The algorithm is generic for $(E,\pi)$ with lower bound $0$. More precisely, let $\bar\ell$ and $\bar\ell'$ be powers of primes $\ell$ and $\ell'$ such that $e\mid (\ell-1),(\ell'-1)$, and write $G_1=G_{\bar\ell}$ and $G_2=G_{\bar\ell'}$, using the construction above. If $\pi:\{1,\dots,e\}\to\Z_{\geq 0}$ is a perversity function then, if $X_i$ and $X_i'$ ($1\leq i\leq e$) are the complexes describing the results of the algorithm applied to $(G_1,\pi)$ and $(G_2,\pi)$ respectively, we have:
\begin{enumerate}
\item for $1\leq j\leq \pi(i)$, the projective module in degree $-j$ for both $X_i$ and $X_i'$ is the projective cover $\Proj{T_\alpha}$ for some $1\leq \alpha\leq e$;
\item the module $H^{-j}(X_i)$ is a uniserial module $U_{\alpha,\beta}$, and this is identified with $H^{-j}(X_i')$;
\item writing $A_i$ for the term in degree $0$ of $X_i$, and $A_i'$ for the term in degree $0$ of $X_i'$, if $\pi(i)$ is even then $A_i$ and $A_i'$ are identified uniserial modules, and if $\pi(i)$ is odd then $\Omega(A_i)$ and $\Omega(A_i')$ are identified uniserial modules.
\end{enumerate}
\end{prop}
\begin{pf} Label the uniserial $kG_1$-modules of length at most $e$ (and hence also the $k'G_2$-modules via identification) $U_{\alpha,\beta}$, as above. Fix $1\leq i\leq e$, and for $kG_1$ and $1\leq j\leq \pi(i)$, we construct the modules $P_j$, $M_j$ and $L_j$, as in the algorithm, so that $P_j$ is the injective hull of $L_{j+1}$, and $M_j$ is the largest submodule of $P_j$, containing $L_{j-1}$, such that $M_j/L_{j-1}$ contains as composition factors only modules $T_\alpha$ where $\pi(\alpha)<j$. For $k'G_2$ we construct the modules $P_j'$, $M_j'$ and $L_j'$ similarly.

We proceed by reverse induction on $j$, starting with the case $j=\pi(i)$. Here, $P_j=\Proj{T_i}$ and $P_j'=\Proj{T_i}$, so (i) of the proposition is true for $j=\pi(i)$. Additionally, $H^{-\pi(i)}(X_i)$ is uniserial of length $r+1$ for some $r\geq 0$, so is the module $U_{i,r+1}$, with radical layers $T_{i-r},T_{i-r+1},\dots,T_i$ (with indices read modulo $e$); this is the largest $r\geq 0$ such that all of $T_{i-r},T_{i-r+1},\dots,T_{i-1}$ have $\pi$-value less than $\pi(i)$. Clearly $r<e$, as the $(e+1)$th socle layer of $\Proj{T_i}$ is $T_i$, which cannot be part of $H^{-\pi(i)}(X_i)$; hence $r$ is independent of the particular exceptionality of the vertex, and so $H^{-\pi(i)}(X_i)$ and $H^{-\pi(i)}(X_i')$ are both $U_{i,r+1}$, proving (ii) for $j=\pi(i)$.

Now let $j$ be less than $\pi(i)$. We notice that, if the top of $H^{-(j+1)}(X_i)$ -- which is the top of $M_{j+1}$ -- is $T_\alpha$ for some $\alpha$, then the projective module in degree $-j$ is $\Proj{T_{\alpha-1}}$; since $T_{\alpha-1}$ was not included in $M_{j+1}$, we must have that $\pi(\alpha-1)\geq j+1$. Since $H^{-(j+1)}(X_i)$ is identified with $H^{-(j+1)}(X_i')$, we see that both $P_j$ and $P_j'$ are $\Proj{T_{\alpha-1}}$, and so (i) is true for $j$. Also, if $P_{j+1}=\Proj{T_\beta}$, then the top of $P_{j+1}$, and hence the top of $L_j$, is $T_\beta$: by the remark just before the start of this subsection, $\pi(\beta)>j$. 

The module $M_j/L_{j-1}$ is uniserial, with radical layers $T_{\beta-s},T_{\beta-s+1},\dots,T_{\beta-1}$ (with indices read modulo $e$), and some $s$, possibly zero; this is the largest $s\geq 0$ such that all of $T_{\beta-s},T_{\beta-s+1},\dots,T_{\beta-1}$ have $\pi$-value less than $j$. Clearly $s<e$, as the $T_{\beta-e}=T_\beta$, and $\pi(\beta)>j$. Hence $H^{-j}(X_i)=U_{\beta-1,s}$; as the top of $L_j'$ is also $T_\beta$, we must also have that $H^{-j}(X_i')=U_{\beta-1,s}$, proving (ii) for this $j$. Hence, by reverse induction, (i) and (ii) hold.

It remains to deal with (iii). We note that $A_i=\Omega^{-1}(M_1)$ and $A_i'=\Omega^{-1}(M_1')$; since all projective modules have dimension $\bar\ell$ and $\bar\ell'$ respectively, $\dim(A_i)+\dim(M_1)=\bar\ell$ and $\dim(A_i')+\dim(M_1')=\bar\ell'$. As the tops of $M_1$ and $M_1'$ are identified simple modules, the socles of $A_i$ and $A_i'$ are identified simple modules; as $A_i$ and $A_i'$ are determined by their dimension and their socle, we need to show that if $\pi(i)$ is even then $\dim A_i=\dim A_i'$, and if $\pi(i)$ is odd then $\dim(\Omega(A_i))=\dim(\Omega(A_i'))$, or equivalently $\dim(M_1)=\dim(M_1')$.

Firstly, $\dim(L_j)+\dim(M_j)=\bar\ell$, and $\dim(M_j)=\dim(L_{j+1})+\dim(H^{-j}(X_i))$; by repeating this calculation, we see that if $\pi(i)-j$ is even, we have
\[ \dim(M_j)=\sum_{\alpha=j}^{\pi(i)} (-1)^{\alpha-j}\dim(H^{-\alpha}(X_i)).\]
If $\pi(i)-1$ is even, so $\pi(i)$ is odd, then $\dim(M_1)=\dim(M_1')$, as the cohomology of $X_i$ and $X_i'$ is the same, yielding (iii) in this case. If $\pi(i)$ is even,
\[ \dim(A_i)=\sum_{\alpha=1}^{\pi(i)} (-1)^{\alpha-j}\dim(H^{-\alpha}(X_i)),\]
and so we get $\dim(A_i)=\dim(A_i')$, as needed for (iii).
\end{pf}

Because of Proposition \ref{prop:genericity}, we can run the algorithm constructing perverse equivalences `generically', at least for cyclic defect groups. In this situation, let $E$ be a cyclic group of order $e$, and $\pi:\{1,\dots,e\}\to\Z_{\geq 0}$ be a perversity function. We say that we \emph{apply the algorithm generically} to $(E,\pi)$ if we apply the algorithm to the pair $(kG_\ell,\pi)$ for some prime $\ell\equiv 1\bmod e$. The data we retrieve are:
\begin{enumerate}
\item \emph{generic complexes} for each $i$, that is, a sequence $(n_{i,1},\dots,n_{i,\pi(i)})$ of $\pi(i)$ integers in $[1,e]$, with $n_{i,j}$ being the label of the projective modules in degree $\pi(i)+1-j$, so that for example $n_{i,1}=i$;
\item \emph{generic cohomology} for each $i$, that is, a sequence $(M_{i,1},\dots,M_{i,\pi(i)})$ of $\pi(i)$ uniserial modules of dimension at most $e$ (these exist for any admissible $\ell$, with $M_{i,j}$ being the module $H^{-(\pi(i)+1-j)}(X_i)$ in the complex $X_i$ associated to $T_i$), and the associated \emph{generic alternating sum of cohomology};
\item a \emph{generic Green correspondent} for each $i$, that is, if $C_i$ is the module in degree $0$, then the generic Green correspondent is the integer pair $(c_{i,1},c_{i,2})$, where $C_i$ has Socle $T_{c_{i,1}}$ and top $T_{c_{i,2}}$. (There is a unique such uniserial module with dimension at most $e$, and a unique such uniersial module with dimension at least $\bar\ell-e$, so this pair, together with the parity of $\pi(i)$, determine $C_i$.)
\end{enumerate}
This generic setup will be needed particularly in Section \ref{ssec:perverseBrauer2}, where we want to compare the outputs of the algorithm when the acting group $E$ has order $e$ and $e-1$; of course there can be no prime $\ell$ such that $e\mid \ell$ and $(e-1)\mid\ell$, so we are forced to work in a generic situation.

\subsection{Perverse Equivalences and Brauer Trees}
\label{ssec:perverseBrauer1}

We continue with notation from previous sections, specialized to the cyclic defect group case: $E$ is a cyclic group of order $e$, acting faithfully on $Z_{\bar\ell}$, and $G_{\bar\ell}=Z_{\bar\ell}\rtimes E$. We label the simple $kG_{\bar\ell}$-modules $T_1,\dots,T_e$, with $T_1$ being the trivial module, as in Section \ref{sec:gensetup}. If $B$ has cyclic defect groups (and recall that $B'$ is its Brauer correspondent) then $B'$ is isomorphic to $kG_{\bar\ell}$; we will in this section describe a specific identification of the simple $B'$-modules and the $T_i$. If $\pi:\{1,\dots,e\}\to\Z_{\geq 0}$ is a perversity function we write $X_i$ for the complex corresponding to $T_i$ when we apply the algorithm to $(kG_{\bar\ell},\pi)$, or after the identification, $(B',\pi)$.

In \cite{rickard1989}, Rickard proved that there is a derived equivalence between any block with cyclic defect group and the block of the normalizer of the defect group; in fact, the equivalence he constructed is perverse, for some bijection between the simple $B$- and $B'$-modules. We will produce this particular bijection and perversity function here, and using Green's walk on the Brauer tree \cite{green1974} we will show that the Green correspondents of the simple $B$-modules are indeed the terms in degree $0$ of the complexes. (In terms of \cite{rickard1989}, the perversity function can easily be extracted, and the bijection is slightly more subtle.) The proof that Rickard's equivalence is perverse is in \cite{chuangrouquierun}, but our proof of the existence of a perverse equivalence does not require either paper.

\medskip

We now make some important remarks about the rest of this section and the next: in Sections \ref{ssec:defnalg} and \ref{ssec:genericity} a perversity function was defined on the simple modules for the algebra $kG_{\bar\ell}$, and in the case of a block $B$, on the Brauer correspondent $B'$. However, in the groups of Lie type, the perversity function is defined for the simple $B$-modules, and must be passed to the simple $B'$-modules via a bijection between the two sets of simples. Technically speaking, a function defined on the simple $B$-modules is not a perversity function; however, since it can be turned into one via a bijection between the simple $B$- and $B'$-modules, which we will always provide, we will often abuse nomenclature somewhat and conflate the two.

In this section a perversity function will be defined on the simple $B$-modules with the help of the Brauer tree, and in the next section we will find it useful to think of the perversity function as defined on the simple $B$-modules, with the bijection between the simple $B$- and $B'$-modules being altered.

\medskip

We begin with a result that will be of great use in computing the degree $0$ terms in the complexes of our putative equivalence. 

\begin{lem}\label{lem:easycomplex} Suppose that $\pi:\{1,\dots,e\}\to\Z_{\geq 0}$ is a perversity function, and that for all $1\leq j<e$ we have $\pi(j+1)-\pi(j)\leq 1$, and $\pi(1)-\pi(e)\leq 1$. Let $\bar\ell$ be a power of $\ell$ such that $e\mid(\ell-1)$, write $G_{\bar\ell}=Z_{\bar\ell}\rtimes Z_e$, and apply the algorithm to $(kG_{\bar\ell},\pi)$, to form complexes $X_i$. The cohomology $H^{-j}(X_i)$ of the complex $X_i$ is zero for $1\leq j<\pi(i)$; in other words, the cohomology is concentrated in degree $-\pi(i)$.
\end{lem}
\begin{pf} We use the notation $M_j$ and $L_j$ introduced in Section \ref{ssec:defnalg} for the complex $X_i$, and write $H_i=H^{-\pi(i)}(X_i)$. It is easy to see, since there is at most one copy of $T_i$ in $H_i$, that $H_i$ is a (uniserial) module of dimension $a$, with $a\leq e$. In addition, as the socle of $L_{\pi(i)-1}$ is the module $T_{i-a}$, and it is not a part of $H_i$, we have that $\pi(i-a-1)\geq \pi(i)$. Finally, the hypothesis $\pi(j+1)-\pi(j)\leq 1$ clearly implies that $\pi(i-j)\geq \pi(i)-j$ and $\pi(i-a-j)\geq \pi(i)-j$ for all $j\geq 1$.

Our aim is to show by downward induction on $j$ for $\pi(i)>j>0$, that $H^{-j}(X_i)=0$, so that $L_j=M_j$, and therefore that $M_j=\Omega^{j-\pi(i)}(M_{\pi(i)})$. We assume that $M_j=L_j=\Omega^{j-\pi(i)}(M_{\pi(i)})$: as the socle of $M_{\pi(i)}$ is $T_i$ and the top of $M_{\pi(i)}$ is $T_{i-a+1}$, the socle of $M_{\pi(i)-j}$ is $T_{i-b}$ if $\pi(i)-j=2b$ and $T_{i-a+1-b}$ if $\pi(i)-j=2b+1$. As $L_{j-1}=P_j/M_j$, we have that the top of $L_{j-1}$ is the socle of $M_j$; hence
\[ \soc(P_{j-1}/L_{j-1})=\begin{cases}T_{i-b-1}&\pi(i)-j=2b\\ T_{i-a+1-b}&\pi(i)-j=2b+1.\end{cases}\]
In order for $H^{j-\pi(i)}(X_i)$ to be non-zero, we must have that the socle $T_\alpha$ of $P_{j-1}/L_{j-1}$ satisfies $\pi(\alpha)<\pi(i)-j$. However, $\pi(i-j)\geq \pi(i)-j$ and $\pi(i-a-j)\geq \pi(i)-j$ for all $j\geq 1$, proving that $H^{j-\pi(i)}(X_i)=0$ for all $\pi(i)>j>0$, as required.
\end{pf}

The point of this is that, if there is no cohomology except at degree $-\pi(i)$, then moving along the complex has the effect of applying $\Omega^{-1}$. We get the following corollary.

\begin{cor}Use the notation of the previous lemma. Write $H_i$ for $H^{-\pi(i)}(X_i)$. The term in degree $0$ of $X_i$ is $\Omega^{-\pi(i)}(H_i)$.
\end{cor}

Since the projective covers of the simple $kG_{\bar\ell}$-modules are all uniserial, the effect of applying $\Omega^i$ is very easy to describe: namely, $\Omega(T_i)$ is the indecomposable module of dimension $\bar\ell-1$ with socle $T_i$, and $\Omega^2(T_i)=T_{i+1}$.

\medskip

There is a natural one-to-one correspondence between non-exceptional characters and simple modules in a block with cyclic defect group, obtained by making a non-exceptional character of valency $1$ in its Brauer tree correspond to the unique simple module to which it is incident, then removing both character and module from the tree and repeating.

\begin{defn}\label{defn:canonicalperv} Let $B$ be a block with cyclic defect group. For $S$ a simple $B$-module, let $f(S)$ denote the length of the path from the exceptional vertex of the Brauer tree of $B$ to the vertex incident to $S$ that is closest to the exceptional vertex; let $r$ be the maximum of the $f(S)$. Write $\pi_0(S)=r-f(S)$, and $\pi_\alpha(S)=\pi_0(S)+\alpha$ for any $\alpha\geq 0$. We call $\pi_0$ the \emph{canonical perversity function}, and $\pi_\alpha$ the \emph{$\alpha$-shifted canonical perversity function}.
\end{defn}

If $\chi$ is a non-exceptional character in $B$, then $\chi(1)$ is congruent to either $1$ or $-1$ modulo $\ell$: we only consider $\alpha$-shifted canonical perversity functions such that $\chi(1)\equiv (-1)^{\pi_\alpha(\chi)}\bmod \ell$, where $\pi_\alpha$ is transferred from the simple $B$-modules to the non-exceptional $B$-characters using the correspondence described just before Definition \ref{defn:canonicalperv}. (Recall the remark at the start of this subsection about perversity functions being defined on $B$ and transferred to $B'$.)

\medskip

We recall Green's walk on the Brauer tree from \cite{green1974}, which can be used to construct the Green correspondents of the simple $B$-modules. Let $\chi$ be a non-exceptional character of maximal distance from the exceptional node: as $\chi$ lies on the boundary of the Brauer tree, $\chi$ has simple reduction modulo $\ell$, say $S$. We define $T_1$ by the statement that $\Omega^{\pi_\alpha(S)}(T_1)$ is the Green correspondent of $S$. (The module $T_1$ is always simple.)
We keep our notation for modules of $B'$, writing $T_i=\Omega^{2i}(T_1)$. We now wish to label the simple $B$-modules as $S_i$, so that, with respect to $\pi_\alpha$, the bijection between simple $B$- and $B'$-modules is $S_i\mapsto T_i$.

Starting from the vertex $\chi_1$ we walk around the edges of the Brauer tree in an anti-clockwise direction, labelling the edges $\delta(1)$, $2$, $\delta(2)$, \dots, $e$, $\delta(e)$, $1$. It is easy to see that every edge is labelled exactly twice, with $i$ and $\delta(j)$ for some $i$ and $j$, so that $\delta$ is a permutation of $\{1,\dots,e\}$. We then rotate the labelling of the edges anti-clockwise $\alpha$ times, so that for example if $\alpha=1$ then the first edge in the sequence is now labelled $1$, not $\delta(1)$.

Write $A_i$ for the simple module whose edge is labelled `$i$'. (Note that each edge is labelled `$i$' and `$\delta(j)$' for some $i$ and $j$.) The Green correspondent of $A_i$ is an indecomposable $B'$-module whose top is $T_i$ and whose socle is $T_{\delta^{-1}(i)}$ \cite[(6.1)]{green1974}, and its dimension is between $1$ and $e$, or between $\bar\ell-1$ and $\bar\ell-e$, depending on whether $\pi_\alpha(A_i)$ is even or odd.

If $\pi_\alpha(A_i)=2c$ is even, write $S_{\delta^{-1}(i)+c}=A_i$, and if $\pi_\alpha(A_i)=2c+1$ is odd, write $S_{i+c}=A_i$. (These indices should be taken modulo $e$.) An important point to notice, and that we will use in the proof of the next theorem, is that if we start from \emph{any} vertex $\chi$ of the Brauer tree and start our walk in the same way, finally rotating the labelling by $\pi_\alpha(\chi)$ rather than $\alpha$, then we get the same labelling of the $A_i$ and the $S_i$, except that the indices are shifted all by the same amount. This observation yields the following lemma.

\begin{lem}\label{lem:Sipi0isfixedpoint} Let $S$ be a simple module whose associated non-exceptional vertex lies on the boundary of the Brauer tree. If $\pi_\alpha(A_i)$ is even then $\delta(i)=i$, and if $\pi_\alpha(A_i)$ is odd then $\delta(i-1)=i$.
\end{lem}
\begin{pf} Travelling on Green's walk alternates between positive and negative vertices (i.e., vertices whose associated character is congruent to $+1$ and $-1$ modulo $\ell$), and also alternates between labels of the form `$i$' and `$\delta(i)$'. Hence a label of the form `$i$' always occurs when moving from a positive to a negative vertex, and `$\delta(i)$' when moving from negative to positive. The result now follows since when encountering a boundary vertex the walk doubles back and labels the same edge on consecutive steps of the walk.
\end{pf}

The most important case is when the Brauer tree is a line, and of course in this case we can be completely explicit: the ordering on the simple modules is to start from one end of the tree and travel to the exceptional node, then to repeat this from the other end.

\begin{lem}\label{lem:Siinline} Suppose that the Brauer tree of a block is a line, with $s$ vertices to the right of the exceptional node and $t$ vertices to the left. Assume that $s\geq t$. Write $\chi$ for the vertex farthest to the right, and start Green's walk from this vertex; let $\alpha\leq 1$ and consider the corresponding $\pi_\alpha$ and $S_i$. The simple modules $S_i$ for $1\leq i\leq s$ are the simple modules, in sequence, starting from $\chi$ and ending at the exceptional node, and for $s+1\leq i\leq t$ the $S_i$ are the simple modules, in sequence starting from the farthest-left vertex and ending at the exceptional node. Moreover, for each $i$, if $\pi_\alpha(A_i)=2c$ then $i+c\leq e+1$ and $\delta^{-1}(i)+c\leq e$, and if $\pi_\alpha(A_i)=2c+1$ then $i+c\leq e$ and $\delta^{-1}(i)+c\leq e$.
\end{lem}
\begin{pf} If $\alpha=0$, then the permutation $\delta$ swaps $i$ and $e+2-i$. For $1<i\leq s/2+1$ we have $\pi_\alpha(A_i)=2i-3$, so that $i+c=\pi_\alpha(A_i)+1\leq e$ and $\delta^{-1}(i)+c=(e+2-i)+(i-2)=e$. If $i>s/2+1+t$ then $\pi_\alpha(A_i)=2(e+1-i)$ so that $\delta^{-1}(i)+c=(e+2-i)+(e+1-i)=\pi_\alpha(A_i)+1$ and $i+c=(e+1-i)+i=e+1$. This proves the lemma for the modules to the right of the exceptional node; for the other side the result is similar.

If $\alpha=1$, then the permutation $\delta$ swaps $i$ and $e+1-i$. For $1<i\leq (s+1)/2$ we have $\pi_\alpha(A_i)=2i-1$, so that $i+c=\pi_\alpha(A_i)\leq e$ and $\delta^{-1}(i)+c=(e+1-i)+(i-1)=e$. If $i>(s+1)/2+t$ then $\pi_\alpha(A_i)=2(e+1-i)$ so that $\delta^{-1}(i)+c=(e+1-i)+(e+1-i)=\pi_\alpha(A_i)$ and $i+c=(e+1-i)+i=e+1$. This again proves the lemma for the modules to the right of the exceptional node; for the other side the result is similarly similar.
\end{pf}

At the moment it is not clear that distinct $A_j$ produce distinct $S_i$ in general. The first step along the way is to show that, if $\alpha\leq 1$ and $\chi$ is a vertex of maximal distance from the exceptional vertex with respect to which the edges are labelled, then as in the case of the line above, the index $i$ of $S_i$ need not be taken modulo $e$, i.e., that $i+c\leq e$ and $\delta^{-1}(i)+c\leq e$.

\begin{lem}\label{lem:Sismall} Let $\chi$ be a vertex of maximal distance from the exceptional node, and let the $A_i$ and $\delta$ be as constructed above. Let $\alpha=0$ or $\alpha=1$. If $\pi_\alpha(A_i)=2c$ is even then $i+c\leq e+1$ and $\delta^{-1}(i)+c\leq e$, and if $\pi_\alpha(A_i)=2c+1$ is odd then $i+c\leq e$ and $\delta^{-1}(i)+c\leq e$.
\end{lem}
\begin{pf} If one removes from the Brauer tree any edge that does not lie on the unique line connecting the edges $A_1$, $A_i$ and the exceptional node, then $e$ reduces by $1$ and $i$ and $\delta^{-1}(i)$ either remain the same or at least one is reduced by $1$, with the $\pi_\alpha$-function remaining constant. Hence, if $\psi$ denotes the vertex incident to $A_i$ that is farther from the exceptional node than the other, then the case where $i+c-e$ is maximal is where the Brauer tree is a line, with $\chi$ as one of the boundary vertices and either the exceptional node or $\psi$ as the other. This case is done in Lemma \ref{lem:Siinline}.
\end{pf}

Having proved something about the $S_i$, in every case, we can proceed with an inductive proof of the well definedness of the $S_i$, together with enough facts about $\pi_\alpha(S_i)$ to get a perverse equivalence between $B$ and $B'$.

\begin{thm}\label{thm:Greenbijectionwelldefined} Let $\pi_\alpha$, the $S_i$ and $T_i$ be as above, and let $C_i$ denote the Green correspondent of $S_i$ in $B'$.
\begin{enumerate}
\item The module $S_i$ is well defined; i.e., for a given $i$, there is a unique $A_j$ such that in the definition above we get $A_j=S_i$. Consequently, the $S_i$ form a complete set of simple $B$-modules.
\end{enumerate}
We may now produce a perversity function on the $T_i$ by defining $\pi_\alpha(i):=\pi_\alpha(S_i)$.
\begin{enumerate}
\setcounter{enumi}{1}
\item We have that $\pi_\alpha(i+1)-\pi_\alpha(i)\leq 1$ for all $1\leq i\leq e$ (where $e+1$ and $1$ are identified).
\item If $\pi_\alpha(A_i)=2c$ is even, then $\pi_\alpha(S_j)<\pi_\alpha(A_i)$ for $\delta^{-1}(i)+c>j\geq i+c$, and $\pi_\alpha(S_{i+c-1})\geq \pi_\alpha(A_i)$. If $\pi_\alpha(A_i)=2c+1$ is odd, then $\pi_\alpha(S_j)<\pi_\alpha(A_i)$ for $i+c>j\geq \delta^{-1}(i)+c+1$, and $\pi_\alpha(S_{\delta^{-1}(i)+c})\geq \pi_\alpha(A_i)$.
\item Denote by $X_i$ is the complex associated to $T_i$ when one applies the algorithm to $(B',\pi_\alpha)$. For $j>0$ we have $H^{-j}(X_i)=0$ unless $j=\pi_\alpha(S_i)$, and the $0$th term of $X_i$ is the module $C_i$. Consequently, there is a perverse equivalence between $B$ and $B'$, with perversity function $\pi_\alpha(-)$ and bijection given by $S_i\mapsto T_i$.
\end{enumerate}
\end{thm}
\begin{pf} We will begin by proving (i), (ii) and (iii), and then prove that (iv) follows from these.

Firstly, we show that the result holds for a given $\alpha$ if and only if it holds for $\alpha+2$. For this, we simply note that replacing $\alpha$ by $\alpha+2$ has the effect of cycling the $A_i$ (i.e., replacing $A_i$ by $A_{i+1}$), doing the same to the $T_i$, and increasing the $c$ used to relate the $A_i$ and the $S_i$ by $1$, so that $S_i$ is replaced by $S_{i+2}$. Since the $S_i$ are all shifted by $2$, this means that (i), (ii) and (iii) all hold for $\alpha$ if and only if they hold for $\alpha+2$.

We can therefore assume that $\alpha=0$ or $\alpha=1$, and in particular $A_1=S_1$ has Green correspondent either $T_1$ (if $\alpha=0$) or $\Omega(T_1)$ (if $\alpha=1$). Our plan is to remove the edge labelled by $S_1$ and use induction on the smaller Brauer tree, since the construction of the $A_i$ and $S_i$ do not require that the tree is a Brauer tree for a particular block. Notice that in the base case, where there is one simple module, all parts are true.

\medskip

We first assume that $\alpha=0$: as the edge corresponding to $A_1$ is labelled by `$1$' and `$\delta(1)$', removing it results in a Brauer tree with one fewer edge, and with edges labelled in a Green's walk from $2$ to $e$, instead of from $1$ to $e-1$.
By induction, if we subtract $1$ from the labels of these edges, we get a labelling $A_i'$ and $S_i'$, with $A_i'=A_{i+1}$ for $1\leq i\leq e-1$, and the $S_i'$ are well defined, with $S_i'=S_{i+1}$. Lemma \ref{lem:Sismall} now proves that the $S_i$ are well defined, proving (i).

For (ii), we notice that the relative orders of the $S_i$ and $S_i'$ are unchanged except for the insertion of $S_1$, so that (ii) is valid except possibly for $i=e$ and $i=1$. If $i=e$ then, since $\pi_\alpha(S_1)=0$ we clearly have $\pi_\alpha(S_{i+1})-\pi_\alpha(S_i)\leq 1$ for $i=e$.

If $i=1$ then we need to locate $S_2$. If $\pi_\alpha(S_2)\leq 1$ then the result holds, so $\pi_\alpha(S_2)\geq 2$, in which case $c\geq 1$. As $S_2=A_{j+c}$ or $A_{\delta^{-1}(j)+c}$, and $c,j\geq 1$, we get that $c=1$ and hence either $j=1$ or $\delta^{-1}(j)=1$; since $\delta(1)=1$, this implies $j=1$ in either case, a contradiction as $A_1=S_1$. Hence $\pi_\alpha(S_2)\leq 1$ and so (ii) holds.

It remains to show (iii). As with the previous part, the only possible problem is when we reintroduce $A_1=S_1$, for which $\pi_\alpha(S_1)=0$. For a given $i>1$, the only way that $\pi_\alpha(S_1)$ cannot satisfy $\pi_\alpha(S_j)<\pi_\alpha(A_i)$ is if $\pi_\alpha(A_i)=0$, in which case Lemma \ref{lem:Sipi0isfixedpoint} proves that the requirement is vacuous for this $i$. Hence we need to show the remaining requirement: if $\pi_\alpha(A_i)=2c$ is even then $i+c-1\leq e$ by Lemma \ref{lem:Sismall}, so that $S_{i+c-1}\neq S_1$ unless $i=2$ and $c=0$, in which case $\pi_\alpha(S_1)=\pi_\alpha(S_2)$, and if $\pi_\alpha(A_i)=2c+1$ then for the same reason $S_{\delta^{-1}(i)+c}\neq S_1$, this time with no exception. Hence by the inductive hypothesis the inequality does hold, and we complete the proof of (i), (ii) and (iii) when $\alpha=0$.

\medskip

Now assume that $\alpha=1$. In this case, removing the edge corresponding to $A_1=S_1$ again results in a Brauer tree with $e-1$ edges, but this time, the label `$\delta(1)$' must be replaced by `$\delta(e)$'; hence the resulting permutation on $\{2,\dots,e\}$ is simply $(1,e)\delta$. By induction, subtracting $1$ from the labelling produces $A_i'=A_{i+1}$, a permutation $\delta'$ on $\{1,\dots,e-1\}$ such that $\delta'^{-1}(i-1)=\delta^{-1}(i)-1$ unless $\delta^{-1}(i)=1$, in which case $\delta'^{-1}(i-1)=e-1$, and well-defined $S_i'$ with $1\leq i\leq e-1$.

If $\pi_\alpha(A_i)=\pi_\alpha(A_{i-1}')=2c+1$ is odd, then $S_{i+c}=S_{i+c-1}'$, and if $\pi_\alpha(A_i)=\pi_\alpha(A_{i-1}')=2c$ is even then $S_{\delta^{-1}(i)+c}=S_{\delta'^{-1}(i-1)+c}'$ unless $\delta^{-1}(i)=1$, in which case $S_{1+c}=S_{e-1+c}'$. Note that $A_1=S_1$: if $\delta^{-1}(i)=1$ then the edges labelled $i$ and $1$ share a vertex, so that $\pi_\alpha(A_i)$ differs from $\pi_\alpha(A_1)=1$ by at most $1$. Since $\pi_\alpha(A_1)$ is minimal, we have $\pi_\alpha(A_i)=1$ or $\pi_\alpha(A_i)=2$. Since we require $\pi_\alpha(A_i)$ to be even, we have $c=1$, so that $A_i=S_2=S_1'$, and hence $S_{i+1}=S_i'$ for all $1\leq i\leq e-1$. This proves (i).

For (ii), the same argument as in case $\alpha=0$ means we only need to investigate whether $\pi_\alpha(S_2)\leq \pi_\alpha(S_1)+1$. Let $j$ be such that $A_j=S_2$; if $\pi_\alpha(A_j)=2c+1$ is odd then since $j+c\leq e$ we have that $j=2$ and $c=0$, so that $\pi_\alpha(S_2)=\pi_\alpha(S_1)$. On the other hand, if $\pi_\alpha(A_j)=2c$ is even then $\delta^{-1}(j)+c=2$ then $c\leq 1$ (with $c=1$ implying that $\delta(1)=j$), in which case $\pi_\alpha(S_2)=\pi_\alpha(S_1)+1$, as $\pi_\alpha(A_j)\neq 0$. This completes the proof of (ii).

As with the case where $\alpha=0$, the only possible problem is for $A_1=S_1$, for which $\pi_\alpha(A_1)=1$ is minimal. For a given $i>1$, the only way that $\pi_\alpha(S_1)$ cannot satisfy $\pi_\alpha(S_j)<\pi_\alpha(A_i)$ is if $\pi_\alpha(A_i)=1$, in which case Lemma \ref{lem:Sipi0isfixedpoint} proves that the requirement is vacuous for this $i$. Hence we need to show the remaining requirement: if $\pi_\alpha(A_i)=2c$ is even then $i+c-1\leq e$ by Lemma \ref{lem:Sismall}, so that $S_{i+c-1}\neq S_1$ unless $i=2$ and $c=0$, in which case $\pi_\alpha(S_1)=1>0=\pi_\alpha(A_i)$ (and in any case, there is no such $A_i$), and if $\pi_\alpha(A_i)=2c+1$ then for the same reason $S_{\delta^{-1}(i)+c}\neq S_1$ unless $\delta^{-1}(i)=1$ and $c=0$, in which case $\pi_\alpha(S_i)=\pi_\alpha(A_i)$. Hence by the inductive hypothesis the inequality does hold, and we complete the proof of (i), (ii) and (iii) when $\alpha=1$ as well.

\medskip

It remains to prove that the first three parts imply the existence of a perverse equivalence, for which we require Green's walk. Let $1\leq i\leq e$, consider $A_i=S_{i+c}$ (assume that $\pi_\alpha(A_i)=2c+1$ is odd, as the even case is exactly the same) and the complex $X_j$ of $B'$-modules corresponding to $S_{i+c}$. By (iii) we have that $H^{-\pi(i)}(X_i)$ is the uniserial module $M_i$ of dimension at most $e$, with socle $T_{i+c}$ and top $T_{\delta^{-1}(i)+c+1}$. By (ii) and Lemma \ref{lem:easycomplex}, the $0$th term of the complex $X_i$ is the module $\Omega^{-\pi(i)}(M_i)$, which is a uniserial module of dimension at least $\bar\ell-e$ with socle $T_{\delta^{-1}(i)}$ and top $T_i$; this is the Green correspondent of $A_i=S_{i+c}$, and hence we get a perverse equivalence, as needed.
\end{pf}

The bijection given when $\alpha=0$ is called the \emph{canonical bijection}, and for a given $\alpha$ it is called the \emph{$\alpha$-shifted canonical bijection}. We summarize the results of this section for future reference.

\begin{cor}\label{cor:standardperversecyclic} Let $B$ be a block of $kG$ with a cyclic defect group $D$, and let $B'$ be its Brauer correspondent in $k\Norm_G(D)$. The $\alpha$-shifted canonical perversity function, together with the $\alpha$-shifted canonical bijection, yields a perverse equivalence between $B$ and $B'$.
\end{cor}

\subsection{A Family of Perverse Equivalences}
\label{ssec:perverseBrauer2}

We will describe a family of perverse equivalences for blocks with cyclic defect groups: by varying the perversity function in a natural way, we get infinitely many different perverse equivalences, for \emph{some} bijection between the simple modules. We use the canonical perversity function and bijection from the previous section. As in the previous section, the canonical ordering on $B'$ is the ordering where $T_i$ is the $i$th radical layer of the projective cover of $T_1$, for all $1\leq i\leq e$. Therefore, if the exceptionality of the vertex of the Brauer tree is $1$, then the projective cover of $T_1$ has radical layers
\[ 1/2/3/4/\cdots/e/1.\]
As in the previous section, we extend the perversity function given in Definition \ref{defn:canonicalperv} to an arbitrary Brauer tree algebra, as since we will again be proceeding by induction on the number of vertices, we need to consider Brauer trees that do not necessarily come from groups. We will also be pursuing the same proof method as in the previous section -- removing a single simple module and using induction -- and so we need to compare results of the algorithm when there are $e$ and $e-1$ simple modules; we introduced the idea of applying the algorithm generically in Section \ref{ssec:genericity} for precisely this purpose.

\medskip

There are two situations that we are interested in: the first is comparing two different perversity functions, $\pi$ and $\pi'$, and we want to know whether, when we apply the algorithm generically to $(E,\pi)$ and $(E,\pi')$, if the sets generic alternating sums of cohomology of generic Green correspondents are identical; the second is where we have a block $B$ of a finite group and we have a perversity function $\pi$ on the simple $B$-modules, and we wish to know if there is a bijection between the simple $B$- and $B'$-modules that yields a perverse equivalence, with perversity function $\pi$. (Recall the remark made at the start of the previous section about defining perversity functions on simple $B$-modules.)

Let $E$ be the cyclic group of order $e$, and let $\pi$ and $\pi'$ be perversity functions from $\{1,\dots,e\}$ to $\Z_{\geq 0}$. We say that the perversity functions $\pi$ and $\pi'$ are \emph{algorithmically equivalent} if there exists some permutation $\rho\in\Sym(e)$ such that, for all $1\leq i\leq e$, if one applies the algorithm generically to $(E,\pi)$ and $(E,\pi')$, and consider the modules $T_i$ and $T_{\rho(i)}$ respectively, then the generic alternating sums of cohomology and generic Green correspondents of $T_i$ under $(E,\pi)$ and of $T_{\rho(i)}$ under $(E,\pi')$ are identical, and $\pi(i)\equiv \pi'(\rho(i))\bmod 2$. (Notice that we do not need to include the group $E$ in the definition because $|E|$ is the size of the domain of $\pi$ and $\pi'$.)

To turn the second situation described into a version of the first, let $\pi'$ be the $\alpha$-shifted canonical perversity function on $B$ for a suitable $\alpha$, passed through the $\alpha$-shifted canonical bijection between $B$ and $B'$; we ask whether there is a function $\bar\pi:\{1,\dots,,e\}\to\Z_{\geq 0}$ taking the same values as $\pi$ and such that $\bar\pi$ and $\pi'$ are algorithmically equivalent.

\bigskip

The bulk of this section will be spent proving that, given any perversity function $\pi$, one may construct another perversity function $\pi'$ that is algorithmically equivalent to $\pi$ via a certain map $\rho\in\Sym(e)$ and such that $\pi'(i)=\pi(i)$ if $\rho(i)=i$, and $\pi'(\rho(i))=\pi(i)+2$ otherwise; furthermore the map $\rho$ is easy to describe, being a single cycle on the non-fixed points.

As usual, let $B$ be a block of a finite group and let $B'$ be its Brauer correspondent. There are two obvious conditions that a perversity function $\pi$ \textbf{on the simple $B$-modules} must satisfy if there is to be a bijection from the simple $B$- and $B'$-modules $\sigma$ that yields a perverse equivalence between $B$ and $B'$:
\begin{enumerate}
\item if $S_j$ lies on a path between $S_i$ and the exceptional node, then $\pi(S_j)>\pi(S_i)$.
\item $\chi_i(1)\equiv (-1)^{\pi(S_i)}\bmod \ell$, where $\chi_i$ is the non-exceptional character in bijection with the simple module $S_i$, as described just before Definition \ref{defn:canonicalperv};
\end{enumerate}
The first condition is simply because, ordered by the $\pi$-function, the decomposition matrix of $B$ must be lower triangular; the second condition is required by the perfect isometry induced by the perverse equivalence (see Theorem \ref{thm:bijectionwithsigns}). The main result of this section is that these are the \emph{only} restrictions on $\pi$. This therefore classifies all perverse equivalences between a block with cyclic defect groups and its Brauer correspondent.

In order to prove this, we first prove Theorem \ref{thm:perversechange}, which works directly with perversity functions. We then interpret it in Theorem \ref{thm:perversebijchange} in the language above, of perversity functions on the simple $B$-modules and a bijection with the simple $B'$-modules. This leads to Theorem \ref{thm:allperversecyclic}, which proves the asserted classification above.

\medskip

The proof of Theorem \ref{thm:perversechange} is a similar approach to the main result of the previous section, removing a vertex of degree $1$ from the Brauer tree and using induction; we extract part of the inductive step into the next technical lemma.

\begin{lem}\label{lem:inductivegc} Let $\pi:\{1,\dots,e\}\to\Z_{\geq 0}$ be a perversity function, and apply the algorithm generically to $\pi$, yielding generic complexes $(n_{i,j})$, generic cohomology $(M_{i,j})$ and generic Green correspondents $(c_{i,1},c_{i,2})$. Suppose that the generic cohomology associated to $T_e$ is the sequence $(T_e,0,\dots,0)$, i.e., $M_{e,1}=T_e$ and $M_{e,j}=0$ for $j>1$.  Let $\bar\pi$ be the restriction of $\pi$ to the subset $\{1,\dots,e-1\}$, and apply the algorithm generically to $\bar\pi$, yielding $(\bar n_{i,j})$, $(\bar M_{i,j})$ and $(\bar c_{i,1},\bar c_{i,2})$ analogously. Fix $a\neq e$.
\begin{enumerate}
\item The sequences $(n_{a,j})$ and $(\bar n_{a,j})$ are identical if and only if neither $c_{a,1}$ nor $c_{a,2}$ is $e-x$ for some $0\leq x<\pi(e)/2$.
\item If $c_{a,1}=e-x$ for some $0\leq x<\pi(e)/2$, then:
\begin{enumerate}
\item $\bar c_{a,1}=e-x-1$;
\item for all $j\leq x$ we have $n_{a,2j}=e-x+j$ and $\bar n_{a,2j}=e-x+j-1$.
\end{enumerate}
\item If $c_{a,2}=e-x$ for some $0\leq x<(\pi(e)-1)/2$, then:
\begin{enumerate}
\item $\bar c_{a,2}=e-x-1$;
\item for all $j\leq x$ we have $n_{a,2j-1}=e-x+j$ and $\bar n_{a,2j-1}=e-x+j-1$.
\end{enumerate}
\item The multiplicity of each $T_\alpha$, $\alpha\neq e$, in the modules $M_{a,j}$ and $\bar M_{a,j}$ is the same, for each $j>0$.
\item Write $\pi(e)=2c+\delta$ where $\delta\in\{0,1\}$. We have $(c_{e,1},c_{e,2}=(e-c-\delta,e-c)$.
\end{enumerate}
\end{lem}
\begin{pf} In order to prove this, we choose primes $\ell$ and $\bar\ell$ such that $\ell\equiv 1\bmod e$ and $\bar\ell\equiv 1\bmod (e-1)$, form the groups $G_\ell=Z_\ell\rtimes Z_e$ and $G_{\bar\ell}=Z_{\bar\ell}\rtimes E$, let $k$ be of characteristic $\ell$ and $\bar k$ of characteristic $\bar \ell$, and apply the algorithm to $(kG_\ell,\pi)$ -- yielding complexes $X_i$ associated to $T_i$ -- and apply the algorithm to $(\bar k G_{\bar\ell},\bar\pi)$ -- yielding complexes $\bar X_i$ associated to $T_i$. We apply Proposition \ref{prop:genericity} repeatedly to pass between the generic and particular cases. Fixing $a\neq e$, write $P_j$ for the projective in degree $-j$ of the complex $X_a$, and similarly for $\bar P_j$ and $\bar X_a$.

Firstly, since $X_e$ has no cohomology outside of the module $T_e$, we must have that $\pi(e-j)\geq \pi(e)-2j+2$, for all $1\leq j\leq\pi(e)/2$. We prove (iv) while we prove (i) and (ii), noting that if $P_j$ and $\bar P_j$ are covers of the same simple modules for all $j$ then clearly the multiplicities of all simple modules $T_x$ ($1\leq x\leq e-1$) are the same in $H^{-j}(X_a)$ and $H^{-j}(\bar X_a)$, so we assume that this is not the case.

Consider the difference between the terms of $X_a$ and $\bar X_a$: use the notation $M_j$ and $\bar M_j$ as in Section \ref{ssec:defnalg}. In all degrees $-j$ with $j\geq\pi(e)$, the projective modules $P_j$ and $\bar P_j$ are labelled by the same simple module, since the module $T_e$ is taken in cohomology whenever it has the opportunity to be. In particular, the multiplicity of $T_x$ in $H^{-j}(X_a)$ and $H^{-j}(\bar X_a)$ coincide for $1\leq x\leq e-1$, proving (iv) for $j>\pi(e)$.

For $j<\pi(e)$ however, the modules $P_j$ and $\bar P_j$ will not coincide if $\bar P_j=\Proj{T_{e-1}}$ for some $j$, in which case $P_j=\Proj{T_e}$. Let $-\alpha$ be the lowest degree for which $P_\alpha\neq\bar P_\alpha$ (note $\alpha<\pi(e)$), so that $P_\alpha=\Proj{T_e}$ and $\bar P_\alpha=\Proj{T_{e-1}}$. We also see that therefore $H^{-j}(X_a)$ and $H^{-j}(\bar X_a)$ coincide for $\alpha<j\leq \pi(e)$, proving (iv) in this range as well. An easy induction yields that $H^{-\alpha+2j-1}(X_a)=0$ and hence $P_{\alpha-2j}=\Proj{T_{e-j}}$, since $\pi(e-j)\geq \pi(e)-2j+2$. What this proves is that, if $\alpha$ is even then the socle of the degree $0$ term of $X_a$ is $T_{e-\alpha/2}$ and that of $\bar X_a$ is $T_{e-\alpha/2-1}$, whereas if $\alpha$ is odd then the top of the degree $0$ term of $X_a$ is $T_{e-(\alpha-1)/2}$ and that of $\bar X_a$ is $T_{e-(\alpha+1)/2}$. This completes the proof of one direction of (i) and, assuming the other direction, proves all parts of (ii). It remains to prove (iv) for $j\leq \alpha$. In this case, for every other $\alpha-j$ even, $H^{-j}(X_a)=H^{-j}(\bar X_a)=0$ so the result holds there, and for $\alpha-j$ odd, the only time $H^{-j}(X_a)$ and $H^{-j}(\bar X_a)$ can differ is if $T_e$ is not taken in cohomology, in which case $T_{e-1}$ would not be either (as $\pi(e)\leq \pi(e-1)$, and so again the cohomologies coincide, completing the proof of (iv).

We now prove the converse for (i); let $0\leq x<\pi(e)/2$, and suppose that $a\neq e$ is such that the socle of the degree $0$ term of $X_a$ is $T_{e-x}$. We reverse the algorithm, and prove that we must have that $P_{2x}=\Proj{T_e}$, which will be enough. By the condition on the degree $0$ term, we see that $M_1$ has top $T_{e-x+1}$. Using the observation at the start of this proof, $\pi(e-x+1)\geq \pi(e)-2(x-1)+2>1$, so that $T_{e-x+1}$ cannot lie in $H^{-1}(X_a)$. Thus $L_1=M_1$ and $H^{-1}(X_a)=0$, so that $P_2=\Proj{T_{e-x+1}}$. A simple induction shows that, given that $\pi(e-x+j)\geq \pi(e)-2(x-j)+2>2j+2$, $L_{2j-1}=M_{2j-1}$ has top $T_{e-x+j}$, $H^{-2j-1}(X_a)=0$ and $P_{2j}=\Proj{T_{e-x+j}}$ until $j=x$, at which point $P_{2x}$ and $\bar P_{2x}$ will differ and the converse of (i) holds for the socle.

The proof of the result for the top, i.e., (iii), is similar and omitted.

The statement (v) is simply determining the socle and top of $\Omega^{-\pi(e)}(T_e)$.
\end{pf}

Call a subset $I\subs \{1,\dots,e\}$ \emph{cohomologically closed} for $\pi$ if, when we apply the algorithm generically to $\pi$, yielding generic cohomology $(M_{i,j})$, and whenever $x\in I$ and $T_x$ appears in the cohomology $M_{i,j}$ for some $i$ and $j$, then $i\in I$. (In other words, if $T_x$ appears in the cohomology of the complex corresponding to $T_i$ and $x\in I$, then $i$ is also in $I$.) This concept has a very natural interpretation for the canonical perversity function and bijection: if $J$ is a subset of the simple $B$-modules such that, if $S_j\in J$ and $S_i$ lies on the path between $S_j$ and the exceptional node, then $S_i\in I$, then the image of $J$ under the canonical bijection is cohomologically closed. It is this interpretation that the reader should bear in mind, especially when it is used in Theorem \ref{thm:allperversecyclic}.

\begin{thm}\label{thm:perversechange} Let $\pi$ be a perversity function. Let $I=\{x_1,\dots,x_m\}$, with $x_i<x_{i+1}$, and $I_1\subs I$ be cohomologically closed subsets for $\pi$. Define $\rho\in\Sym(e)$ by $\rho(i)=i$ for $i\not\in I$ and $\rho(x_i)=x_{i+1}$ (cycling indices modulo $m$), and let $\pi'(-)$ be defined by $\pi'(i)=\pi(i)$ for $i\notin I$, and $\pi'(\rho(i))=\pi(i)+2$ for $i\in I$. The functions $\pi$ and $\pi'$ are algorithmically equivalent, via the permutation $\rho\in\Sym(e)$. Moreover, the sets $I$ and $I_1$ are cohomologically closed for $\pi'$.
\end{thm}
\begin{pf} If $I=\{1,\dots,e\}$ then the same argument as at the start of the proof of Theorem \ref{thm:Greenbijectionwelldefined} yields the result, so we may assume that $I\subset \{1,\dots,e\}$. One may apply a cyclic permutation to the function $\pi$ so that $\pi(e)$ is minimal subject to not being in $I$; applying the algorithm generically to a cyclic permutation of $\pi$ results is the same generic objects, but with the integers and module labels cyclically permuted. In particular, when applying the algorithm generically to $\pi$, the generic cohomology corresponding to $T_e$ is $(T_e,0,\dots,0)$. Let $i\notin I$, and consider the generic complexes -- $(n_{i,j})$ and $(n_{i,j}')$ -- and cohomologies -- $(M_{i,j})$ and $M_{i,j}')$ -- of $T_i$ with respect to $\pi$ and $\pi'$; we claim that the generic complexes and cohomologies are identical.

To see this, firstly note that the $M_{i,j}$ only have composition factors $T_x$ for $x\notin I$ by the definition of cohomological closure. We need to prove that the $M_{i,j}'$ only involve $T_x$ for $x\notin I$ as well, for then $M_{i,j}=M_{i,j}'$ for all $j$ as the $\pi$- and $\pi'$-functions coincide outside of $I$; as the cohomology is identical, $n_{i,j}=n_{i,j}'$ for all $j$, and we have proved algorithmic equivalence for $i\notin I$.

Let $j$ be minimal subject to $M_{i,j}\neq M_{i,j}'$; then $M_{i,j}'$ contains some $T_{\rho(x)}$ for $x\in I$, and $\pi'(\rho(x))=\pi(x)+2<\pi(i)-j+1$ but $\pi(\rho(x))\geq \pi(i)-j+1$. (As $\pi'(\rho(x))\geq 2$ we must have that $j<\pi(i)-1$.) Choose $\rho(x)$ with this property so that $T_{\rho(x)}$ is in the smallest socle layer of $M_{i,j}'$: we see that $M_{i,j}$ is a submodule of $M_{i,j}'$, and the socle of $M_{i,j}'/M_{i,j}$ is $T_{\rho(x)}$, so that $n_{i,j+1}=\rho(x)$. By the description of the algorithm, the module $M_{i,j+2}$, if it is non-zero, has socle $T_{\rho(x)-1}$.

Notice that for $y$ such that $\rho(x)>y>x$ (with $e$ and $0$ identified), $y\notin I$, by the definition of $\rho$. If $\pi(y)<\pi(i)-j-1=\pi(i)-(j+2)+1$ for all $\rho(x)>y>x$, then $M_{i,j+2}$ contains as a composition factor each $T_y$; however, in this case, since $\pi(x)=\pi(\rho(x))-2<\pi(i)-j-1$, we see that $M_{i,j+2}$ would also contain $T_x$ as a composition factor, a contradiction since $T_x$ cannot be a composition factor by the second paragraph. From those $y$ such that $\pi(y)\geq \pi(i)-j-1$, choose $y$ so that $T_y$ is in the highest socle layer in the uniserial module with socle $T_{\rho(x)}$ and top $T_x$: we now claim that $T_x$ is a composition factor in the module $M_{y,1}$, which is a contradiction as $x\in I$ and $y\notin I$. This can easily be seen as $\pi(y)$ is greater than all $z$ from $y$ to $x$, including $x$, and this final contradiction completes the proof. Thus the generic complexes and generic cohomology for $i\notin I$ are the same for $\pi$ and $\pi'$.

\medskip

Let $i\in I$, and write $\bar\pi$ and $\bar\pi'$ for the restrictions of $\pi$ and $\pi'$ to $\{1,\dots,e-1\}$. The functions $\bar\pi$ and $\bar\pi'$ are algorithmically equivalent by the restriction of $\rho$ to $\Sym(e-1)$ by induction (as $\rho(e)=e$). Hence the generic Green correspondents for $T_i$ for $\bar\pi$ and for $T_{\rho(i)}$ for $\bar\pi'$ are identical; however, by Lemma \ref{lem:inductivegc}, we can construct the generic Green correspondent for $T_i$ for $\pi$ from that of $\bar\pi$, and similarly for $T_{\rho(i)}$ and $\bar\pi'$ and $\pi'$. This means that the generic Green correspondents for $T_i$ for $\pi$ and for $T_{\rho(i)}$ for $\pi'$ are identical as well. Hence the second criterion of being algorithmically equivalent -- that the generic Green correspondents match up -- is true for all $1\leq i\leq e$.

In addition, Lemma \ref{lem:inductivegc}(iii) states that the coefficient of $T_j$ ($1\leq j\leq e-1$) in the generic alternating sums of cohomologies of $\pi$ for $T_i$ and $\bar\pi$ for $T_i$ coincide, and similarly for $\pi'$ for $T_{\rho(i)}$ and $\bar\pi'$ for $T_{\rho(i)}$. However, as $\bar\pi$ for $T_i$ and $\bar\pi'$ for $T_{\rho(i)}$ have the same generic alternating sum of cohomologies, this means that so must $\pi$ for $T_i$ and $\pi'$ for $T_{\rho(i)}$, except possibly for the multiplicity of $T_e$. However, the multiplicity of $T_e$ is determined by the generic Green correspondent and the multiplicities of the other $T_j$ in the alternating sum of cohomology, and since these are the same for $X_i$ and $X_i'$, the multiplicity of $T_e$ in the alternating sum of cohomologies must also be the same. (To see this statement, use the same method of proof as that of the end of Proposition \ref{prop:genericity}.)

This proves that $\pi$ and $\pi'$ are algorithmically equivalent, as claimed.

Finally we prove that $I_1$ is cohomologically closed for $\pi'$. Using Lemma \ref{lem:inductivegc}(iv) we see that $I_1$ is cohomologically closed for $\bar\pi$, so by induction $I_1$ is cohomologically closed for $\bar\pi'$. Another application of Lemma  \ref{lem:inductivegc}(iv), together with the fact that the generic cohomology of the module $T_e$ for $\pi'$ is the sequence $(T_j,0,\dots,0)$, proves that $I_1$ is cohomologically closed for $\pi'$, as needed.
\end{pf}

We now translate this theorem into a statement about producing a new perverse equivalence between two blocks from an old one. We simply state this theorem, as it is merely a rewriting of the previous result.

\begin{thm}\label{thm:perversebijchange} Let $\pi$ be a perversity function on the simple $B$-modules $S_1,\dots,S_e$, and order the simple $B'$-modules $T_1,\dots,T_e$ in accordance with Section \ref{sec:gensetup}. Let $\sigma:\{1,\dots,e\}\to\{1,\dots,e\}$ be a bijection from the $S_i$ to the $T_i$, such that there is a perverse equivalence from $B$ to $B'$ with perversity function $\bar\pi(T_{\sigma(i)}):=\pi(S_i)$ and bijection $\sigma$. Let $I=\{x_1,\dots,x_m\}$, with $x_i<x_{i+1}$, and $I_1$ be cohomologically closed subsets for $\bar\pi$, with $I_1\subs I$, and let $J$ be the preimage of $I$ under $\sigma$.

Define $\pi'(S_i)=\pi(S_i)$ for $i\notin J$ and $\pi'(S_i)=\pi(S_i)+2$ for $i\in J$. Let $\rho\in\Sym(e)$ be defined by $\rho(i)=i$ for $i\notin I$ and $\rho(x_i)=x_{i+1}$ (cycling indices modulo $m$). Finally, define $\bar\pi'(T_{\rho(\sigma(i))}):=\pi'(S_i)$. There is also a perverse equivalence from $B$ to $B'$ with perversity function $\bar\pi'$ and bijection $\rho\circ\sigma$.
\end{thm}

Starting from the canonical perversity function and canonical bijection, we can therefore add $2$ to the perversity function for any collection of simple modules as long as we also do it to ones on a path to the exceptional node. The main theorem of this section is the result of allowing repeated uses of the previous theorem.

\begin{thm}\label{thm:allperversecyclic} Let $B$ be a block of $kG$ with a cyclic defect group $D$, and let $B'$ be its Brauer correspondent in $k\Norm_G(D)$. Let $\pi(-)$ be a $\Z_{\geq 0}$-valued function on the set $\{S_1,\dots,S_e\}$ of simple $B$-modules such that:
\begin{enumerate}
\item if $S_i$ and $S_j$ share a non-exceptional vertex in the Brauer tree of $B$, with $S_j$ closer to the exceptional vertex than $S_i$, then $\pi(S_j)-\pi(S_i)$ is positive;
\item if $\chi_i$ is the non-exceptional ordinary character associated to $S_i$, then $\chi(1)\equiv (-1)^{\pi(S_i)}\bmod \ell$.
\end{enumerate}
There is a bijection between the simple $B$- and $B'$-modules such that, via this bijection, there is a perverse equivalence from $B$ to $B'$ with $\pi$ as perversity function.
\end{thm}
\begin{pf} Write $\pi_0(S_i)$ for the appropriate $\alpha$-shifted perversity function with $\alpha\in\{0,1\}$, and let $\sigma_0:\{1,\dots,e\}\to\{1,\dots,e\}$ be the $\alpha$-shifted canonical bijection from the $S_i$ to the $T_i$. Let $m$ be the smallest even non-negative integer such that $\pi(S_i)+m\geq \pi_0(S_i)$ for all $i$, and let $\pi'(S_i)=\pi(S_i)+m$. By Theorem \ref{thm:perversebijchange}, the claimed result holds for $\pi$ if and only if it holds for $\pi'$, so we may replace $\pi$ by $\pi'$ and assume that $\pi(S_i)\geq \pi_0(S_i)$ for all $i$. (Notice that $\pi(S_i)-\pi_0(S_i)$ is even by the second hypothesis.) For each $j\geq 1$, let $J_j$ denote the set of all $1\leq i\leq e$ such that $\pi'(S_i)-\pi_0(S_i)\geq 2j$, and note that $J_j\subs J_{j-1}$. By the first hypothesis on $\pi$, the images of the $J_j$ under $\sigma_0$ are all cohomologically closed with respect to $\pi_0$. Suppose that $J_n\neq\emptyset$ but $J_{n+1}=\emptyset$.

Inductively for $j\geq 1$, write $\pi_j$ and $\sigma_j$ for the perversity function and bijection $\sigma_j:\{1,\dots,e\}\to\{1,\dots,e\}$ that results when applying Theorem \ref{thm:perversebijchange} with the (cohomologically closed) set $I_j=\sigma_{j-1}(J_j)$ and the perversity function $\pi_{j-1}$, with associated bijection $\sigma_{j-1}$, noting that at each stage all subsets $\sigma_j(J_x)$ with $x\geq j$ stay cohomologically closed with respect to $\pi_j$. Clearly, $\pi_n=\pi$ and $\pi_j$ and $\pi_{j-1}$ are algorithmically equivalent via $\sigma_j\circ\sigma_{j-1}^{-1}$ for all $j$, hence the result is proved.
\end{pf}

We will show in later sections that the perversity function on blocks with cyclic defect group, for groups of Lie type, does satisfy the hypotheses of this corollary in the cases where the Brauer tree is known.

\section{The Combinatorial Brou\'e Conjecture}
\label{sec:combbroueconj}

In this section we give a complete description of the combinatorial Brou\'e conjecture for unipotent blocks with cyclic defect groups, and give an outline of how to prove it for all blocks where the Brauer tree is known, which is all but two unipotent blocks for $E_8$ at this stage.

In order to give a perverse equivalence between a block and its Brauer correspondent, we need a perversity function and a bijection. The perversity function is given by $\pi_{\kappa/d}(-)$ applied to the unipotent $B$-characters; we need to provide a bijection between the simple $B$- and $B'$-modules, and also a bijection between the simple $B$-modules and unipotent $B$-characters. This latter bijection was given in the previous section -- we associate a vertex of valency $1$ on the Brauer tree of $B$ to its incident edge, remove both, and repeat the process -- but we repeat it below in a more general setting for all unipotent blocks. We now produce a bijection between simple $B$-modules and simple $B'$-modules.

In Section \ref{sec:cyclohecke} we recall the definition of a cyclotomic Hecke algebra. This is a deformation of the group algebra of the cyclic group $Z_e$ (in our case, being deformations of group algebras of any complex reflection group in general) that take parameters of the form $u_i=\omega_iq^{v_i}$ for $1\leq i\leq e$, where $v_i$ is a semi-integer and $\omega_i$ is a root of unity. These parameters are defined up to a global multiplication by any root of unity and any power of $q$.

To each parameter $u_i$ one can associate a generic degree, given in (\ref{eq:relativedegree}). For a given unipotent block $B$ with cyclic defect groups, it was proved in \cite{bmm1993} that there is a collection of parameters $u_1,\dots,u_e$ such that the generic degrees of the $u_i$ are the degrees of the unipotent characters in the block $B$, up to a global scaling factor of a polynomial. Using Lemma \ref{lem:scalingparams}, we see that, up to scaling by a power of $q$, the $v_i$ satisfy $v_i=-\aA(\chi_i)/e$, where $\chi_i$ is the unipotent character whose degree is the relative degree (up to scaling again) associated to the parameter $u_i$.

\if\withapp1 For the root $\omega_i$, if the Brauer tree is a line -- in particular if $G$ is a classical group -- then $\omega_i=\pm1$, with all parameters corresponding to characters on one side of the exceptional node having the same sign, so with the power of $q$ given above, this completely determines the parameters in this case. For the unipotent blocks of exceptional groups, a case-by-case description of the parameters is given in \cite{bmm1993} with some cases missing, and in the appendix here for all cases. Thus there is a bijection between unipotent ordinary characters of $B$ and parameters of the cyclotomic Hecke algebra. Finally, recall that the decomposition matrix for $B$ is conjecturally lower triangular in all cases, with the top square consisting of unipotent characters (and of course is for Brauer trees): this produces a natural bijection between the simple $B$-modules and the ordinary unipotent $B$-characters.

\else
For the root $\omega_i$, if the Brauer tree is a line -- in particular if $G$ is a classical group -- then $\omega_i=\pm1$, with all parameters corresponding to characters on one side of the exceptional node having the same sign, so with the power of $q$ given above, this completely determines the parameters in this case. For the unipotent blocks of exceptional groups, a case-by-case description of the parameters is given in \cite{bmm1993} with some cases missing, and on the author's website for all cases. Thus there is a bijection between unipotent ordinary characters of $B$ and parameters of the cyclotomic Hecke algebra. Finally, recall that the decomposition matrix for $B$ is conjecturally lower triangular in all cases, with the top square consisting of unipotent characters (and of course is for Brauer trees): this produces a natural bijection between the simple $B$-modules and the ordinary unipotent $B$-characters.
\fi

The collection $\omega_iq^{-v_i}$ (note the minus sign), upon the substitution $q\mapsto \zeta$, produces a complete set of $e$th roots of unity (they also do so without the minus sign, but it is these ones that we want). Furthermore, the Brauer tree of $B'$, the Brauer correspondent, is a star embedded in $\C$, with exceptional node at $0$ and $e$ edges -- corresponding to simple $B'$-modules -- equally spaced around $0$. In order to achieve a bijection between the simple $B'$-modules and the $e$th roots of unity, we need to determine the exact embedding of the Brauer tree in $\C$, i.e., the rotational position of the star.

In order to do this, we choose $\chi$ a unipotent character with minimal $\pi_{\kappa/d}$-function in the block $B$. This has simple reduction modulo $\ell$, since it must lie on the boundary of the Brauer tree, so write $S$ for the simple $B$-module to which it corresponds. By Green's walk on the Brauer tree either the Green correspondent $T$ is  simple, if $\pi_{\kappa/d}(\chi)$ is even, or $\Omega(T)$ is simple, if $\pi_{\kappa/d}(\chi)$ is odd. In either case, the Green correspondent $T$ lies on the doubled Brauer tree (see the end of Section \ref{sec:gensetup}), and we embed the Brauer tree of $B'$ in such a way so that $T$ is at position $\omega_\chi$. This fixes the rotation of the Brauer tree, and completes the description of the bijection; in particular, this allows us to pass the $\pi_{\kappa/d}$-function to the simple $B'$-modules, so we may apply the algorithm to $(B',\pi_{\kappa/d})$.

\begin{conj}\label{conj:combbroue} Let $B$ be a unipotent block of $kG$ with cyclic defect groups. The perversity function $\pi_{\kappa/d}$ given above, and the bijection between the simple $B$- and $B'$-modules, induce a perverse equivalence between $B$ and $B'$.
\end{conj}

In this paper we will prove this conjecture whenever the Brauer tree is known; our task, therefore, is twofold:
\begin{enumerate}
\item prove that the perversity function $\pi_{\kappa/d}$ satisfies the first condition in Theorem \ref{thm:allperversecyclic} (as it is known to satisfy the second by Theorem \ref{thm:bijectionwithsigns});
\item prove that the associated bijection described recursively in Theorem \ref{thm:perversebijchange} matches the bijection given above.
\end{enumerate}
The first task will be done case by case for the exceptional groups, but for the classical groups we go via the cyclotomic Hecke algebra. In Section \ref{sec:cyclohecke} we prove a result, Proposition \ref{prop:increasingperv}, which states that, if the Brauer tree is a line, then the $\pi_{\kappa/d}$-function increases towards the exceptional node if and only if the exponent of the $q$-part of the parameter of the cyclotomic Hecke algebra \emph{decreases}, in other words, the quantity $\aA(-)$ increases towards the exceptional node.

We therefore prove the following theorem over the course of Section \ref{sec:BrauerMinPi}, using the standard combinatorial devices of partitions and symbols introduced in Section \ref{sec:combclassgroups}.

\begin{thm}\label{thm:piincreasetoexc} Let $B$ be a unipotent block of $kG$ with cyclic defect groups, whose Brauer tree is a line. If $\chi$ and $\psi$ are two unipotent characters in $B$, with $\psi$ appearing on a minimal path from $\chi$ to the exceptional node, then $\aA(\chi)<\aA(\psi)$.
\end{thm}

This proves that the $\pi_{\kappa/d}$-function induces a perverse equivalence with \emph{some} bijection, completing the first objective given above. The second objective itself splits into two parts: the first is to prove that a unipotent character with minimal $\pi_{\kappa/d}$-function (amongst those of its block) has the correct image under the bijection; the second is to prove that the relative positions of the images of all unipotent characters are correct, i.e., that the bijection is correct up to a rotation of the Brauer tree of $B'$. Of course, combining these two statements yields that the bijection is correct, and proves the combinatorial Brou\'e conjecture.

\medskip

To prove that a unipotent character with minimal $\pi_{\kappa/d}$-function has the correct image, we note that using the bijection described above, $\Omega^{\pi_{\kappa/d}(S)}(T)$ is a simple $B'$-module, and its position on the Brauer tree has argument $\arg(\omega_\chi)+\pi_{\kappa/d}(S)\cdot \pi/e$. On the other hand, evaluating $\omega_\chi q^{\aA(\chi)}$ at $q=\zeta$ yields a root of unity with argument $\arg(\omega_\chi)+\aA(\chi)/e\cdot 2\pi\kappa/d$. We therefore need to prove the following theorem. 

\begin{thm}\label{thm:minpivalue} A unipotent character $\chi$ in $B$ with $\pi_{\kappa/d}(\chi)$ minimal satisfies $\pi_{\kappa/d}(\chi)=2\kappa\aA(\chi)/d$.
\end{thm}

By Lemma \ref{lem:add2pi}, if one moves from $\kappa$ to $\kappa+d$, the change in the $\pi$-function satisfies
\[ \pi_{(\kappa+d)/d}(\chi)-\pi_{\kappa/d}(\chi)=2(A(\chi)-A(\bo\lambda)),\]
whereas the theorem says it should be $2\aA(\chi)$. This yields the following corollary.

\begin{cor}\label{cor:minpiaiszero} A unipotent character $\chi$ in $B$ with $\pi_{\kappa/d}(\chi)$ minimal satisfies $a(\chi)=a(\bo\lambda)$, and $\pi_{\kappa/d}(\chi)=2\kappa(A(\chi)-A(\bo\lambda))/d$.
\end{cor}

In fact, while the method of proof for the classical groups is as above, for exceptional groups Corollary \ref{cor:minpiaiszero} is established first, and then that $\pi_{\kappa/d}(\chi)=2\kappa(A(\chi)-A(\bo\lambda))/d$ for $\kappa<d$, which yields Theorem \ref{thm:minpivalue} for all $\kappa$; we prove this theorem in Section \ref{sec:combclassgroups} as well.

The last part is to prove the statement about the relative position of the images in the bijection: for blocks whose Brauer tree is a line -- so the parameters all have roots of unity $\pm 1$ -- this is performed using the cyclotomic Hecke algebra. We introduce a combinatorial procedure in Section \ref{sec:perturbHecke} called \emph{perturbation}, and a generalization of the cyclotomic Hecke algebra associated to the principal $\Phi_d$-block for $d$ the Coxeter number, called the \emph{Coxeter Hecke algebra}. Perturbing a cyclotomic Hecke algebra involves replacing the parameter $q^a$ with lowest exponent by another $q^{a+d}$, and then reordering the parameters in order of decreasing exponent. (There are two other types of perturbation, involving replacing $-q^b$ by $-q^{b+d}$, and $q^a$ by $-q^{a+d/2}$ and $-q^b$ by $q^{b+d/2}$, where $-q^b$ is the negative parameter with lowest exponent.) Because of the conditions placed upon these three types of perturbations, only one is allowed for any cyclotomic Hecke algebra.

Given a cyclotomic Hecke algebra, the generic degrees and parameter specialization give a perversity function and bijection with roots of unity. The main result of Section \ref{sec:perturbHecke} is the statement that perturbing the cyclotomic Hecke algebra induces changes in both the perversity function and bijection, and these are identical to that given in Theorem \ref{thm:perversebijchange} for adding $2$ to the $\pi$-function associated to certain simple modules and cycling their images under the bijection. Thus when checking if the bijection induced by parameter specialization is consistent with the perversity function, we may perturb the cyclotomic Hecke algebra as often as we like before checking this.

Finally, in Section \ref{sec:CoxHecke} we prove that repeated perturbation eventually results in a Coxeter Hecke algebra, and prove that in this case the bijection \emph{is} consistent with the perversity function, finally proving the combinatorial Brou\'e conjecture whenever the Brauer tree is a line, in particular for classical groups.

\medskip
\if\withapp1
For blocks of exceptional groups whose Brauer tree is not a line, we unfortunately do not have a general method like the perturbation of the cyclotomic Hecke algebra. (It should exist, but developing the theory is currently outside of our understanding.) We instead resort to a case-by-case check, which is performed for three representative blocks, and we relegate the list of all blocks with cyclic defect groups for exceptional groups to the appendix; this gives Brauer trees and parameters in every case, and completes the proof of the combinatorial Brou\'e conjecture for all unipotent blocks whose Brauer tree is known.

\else
For blocks of exceptional groups whose Brauer tree is not a line, we unfortunately do not have a general method like the perturbation of the cyclotomic Hecke algebra. (It should exist, but developing the theory is currently outside of our understanding.) We instead resort to a case-by-case check, which is performed for three representative blocks, and we relegate the list of all blocks with cyclic defect groups for exceptional groups to the author's website; this gives Brauer trees and parameters in every case, and completes the proof of the combinatorial Brou\'e conjecture for all unipotent blocks whose Brauer tree is known.
\fi

\section{Evaluating $\pi_{\kappa/d}$}
\label{sec:evalpi}

This section contains some calculations of the $\pi_{\kappa/d}$-function needed for evaluating it on character degrees of classical groups. We assume that $d\geq 2$, as if $d=1$ and $f=f(q)$ is a polynomial such that $f(1)\neq 0$, then $\pi_{\kappa/d}(f)=2A(f)$, so this case is easy.

\begin{prop} Let $i$ and $j$ be integers with $i>j$. We have
\[ \pi_{\kappa/d}(q^i-q^j)=\frac{\kappa(i+j)}{d}+\left\lfloor \frac{\kappa(i-j)}{d}\right\rfloor+\frac{1}{2},\]
and
\[ \pi_{\kappa/d}(q^i+q^j)=\frac{\kappa(i+j)}{d}+\left\lfloor\frac{2\kappa(i-j)}{d}\right\rfloor-\left\lfloor\frac{\kappa(i-j)}{d}\right\rfloor.\]
\end{prop}
\begin{pf} Suppose that $j=0$. Then $\pi_{\kappa/d}(q^i-1)$ is the sum of $\kappa/d\cdot A(q^i-1)$, namely $\kappa i/d$, the number of $i$th roots of unity of positive argument less than $2\pi \kappa/d$ -- of which there are $\lfloor \kappa i/d\rfloor$ -- and $1/2$ (for the single root at $1$); this gives the result. The general case easily follows since $q^i-q^j=q^j(q^{i-j}-1)$. For the second equality, we have $q^i+q^j=(q^{2i}-q^{2j})/(q^i-q^j)$, so that 
\begin{align*} \pi_{\kappa/d}(q^i+q^j)&=\(\frac{\kappa(2i+2j)}{d}+\left\lfloor \frac{2\kappa(i-j)}{d}\right\rfloor+\frac{1}{2}\)-\(\frac{\kappa(i+j)}{d}+\left\lfloor \frac{\kappa(i-j)}{d}\right\rfloor+\frac{1}{2}\)
\\ &=\frac{\kappa(i+j)}{d}+\left\lfloor\frac{2\kappa(i-j)}{d}\right\rfloor-\left\lfloor\frac{\kappa(i-j)}{d}\right\rfloor.\end{align*}
\end{pf}

This yields the following proposition in an obvious way, which deals with the effect on the second term in the numerator for the character degrees for $\GL_n(q)$, which we will see in Section \ref{sec:combclassgroups}. (We also include a case that will be needed for symplectic and orthogonal groups.)

\begin{prop}\label{prop:diffsinBdslinear} Let $i$ and $j$ be integers, and let $d\geq 2$ be an integer. Write $\kappa(j-i)=ad+b$, where $0\leq b<d$. We have that
\[\pi_{\kappa/d}(q^{i+d}-q^j)-\pi_{\kappa/d}(q^i-q^j)=\begin{cases}2\kappa & i-j>0\\ 2(\kappa-a)-1 & -d<i-j<0 \\ 0 & i-j<-d\end{cases}.\]
Now write $\kappa(j-i)=ad+b-d/2$, where $0\leq b<d$. We have that
\[\pi_{\kappa/d}(q^{i+d}+q^j)-\pi_{\kappa/d}(q^i+q^j)=\begin{cases} 2\kappa& i-j>-d/2\kappa\\2(\kappa-a)+\delta_{0,b}&-d+d/2\kappa\leq i-j\leq -d/2\kappa\\0&i-j<-d+d/2\kappa\end{cases}\]
\end{prop}
\begin{pf} For the first equation, the only case needing comment is when $0>i-j>-d$, in which case we have
\begin{align*} \pi_{\kappa/d}(q^{i+d}-q^j)-\pi_{\kappa/d}(q^j-q^i)&=\frac{\kappa(i+j+d)}{d}-\frac{\kappa(i+j)}{d}+\left\lfloor \frac{\kappa(i-j+d)}{d}\right\rfloor-\left\lfloor \frac{\kappa(j-i)}{d}\right\rfloor
\\ &=2\kappa-\left(\left\lfloor \frac{\kappa(j-i)}{d}\right\rfloor-\left\lfloor \frac{\kappa(i-j)}{d}\right\rfloor\right)
\\ &=2\kappa-2\left\lfloor \frac{\kappa(j-i)}{d}\right\rfloor-\begin{cases}1 & b\neq 0\\ 0& b=0\end{cases}.
\end{align*}
(The last equality relies upon the simple statement that for $a>0$, $\lfloor -a\rfloor=-\lfloor a\rfloor$ if $a\in\Z$ and $\lfloor -a\rfloor=-\lfloor a\rfloor-1$ otherwise.) Of course, since $(\kappa,d)=1$ and $0<j-i<d$, $\kappa(j-i)$ cannot be divisible by $d$, so $b\neq 0$. For the second equation the same statement about the case needing comment holds, and here we have
\begin{align*} \pi_{\kappa/d}(q^{i+d}+q^j)-\pi_{\kappa/d}(q^j+q^i)&=\frac{\kappa(i+j+d)}{d}-\frac{\kappa(i+j)}{d}+\left\lfloor \frac{2\kappa(i-j+d)}{d}\right\rfloor-\left\lfloor \frac{2\kappa(j-i)}{d}\right\rfloor
\\ &\qquad\quad-\left\lfloor \frac{\kappa(i-j+d)}{d}\right\rfloor-\left\lfloor \frac{\kappa(j-i)}{d}\right\rfloor
\\ &=2\kappa-\left(\left\lfloor \frac{2\kappa(j-i)}{d}\right\rfloor-\left\lfloor \frac{2\kappa(i-j)}{d}\right\rfloor\right)+\left(\left\lfloor \frac{\kappa(j-i)}{d}\right\rfloor-\left\lfloor \frac{\kappa(i-j)}{d}\right\rfloor\right)
\\ &=2\kappa-2\left\lfloor \frac{2\kappa(j-i)}{d}\right\rfloor+2\left\lfloor \frac{\kappa(j-i)}{d}\right\rfloor+\begin{cases}1 & b=0\\ 0& b\neq 0\end{cases}.
\end{align*}
\end{pf}

When working with unitary groups our integers $d$ and $e$ (i.e., where $\ell\mid\Phi_d$ and there are $e$ unipotent characters in $B$) satisfy $e=d$ if $4\mid d$, $e=2d$ if $d$ is odd, and $e=d/2$ otherwise. Evaluating $\pi_{\kappa/d}((-q)^{i+e}-(-q)^j)-\pi_{\kappa/d}((-q)^i-(-q)^j)$ is much more complicated than the previous proposition, and so we will content ourselves with simply giving the special cases that we need, namely $i>j$ and $j=i+1$. These particular cases follow a similar pattern to the previous proposition, and so the proof is omitted.

\begin{prop}\label{prop:diffsinBdsunitary} Let $d\geq 2$ and $\kappa\geq 1$ be coprime integers. Write $e=d$ if $4\mid d$, $e=2d$ if $d$ is odd and $e=d/2$ otherwise. If $i>j$ are non-negative integers, then
\[ \pi_{\kappa/d}\left((-q)^{i+e}-(-q)^j\right)-\pi_{\kappa/d}\left((-q)^i-(-q)^j\right)=2\kappa\frac ed,\]
and
\[ \pi_{\kappa/d}\left((-q)^e+q\right)-\pi_{\kappa/d}(q+1)=\frac{2\kappa e}{d}-2\left\lfloor\frac{2\kappa}{d}\right\rfloor+2\left\lfloor\frac{\kappa}{d}\right\rfloor.\]
This latter quantity is positive unless $d=2$, in which case it is $-1$.
\end{prop}

Using Propositions \ref{prop:diffsinBdslinear} and \ref{prop:diffsinBdsunitary}, we can prove the next difference, which is necessary when evaluating $\pi_{\kappa/d}$ on character degrees for linear and unitary groups.

\begin{prop}\label{prop:diffbetweenprodslun} Let $d\geq 2$ and $\kappa\geq 1$ be coprime integers. Write $e=d$ if $4\mid d$, $e=2d$ if $d$ is odd and $e=d/2$ otherwise. We have
\[  \pi_{\kappa/d}\(\prod_{i=n+1}^{n+d} (q^i-1)\)-\pi_{\kappa/d}\(\prod_{i=m+1}^{m+d} (q^i-1)\)=2\kappa(n-m),\]
and
\[ \pi_{\kappa/d}\(\prod_{i=n+1}^{n+e} ((-q)^i-1)\)-\pi_{\kappa/d}\(\prod_{i=m+1}^{m+e} ((-q)^i-1)\)=2\kappa(n-m)\frac ed.\]
\end{prop}
\begin{pf} Notice that $\pi_{\kappa/d}(q^{n+d}-1)-\pi_{\kappa/d}(q^{n}-1)=2\kappa$ and
\[ \pi_{\kappa/d}\((-q)^{n+e}-1\)-\pi_{\kappa/d}\((-q)^n-1\)=\frac{2e\kappa}{d}:\]
hence if $m=n-1$ the result holds. For general $m$ it is an obvious induction.
\end{pf}

For the symplectic and orthogonal groups, as well as Proposition \ref{prop:diffsinBdslinear} we need to deal with polynomials like $(q^{2i}-1)$.

\begin{prop}\label{prop:diffbetweenprodsorthog} Let $d\geq 2$ and $\kappa\geq 1$ be coprime integers. Write $d'=d$ if $d$ is odd and $d'=d/2$ if $d$ is even. We have
\[  \pi_{\kappa/d}\(\prod_{i=n+1}^{n+d'} (q^{2i}-1)\)-\pi_{\kappa/d}\(\prod_{i=m+1}^{m+d'} (q^{2i}-1)\)=4\kappa(n-m)\frac{d'}d.\]
\end{prop}
\begin{pf} Notice that $\pi_{\kappa/d}(q^{2(n+d')}-1)-\pi_{\kappa/d}(q^{2n}-1)=4\kappa d'/d$ and hence if $m=n-1$ the result holds. For general $m$ it is an obvious induction.
\end{pf}

Finally, we will have to take so-called cohooks when $d$ is even, and this interchanges plus and minus.

\begin{prop}\label{prop:diffsinBdsorthog} Let $d\geq 2$ and $\kappa\geq 1$ be coprime integers. Assume that $d$ is even and write $d'=d/2$. Let $i$ and $j$ be integers, and write $\kappa(j-i)=ad+b-d/2$, with $0\leq b<d$. We have
\[ \pi_{\kappa/d}(q^{i+d'}-q^j)-\pi_{\kappa/d}(q^i+q^j)=\begin{cases}\kappa & i-j>0\\\kappa-2a+\delta_{b,0}&0<j-i<d' \\0 & j-i>d'\end{cases}.\]
Writing $\kappa(j-i)=ad+b$ with $0\leq b<d$, we have
\[ \pi_{\kappa/d}(q^{i+d'}+q^j)-\pi_{\kappa/d}(q^i-q^j)=\begin{cases}\kappa & i-j>0\\\kappa-2a-1& 0<j-i<d'\\ 0& j-i>d'.\end{cases}\]
\end{prop}
\begin{pf} Firstly, since $d$ is even $\kappa$ must be odd. This means that, for any integer $b$,
\[ \left\lfloor\frac{b}{d}\right\rfloor+\left\lfloor\frac{b+\kappa d'}{d}\right\rfloor-\left\lfloor\frac{2b}{d}\right\rfloor+\frac12=\frac \kappa2;\]
the rest of the first part is an easy calculation of the same type as the previous propositions. The second property is similar to previous statements and its proof is omitted.
\end{pf}

\section{Cyclotomic Hecke Algebras}
\label{sec:cyclohecke}

Cyclotomic Hecke algebras were first introduced in \cite{brouemalle1993}: in some sense they parametrize the unipotent characters belonging to a given unipotent block $B$ in a group of Lie type. The general definition involves a complex reflection group, but since we are only concerned about blocks with cyclic defect group, our complex reflection group is the cyclic group $Z_e$ and so the definitions are much easier.

\begin{defn} Let $e\geq 1$ be an integer, and let $\b u=(u_1,\dots,u_e)$ be a sequence of transcendentals over $\Z$. The \emph{cyclotomic Hecke algebra} $\mc H(Z_e,\b u)$ is the algebra
\[ \frac{\Z[\b u,T]}{\left(\,(T-u_1)(T-u_2)\ldots(T-u_e)\,\right)}.\]
The \emph{relative degree} associated to $u_i$ is, up to sign,
\begin{equation} \prod_{j\neq i} \frac{u_j}{u_i-u_j}.\label{eq:relativedegree}\end{equation}
\end{defn}

Notice that, by scaling $T$, we can replace the parameters $u_i$ with $\alpha u_i$ for any element $\alpha\in\Z[\b u]$; in \cite{brouemalle1993}, the authors use this to set $u_1=1$, but we will not do this here for reasons that will become clear later.

For a particular group of Lie type $G$ and unipotent block $B$ of $kG$ with cyclic defect groups, to produce the cyclotomic Hecke algebra of $B$ we need specializations of the parameters $u_i$. These are of the form $u_i\mapsto \omega_i q^{v_i}$, where $q$ is another transcendental, $\omega_i$ is a root of unity (at most a sixth root in fact if $G$ is not of Suzuki or Ree type) and the $v_i$ are rationals (in fact semi-integers). In \cite{brouemalle1993}, it was proved that there is a choice for the $v_i$ and $\omega_i$ such that the relative degrees associated with the $u_i$, multiplied by $\Deg(R_{\b T}^{\b G}(\bo\lambda))$, are equal to the generic degrees of the unipotent characters in $B$.

With this information, it is easy to reconstruct the exponents $v_i$ in the specialized parameters $\omega_i q^{v_i}$; this lemma is well known, and we reproduce it here for completeness. As usual, if $\psi$ is a unipotent character, let $a(\psi)$ denote the power of $q$ dividing the generic degree of $\psi$ (as a polynomial in $q$) and $A(\psi)$ denote the degree of the generic degree of $\psi$.

\begin{lem}\label{lem:scalingparams} Let $\chi_1,\dots,\chi_e$ be the unipotent characters in $B$. If $\mc H$ denotes the cyclotomic Hecke algebra of $B$ then, up to scaling, the specialized parameters $\omega_i q^{v_i}$ satisfy
\begin{equation} v_i=-\aA(\chi_i)/e=-\left(\frac{a(\chi_i)+A(\chi_i)}{e}-\frac{a(\bo\lambda)+A(\bo\lambda)}{e}\right),\label{eq:relativedegsep}\end{equation}
where we recall that $(\L,\bo\lambda)$ is a $d$-cuspidal pair for $B$.
\end{lem}
\begin{pf} By scaling, we can assume the result for $i=1$. Notice that the quotient of the relative degree for $\chi_i$ by that of $\chi_1$ is
\[ \frac{u_1}{u_i} \prod_{j\neq 1,i} \frac{u_1-u_j}{u_i-u_j}=\omega_1\omega_i^{-1} q^{v_1-v_i}\prod_{j\neq 1,i} \frac{\omega_1q^{v_1}-\omega_jq^{v_j}}{\omega_iq^{v_i}-\omega_jq^{v_j}}.\]
Notice that $a(-)$ and $A(-)$ are both homomorphisms from the multiplicative monoid of polynomials over $\C$ in rational powers of $q$ (without the zero polynomial) to the rationals under addition, and so to evaluate $a(f)+A(f)$ it suffices to do so on each factor of $f$. Clearly $a(\omega_iq^{v_i}-\omega_j q^{v_j})+A(\omega_iq^{v_i}-\omega_j q^{v_j})=v_i+v_j$, and so we get that $a(-)+A(-)$, applied to the quotient of specialized relative degrees, is
\[ 2(v_1-v_i)+(e-2)(v_1-v_i)=e(v_1-v_i).\]
Since $\bo\lambda$ and $\Deg(R_{\b T}^{\b G}(\bo\lambda))$ are the same for $\chi_i$ and $\chi_1$, we get $a(\chi_i)-a(\chi_1)+A(\chi_i)-A(\chi_1)=e(v_1-v_i)$, which is consistent with  (\ref{eq:relativedegsep}), as needed.
\end{pf}

In the case of classical groups, the signs $\omega_i$ are simply $\pm 1$, whereas for exceptional groups $\omega_i$ can have order up to $12$ for non-real characters. For exceptional groups, however, there is a finite list of possible cyclotomic Hecke algebras to construct, and we will simply consider each one in turn. For the classical groups however, we need to develop a general theory.

In what follows we let $\mc H$ be a cyclotomic Hecke algebra with specialized parameters $\omega_i q^{v_i}$, where $\omega_i=\pm 1$ and $v_i$ is a semi-integer. We introduce the definitions formally now. (Note that these definitions and notation are non-standard.)

\begin{defn} Let $\mc H=\mc H(Z_e,\b u)$ be a cyclotomic Hecke algebra, with specialization $u_i\mapsto \omega_i q^{v_i}$, with $\omega_i$ a root of unity in $\C$ and $v_i$ a rational. We say that $\mc H$ \emph{has type $(s,t)$ and ambiance $d$} if
\begin{enumerate}
\item $e=s+t$,
\item $\omega_1,\dots,\omega_s=1$, $\omega_{s+1},\dots,\omega_e=-1$,
\item $v_i>v_j$ for $1\leq i<j\leq s$ and $s+1\leq i<j\leq e$, and
\item if we evaluate $q$ at a primitive $d$th root of unity $\zeta$, the set of $\omega_i\zeta^{v_i}$ form a complete set of $e$th roots of unity (up to a global multiplication by a root of unity).
\end{enumerate}
We write $\chi_i$ for the relative degree associated to $u_i$ for $1\leq i\leq s$, and $\psi_i$ for the relative degree associated to $u_{s+i}$ for $1\leq i\leq t$. Similarly, we write $a_i=v_i$ for $1\leq i\leq s$ and $b_i=u_{s+i}$ for $1\leq i\leq t$.
\end{defn}

Hence the relative degrees of $\mc H$ are $\chi_1,\dots,\chi_s$ and $\psi_1,\dots,\psi_t$. As an example, if the specialized parameters are (in order) $1$, $q^{-2}$, $-q$, $-q^3$ then $\mc H$ has type $(2,2)$. There is an ordering of the unipotent characters (and hence the specialized parameters) of a unipotent $\Phi_d$-block (with cyclic defect groups) in any classical group such that the associated cyclotomic Hecke algebra has type $(s,t)$ and ambiance $d$ for some $s,t$ with $s+t=e$ being the number of unipotent $B$-characters.

\medskip

At this juncture we will summarize the ideas behind this definition, which should help the reader follow the rest of the proof. If $B$ is a unipotent block with cyclic defect groups in a classical group, then the Brauer tree of $B$ is a line, with $s$ vertices on one side of the exceptional node, all of whose associated unipotent characters have parameters $+q^{a_i}$, and $t$ nodes on the other side, all of whose associated characters have parameters $-q^{b_i}$, with the labelling so that $q^{a_1}$ and $-q^{b_1}$ label the vertices of degree $1$, as the diagram below suggests.

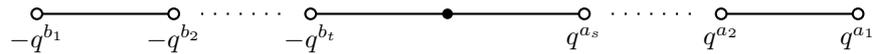
\begin{center}\begin{tikzpicture}[thick,scale=1.8]
\draw (0,-0.18) node {$-q^{b_1}$};
\draw (1,-0.18) node {$-q^{b_2}$};
\draw (2,-0.18) node {$-q^{b_t}$};
\draw (4,-0.18) node {$q^{a_s}$};
\draw (5,-0.18) node {$q^{a_2}$};
\draw (6,-0.18) node {$q^{a_1}$};

\draw (0,0) -- (1,0);
\draw (2,0) -- (4,0);
\draw (5,0) -- (6,0);
\draw[loosely dotted] (1.2,0) -- (1.8,0);
\draw[loosely dotted] (4.2,0) -- (4.8,0);

\draw (6,0) node [draw] (l2) {};
\draw (5,0) node [draw] (l2) {};
\draw (2,0) node [draw] (l2) {};
\draw (4,0) node [draw] (l2) {};
\draw (3,0) node [fill=black!100] (ld) {};
\draw (1,0) node [draw] (l4) {};
\draw (0,0) node [draw] (l4) {};
\end{tikzpicture}\end{center}
The ambiance, $d$, is the order of $q$ modulo $\ell$, i.e., so that the defect groups lie inside a $\Phi_d$-torus. For classical groups, we have that $e=d$, $e=2d$ or $e=d/2$, as we will see later.

If $d$ is the Coxeter number then the parameters are consecutive, in the sense that $a_{i+1}=a_i-1$ and similarly for the $b_i$, with $a_1=1$ corresponding to the trivial character; we will define a Coxeter Hecke algebra to be a generalization of this case. For this case it is easy to compute the $\pi_{\kappa/d}$-function, and we will show that the combinatorial form of Brou\'e's conjecture holds in these cases, by showing that the bijection induced is consistent with the $\pi_{\kappa/d}$-function.

We then consider an arbitrary cyclotomic Hecke algebra $\mc H$ with type $(s,t)$ and ambiance $d$, and `perturb' the specializations of the parameters one by one until we reach a Coxeter Hecke algebra. By keeping track of the changes to the positions of the parameters (recall that we maintain an ordering on them) and their associated $\pi_{\kappa/d}$-functions, we show that these two movements are consistent with those given in Theorem \ref{thm:perversebijchange}. This will prove combinatorial Brou\'e's conjecture for an arbitrary unipotent block $B$ with cyclic defect groups whose Brauer tree is a line, provided we can prove Theorems \ref{thm:piincreasetoexc} and \ref{thm:minpivalue}.

The next proposition proves that Theorem \ref{thm:piincreasetoexc} is equivalent to the statement that $\aA(-)$ increases towards the exceptional node, so that the parameter associated to a given character lines up in the way the diagram above suggests. The need for Theorem \ref{thm:minpivalue} arises from Theorem \ref{thm:perversebijchange}, where it is seen that subtracting $2$ from $\pi_{\kappa/d}(\psi)$ for all $\psi$ results in rotating the bijection by $2\pi/e$. If $\pi_{\kappa/d}$ is exactly twice $\kappa\cdot\aA(\chi)/d$ then subtracting $\pi_{\kappa/d}(\psi)$ from all $\psi$ makes $\psi$ in bijection with its Green correspondent, which is consistent with the case where $\pi_{\kappa/d}(\psi)=0$. This will be explained in more detail later, but this brief explanation should suffice to have an idea of the direction we will take.

\medskip

We now prove an important proposition about the $\pi_{\kappa/d}$ function on the relative degrees of such cyclotomic Hecke algebras.

\begin{prop}\label{prop:increasingperv} Let $\mc H$ be a cyclotomic Hecke algebra of type $(s,t)$ and ambiance $d$. For $1\leq i<s-1$ we have $\pi_{\kappa/d}(\chi_{i+1})>\pi_{\kappa/d}(\chi_i)$, and for $1\leq j\leq t-1$ we have $\pi(\psi_{i+1})>\pi(\psi_i)$.
\end{prop}
\begin{pf}Firstly, scale the parameters so that $d\mid a_{i+1}$. For any positive rational $x$, write $\bar x$ for the remainder upon division by $d$, so that $0\leq \bar x\leq d-1$. Define a positive parameter $q^{a_j}$ to be \emph{problematic} if $a_j>a_i$ and $\overline{\kappa a_i}>\overline{\kappa a_j}>\overline{\kappa a_{i+1}}=0$, and define a negative parameter $-q^{b_j}$ to be \emph{problematic} if $b_j>a_{i+1}$ and $\overline{\kappa a_i}>\overline{\kappa b_j+d/2}>\overline{\kappa a_{i+1}}=0$. Write $\kappa(a_i-a_{i+1})=\alpha+d\gamma$, where $\alpha=\overline{\kappa a_i}$. Notice that, since evaluation of $q$ at a primitive $d$th root of $1$ yields a bijection between the parameters and all $e$th roots of $1$, the number $z$ of problematic parameters is at most $(\alpha-1)e/d$. Write $z_+$ for the number of positive problematic parameters, and $z_-$ for the number of negative problematic parameters.

We firstly note that
\begin{equation} \frac{\chi_{i+1}(1)}{\chi_i(1)}=q^{a_i-a_{i+1}}\frac{\D \mathop{\prod_{j=1}^s}_{j\neq i}(q^{a_i}-q^{a_j})}{\D \mathop{\prod_{j=1}^s}_{j\neq i+1}(q^{a_{i+1}}-q^{a_j})}\cdot\frac{\D \prod_{j=1}^t(q^{a_i}+q^{b_j})}{\D\prod_{j=1}^t(q^{a_{i+1}}+q^{b_j})}.\label{eqcharquotient}\end{equation}
We apply the $\pi_{\kappa/d}$-function to this quotient. The first term clearly gives $2\kappa(a_i-a_{i+1})/d$, and the second term in (\ref{eqcharquotient}) yields
\begin{equation} \frac{\kappa(s-2)(a_i-a_{i+1})}{d}+\(\sum_{j=1}^s\left\lfloor\frac{\kappa |a_i-a_j|}{d}\right\rfloor-\left\lfloor\frac{\kappa |a_{i+1}-a_j|}{d}\right\rfloor\)\label{eqsecondterm}\end{equation}
Consider the sum in (\ref{eqsecondterm}) above: for a given $j$, notice that this expression is non-negative if $a_{i+1}>a_j$, and if $a_j>a_i$ we see that it is $-\gamma-1$ if $q^{a_j}$ is problematic, and $-\gamma$ otherwise. Hence (\ref{eqcharquotient}) is at least $\kappa(s-2)(a_i-a_{i+1})/d-\gamma(i-1)-z_+\geq \kappa(s-2)(a_i-a_{i+1})/d-\gamma(s-2)-z_+$.

The third term in (\ref{eqcharquotient}) yields
\begin{equation} \frac{\kappa t(a_i-a_{i+1})}{d}+\(\sum_{j=1}^t\left\lfloor\frac{2\kappa|a_i-b_j|}{d}\right\rfloor-\left\lfloor\frac{2\kappa|a_{i+1}-b_j|}{d}\right\rfloor-\left\lfloor\frac{\kappa |a_i-b_j|}{d}\right\rfloor+\left\lfloor\frac{\kappa |a_{i+1}-b_j|}{d}\right\rfloor\)\label{eqthirdterm}\end{equation}
Consider the sum in (\ref{eqthirdterm}) above: for a given $j$, as before, if $a_{i+1}>b_j$ then the expression is non-negative, so we may assume that $b_j>a_{i+1}$. Write $\kappa(b_j-a_{i+1})=\delta d+\beta$ where $\bar b_j=\beta$, and recall that $\kappa(a_i-a_{i+1})=\gamma d+\alpha$. We first deal with the case where $b_j>a_i$. We have
\[ \left\lfloor\frac{2\kappa|a_i-b_j|}{d}\right\rfloor-\left\lfloor\frac{2\kappa|a_{i+1}-b_j|}{d}\right\rfloor-\left\lfloor\frac{\kappa |a_i-b_j|}{d}\right\rfloor+\left\lfloor\frac{\kappa |a_{i+1}-b_j|}{d}\right\rfloor=-\gamma+\left\lfloor\frac{2(\beta-\alpha)}{d}\right\rfloor-\left\lfloor\frac{2\beta}{d}\right\rfloor-\left\lfloor\frac{\beta-\alpha}{d}\right\rfloor+\left\lfloor\frac{\beta}{d}\right\rfloor.\]
The last term is always $0$. For the rest of the terms, we have (noting that $\beta-\alpha$ cannot be equal to $\pm d/2$)
\[ \left\lfloor\frac{2(\beta-\alpha)}{d}\right\rfloor=\begin{cases} 1& d/2<\beta-\alpha,\\0& 0<\beta-\alpha<d/2,\\1& -d/2<\beta-\alpha<0,\\1& \beta-\alpha<-d/2;\end{cases} \quad -\left\lfloor\frac{2\beta}{d}\right\rfloor=\begin{cases} -1 & \beta>d/2,\\0 & \beta<d/2;\end{cases}\quad -\left\lfloor\frac{\beta-\alpha}{d}\right\rfloor=\begin{cases}1&\beta-\alpha>0,\\0&\beta-\alpha<0.\end{cases}\]
The sum of all of these becomes
\[ -\gamma+\left\lfloor\frac{2(\beta-\alpha)}{d}\right\rfloor-\left\lfloor\frac{2\beta}{d}\right\rfloor-\left\lfloor\frac{\beta-\alpha}{d}\right\rfloor+\left\lfloor\frac{\beta}{d}\right\rfloor=\begin{cases}-\gamma&\beta-\alpha>d/2,\\-\gamma&\beta<d/2,\beta-\alpha>-d/2,\\-\gamma-1&\beta>d/2,\beta-\alpha<d/2,\\-\gamma-1&\beta-\alpha<-d/2.\end{cases}\]
We see that this sum is $-\gamma-1$ if $-q^{b_j}$ is problematic and $-\gamma$ otherwise. Finally, if $a_i>b_j>a_{i+1}$ then
\begin{align*} \left\lfloor\frac{2(a_i-b_j)}{d}\right\rfloor&-\left\lfloor\frac{2(b_j-a_{i+1})}{d}\right\rfloor-\left\lfloor\frac{a_i-b_j}{d}\right\rfloor+\left\lfloor\frac{b_j-a_{i+1}}{d}\right\rfloor
\\ &=2\delta-\gamma+\left\lfloor\frac{2(\alpha-\beta)}{d}\right\rfloor-\left\lfloor\frac{2\beta}{d}\right\rfloor-\left\lfloor\frac{\alpha-\beta}{d}\right\rfloor+\left\lfloor\frac{\beta}{d}\right\rfloor=\begin{cases} 2\delta-\gamma+1&\alpha-\beta>d/2,\\2\delta-\gamma&\beta,\alpha-\beta<d/2,\\2\delta-\gamma-1&\beta>d/2;\end{cases}\end{align*}
we have $2\delta-\gamma\geq-\gamma$, and so this expression is at least $-\gamma-1$ if $\beta>d/2$ -- so that $-q^{b_j}$ is problematic -- and at least $-\gamma$ otherwise. Hence (\ref{eqthirdterm}) is at least $t\kappa(a_i-a_{i+1})/d-\gamma t-z_-$.

Adding these three contributions, we see that
\[ \pi_{\kappa/d}(\chi_{i+1})-\pi_{\kappa/d}(\chi_i)\geq \frac{\kappa e(a_i-a_{i+1})}{d}-\gamma(e-2)-z\geq \frac ed(d\gamma+\alpha)-\gamma(e-2)-\frac ed(\alpha-1)=2\gamma+1>0.\]
Hence $\pi_{\kappa/d}(\chi_{i+1})>\pi_{\kappa/d}(\chi_i)$, as needed.

\medskip

The proof that $\pi_{\kappa/d}(\psi_{i+1})>\pi_{\kappa/d}(\psi_i)$ is similar.
\end{pf}

This shows that the $\pi_{\kappa/d}$-function increases towards the exceptional node if and only if the $a_i$ and $b_i$ decrease (as they are negative) towards the exceptional node; the $a_i$ and $b_i$ are much easier to calculate than the $\pi_{\kappa/d}$-function, and we will prove this statement in Section \ref{sec:BrauerMinPi}.

\section{Combinatorics for Classical Groups}
\label{sec:combclassgroups}

The purpose of this section is to introduce the combinatorial objects needed for discussion of unipotent characters of classical groups, and then describe the degrees and distribution into blocks for unipotent characters.

\subsection{Partitions and Symbols}
\label{ssec:partssymbs}

In this section we introduce partitions and symbols. Much of this is well known and we summarize it briefly here, both to fix notation and for the reader's convenience.

We often identify a partition with its Young diagram, and talk of boxes for a partition. A \emph{hook} of a partition (really, a Young diagram) consists of a box $x$, all boxes below $x$, and all boxes to the right of $x$; if this is $t$ boxes in total, and there are $i$ boxes below $x$ or equal to $x$, then this hook is a \emph{$t$-hook} (or of length $t$) of \emph{leg length $i$}. \emph{Removing} a $t$-hook consists of deleting all boxes in a hook, and then pushing the boxes that were below and right of the hook up and to the left, creating a new partition.

If $\lambda=(\lambda_1,\lambda_2,\dots,\lambda_a)$ is a partition of $n$ (with $\lambda_i\geq \lambda_{i+1}>0$ being the \emph{parts}), the \emph{first-column hook lengths} of $\lambda$ is the set $X=\{x_1,\dots,x_a\}$, where $x_i=\lambda_i+a-i$, i.e., the lengths of the hooks of the boxes in the far-left column. It is easy to see that the set of all partitions (including the empty partition) is in bijection with the set of all finite subsets of $\Z_{>0}$, via sending a partition to its set of first-column hook lengths.

A \emph{$\beta$-set} is a finite subset of $\Z_{\geq 0}$. We introduce an equivalence relation on all such sets, generated by $X\sim X'$ if $X'=\{0\}\cup\{x+1\,:\,x\in X\}$. The \emph{rank} of $X$ is the quantity $\sum_{x\in X} x-a(a-1)/2$, where $a=|X|$. Notice that the rank is independent of the representative of the equivalence class of $\beta$-set; indeed, if we take the unique representative $X$ with $0\notin X$, then the rank of $X$ is the size of the partition $\lambda$ whose first-column hook lengths are $X$. We tend to order the elements of a $\beta$-set $X=\{x_1,\dots,x_a\}$ so that $x_i>x_{i+1}$.

If $X=\{x_1,\dots,x_a\}$ is a $\beta$-set, then the act of removing a $t$-hook is simply replacing some $x_i$ by $x_i-t$ (where $x_i-t\notin X$), and similarly adding a $t$-hook to $X$ involves replacing some $x_i$ by $x_i+t$ (assuming that $x_i+t\notin X$). The \emph{$t$-core} of $X$ is the $\beta$-set obtained by removing all possible $t$-hooks.

The $\beta$-sets of partitions can be more easily understood on the abacus. If $t$ is a positive integer, the \emph{$t$-abacus} is a diagram consisting of $t$ columns, or \emph{runners}, labelled $0,\dots,t-1$ from left to right. Starting with $0$ at the top of the left-most runner, we place all non-negative integers on the runners of the abacus, first by moving across the runners left to right, then moving down the runners, as below.
\begin{center}\begin{tabular}{|c|c|c|c|c|}
\hline 0&1&2&3&4
\\ 5&6&7&8&9
\end{tabular}\end{center}
If $X$ is a $\beta$-set, it can be represented on the $t$-abacus by placing a \emph{bead} at position $i$ whenever $i\in X$.

The act of adding or removing a $t$-hook is very easy to describe on the abacus: it consists of moving a bead one place on its runner, down or up respectively. The $t$-core of $X$ is obtained by moving all beads on the $t$-abacus as far upwards as possible.

\medskip

A \emph{symbol} is an unordered pair $\lambda=\{X,Y\}$ of subsets of $\Z_{\geq 0}$. We will write $X=\{x_1,\dots,x_a\}$ with $x_i>x_{i+1}$, and $Y=\{y_1,\dots,y_b\}$ with $y_i>y_{i+1}$. We introduce an equivalence relation on the set of symbols, which is generated by the relation that $\{X,Y\}\sim \{X',Y'\}$ if $X'=\{0\}\cup\{x+1\,:\,x\in X\}$ and $Y'=\{0\}\cup\{y+1\,:\,y\in Y\}$. If $X=Y$ then the symbol is \emph{degenerate}, and otherwise is \emph{non-degenerate}.

The \emph{defect} of $\lambda=\{X,Y\}$ is the quantity $|a-b|$, and the \emph{rank} of $\lambda$ is the quantity $\sum_{x\in X} x+\sum_{y\in Y} y-\lfloor(a+b-1)^2/4\rfloor$. Notice that equivalent symbols have the same defect and rank.

Let $\lambda=\{X,Y\}$ be a symbol. Adding a \emph{$t$-hook} to $\lambda$ involves adding $t$ to one of the elements of either $X$ or $Y$ to get another symbol $\mu$. Adding a \emph{$t$-cohook} to $\lambda$ involves adding $t$ to one of the elements of $X$ and transferring it to $Y$, or vice versa, to get another symbol $\mu$. By removing all $t$-hooks we get the \emph{$t$-core}, and by removing all $t$-cohooks we get the \emph{$t$-cocore}. Adding a $t$-hook does not change the defect of a symbol, but adding a $t$-cohook adds or subtracts $2$.

(If one envisages a symbol as a pair of $\beta$-sets, adding a $t$-hook is simply adding a $t$-hook on the abacus of one of the $\beta$-sets; a $t$-cohook is less easy to visualize.)

\subsection{Unipotent Characters and Blocks for Classical Groups}

In this section we describe the unipotent characters for the classical groups and their distribution into blocks.

\medskip

Let $G=\GL_n(q)$ for some $n$ and $q$. We describe briefly the unipotent characters and blocks of $\GL_n(q)$, as discussed in \cite{fongsrin1982}. The unipotent characters of $\GL_n(q)$ are labelled by partitions $\lambda$ of $n$, or equivalently $\beta$-sets of rank $n$ (up to equivalence). Let $X=\{x_1,\dots,x_a\}$ (with $x_i>x_{i+1}$) be a $\beta$-set of rank $n$, and let $\lambda$ be its corresponding partition. If $\chi_\lambda$ is the unipotent character of $\GL_n(q)$ corresponding to $\lambda$, then
\begin{equation} \chi_\lambda(1)=\frac{\D\(\prod_{i=1}^n (q^i-1)\)\(\prod_{1\leq i<j\leq a}(q^{x_i}-q^{x_j})\)}{\D \(\vphantom{\prod_{i=1}^a}q^{\binom{a-1}{2}+\binom{a-2}{2}+\cdots}\)\(\prod_{i=1}^a\prod_{j=1}^{x_i}(q^j-1)\)}.\label{eq:lineardegrees}\end{equation}
(Later we will refer to the `first' and `second' terms of the numerator and denominator of this equation: these have the obvious meanings.)

It is easy to see that $\chi_\lambda(1)$ does not depend on the choice of $\beta$-set $X$ representing $\lambda$. Two $\beta$-sets $X$ and $Y$, with partitions $\lambda$ and $\mu$, have the same $d$-core if and only if the corresponding unipotent characters, $\chi_\lambda$ and $\chi_\mu$, lie in the same $\ell$-block of $G$: the $d$-cuspidal pair for that block has character labelled by the $d$-core of $\lambda$.

\medskip

Let $G=\GU_n(q)$ for some $n$ and $q$, and write $d$ and $e$ for the multiplicative orders of $q$ and $-q$ respectively modulo $\ell$; then $e=d$ if $4\mid d$, $e=2d$ if $d$ is odd and $e=d/2$ otherwise. As with the linear groups, the facts about unipotent characters and blocks that we need are taken from \cite{fongsrin1982}. The unipotent characters of $\GU_n(q)$ are similar to those of $\GL_n(q)$, in that they are again associated to partitions of $n$. If $\chi_\lambda$ is the unipotent character of $\GL_n(q)$ associated to $\lambda$ and $\phi_\lambda$ is the unipotent character of $\GU_n(q)$ associated to $\lambda$, then the degree of $\phi_\lambda$ is obtained from that of $\chi_\lambda$ by replacing $q$ with $(-q)$ (with a sign change if this makes the character degree negative). In the expansion of $\phi_\lambda(1)$ into powers of $q$ and cyclotomic polynomials, this has the effect of replacing $\Phi_r$ with $\Phi_{2r}$ and vice versa, whenever $r$ is odd.

The structure of the $\ell$-blocks of $G$ is similar as well: these are parametrized by $e$-cores, and two unipotent characters $\phi_\lambda$ and $\phi_\mu$ lie in the same $\ell$-block of $G$ if and only if $\lambda$ and $\mu$ have the same $e$-core: the $d$-cuspidal pair for that block has character labelled by the $e$-core of $\lambda$.

\medskip

For classical groups of types $B$, $C$ and $D$, the unipotent characters for a group of Lie type of rank $n$ are parametrized by symbols $\Lambda=\{X,Y\}$ of rank $n$, with each non-degnerate symbol parametrizing one character and a degenerate one parametrizing two. Let $X=\{x_1,\dots,x_a\}$ and $Y=\{y_1,\dots,y_b\}$, with $x_i>x_{i+1}$ and $y_i>y_{i+1}$, and let $\delta$ be the defect of $\Lambda$, the quantity $|a-b|$. The symbols of odd defect and a given rank $n$ parametrize the unipotent characters of the groups of type $B_n$ and $C_n$, whereas the symbols of defect divisible by $4$ correspond to unipotent characters of the groups of type $D_n$ (with two unipotent characters corresponding to each degenerate symbol), and symbols of defect congruent to $2$ modulo $4$ correspond to unipotent characters of the groups of type ${}^2\!D_n$.

\medskip

In the case of $B_n$ and $C_n$, if $\chi_\Lambda$ is the unipotent character corresponding to the symbol $\Lambda$ (which has odd defect), then
\begin{equation} \chi_\Lambda(1)=\frac{\D\(\prod_{i=1}^n (q^{2i}-1)\)\(\prod_{1\leq i<j\leq a}(q^{x_i}-q^{x_j})\)\(\prod_{1\leq i<j\leq b}(q^{y_i}-q^{y_j})\)\(\prod_{i,j}(q^{x_i}+q^{y_j})\)}{\D 2^{(a+b-1)/2} q^{\binom{a+b-2}{2}+\binom{a+b-4}{2}+\cdots}\(\prod_{i=1}^a\prod_{j=1}^{x_i}(q^{2j}-1)\)\(\prod_{i=1}^b\prod_{j=1}^{y_i}(q^{2j}-1)\)}.\label{eq:sympdegrees}\end{equation}
As with the linear and unitary groups, this degree is invariant under the equivalence relation on symbols.

\medskip

In type $D_n$, so $G=(\CSO_{2n}^+)^0(q)$, if $\chi_\Lambda$ is the (or `a' if $\Lambda$ is degenerate) unipotent character corresponding to the symbol $\Lambda$ (which has defect divisible by $4$), then
\begin{equation} \chi_\Lambda(1)=\frac{\D (q^n-1)\(\prod_{i=1}^{n-1}(q^{2i}-1)\)\(\prod_{1\leq i<j\leq a}(q^{x_i}-q^{x_j})\)\(\prod_{1\leq i<j\leq b}(q^{y_i}-q^{y_j})\)\(\prod_{i,j}(q^{x_i}+q^{y_j})\)}{\D 2^c q^{\binom{a+b-2}{2}+\binom{a+b-4}{2}+\cdots}\(\prod_{i=1}^a\prod_{j=1}^{x_i}(q^{2j}-1)\)\(\prod_{i=1}^b\prod_{j=1}^{y_i}(q^{2j}-1)\)},\label{eq:orthdegrees1}\end{equation}
where $c=\lfloor (a+b-1)/2\rfloor$ if $X\neq Y$, and $a$ if $X=Y$. Again, this degree is invariant under the equivalence relation on symbols.

\medskip

In type $^2\!D_n$, so $G=(\CSO_{2n}^-)^0(q)$, if $\chi_\Lambda$ is the unipotent character corresponding to the symbol $\Lambda$ (which has even defect not divisible by $4$), then
\begin{equation} \chi_\Lambda(1)=\frac{\D (q^n+1)\(\prod_{i=1}^{n-1}(q^{2i}-1)\)\(\prod_{1\leq i<j\leq a}(q^{x_i}-q^{x_j})\)\(\prod_{1\leq i<j\leq b}(q^{y_i}-q^{y_j})\)\(\prod_{i,j}(q^{x_i}+q^{y_j})\)}{\D 2^c q^{\binom{a+b-2}{2}+\binom{a+b-4}{2}+\cdots}\(\prod_{i=1}^a\prod_{j=1}^{x_i}(q^{2j}-1)\)\(\prod_{i=1}^b\prod_{j=1}^{y_i}(q^{2j}-1)\)},\label{eq:orthdegrees2}\end{equation}
where $c=(a+b-2)/2$. This degree is also invariant under the equivalence relation on symbols.

\medskip

In all of these groups, two unipotent characters lie in the same $\ell$-block of their respective group if and only if the corresponding symbols have the same $d$-core if $d$ is odd, and $d/2$-cocore if $d$ is even.

\medskip

\section{Brauer Trees and the Minimal $\pi_{\kappa/d}$-Value}
\label{sec:BrauerMinPi}
We first go through the classical groups type by type; in all cases, we associate to the block $B$ either a partition $\lambda$ or a symbol $\Lambda$. We give the description of the Brauer tree, and from this it is easy to describe the parameters of the cyclotomic Hecke algebra, from \cite[Section 2]{brouemalle1993}: the sign $\omega_\chi$ for all characters $\chi$ is $+1$ on one side of the exceptional node and $-1$ on the other, and the power of $q$ is $-\aA(\chi)/e$. We prove that the quantity $\aA(\chi)$ increases towards the exceptional node (as needed for Theorem \ref{thm:piincreasetoexc} using Proposition \ref{prop:increasingperv}) and finally prove that Theorem \ref{thm:minpivalue} is satisfied.

We then consider the exceptional groups, giving a table of those unipotent characters with minimal $\pi_{\kappa/d}$-value.

\subsection{Linear Groups}

Let $n$ be a positive integer, let $q$ be a prime power, let $\ell$ be a prime, and write $d$ for the multiplicative order of $q$ modulo $\ell$. Let $B$ be an $\ell$-block of $G=\GL_{n+d}(q)$ with a cyclic defect group, with $d$-core a partition $\lambda$ of $n$; let $X=\{x_1,\dots,x_a\}$ (with $x_i>x_{i+1}$) be a $\beta$-set corresponding to $\lambda$. We will compute the function $\pi(-)$ for the unipotent characters in $B$. There are $d$ unipotent characters $\chi_\mu$, each with $\lambda$ as $d$-core and $|\mu|-|\lambda|=d$; by choosing $X$ sufficiently large, we have the subset $X'=\{x_{i_1},\dots,x_{i_d}\}$ of $X$ consisting of those $d$ integers such that $x_{i_j}+d\notin X$ (i.e., they represent the possible $d$-hooks that may be added), and order them so that $x_{i_j}>x_{i_{j+1}}$. Notice that if one adds $d$ to $x_{i_j}$, then $j$ is the leg length of the corresponding $d$-hook added to $\lambda$.

Label the unipotent characters $\chi_1,\dots,\chi_d$ in $B$ by $\chi_j$ having partition with $x_{i_j}$ incremented by $d$. By \cite{fongsrin1984}, the Brauer tree of a block $B$, with $d$-core $\lambda$, is a line, with the exceptional vertex at the end, $\chi_d$ adjacent to it, and $\chi_i$ adjacent to $\chi_{i+1}$, as in the following diagram.
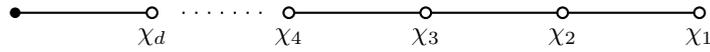
\tikzstyle{every node}=[circle, fill=black!0,
                        inner sep=0pt, minimum width=4pt]
\begin{center}\begin{tikzpicture}[thick,scale=1.8]
\draw \foreach \x in {0,1,2,4}{
(-\x,0) -- (-\x-1,0)};
\draw[loosely dotted] (-3.2,0) -- (-3.8,0);
\draw (0,0) node [draw,label=below:$\chi_{1}$] (l0) {};
\draw (-1,0) node [draw,label=below:$\chi_{2}$] (l1) {};
\draw (-2,0) node [draw,label=below:$\chi_{3}$] (l2) {};
\draw (-3,0) node [draw,label=below:$\chi_{4}$] (l3) {};
\draw (-4,0) node [draw,label=below:$\chi_{d}$] (l5) {};
\draw (-5,0) node [fill=black!100] (ld) {};
\end{tikzpicture}\end{center}

We first want to prove that the $\pi_{\kappa/d}$-function increases towards the exceptional node, using Proposition \ref{prop:increasingperv}.

\begin{prop}\label{prop:increasingGL} We have that $\aA(\chi_{j+1})>\aA(\chi_j)$.
\end{prop}
\begin{pf} Write $x_{i_j}=x_\alpha$ and $x_{i_{j+1}}=x_\beta$, so that $\alpha<\beta$. We have, using (\ref{eq:lineardegrees}),

\[ \frac{\chi_{j+1}(1)}{\chi_j(1)}=\frac{\D\prod_{i=x_\alpha+1}^{x_\alpha+d}(q^i-1)}{\D\prod_{i=x_\beta+1}^{x_\beta+d}(q^i-1)} \cdot \frac{\D\mathop{\prod_{i=1}^a}_{i\neq \beta} (q^{x_\beta+d}-q^{x_i})}{\D\mathop{\prod_{i=1}^a}_{i\neq \beta} (q^{x_\beta}-q^{x_i})}\cdot \frac{\D\mathop{\prod_{i=1}^a}_{i\neq \alpha} (q^{x_\alpha}-q^{x_i})}{\D\mathop{\prod_{i=1}^a}_{i\neq \alpha} (q^{x_\alpha+d}-q^{x_i})}.\]
Clearly, evaluating $\aA(-)$ on the first quotient yields $d(x_\alpha-x_\beta)$, and evaluating it on the second and third terms give $(a-1)d$ and $-(a-1)d$ respectively, so that
\[ \aA(\chi_{j+1})-\aA(\chi_j)=d(x_\alpha-x_\beta)>0,\]
as needed.
\end{pf}

We now consider the unipotent character with minimal $\pi_{\kappa/d}$-function; by Proposition \ref{prop:increasingGL} this must be $\chi_1$. We have, using (\ref{eq:lineardegrees}),

\[ \frac{\chi_1(1)}{\chi_\lambda(1)}=\frac{\D\prod_{i=n+1}^{n+d}(q^i-1)}{\D\prod_{i=x_1+1}^{x_1+d}(q^i-1)} \cdot \frac{\D\prod_{i=2}^a (q^{x_1+d}-q^{x_i})}{\D\prod_{i=2}^a (q^{x_1}-q^{x_i})}.\]
Applying the $\pi_{\kappa/d}$-function to the first quotient yields $2\kappa(n-x_1)$ by Proposition \ref{prop:diffbetweenprodslun}, and to the second quotient yields $2\kappa(a-1)$ by Proposition \ref{prop:diffsinBdslinear}. Hence
\[ \pi_{\kappa/d}(\chi_1)=2\kappa(n-\lambda_1),\]
as $\lambda_1=x_1-a+1$. On the other hand,
\[ \aA(\chi_1)=(n-x_1)d+(a-1)d=(n-\lambda_1)d,\]
so that $\pi_{\kappa/d}(\chi_1)=2\kappa\aA(\chi_1)/d$, as claimed by Theorem \ref{thm:minpivalue}.

\subsection{Unitary Groups}

 Let $n$ be a positive integer, let $q$ be a prime power, let $\ell\mid |G|$ be a prime, and write $d$ and $e$ for the multiplicative orders of $q$ and $-q$ respectively modulo $\ell$; then $e=d$ if $4\mid d$, $e=2d$ if $d$ is odd and $e=d/2$ otherwise. Let $G=\GU_{n+e}(q)$, and let $B$ be an $\ell$-block of $G$ with cyclic defect group.

We use the description of the Brauer trees from \cite{fongsrin1990}. Let $\lambda$ be an $e$-core of size $n$ and let $X$ be a $\beta$-set corresponding to $\lambda$. Let $X'$ denote the subset of $X$ consisting of all $x\in X$ such that $x+e\notin X$, as in the case of $\GL_n(q)$. By replacing $X$ with an equivalent $\beta$-set, we have $|X'|=e$. Divide $X'$ into $Y$ and $Z$, where $Y$ consists of all even elements of $X'$, and $Z$ consists of all odd elements of $X'$, with the ordering on $Y=\{y_1,\dots,y_a\}$ and $Z=\{z_1,\dots,z_b\}$ given by $y_i>y_{i+1}$ and $z_i>z_{i+1}$, as with $X$. Let $\sigma_i$ be the character of $\GU_{n+e}(q)$ obtained by replacing $y_i$ with $y_i+e$, and similarly let $\tau_i$ be the character obtained by replacing $z_i$ with $z_i+e$. The Brauer tree is as follows.

\tikzstyle{every node}=[circle, fill=black!0,
                        inner sep=0pt, minimum width=4pt]
\begin{center}\begin{tikzpicture}[thick,scale=1.8]
\draw \foreach \x in {0,1,3,4,6,7}{
(-\x,0) -- (-\x-1,0)};
\draw[loosely dotted] (-2.2,0) -- (-2.8,0);
\draw[loosely dotted] (-5.2,0) -- (-5.8,0);
\draw (0,0) node [draw,label=below:$\sigma_{1}$] (l0) {};
\draw (-1,0) node [draw,label=below:$\sigma_{2}$] (l1) {};
\draw (-2,0) node [draw,label=below:$\sigma_{3}$] (l2) {};
\draw (-3,0) node [draw,label=below:$\sigma_{a}$] (l5) {};
\draw (-4,0) node [fill=black!100] (ld) {};
\draw (-8,0) node [draw,label=below:$\tau_{1}$] (l0) {};
\draw (-7,0) node [draw,label=below:$\tau_{2}$] (l1) {};
\draw (-6,0) node [draw,label=below:$\tau_{3}$] (l2) {};
\draw (-5,0) node [draw,label=below:$\tau_{b}$] (l5) {};
\end{tikzpicture}\end{center}
Notice that if $e$ is even then the two branches of the tree have the same length.

As in the previous section, we firstly prove that the $\pi_{\kappa/d}$-function increases towards the exceptional node, again using Proposition \ref{prop:increasingperv}.

\begin{prop}\label{prop:increasingGU} We have that $\aA(\sigma_{j+1})>\aA(\sigma_j)$ and $\aA(\tau_{j+1})>\aA(\tau_j)$.
\end{prop}
\begin{pf} Write $y_j=x_\alpha$ and $y_{j+1}=x_\beta$, so that $\alpha<\beta$. The degrees $\sigma_j(1)$ and $\sigma_{j+1}(1)$ are obtained from (\ref{eq:lineardegrees}) by replacing $q$ with $-q$ and $d$ with $e$; this does not affect the $\aA$-function, and so the exact same proof as in Proposition \ref{prop:increasingGL} holds. The case of the $\tau_j$ is identical.
\end{pf}

Since there are now two unipotent characters, $\sigma_1$ and $\tau_1$, on the boundary of the Brauer tree, these are the two possibilities for a unipotent character with minimal $\pi_{\kappa/d}$-function. We may suppose without loss of generality that $x_1=y_1$ is even, and so $\sigma_1$ corresponds to adding an $e$-hook of leg length $1$ to $\lambda$. We can calculate its $\pi_{\kappa/d}$-function exactly as in the previous subsection, to get firstly (via (\ref{eq:lineardegrees}) with $-q$ instead of $q$)
\[ \frac{\sigma_1(1)}{\chi_\lambda(1)}=\frac{\D\prod_{i=n+1}^{n+e}((-q)^i-1)}{\D\prod_{i=x_1+1}^{x_1+e}((-q)^i-1)} \cdot \frac{\D\prod_{i=2}^a ((-q)^{x_1+e}-(-q)^{x_i})}{\D\prod_{i=2}^a ((-q)^{x_1}-(-q)^{x_i})}.\]
Applying the $\pi_{\kappa/d}$-function to the first quotient yields $2\kappa(n-x_1)e/d$ by Proposition \ref{prop:diffbetweenprodslun}, and to the second quotient yields $2\kappa(a-1)e/d$ by Proposition \ref{prop:diffsinBdsunitary}. Hence
\[ \pi_{\kappa/d}(\sigma_1)=2\kappa(n-\lambda_1)\frac ed,\]
as $\lambda_1=x_1-a+1$. On the other hand,
\[ \aA(\sigma_1)=(n-x_1)e+(a-1)e=(n-\lambda_1)e,\]
so that again $\pi_{\kappa/d}(\sigma_1)=2\kappa\aA(\sigma_1)/d$, as claimed by Theorem \ref{thm:minpivalue}.

It remains to check that the other unipotent character, $\tau_1$, has a larger $\pi_{\kappa/d}$-value than $\sigma_1$. We get that $z_1=x_\alpha$ for some $\alpha>1$, and in this case
\[ \frac{\tau_1(1)}{\chi_\lambda(1)}=\frac{\D\prod_{i=n+1}^{n+e}((-q)^i-1)}{\D\prod_{i=x_\alpha+1}^{x_\alpha+e}((-q)^i-1)} \cdot \frac{\D\mathop{\prod_{i=1}^a}_{i\neq \alpha} ((-q)^{x_\alpha+e}-(-q)^{x_i})}{\D\mathop{\prod_{i=1}^a}_{i\neq \alpha} ((-q)^{x_\alpha}-(-q)^{x_i})}.\]
As before, applying the $\pi_{\kappa/d}$-function to the first term yields $2\kappa(n-x_\alpha)e/d$, and applying it to the second quotient yields at least $2\kappa(a-\alpha)e/d$ (for each of the $x_i$ with $i>\alpha$), so we have
\[ \pi_{\kappa/d}(\tau_1)\geq 2\kappa(n-x_\alpha+a-\alpha)\frac ed=2\kappa(n-\lambda_\alpha)\frac ed.\]
The only way that $\pi_{\kappa/d}(\tau_1)$ can equal $\pi_{\kappa/d}(\sigma_1)$ is if $\lambda_1=\lambda_\alpha$: since $x_\alpha=z_1$ is the largest odd $\beta$-number, we must have $\alpha=2$ and $x_\alpha=x_1-1$: in this case, if $d=1$ then $B$ is the principal block and the result is clear, and if $d>1$ we have
\[ \pi_{\kappa/d}(\tau_1)=2\kappa(n-x_1+a-1)\frac ed+\left(\pi_{\kappa/d}\left(\frac{((-q)^{x_2+e}-q^{x_1})}{q^{x_1}+q^{x_1-1}}\right)\right).\]
The first term is simply $\pi_{\kappa/d}(\sigma_1)$, and the last term is $\pi_{\kappa/d}\left(((-q)^e-q)/(q+1)\right)$, which is positive by Proposition \ref{prop:diffsinBdsunitary}. Hence $\pi_{\kappa/d}(\tau_1)>\pi_{\kappa/d}(\sigma_1)$, and this completes the proof of Theorem \ref{thm:minpivalue} for unitary groups.

\subsection{Symplectic and Odd-Dimensional Orthogonal Groups}

If $d$ is even, we write $d'=d/2$. Let $G_n$ be one of the groups $\SO_{2n+1}(q)$ and $\CSp_{2n}(q)$. Let $\Lambda=\{X,Y\}$ be a symbol of rank $n$ and odd defect $\delta$, with $X=\{x_1,\dots,x_a\}$ and $Y=\{y_1,\dots,y_b\}$, ordered so that $x_i>x_{i+1}$ and $y_i>y_{i+1}$. Assume that $\Lambda$ is a $d$-core if $d$ is odd, and a $d'$-cocore if $d$ is even.

We start with the case $d$ is odd. Recall that we view $\Lambda$ as a pair of $\beta$-sets: let $X'$ denote the beads of $X$ on the end of their runners of the $d$-abacus, and let $Y'$ denote the beads of $Y$ on the end of their runners of the $d$-abacus. By choosing $\Lambda$ suitably, $|X'|=|Y'|=d$. Write $X'=\{x_1',\dots,x_d'\}$ and $Y'=\{y_1',\dots,y_d'\}$, with $x_i'>x_{i+1}'$ and $y_i'>y_{i+1}'$.

Let $\sigma_1,\dots,\sigma_d$ be the unipotent characters of $G=G_{n+d}$ corresponding to adding $d$ to the elements of $X'$, with $\sigma_i$ coming from $x_i'$; similarly, let $\tau_1,\dots,\tau_d$ be the unipotent characters of $G$ corresponding to adding $d$ to the elements of $Y'$, with $\tau_i$ coming from $y_i'$. In this case the Brauer tree is as follows.
\tikzstyle{every node}=[circle, fill=black!0,
                        inner sep=0pt, minimum width=4pt]
\begin{center}\begin{tikzpicture}[thick,scale=1.8]
\draw \foreach \x in {0,1,3,4,6,7}{
(-\x,0) -- (-\x-1,0)};
\draw[loosely dotted] (-2.2,0) -- (-2.8,0);
\draw[loosely dotted] (-5.2,0) -- (-5.8,0);
\draw (0,0) node [draw,label=below:$\sigma_{1}$] (l0) {};
\draw (-1,0) node [draw,label=below:$\sigma_{2}$] (l1) {};
\draw (-2,0) node [draw,label=below:$\sigma_{3}$] (l2) {};
\draw (-3,0) node [draw,label=below:$\sigma_{d}$] (l5) {};
\draw (-4,0) node [fill=black!100] (ld) {};
\draw (-8,0) node [draw,label=below:$\tau_{1}$] (l0) {};
\draw (-7,0) node [draw,label=below:$\tau_{2}$] (l1) {};
\draw (-6,0) node [draw,label=below:$\tau_{3}$] (l2) {};
\draw (-5,0) node [draw,label=below:$\tau_{d}$] (l5) {};
\end{tikzpicture}\end{center}

We now need to prove, as in the last two sections, that the $\aA$-function increases towards the exceptional node.

\begin{prop}\label{prop:increasingSp1} We have that $\aA(\sigma_{j+1})>\aA(\sigma_j)$ and $\aA(\tau_{j+1})>\aA(\tau_j)$.
\end{prop}

The proof is almost identical to that of Proposition \ref{prop:increasingGL}, and is safely left to the reader.

As with the previous cases, the minimal $\pi_{\kappa/d}$-value must come from either $\sigma_1$ or $\tau_1$. Without loss of generality, $x_1\geq y_1$. This time we get, using (\ref{eq:sympdegrees})
\[ \frac{\sigma_1(1)}{\chi_\Lambda(1)}=\frac{\D\prod_{i=n+1}^{n+d}(q^{2i}-1)}{\D\prod_{i=x_1+1}^{x_1+d}(q^{2i}-1)}\,\frac{\D\prod_{i=2}^a (q^{x_1+d}-q^{x_i})}{\D\prod_{i=2}^a (q^{x_1}-q^{x_i})}\,\frac{\D\prod_{i=1}^b (q^{x_1+d}+q^{y_i})}{\D\prod_{i=1}^b (q^{x_1}+q^{y_i})}.\]
Applying the $\pi_{\kappa/d}$-function to the first quotient yields $4\kappa(n-x_1)$ by Proposition \ref{prop:diffbetweenprodsorthog}, to the second quotient yields $2\kappa(a-1)$ as in the case of $\GL_n(q)$, and to the third quotient yields $2\kappa b$ by Proposition \ref{prop:diffsinBdslinear}. Hence
\[\pi_{\kappa/d}(\sigma_1)=2\kappa(2n-2x_1+a+b-1).\]
We now wish to evaluate $\aA(\sigma_1)$: we get
\[ \aA(\sigma_1)=2d(n-x_1)+d(a-1)+db=(2n-2x_1+a+b-1)d,\]
so that $\pi_{\kappa/d}(\sigma_1)=2\kappa\aA(\sigma_1)/d$, again in line with Theorem \ref{thm:minpivalue}. As with the previous case of the unitary groups, we need to evaluate $\pi_{\kappa/d}(\tau_1)$ as well: in this case,
\[ \frac{\tau_1(1)}{\chi_\Lambda(1)}=\frac{\D\prod_{i=n+1}^{n+d}(q^{2i}-1)}{\D\prod_{i=y_1+1}^{y_1+d}(q^{2i}-1)}\,\frac{\D\prod_{i=1}^a (q^{y_1+d}+q^{x_i})}{\D\prod_{i=1}^a (q^{y_1}+q^{x_i})}\,\frac{\D\prod_{i=2}^b (q^{y_1+d}-q^{y_i})}{\D\prod_{i=2}^b (q^{y_1}-q^{y_i})}.\]
We apply the $\pi_{\kappa/d}$-function to get $4\kappa(n-y_1)$ and $2\kappa(b-1)$ for the first and third quotients: if $y_1\geq x_\alpha$ but $y_1<x_{\alpha-1}$, then the second quotient yields at least $2\kappa(a-\alpha+1)$; as $x_1-x_\alpha\geq \alpha-1$, we get
\[\pi_{\kappa/d}(\tau_1)\geq 2\kappa(2n-2y_1+a-\alpha+1+b-1)\geq 2\kappa(2n-2x_1+\alpha+a+b-2),\]
which can only equal $\pi_{\kappa/d}(\sigma_1)$ if $\alpha=1$, i.e., $x_1=y_1$. In this case $\pi_{\kappa/d}(\sigma_1)$ is actually equal to $\pi_{\kappa/d}(\tau_1)$, and indeed $\aA(\sigma_1)=\aA(\tau_1)$, so Theorem \ref{thm:minpivalue} is verified when $d$ is odd.

\medskip

If $d$ is even, the description of the Brauer tree is very similar to the case where $d$ is odd: let $\Lambda=\{X,Y\}$ be a $d'$-cocore of odd defect $\delta$ and rank $n$, and let $X'$ and $Y'$ denote the subsets of $X$ and $Y$ given by
\[ X'=\{x\in X:x+d'\notin Y\},\qquad Y'=\{y\in Y:y+d'\notin X\}.\]
Assume that $|X|>|Y|$, so that $|X|-|Y|=\delta$. By \cite[(3E)]{fongsrin1990}, we have that $|X'|=d'+\delta$ and $|Y'|=d'-\delta$. Write $X'=\{x_1',\dots,x_{d'+\delta}'\}$, ordered so that $x_i'>x_{i+1}'$, and similarly for $Y'$. If $\sigma_i$ is the unipotent character corresponding to the symbol obtained by adding $d'$-cohook to $x_i'$, and similarly for $\tau_i$ and $y_i'$, then the Brauer tree is as follows.
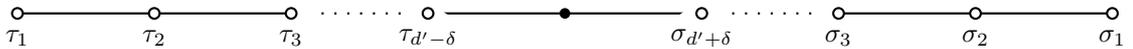
\tikzstyle{every node}=[circle, fill=black!0,
                        inner sep=0pt, minimum width=4pt]
\begin{center}\begin{tikzpicture}[thick,scale=1.8]
\draw \foreach \x in {0,1,3,4,6,7}{
(-\x,0) -- (-\x-1,0)};
\draw[loosely dotted] (-2.2,0) -- (-2.8,0);
\draw[loosely dotted] (-5.2,0) -- (-5.8,0);
\draw (0,-0.18) node{$\sigma_1$};
\draw (-1,-0.18) node{$\sigma_2$};
\draw (-2,-0.18) node{$\sigma_3$};
\draw (-3,-0.18) node{$\sigma_{d'+\delta}$};
\draw (-8,-0.18) node{$\tau_1$};
\draw (-7,-0.18) node{$\tau_2$};
\draw (-6,-0.18) node{$\tau_3$};
\draw (-5,-0.18) node{$\tau_{d'-\delta}$};

\draw (0,0) node [draw] (l0) {};
\draw (-1,0) node [draw] (l0) {};
\draw (-2,0) node [draw] (l0) {};
\draw (-3,0) node [draw] (l0) {};
\draw (-4,0) node [fill=black!100] (ld) {};
\draw (-5,0) node [draw] (l0) {};
\draw (-6,0) node [draw] (l0) {};
\draw (-7,0) node [draw] (l0) {};
\draw (-8,0) node [draw] (l0) {};
\end{tikzpicture}\end{center}
The proof that the $\aA$-function increases towards the exceptional node is essentially identical to that for odd $d$, and is again omitted.

Again, the minimal $\pi_{\kappa/d}$-value must come from either $\sigma_1$ or $\tau_1$. We have
\[ \frac{\sigma_1(1)}{\chi_\Lambda(1)}=\frac{\D\prod_{i=n+1}^{n+d'}(q^{2i}-1)}{\D\prod_{i=x_1+1}^{x_1+d'}(q^{2i}-1)}\,\frac{\D\prod_{i=2}^a (q^{x_1+d'}+q^{x_i})}{\D\prod_{i=2}^a (q^{x_1}-q^{x_i})}\,\frac{\D\prod_{i=1}^b (q^{x_1+d'}-q^{y_i})}{\D\prod_{i=1}^b (q^{x_1}+q^{y_i})}.\]
If $x_1\geq y_1$ then we can use Propositions \ref{prop:diffbetweenprodsorthog} and \ref{prop:diffsinBdsorthog} to get
\[ \pi_{\kappa/d}(\sigma_1)=2\kappa(n-x_1)+\kappa(a-1)+\kappa b=\kappa(2n-2x_1+a+b-1),\qquad \aA(\sigma_1)=(n-x_1+(a+b-1)/2)d,\]
with a similar statement holding for $\tau_1$ in the case where $y_1\geq x_1$. In particular, this character satisfies Theorem \ref{thm:minpivalue}.

It remains to check that $\pi_{\kappa/d}(\sigma_1)>\pi_{\kappa/d}(\tau_1)$ if and only if $x_1<y_1$. Hence we assume that $x_1<y_1$ and compute $\pi_{\kappa/d}(\sigma_1)$ using the above equation. Let $\alpha$ be such that $x_1\geq y_\alpha$ but $x_1<y_{\alpha-1}$. Since $y_1-y_\alpha\geq \alpha-1$ (as in the odd case), we get
\[ \pi_{\kappa/d}(\sigma_1)\geq 2\kappa(n-x_1)+\kappa(a-1)+\kappa(b-\alpha+1)\geq \kappa(2n-2y_1+\alpha+a+b-2)=\pi_{\kappa/d}(\tau_1)-1+\alpha,\]
and so we can only equality between $\pi_{\kappa/d}(\sigma_1)$ and $\pi_{\kappa/d}(\tau_1)$ when $\alpha=1$, i.e., $x_1=y_1$, as seen before. This completes the proof of Theorem \ref{thm:minpivalue} for types $B_n$ and $C_n$.

\subsection{Even-Dimensional Orthogonal Groups}

If $d$ is even, we write $d'=d/2$. Let $G$ be one of the groups $(\CSO_{2n}^\pm)^0(q)$ (see Section \ref{sec:gensetup}). Let $\Lambda=\{X,Y\}$ be a symbol of rank $n$ and even defect $\delta$, with $X=\{x_1,\dots,x_a\}$ and $Y=\{y_1,\dots,y_b\}$, ordered so that $x_i>x_{i+1}$ and $y_i>y_{i+1}$. Assume that $\Lambda$ is a $d$-core if $d$ is odd, and a $d'$-cocore if $d$ is even.

Notice that a degenerate symbol cannot have weight $1$, so all unipotent characters in blocks of weight $1$ are labelled by non-degenerate symbols. Constructing the sets $X'$ and $Y'$ as in the previous section, if $\Lambda$ is non-degenerate then the Brauer tree of the block $B$ is exactly the same as that of the previous section, whereas if $\Lambda$ is degenerate then only one branch of this Brauer tree exists, and there are half as many unipotent characters in $B$ as expected.

Adding a $d$-hook does not alter the defect of a symbol, but adding a $d'$-cohook to a symbol adds or subtracts $2$ from the defect: hence when we add a $d$-cohook we move from a symbol labelling a unipotent character of ${}^\ep\!D_n$ to one labelling ${}^{-\ep}\!D_{n+d'}$. This will of course be relevant when comparing the character degree $\chi_\Lambda$ with that obtained by adding a $d'$-cohook. Recall that $\sigma_1$ is the character obtained by adding a $d$-hook or $d'$-cohook to $x_1$, and similarly for $\tau_1$ and $y_1$.

Firstly, since the equations for character degrees are so similar between types $B_n/C_n$, $D_n$ and ${}^2\!D_n$, the proof that the $\aA$-function increases towards the exceptional node is the same, so omitted; it remains to discuss the minimal $\pi_{\kappa/d}$-function. In the case where $d$ is odd, we get, using (\ref{eq:orthdegrees1}) and (\ref{eq:orthdegrees2})
\[ \frac{\sigma_1(1)}{\chi_\Lambda(1)}=\frac{q^{n+d}\pm1}{q^n\pm 1}\frac{\D\prod_{i=n}^{n+d-1}(q^{2i}-1)}{\D\prod_{i=x_1+1}^{x_1+d}(q^{2i}-1)}\,\frac{\D\prod_{i=2}^a (q^{x_1+d}-q^{x_i})}{\D\prod_{i=2}^a (q^{x_1}-q^{x_i})}\,\frac{\D\prod_{i=1}^b (q^{x_1+d}+q^{y_i})}{\D\prod_{i=1}^b (q^{x_1}+q^{y_i})}.\]
We may assume that $x_1\geq y_1$, in which case a very similar analysis to the symplectic case yields
\[ \pi_{\kappa/d}(\sigma_1)=2\kappa(2n-2x_1+a+b-2),\qquad \aA(\sigma_1)=(2n-2x_1+a+b-2)d,\]
and the same argument as in the last section proves that $\pi_{\kappa/d}(\tau_1)=\pi_{\kappa/d}(\sigma_1)$ if and only if $y_1=x_1$, in which case $\aA(\sigma_1)=\aA(\tau_1)$ and Theorem \ref{thm:minpivalue} follows for $d$ odd.

For $d$ even, we get
\[ \frac{\sigma_1(1)}{\chi_\Lambda(1)}=\frac{q^{n+d'}\pm 1}{q^n\mp1}\frac{\D\prod_{i=n}^{n+d'-1}(q^{2i}-1)}{\D\prod_{i=x_1+1}^{x_1+d'}(q^{2i}-1)}\,\frac{\D\prod_{i=2}^a (q^{x_1+d'}+q^{x_i})}{\D\prod_{i=2}^a (q^{x_1}-q^{x_i})}\,\frac{\D\prod_{i=1}^b (q^{x_1+d'}-q^{y_i})}{\D\prod_{i=1}^b (q^{x_1}+q^{y_i})}\]
and a similar expression for $\tau_1$. If $x_1\geq y_1$ then we get
\[ \pi_{\kappa/d}=\kappa(2n-2x_1+a+b-2),\qquad \aA(\sigma_1)=(2n-2x_1+a+b-2)d',\]
with a similar statement for $\tau_1$ if $y_1\geq x_1$. If $x_1<y_1$ then the same argument as the previous section applies, and so Theorem \ref{thm:minpivalue} is true for the even-dimensional orthogonal groups.

\bigskip

In each section, it was proved that the $\aA$-function increases towards the exceptional node, yielding Theorem \ref{thm:piincreasetoexc} for the classical groups; the rest of the calculations conclude the proof of Theorem \ref{thm:minpivalue} for the classical groups.

\subsection{Exceptional Groups}

In Table \ref{t:minpiexceptionals} we give a complete list of the characters with the minimal $\pi_{\kappa/d}$-value for each unipotent block of each exceptional group; the character degrees are available in \cite{carterfinite} or on GAP. Of course, if the block is the principal block then the character with the smallest $\pi_{\kappa/d}$-value is the trivial character and Theorem \ref{thm:minpivalue} is obvious, so we need only consider non-principal blocks.

In \cite[Table 6.1]{bmm1993} a complete list of the non-principal unipotent blocks with cyclic defect groups for exceptional groups is given along with the $d$-cuspidal pair involved, and the table afterwards gives the unipotent characters in that block. If there is more than one unipotent character then these are listed in brackets.

The easiest way to prove that each of the characters in Table \ref{t:minpiexceptionals} does satisfy the equation $\pi_{\kappa/d}(\chi)=2\kappa\aA(\chi)/d$ is to find a set of polynomials $f$ such that $\phi_{\kappa/d}(f)=\kappa$ for every $d$ and $\kappa<d$, and then note that each of the relative degrees is a product of such polynomials: for example, when $d=3$, $\Phi_1^2\Phi_5$, $\Phi_7$, $\Phi_9$ and $\Phi_{14}$ have this property, and the product of these is the relative degree of $D_4,1$ for $E_7$, so that character satisfies Theorem \ref{thm:minpivalue}. It is already interesting that there is a unipotent character in each block that has the property that its degree $f$ satisfies $\phi_{\kappa/d}(f)=\kappa$, as it is obvious that not every product of cyclotomic polynomials has this property. Of course, the condition that $a(\chi)=a(\bo\lambda)$ for such a character is easy to check by inspection, and so $\pi_{\kappa/d}(\chi)=2\kappa\aA(\chi)/d$ for all $\kappa$, including those greater than $d$.

Either a hand calculation or the use of a computer (which might be wise for $E_8$) verifies in each case that the character in the table has the minimal $\pi_{\kappa/d}$-value in the block and that it satisfies Theorem \ref{thm:minpivalue}. With these calculations, we conclude the proof of Theorem \ref{thm:minpivalue} for all unipotent blocks with cyclic defect groups for all groups of Lie type.

\begin{table}\begin{center}\begin{tabular}{|cccc|cccc|}
\hline
$G$ & $d$ & $\chi$ & $\aA(\chi)/d$ & $G$ & $d$ & $\chi$ & $\aA(\chi)/d$
\\\hline$F_4$ & $2$ & $\phi_{2,4}'$, $\phi_{2,4}''$ & $2$, $2$ & & $10$ & $\phi_{7,1}$, $\phi_{35,4}$ & $3/2$, $3$
\\ ${}^2\!F_4$ & $2$ & ${}^2\!B_2[\psi^3];1$, ${}^2\!B_2[\psi^5];1$ & $4$, $4$ & $E_7$ & $12$ & $\phi_{21,3}$ & $2$
\\ $E_6$ & $2$ & $\phi_{64,4}$ & $6$ & $E_8$ & $1$ & $E_7[\I];1$, $E_7[-\I];1$ & $42$, $42$
\\ & $3$ & $D_4;1$ & $4$ & & $2$ & $\phi_{4096,11}$, $\phi_{4096,12}$ & $21$, $21$
\\ & $4$ & $\phi_{20,2}$ & $3$ & & $3$ & $\phi_{567,6}$, $\phi_{1008,9}$, ($\phi_{1296,13},\phi_{2268,10}$)& $16$, $18$, $20$
\\ & $5$ & $\phi_{6,1}$ & $2$ & & $5$ &$\phi_{35,2}$, ($\phi_{210,4},\phi_{160,7}$) & $8$, $12$
\\ ${}^2\!E_6$ & $1$ & ${}^2\!A_5;1$ & $12$& & $6$ & $\phi_{567,6}$, $\phi_{1008,9}$, ($\phi_{972,12}$, $D_4;\phi_{9,2}$) & $8$, $9$, $10$
\\ & $4$ & $\phi_{4,1}$ & $3$ & & $7$ & $\phi_{8,1}$ & $4$
\\ & $6$ & $\phi_{8,3}'$ & $2$  & & $8$ & $\phi_{8,1}$, ($\phi_{1344,8}$, $D_4;\phi_{4,1}$) &  $3$, $9$
\\ & $10$ & $\phi_{2,4}'$ & $1$& && $\phi_{84,4}$, $D_4;\phi_{1,0}$,  $\phi_{112,13}$, $\phi_{28,8}$ & $6$, $6$, $6$, $6$
\\ $E_7$ & $1$ & $E_6[\theta];1$, $E_6[\theta^2];1$ & $18$, $18$& & $9$ & $\phi_{8,1}$, ($\phi_{28,8},\phi_{112,3}$) & $3$, $6$
\\ & $2$ & $E_6[\theta];1$, $E_6[\theta^2];1$ & $9$, $9$  & & $10$ & $\phi_{35,2}$, ($\phi_{50,8},D_4;\phi_{2,4}'$) & $4$, $6$
\\ & $3$ & $\phi_{27,2}$, $D_4;1$, $\phi_{189,7}$ & $6$, $8$, $10$ & & $12$ &$\phi_{8,1}$, $\phi_{84,4}$, $\phi_{28,8}$, ($\phi_{448,9},E_6[\theta^i];1$) & $2$, $4$, $4$, $6$
\\ & $5$ & $\phi_{7,1}$, $(\phi_{21,6},\phi_{56,3})$ & $3$, $6$ & & $14$ & $\phi_{8,1}$ & $2$
\\ & $6$ & $\phi_{27,2}$, $\phi_{56,3}$, $\phi_{189,7}$ & $3$, $4$, $5$ && $18$ & $\phi_{8,1}$, ($\phi_{84,4},D_4;\phi_{1,0}$) & $3/2$, $3$
\\ & $8$ & $\phi_{7,1}$, $\phi_{27,2}$, $\phi_{21,3}$ & $2$, $3$, $3$&&&&
\\ \hline
\end{tabular}\end{center}
\caption{Unipotent characters with minimal $\pi_{\kappa/d}$-value for non-principal blocks with cyclic defect groups for exceptional groups}
\label{t:minpiexceptionals}
\end{table}

\section{Perturbing the Cyclotomic Hecke Algebra}
\label{sec:perturbHecke}

The aim of this section is to produce a method by which one (specialized) cyclotomic Hecke algebra can be turned into another. This proceeds by \emph{perturbing} the parameters, for example replacing $q^a$ by $q^{a+d}$. The eventual aim is to perturb all cyclotomic Hecke algebras into ones with specialized parameters very similar to that of the Coxeter case, and then prove directly that there is a perverse equivalence in this case.

By replacing one parameter by another, we will alter both the $\pi_{\kappa/d}$-function and the ordering on the parameters: the aim is to show that the $\pi_{\kappa/d}$-function reduces by $2$ for a certain set of parameters, and that they are cycled to reorder them in accordance with Theorem \ref{thm:perversebijchange}.

Before we define perturbations of cyclotomic Hecke algebras, we need to reduce to the case $\kappa=1$. This will be performed in the next lemma and proposition.

\begin{lem}\label{lem:equalroots} Let $f(x)\in\C[x]$ be a polynomial. Let $\sigma\in (0,2)$, and let $\kappa$ be a positive integer. We have that $|\Arg_{\kappa\sigma}(f(x))|=|\Arg_\sigma(f(x^\kappa))|$.
\end{lem}
\begin{pf} Let $\xi$ be any non-zero complex number with argument $\lambda\in(0,2\pi]$, and write $\xi_1,\dots,\xi_\kappa$ for the $\kappa$th roots of $\xi$: the elements of $\Arg_{\kappa\sigma}(\xi)$ are $\lambda+2\pi j$ for $0\leq j\leq \alpha$ for some $\alpha$, and the elements of the union of the sets $\Arg_{\sigma}(\xi_i)$ for $1\leq i\leq \kappa$ are $\lambda/\kappa+2\pi/\kappa j$ for $0\leq j\leq \alpha$. Therefore
\[ |\Arg_{\kappa\sigma}(\xi)|=\sum_{i=1}^\kappa |\Arg_\sigma(\xi_i)|.\]
Since the roots of $f(x^k)$ are the $\xi_i$ for $\xi$ running over all roots of $f(x)$, we see that $|\Arg_{\kappa\sigma}(f(x))|=|\Arg_\sigma(f(x^\kappa))|$, as claimed.
\end{pf}

Equipped this this easy lemma, we prove the next proposition.

\begin{prop}\label{prop:kgoesto1} Let $\mc H_1=\mc H(Z_e,\b u)$ be a cyclotomic Hecke algebra with specialization of parameters $u_i\mapsto \omega_i q^{v_i}$, and let $\mc H_2$ be the same algebra with specialization of parameters $u_i\mapsto \omega_i q^{\kappa v_i}$. Let $\zeta_0=\e^{2\pi\I/d}$. Let $\rho_i$ be the relative degree associated to $u_i$ in $\mc H_1$, and $\sigma_i$ the same for $\mc H_2$.
\begin{enumerate}
\item We have $\pi_{\kappa/d}(\rho_i)=\pi_{1/d}(\sigma_i)$.
\item For $\mc H_1$ and $q\mapsto \zeta$, and for $\mc H_2$ and $q\mapsto \zeta_0$, the induced bijections between the parameters and the $e$th roots of unity are identical.
\end{enumerate}
\end{prop}
\begin{pf} By Lemma \ref{lem:equalroots}, $|\Arg_{2\kappa/d}(\rho_i)|=|\Arg_{2/d}(\sigma_i)|$, and clearly the multiplicities of $1$ as a root of $\rho_i$ and $\sigma_i$ are identical. Finally, we have $a(\sigma_i)=\kappa\cdot a(\rho_i)$ and $A(\sigma_i)=\kappa\cdot A(\rho_i)$, so we therefore have $\pi_{\kappa/d}(\rho_i)=\pi_{1/d}(\sigma_i)$, as needed.

The second part is clear.
\end{pf}

As an aside, the appropriate analogue of this proposition holds for all cyclotomic Hecke algebras: the notion of perturbing cyclotomic Hecke algebras can be extended to the non-cyclic case, and will be dealt with in a later paper in this series.

Hence if we can show, for an arbitrary cyclotomic Hecke algebra and $\kappa=1$, that the corresponding bijection is that of Theorem \ref{thm:allperversecyclic}, then it is true for all $\kappa$. From now on in the proof for classical groups we will assume that $\kappa=1$.

We will take an arbitrary cyclotomic Hecke algebra of type $(s,t)$ and ambiance $d$ and perform one of three operations on the parameters:
\begin{enumerate}
\item replace $q^{a_s}$ by $q^{a_s+d}$ and rearrange parameters as required (called a \emph{$+$-type perturbation});
\item replace $-q^{b_t}$ by $-q^{b_t+d}$ and rearrange parameters as required (called a \emph{$-$-type perturbation});
\item replace $q^{a_s}$ by $-q^{a_s+d/2}$ and $-q^{b_t}$ by $q^{a_s+d/2}$ and rearrange parameters as required (called a \emph{$\pm$-type perturbation}).
\end{enumerate}
It is clear that these replacements preserve the property of being a cyclotomic Hecke algebra of type $(s,t)$ and ambiance $d$.

A $+$-type perturbation is only allowed if $a_s+d/2<b_t$, a $-$-type perturbation is only allowed if $b_t+d/2<a_s$, and a $\pm$-type perturbation is only allowed if neither a $+$-type nor a $-$-type perturbation is allowed; notice that this means that there is exactly one way to perturb a given cyclotomic Hecke algebra. Since we have defined a canonical ordering on the parameters of a cyclotomic Hecke algebra of type $(s,t)$, we see that a perturbation permutes some of the parameters. The set of \emph{permuted parameters} is a subset of the parameters of the cyclotomic Hecke algebra, and can also be thought of as a subset of the $\chi_i$ and $\psi_i$: because of the canonical ordering on the parameters, we may compare the sets of permuted parameters for different cyclotomic Hecke algebras of the same type, even though their parameters might be different.

To evaluate the difference in the $\pi_{\kappa/d}$-function between two cyclotomic Hecke algebras with different parameters, we will have normalize the relative degrees in some way; choose a parameter that is not permuted for this, noting that one always exists unless the perturbation permutes all parameters, which we will prove in Theorem \ref{thm:perturbationchain} results in an isomorphic algebra.

If the $a_i$ and $b_i$ are all integers then we know how to evaluate the $\pi_{\kappa/d}$-function on polynomials of the form $q^{a_i}-q^{a_j}$ and $q^{a_i}+q^{b_j}$, so firstly we also need to reduce to the case where all of the $a_i$ and $b_i$ are integers, and secondly if $d$ is odd then a $\pm$-perturbation will reintroduce fractional parameters, so we will need to know that we do not need $\pm$-type perturbations if $d$ is odd, i.e., the exceptional node has valency $1$ in such cases.

From Section \ref{sec:BrauerMinPi} we know the structure of the Brauer tree for classical groups, and we know that one of the following holds:
\begin{enumerate}
\item $d=e$ is arbitrary, $t=0$ ($\GL_n$, all $d$; $\Sp_{2n}$ and $(\CSO_{2n}^\pm)^0$, $d$ odd, $\Lambda$ degenerate);
\item $d$ is odd, $e=2d$, $s=t=e$ ($\GU_n$, $d$ odd; $\Sp_{2n}$ and $(\CSO_{2n}^\pm)^0$, $d$ odd, $\Lambda$ non-degenerate);
\item $d$ is even, $e=d/2$ ($\GU_n$, $d/2$ odd; $\Sp_{2n}$ and $(\CSO_{2n}^\pm)^0$, $d$ even, $\Lambda$ degenerate);
\item $d$ is even, $e=d$ ($\GU_n$, $4\mid d$, $\Sp_{2n}$ and $(\CSO_{2n}^\pm)^0$, $d$ even, $\Lambda$ non-degenerate).
\end{enumerate}

Since $\pi_{\kappa/d}(\chi)$ is an integer, and differs from $(a(\chi)+A(\chi))/d$ by a semi-integer, unless $e=2d$ we must have that $(a(\chi)+A(\chi))/e$ is a semi-integer. If $d$ is odd, $(a(\chi)+A(\chi))/d$ is a semi-integer  so must be an integer, and again $(a(\chi)+A(\chi))/e$ is a semi-integer. (Notice that in the case where $d$ is odd and $d=e$ this proves that $(a(\chi)+A(\chi))/d$ is an integer.)

Lemma \ref{lem:equalroots} allows us to move from semi-integers to integers, and from $d$ odd to $d$ even, easily.

\begin{prop} Let $\mc H_1=\mc H(Z_e,\b u)$ be a cyclotomic Hecke algebra with specialization of parameters $u_i\mapsto \omega_i q^{v_i}$, and let $\mc H_2$ be the same algebra with specialization of parameters $u_i\mapsto \omega_i q^{2v_i}$. Let $\zeta_1=\e^{2\pi\I/d}$ and $\zeta_2=\e^{\pi\I/d}$. Let $\rho_i$ be the relative degree associated to $u_i$ in $\mc H_1$, and $\sigma_i$ the same for $\mc H_2$.
\begin{enumerate}
\item We have $\pi_{1/d}(\rho_i)=\pi_{1/2d}(\sigma_i)$.
\item For $\mc H_1$ and $q\mapsto \zeta_1$, and for $\mc H_2$ and $q\mapsto \zeta_2$, the induced bijections between the parameters and the $e$th roots of unity are identical.
\end{enumerate}
\end{prop}

The proof is the same as that of Proposition \ref{prop:kgoesto1}, and omitted. Using this proposition in the same way as Proposition \ref{prop:kgoesto1}, we may assume that all of the $a_i$ and $b_i$ are integers, and that if we have to make a $\pm$-type perturbation then $d$ is even. For the rest of this section we will do this without further comment.

\begin{prop}\label{prop:plustypeperturbation} Suppose that $\kappa=1$. Let $\mc H'$ be the cyclotomic Hecke algebra obtained from $\mc H$ by applying a $+$-type perturbation that does not permute all parameters. Let $\alpha$ be the number of $i<s$ such that $a_s<a_i<a_s+d$, so that $a_{s-\alpha-1}>a_s+d>a_{s-\alpha}$. Write $\chi_i$ and $\psi_i$ for the characters of $\mc H$, and $\chi_i'$ and $\psi_i'$ for the characters of $\mc H'$. 
\begin{enumerate}
\item We have that $\pi_{\kappa/d}(\psi_i')=\pi_{\kappa/d}(\psi_i)$ for all $i$, $\pi_{\kappa/d}(\chi_i')=\pi_{\kappa/d}(\chi_i)$ for $1\leq i\leq s-\alpha-1$, and $\pi_{\kappa/d}(\chi_i')=\pi_{\kappa/d}(\chi_i)-2$ otherwise (when the $\pi_{\kappa/d}$-function is suitably normalized by an unpermuted parameter).
\item Using the ordering on the permuted parameters inherited from the ordering on all parameters, they are permuted so that the $i$th permuted parameter of $\mc H'$ is the $(i-1)$th permuted parameter of $\mc H$.
\end{enumerate}
\end{prop}
\begin{pf} As usual, $\chi_i$ and $\psi_i$ are the characters of $\mc H$, and write $\chi_i'$ and $\psi_i'$ for the characters of $\mc H'$, ordered in the standard way. We choose $\psi_1$ to normalize our $\pi_{\kappa/d}$-functions, noting that if there are no negative parameters then our statement about $\psi_i$ is vacuous anyway.

We have
\[ \left.\frac{\psi_i'(1)}{\psi_1'(1)}\right/\frac{\psi_i(1)}{\psi_1(1)}=\frac{(q^{b_1}+q^{a_s+d})}{(q^{b_1}+q^{a_s})}\cdot \frac{(q^{b_i}+q^{a_s})}{(q^{b_i}+q^{a_s+d})}.\]
Applying the $\pi_{\kappa/d}$-function to this expression yields $0$ as $a_s+d/2<b_t$, by Proposition \ref{prop:diffsinBdslinear}, and so $\pi_{\kappa/d}(\psi_i')=\pi_{\kappa/d}(\psi_i)$ (suitably normalized), as claimed. For $\chi_i$ ($i<s-\alpha$), we have
\[ \left.\frac{\chi_i'(1)}{\psi_1'(1)}\right/\frac{\chi_i(1)}{\psi_1(1)}=\frac{(q^{b_1}-q^{a_s+d})}{(q^{b_1}-q^{a_s})}\cdot \frac{(q^{a_i}-q^{a_s})}{(q^{a_i}-q^{a_s+d})}.\]
Applying the $\pi_{\kappa/d}$-function to this expression yields $0$ if $a_i>a_s+d$, again by Proposition \ref{prop:diffsinBdslinear}, and so we get $\pi_{\kappa/d}(\chi_i')=\pi_{\kappa/d}(\chi_i)$ in this case. (If there are no negative parameters, we can simply use $\chi_1$ to normalize instead, with the same outcome.)

If $s-\alpha<i\leq s$, then $(\chi_i'(1)/\psi_1(1))/(\chi_{i-1}(1)/\psi_1(1))$ satisfies the same equation as above, except this time $a_i<a_s+d$, so that term contributes $-1$. We therefore see that $\pi_{\kappa/d}(\chi_i')=\pi_{\kappa/d}(\chi_{i-1})-1$ for $i\geq s-\alpha$. Finally, we have
\[ \frac{\chi_{s-\alpha}'(1)}{\chi_s(1)}=\frac{q^{a_s}}{q^{a_s+d}}\prod_{j=2}^{s-1}\frac{(q^{a_s}-q^{a_j})}{(q^{a_s+d}-q^{a_j})}\prod_{j=1}^t\frac{(q^{a_s}+q^{b_j})}{(q^{a_s+d}+q^{b_j})}\]
The first quotient contributes $-2$, the third quotient contributes $0$, and the second quotient contributes $-1$ for each $i$ with $a_s<a_i<a_s+d$; therefore $\pi(\chi_{s-\alpha}')=\pi(\chi_s)-\alpha-2$.

Finally, by Proposition \ref{prop:increasingperv}, we must have that $\pi_{\kappa/d}(\chi_{s-\alpha}')<\pi_{\kappa/d}(\chi_{s-\alpha+1}')$, i.e., $\pi_{\kappa/d}(\chi_s)-\pi_{\kappa/d}(\chi_{s-\alpha})<\alpha+1$. However, since $\pi_{\kappa/d}(\chi_i)-\pi_{\kappa/d}(\chi_{i-1})\geq 1$, we must have that $\pi_{\kappa/d}(\chi_s)-\pi_{\kappa/d}(\chi_{s-\alpha})=\alpha$, and for $s-\alpha\leq i<s$, $\pi_{\kappa/d}(\chi_i)-\pi_{\kappa/d}(\chi_{i-1})=1$. Thus $\pi_{\kappa/d}(\chi_i)-\pi_{\kappa/d}(\chi_i')=2$ for $s-\alpha\leq i\leq s$, as claimed.
\end{pf}

Notice that the permutation of the parameters and change in $\pi_{\kappa/d}$-function is exactly that of Theorem \ref{thm:perversebijchange}.

By multiplying all parameters by $-1$ we change a $+$-type perturbation into a $-$-type perturbation, and so the same result -- that the change in $\pi_{\kappa/d}$-function and bijection is compatible with Theorem \ref{thm:perversebijchange} -- holds.

The last of the perturbations is the $\pm$-type one, where we alter both halves of the Brauer tree.

\begin{prop}\label{prop:plusminustypeperturbation} Suppose that $\kappa=1$. Let $\mc H'$ be the cyclotomic Hecke algebra obtained from $\mc H$ by applying a $\pm$-type perturbation that does not permute all parameters. Let $\alpha$ be the number of $i<s$ such that $a_i<b_t+d/2$, so that $a_{s-\alpha-1}>b_t+d/2>a_{s-\alpha}$, and let $\beta$ be the number of $i<t$ such that $b_i<a_s+d/2$, so that $b_{t-\beta-1}>a_s+d/2>b_{t-\beta}$. Write $\chi_i$ and $\psi_i$ for the characters of $\mc H$, and $\chi_i'$ and $\psi_i'$ for the characters of $\mc H'$. 
\begin{enumerate}
\item We have that $\pi_{\kappa/d}(\chi_i')=\pi_{\kappa/d}(\chi_i)$ for $1\leq i\leq s-\alpha-1$ and $\pi_{\kappa/d}(\psi_i')=\pi_{\kappa/d}(\psi_i)$ for all $1\leq i\leq t-\beta-1$. We have $\pi_{\kappa/d}(\chi_i')=\pi_{\kappa/d}(\chi_i)-2$ and $\pi_{\kappa/d}(\psi_i')=\pi_{\kappa/d}(\psi_i)-2$ otherwise.
\item Using the ordering on the permuted parameters inherited from the ordering on all parameters, they are permuted so that the $i$th permuted parameter of $\mc H'$ is the $(i-1)$th permuted parameter of $\mc H$.
\end{enumerate}
\end{prop}
\begin{pf} As in the previous case, $\chi_i$ and $\psi_i$ are the characters of $\mc H$, and write $\chi_i'$ and $\psi_i'$ for the characters of $\mc H'$, ordered in the standard way. We choose $\chi_1$ to normalize our $\pi_{\kappa/d}$-functions, noting that not all parameters are permuted, so by changing sign if necessary we may assume that $q^{a_1}$ is not moved.

If $i<s-\alpha$ then we have
\[ \left.\frac{\chi_i'(1)}{\chi_1'(1)}\right/\frac{\chi_i(1)}{\chi_1(1)}=\frac{(q^{a_1}+q^{a_s+d/2})}{(q^{a_1}-q^{a_s})}\cdot\frac{(q^{a_s}-q^{a_i})}{(q^{a_s+d/2}+q^{a_i})}\cdot \frac{(q^{b_t+d/2}-q^{a_1})}{(q^{b_t}+q^{a_1})}\cdot\frac{(q^{b_t}+q^{a_i})}{(q^{b_t+d/2}-q^{a_i})}.\]
By Proposition \ref{prop:diffsinBdsorthog}, since $a_i>b_t+d/2$ (and hence certainly $a_1>b_t+d/2$) the third and fourth terms contribute $0$ to the $\pi_{\kappa/d}$-function, and as $\kappa=1$ we also see that the first two terms contribute $0$ as well. Hence (i) holds for these characters. The proof of (i) for $\psi_i$ with $1\leq i\leq t-\beta-1$ is similar.

We now assume that $s-\alpha<i\leq s$, and compare $\chi_i'(1)$ with $\chi_{i-1}(1)$. The exact same formula for the quotient as above holds, but this time the $\pi_{\kappa/d}$-function will not evaluate to $0$. It still will on the first two terms, since $a_i>a_s$ and $a_1>a_s$, and as $a_1>b_t+d/2$, the third term still evaluates to $0$. However, the fourth term evaluates to $-1$, so that $\pi_{\kappa/d}(\chi_i')=\pi_{\kappa/d}(\chi_{i-1})-1$, as in Proposition \ref{prop:plustypeperturbation}. As in the previous paragraph, the proof that $\pi_{\kappa/d}(\psi_i')=\pi_{\kappa/d}(\psi_{i-1})-1$ for $t-\beta<i\leq t$ is similar.

It remains to compare $\chi_s$ with $\psi_{t-\beta}'$, and $\psi_t$ with $\chi_{s-\alpha}'$; we obtain
\[ \left.\frac{\psi_{t-\beta}'(1)}{\chi_1'(1)}\right/\frac{\chi_s(1)}{\chi_1(1)}=\frac{q^{a_s}}{q^{a_s+d/2}}\cdot\frac{(q^{a_1}-q^{b_t+d/2})}{(q^{a_1}+q^{b_t})}\cdot\frac{(q^{a_s}+q^{b_t})}{(q^{a_s+d/2}+q^{b_t+d/2})}\cdot\prod_{j=2}^{s-1}\frac{(q^{a_s}-q^{a_j})}{(q^{a_s+d/2}+q^{a_j})}\cdot\prod_{j=1}^{t-1}\frac{(q^{a_s}+q^{b_j})}{(q^{a_s+d/2}-q^{b_j})}.\]

Evaluating this quotient with the $\pi_{\kappa/d}$-function, and again using Proposition \ref{prop:diffsinBdsorthog}, we get a contribution of $-1$ from the first and third quotients, $0$ from the second and fourth quotients, and a contribution of $-1$ from each $j<t$ such that $a_s+d/2>b_j$, i.e., $\beta$. Therefore $\pi_{\kappa/d}(\psi_{t-\beta}')=\pi_{\kappa/d}(\chi_s)-\beta-2$, and similarly $\pi_{\kappa/d}(\chi_{s-\alpha}')=\pi_{\kappa/d}(\psi_t)-\alpha-2$.

We now combine these inequalities, together with the obvious inequalities $\pi_{\kappa/d}(\chi_{s-\alpha})\leq \pi_{\kappa/d}(\chi_s)-\alpha$ and $\pi_{\kappa/d}(\psi_{t-\beta})\leq \pi_{\kappa/d}(\psi_t)-\beta$, to get
\[ \pi_{\kappa/d}(\psi_t)-\alpha-1=\pi_{\kappa/d}(\chi_{s-\alpha}')+1<\pi_{\kappa/d}(\chi_{s-\alpha+1}')+1=\pi_{\kappa/d}(\chi_{s-\alpha})\leq \pi_{\kappa/d}(\chi_s)-\alpha,\]
so that $\pi_{\kappa/d}(\psi_t)\leq \pi_{\kappa/d}(\chi_s)$. Using the other inequalities we get $\pi_{\kappa/d}(\chi_s)\leq\pi_{\kappa/d}(\psi_t)$, and so $\pi_{\kappa/d}(\chi_s)=\pi_{\kappa/d}(\psi_t)$, and as in Proposition \ref{prop:plustypeperturbation}, $\pi_{\kappa/d}(\chi_s)-\pi_{\kappa/d}(\chi_{s-\alpha})=\alpha$ and for $s-\alpha\leq i<s$, $\pi_{\kappa/d}(\chi_i)-\pi_{\kappa/d}(\chi_{i-1})=1$, with similar statements for the $\psi_i$. Thus $\pi_{\kappa/d}(\chi_i)-\pi_{\kappa/d}(\chi_i')=2$ for $s-\alpha\leq i\leq s$, as claimed, and similarly for $\psi_i$, completing the proof of (i).

To see (ii), we must prove that the permuted parameters of $\mc H$ are, upon evaluation $q\mapsto\zeta$, in sequence 
\[ q^{a_{s-\alpha}},\dots,q^{a_{s-1}},q^{a_s},-q^{b_{t-\beta}},\dots,-q^{b_{t-1}},-q^{b_t},\]
for moving from $\mc H$ to $\mc H'$ results in cycling these parameters, as required by (ii). Clearly the sequence is correct on the $q^{a_i}$ and $-q^{b_i}$, so we must show that when we evaluate $q\mapsto \zeta$, none of the $-q^{b_i}$ evaluates to roots of unity in between $q^{a_j}$ and $q^{a_{j+1}}$ for some $s-\alpha\leq j\leq s-1$, in other words, that we cannot have that $a_j<b_i\pm d/2<a_{j+1}$. However, $b_i+d/2\geq b_t+d/2>a_{s-\alpha}>a_{j+1}$, and $a_j>a_s>b_{t-\alpha}-d/2>b_i-d/2$, so this set of inequalities cannot occur, and the sequence is correct. This proves (ii), and completes the proposition.
\end{pf}

Again, the permutation of the parameters and change in $\pi_{\kappa/d}$-function is exactly that of Theorem \ref{thm:perversebijchange}. In the next section we will compose these perturbations, and the stage at which no more perturbations are possible will be called a \emph{Coxeter Hecke algebra}, since it resembles the cyclotomic Hecke algebra corresponding to the Coxeter torus. It will be easy to prove that the $\pi_{\kappa/d}$-function and bijection are the canonical perversity function and bijection in the Coxeter case, and since repeated perturbations of parameters produce changes as described in Theorem \ref{thm:perversebijchange}, the bijection implied by specialization of parameters $q\mapsto\zeta$ will be the bijection required by Theorem \ref{thm:allperversecyclic}.

\section{Coxeter Hecke Algebras}
\label{sec:CoxHecke}
We continue to use our reductions of the previous section, namely that $\kappa=1$, that $d$ is even, and that $d/e$ is an integer. In cyclotomic Hecke algebras of finite groups of Lie type associated with the Coxeter number and $\kappa=1$, the eigenvalues are consecutive roots of unity. The next definition is the natural generalization of this.

\begin{defn} The \emph{Coxeter Hecke algebra} $\mc H_c$ of type $(s,t)$ and ambiance $d$ is the cyclotomic Hecke algebra of type $(s,t)$ and ambiance $d$ with the specialization of parameters
\[ 1,q^{-\ep},\dots,q^{-(s-1)\ep},-q^{-(s-t)\ep/2},-q^{-(s-t)\ep/2-\ep},\dots,-q^{-(s+t)\ep/2+\ep},\]
where $\ep=d/e$ with $s+t=e$.
\end{defn}

In other words, a Coxeter Hecke algebra -- since the parameters are defined only up to global shift -- consists of parameters whose exponents are in arithmetic progression with difference $\ep$, and such that the exponents of the positive and negative powers have the same arithmetic mean. This definition can be made without our restrictions on $d$ and $d/e$; however if $d/e$ is an integer and $d$ is even then all of the $a_i$ and $b_i$ in this definition are integers.

We now find the $\pi_{\kappa/d}$-function associated to a Coxeter Hecke algebra. We will assume that $s\geq t$, simply so we can take the $\pi_{\kappa/d}$-function relative to $\chi_1$; of course, we can take the $\pi_{\kappa/d}$-function relative to $\psi_1$ if $t>s$.

\begin{prop}\label{prop:Coxeterpifunction} Let $\mc H$ be the Coxeter Hecke algebra of type $(s,t)$ and ambiance $d$, and assume that $s\geq t$. The $\pi_{\kappa/d}$-function on the characters of $\mc H$ is the canonical perversity function on the Brauer tree of the line with exceptional node so that the two branches have lengths $s$ and $t$; in other words, $\pi_{\kappa/d}(\chi_i)=i-1$ and $\pi_{\kappa/d}(\psi_i)=s-t-1+i$.
\end{prop}
\begin{pf} Multiply the parameters by $q^{(s-1)\ep}$ so that all powers are non-negative. All terms involved are of the form $q^\alpha-q^\beta$, where $\alpha$ and $\beta$ lie in the range $\{0,\dots,\ep(s-1)\}$, and $q^\alpha+q^\beta$, where $\alpha\in\{0,\dots,\ep(s-1)\}$ and $\beta\in\{\ep(s-t)/2,\dots,\ep(s+t)/2-\ep\}$. In either case, all cyclotomic polynomials $\Phi_x$ that appear satisfy $x<d$, so that $\pi_{\kappa/d}(\Phi_x)=\deg(\Phi_x)/d$. In particular, this means that $\pi_{\kappa/d}(\chi_i)$ is simply $(A(\chi_i)+a(\chi_i))/d$ plus half the multiplicity of $1$ as a zero of $\chi_i$, and similarly for $\psi_i$. Since $a(q^\alpha\pm q^\beta)+A(q^\alpha\pm q^\beta)=\alpha+\beta$, it is easy to evaluate this for a relative degree.

Normalize with respect to $\chi_1$. We have, writing $\gamma=e/2$,
\[ \frac{\chi_i(1)}{\chi_1(1)}=\frac{q^{\ep(s-1)}}{q^{\ep(s-i)}}\cdot \frac{\D \prod_{j=2}^s (q^{\ep(s-1)}-q^{\ep(s-j)})}{\D \mathop{\prod_{j=1}^s}_{j\neq i} (q^{\ep(s-i)}-q^{\ep(s-j)})}\cdot\frac{\D \prod_{j=1}^t (q^{\ep(s-1)}+q^{\ep(\gamma-j)})}{\D \prod_{j=1}^t (q^{\ep(s-i)}+q^{\ep(\gamma-j)})}.\]
Using the above observation, the sum of the $a$- and $A$-functions on each of these quotients is $2\ep(i-1)$, $\ep(i-1)(s-2)$ and $\ep(i-1)t$, yielding $(i-1)d$. Since there are equal numbers of $\Phi_1$-terms on top and bottom of the quotient, we get that $\pi_{\kappa/d}(\chi_i)=(i-1)$, as claimed.

For $\psi_i$, we get
\[ \frac{\psi_i(1)}{\chi_1(1)}=\frac{q^{\ep(s-1)}}{q^{\ep(\gamma-i)}}\cdot \frac{\D \prod_{j=2}^s (q^{\ep(s-1)}-q^{\ep(s-j)})}{\D \prod_{j=1}^s (q^{\ep(\gamma-i)}+q^{\ep(s-j)})}\cdot\frac{\D \prod_{j=1}^t (q^{\ep(s-1)}+q^{\ep(\gamma-j)})}{\D \mathop{\prod_{j=1}^t}_{j\neq i} (q^{\ep(\gamma-i)}-q^{\ep(\gamma-j)})}.\]
This time there are $(s-1)$ copies of $\Phi_1$ on the top and $(t-1)$ copies of $\Phi_1$ on the bottom, contributing $(s-t)/2=s-\gamma$ to $\pi_{\kappa/d}(\psi_i)$. The $a$- and $A$-functions yield $2\ep(s-1)+2\ep(i-\gamma)$, $\ep(s-1)(s-2)+s\ep(i-\gamma)$ and $t\ep(s-1)+\ep(t-2)(i-\gamma)$, whose sum is $d(s-1+i-\gamma)$, and so 
\[ \pi_{\kappa/d}(\psi_i)=(s-1+i-\gamma)+s-\gamma=s-t-1+i,\]
as needed.
\end{pf}

It is easy to see that the ordering on the simple modules in the Coxeter Hecke algebra is the canonical ordering, and so the $\pi_{\kappa/d}$-function and ordering are compatible in this case.

Our main result is that, given an arbitrary cyclotomic Hecke algebra with our restrictions on $\kappa$, $d$ and $d/e$, repeated perturbation of the parameters eventually reduces it to a Coxeter Hecke algebra. The next result shows that perturbations are nested, i.e., the set of parameters that they permute gets larger: these will become the cohomologically closed sets $I_j$ that we used in the proof of Theorem \ref{thm:allperversecyclic}. Recall that we have no choice about the perturbations that we apply, and so we will simply say `apply a perturbation'. 

\begin{prop}\label{prop:inclusionofpermutation} Let $\mc H$ be a cyclotomic Hecke algebra of type $(s,t)$ and ambiance $d$, with parameters $q^{a_1},\dots,a^{a_s}$ and $-q^{b_1},\dots,-q^{b_t}$. Apply a perturbation on $\mc H$ to produce the algebra $\mc H'$, with the set $I$ of parameters being permuted. Apply a perturbation to $\mc H'$ to get $\mc H''$, with set $I'$ of permuted parameters. We have $I\subs I'$.
\end{prop}

This proposition is a trivial consequence of the definition of perturbations, together with the observation that, if the first perturbation applied is of $\pm$-type then so is the second one.

The main aim of all of the definitions and results of the last section is the following theorem.

\begin{thm}\label{thm:perturbationchain} Let $\mc H_0$ be a cyclotomic Hecke algebra of type $(s,t)$ and ambiance $d$. Write $e=s+t$ and $\ep=d/e$. Inductively we perturb the algebra $\mc H_i$ to produce a new algebra $\mc H_{i+1}$. Assume that $s\geq t$.
\begin{enumerate}
\item There exists $n$ such that $\mc H_n$ and $\mc H_{n+1}$ have the same parameters (recall that parameters are only defined up to a global shift by a power of $q$). The algebra $\mc H_n$ is a Coxeter Hecke algebra.
\end{enumerate}
Let $n$ denote the smallest such number.
\begin{enumerate}
\item[(ii)] Write $I_j$ for the set of permuted parameters of $\mc H_j$. We have a chain
\[ I_1\subs I_2\subs \cdots \subs I_{n-1}\]
of \emph{proper} subsets of $\{1,\dots,e\}$. Let $\chi_i$ and $\psi_i$ denote the relative degrees of $\mc H_0$, normalized by $\chi_1$. For a given $i$, let $f(\chi_i)$ denote the largest $j$ such that $q^{a_i}\in I_j$, and similarly for $f(\psi_i)$. We have that $\pi_{\kappa/d}(\chi_i)=2f(\chi_i)+(i-1)$ and $\pi_{\kappa/d}(\psi_i)=2f(\psi_i)+(s-t+i-1)$.
\end{enumerate}
\end{thm}
\begin{pf} Let $\mc H$ be an arbitrary cyclotomic Hecke algebra with type $(s,t)$ and ambiance $d$, and write $q^{a_1},\dots,q^{a_s}$ and $-q^{b_1},\dots,-q^{b_t}$ for its parameters. Suppose firstly that perturbing $\mc H$ results in all parameters of $\mc H$ being permuted; we will prove that $\mc H$ is a Coxeter Hecke algebra and is isomorphic to its perturbation.

Suppose that $t>0$; since all parameters are permuted in the perturbation, we must have that $-q^{a_s+d/2}$ is the largest negative parameter, so that $a_s+d/2>b_1$, and similarly $b_t+d/2>a_1$. Clearly, since under evaluation $q\mapsto \zeta$ the parameters map to distinct $e$th roots of unity, we must have that $a_i\geq a_s+(s-i)\ep$, and similary $b_i\geq b_t+(t-i)\ep$, and since $b_t+d/2$ and $a_1$ must differ by a multiple of $\ep$, we also get $a_s+d/2\geq b_1+\ep$ and $b_t+d/2\geq a_1+\ep$. Combining these last four inequalities gives
\[ a_s+d\geq b_1+d/2+\ep\geq b_t+t\ep+d/2\geq a_1+(t+1)\ep\geq a_s+e\ep=a_s+d.\]
Thus all of the inequalities are actually equalities, and so the $a_i$ and $b_i$ are consecutive multiples of $\ep$; it is easy to see that $\mc H$ is a Coxeter Hecke algebra, and after permutation of parameters the perturbation is simply multiplying all entries by $q^\ep$, hence an isomorphism. If $t=0$ then the same argument, this time with a $+$-type perturbation, proves that $\mc H$ is also a Coxeter Hecke algebra.

Thus we need to prove that there exists $n$ such that $I_n$ is the set of all parameters of $\mc H_n$. For this, we simply note that, if $\pm q^x$ is some parameter of $\mc H_0$, then repeated perturbations increase the exponent belonging to the smallest parameter, or parameters in the case of $\pm$-type perturbations, so that eventually the smallest parameter will be within $d/2$ of $x$, so that $\pm q^x$ is permuted. Similarly, eventually all parameters are permuted for some $n$, completing the proof of (i).

For (ii), the inclusion of subsets follows from Proposition \ref{prop:inclusionofpermutation}, and the statement about the $\pi_{\kappa/d}$-function of the $\chi_i$ and $\psi_i$ follows from the calculations of the $\pi_{\kappa/d}$-functions in Propositions \ref{prop:plustypeperturbation} and \ref{prop:plusminustypeperturbation}, together with the computation of the $\pi_{\kappa/d}$-function of the Coxeter Hecke algebra in Proposition \ref{prop:Coxeterpifunction}.
\end{pf}

(As before, the restriction that $s\geq t$ can be removed, with the $\pi_{\kappa/d}$-function altered in the obvious way if $t>s$.)

With the results that we have collated so far we are able to produce the proof of Theorem \ref{thm:combbroueconj} for the classical groups, and in fact any unipotent block whose Brauer tree is a line. By Theorems \ref{thm:piincreasetoexc} and \ref{thm:allperversecyclic} there is a perverse equivalence between $B$ and $B'$, and we must show that the bijection is as suggested in Theorem \ref{thm:combbroueconj}. This bijection, up to a rotation of the Brauer tree of the Brauer correspondent $B'$, is correct for the Coxeter case as we have seen in this section; by Theorem \ref{thm:perturbationchain} every cyclotomic Hecke algebra of type $(s,t)$ can be perturbed into a Coxeter Hecke algebra, and by Propositions \ref{prop:plustypeperturbation} and \ref{prop:plusminustypeperturbation}, the alterations to the $\pi_{\kappa/d}$-function and the bijection are consistent with that required from Theorem \ref{thm:perversebijchange}. Since any unipotent block whose Brauer tree is a line has a cyclotomic Hecke algebra of type $(s,t)$, this proves that the bijection given by combinatorial Brou\'e's conjecture is correct up to a rotation of the Brauer tree of $B'$. Finally, we consider a unipotent character $\chi$ with minimal $\pi_{\kappa/d}$-function: if $S$ denotes the associated simple $B$-module and $M$ its Green correspondent in $B'$, then the simple $B'$-module in bijection with $S$ in the perverse equivalence is $\Omega^{\pi_{\kappa/d}(\chi)}(M)$, which is at position $\omega_\chi \zeta^{\aA(\chi)/e}$, since $M$ is in position $\omega_\chi$. This proves that the bijection suggested by combinatorial Brou\'e's conjecture agrees with the actual bijection at a module, and hence they are the same. This completes the proof.

\section{The Exceptional Groups}
\label{sec:exceptionals}

We consider four examples: $d=3$ for $G_2(q)$, $d=12''$ for ${}^2G_2(q)$, $d=24'$ for ${}^2\!F_4(q)$, and $d=14$ for $E_7(q)$. For exceptional groups of Lie type, the Brauer tree is nearly always either a line or a line with one pair of non-real vertices. The first and last case we consider is of this latter type.

It is fairly easy to prove the combinatorial Brou\'e conjecture for $\kappa<d$, but to prove it for all $\kappa$ we need to know something about the $A$-function. Let $X_{\geq i}$ denote the set of unipotent characters $\chi$ of $B$ for which $A(\chi)-A(\bo\lambda)$ is at least $i$, and let $c_1,\dots,c_n$ denote those integers for which $|X_{\geq c_i}|>|X_{\geq c_i+1}|$. By Lemma \ref{lem:add2pi}, when moving from $\kappa$ to $\kappa+d$, one adds $2A(\chi)$ to $\pi_{\kappa/d}(\chi)$. If the parameters are $\omega_i q^{v_i}$, and $v_i$ is an integer, then the root of unity obtained by specialization $q\mapsto \e^{2\pi\I\kappa/d}$ does not change upon replacement of $\kappa$ by $\kappa+d$, so the bijections for the perverse equivalences with $\pi_{\kappa/d}$ and $\pi_{(\kappa+d)/d}$ are the same. (Recall that the $v_i$ are semi-integers, and we give an example where they are not integers in $E_7$ and $d=14$.) Examining Theorem \ref{thm:perversebijchange}, it is easy to see that if we add $2n$ to the $\pi$-function of exactly $n$ of the simple modules, the bijection stays the same, since the change in bijection is applying an $n$-cycle.

We now see what we need to know in order to prove that there is a perverse equivalence with $\pi_{(\kappa+d)/d}$ as perversity function, given that there is one with $\pi_{\kappa/d}$ as one, and that they have the same bijection:
\begin{enumerate}
\item the sets $X_{\geq c_i}$ are cohomologically closed with respect to the canonical perversity function on the simple $B$-modules, and there exists $j$ such that
\[ J_j \subs X_{\geq c_i}\subs J_{j-1};\]
\item the size of $X_{\geq c_i}$ must be divide $c_i-c_{i-1}$;
\item all powers of $q$ in the parameters of the cyclotomic Hecke algebra are integral.
\end{enumerate}
The first condition means that the sets $J_j$ constructed in the proof of Theorem \ref{thm:allperversecyclic}, applied to the function $\pi_{(\kappa+d)/d}$, are firstly cohomologically closed, and secondly are the same as those of $\pi_{\kappa/d}$, together with $c_i-c_{i-1}$ copies of $X_{\geq c_i}$; the second condition implies that the set $X_{\geq c_i}$ appears a multiple of $|X_{\geq c_i}|$ times, and so the bijection remains unchanged; the third condition is trivial to check when it holds.

In each of the first three cases we give the $A$-function, and the reader may note that it does satisfy these three conditions. The final example does not satisfy the third condition, but does satisfy the first. The modification needed to the second condition to take account of the semi-integrality of the parameter powers is intuitive, and we detail it in Section \ref{ssec:e7}

\subsection{$G_2(q)$, $d=3$}

Here there is a single unipotent block, the principal block, and it has six unipotent characters, so that the cyclotomic Weyl group is $Z_6$. Hence substituting $\zeta=\e^{2\kappa\pi\I/3}$ to the parameters should produce the set of $6$th roots of unity. We give the table below, ordered so that substitution $q\mapsto\e^{2\pi\I/3}$ (i.e., $\kappa=1$) gives the $6$th roots of unity in order.
\begin{center}\begin{tabular}{lcccc}
\hline Character & $A(-)$ &  $\omega_i q^{aA/e}$ & $\kappa=1$ & $\kappa=2$
\\\hline $\phi_{1,0}$ & $0$ & $q^{0}$ & $0$ & $0$
\\ $G_2[\theta^2]$ & $5$ & $-\theta q$ & $3$ & $7$
\\ $\phi_{2,2}$ & $5$ & $q$ & $3$ & $7$
\\ $G_2[\theta]$ & $5$ & $-\theta^2 q$ & $3$ & $7$
\\ $\phi_{1,6}$ & $6$ & $q^{2}$ & $4$ & $8$
\\ $G_2[1]$ & $5$ & $-q$ & $4$ & $6$
\\\hline\end{tabular}\end{center}

We now give the Brauer tree of this block, taken from \cite{shamash1989}, with the $\pi_{\kappa/d}$-function in the case $\kappa=1$ attached.

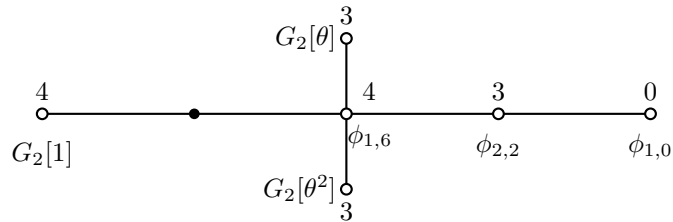
\begin{center}\begin{tikzpicture}[thick,scale=2]
\draw (0,0.5) -- (-4,0.5);
\draw (-2,0) -- (-2,1);

\draw (-0,0.5) node [draw,label=below:$\phi_{1,0}$] (l0) {};
\draw (-1,0.5) node [draw,label=below:$\phi_{2,2}$] (l1) {};
\draw (-4,0.5) node [draw,label=below:${G_2[1]}$] (l2) {};
\draw (-3,0.5) node [fill=black!100] (ld) {};
\draw (-2,0) node [draw,label=left:${G_2[\theta^2]}$] (l4) {};
\draw (-2,0.5) node [draw,label=below right:$\phi_{1,6}$] (l4) {};
\draw (-2,1) node [draw,label=left:${G_2[\theta]}$] (l4) {};
\draw (-0,0.65) node{$0$};
\draw (-1,0.65) node{$3$};
\draw (-1.85,0.65) node{$4$};
\draw (-4,0.65) node{$4$};
\draw (-2,1.15) node{$3$};
\draw (-2,-0.15) node{$3$};
\end{tikzpicture}\end{center}
The canonical ordering here is $\phi_{1,0},\phi_{2,2},G_2[\theta],\phi_{1,6},G_2[1],G_2[\theta^2]$, and all characters apart from the trivial $\phi_{1,0}$ have had $2$ added to them, so the new ordering should be $\phi_{1,0},G_2[\theta^2],\phi_{2,2},G_2[\theta],\phi_{1,6},G_2[1]$ according to Theorem \ref{thm:perversebijchange}. This is the ordering given in the table, and so the combinatorial form of Brou\'e's conjecture is verified in this case. The case $\kappa=2$ is similar, but with more applications of Theorem \ref{thm:perversebijchange}.

\subsection{${}^2G_2(q)$, $d=12'$}

The unipotent characters of ${}^2G_2(q)$ are given in \cite{carterfinite}, but here we use a slightly different notation according to the eigenvalue of the Frobenius, and a different definition of $\Phi_{12}'$ consistent with the cases of Suzuki and big Ree groups; $\Phi_{12}'$ here is defined as $(q-\xi^5)(q-\xi^7)$ (this is $\Phi_{12}''$ in \cite{carterfinite}), where $\xi=\e^{2\pi\I/12}$, so that this is the case $d=12$ and $\kappa=5,7$. Here there is a single unipotent block, the principal block, and it again has six unipotent characters, so that the cyclotomic Weyl group is $Z_6$. Hence substituting $\zeta=\e^{2\kappa\pi\I/12}$ (with $\kappa=5,7$) to the parameters should produce the set of $6$th roots of unity. We give the table below, ordered so that substitution $q\mapsto\e^{5\pi\I/6}$ (i.e., $\kappa=5$) gives the $6$th roots of unity in order.
\begin{center}\begin{tabular}{lcccc}
\hline Character & $A(-)$ & $\omega_i q^{aA/e}$ & $\kappa=5$ & $\kappa=7$
\\\hline $\phi_{1,0}$ &  $0$ & $q^{0}$ & $0$ & $0$
\\ ${}^2G_2^{\mathrm{II}}[-\I]$ & $5$ & $-\I q$ & $4$ & $6$
\\ ${}^2G_2[\xi^7]$ & $5$ & $\xi^{11} q$ & $4$ & $6$
\\ ${}^2G_2[\xi^5]$ & $5$ & $\xi q$ & $4$ & $6$
\\ ${}^2G_2^{\mathrm{II}}[\I]$ & $5$ & $\I q$ & $4$ & $6$
\\ $\phi_{1,2}$ & $6$ & $q^2$ & $5$ & $7$
\\\hline\end{tabular}\end{center}
When $\ell\mid\Phi_{12}'$ we get the following tree, determined in \cite{hiss1991}, with the $\pi_{\kappa/d}$-function in the case $\kappa=5$ attached.
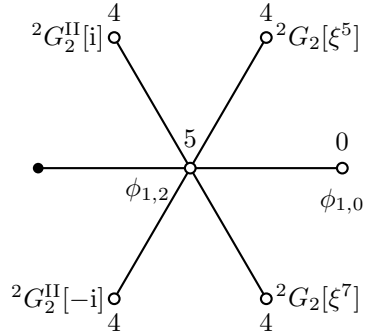
\begin{center}\begin{tikzpicture}[thick,scale=2]
\draw (-1.28,0.85) node{$\phi_{1,2}$};
\draw (0,1) -- (-2,1);
\draw (-0.5,0.1339) -- (-1.5,1.866);
\draw (-0.5,1.866) -- (-1.5,0.1339);
\draw (0,1) node [draw,label=below:$\phi_{1,0}$] (l0) {};
\draw (-1,1) node [draw] (l2) {};
\draw (-2,1) node [fill=black!100] (ld) {};
\draw (0,1.18) node{$0$};
\draw (-1,1.2) node{$5$};

\draw (-0.5,1.866) node [draw,label=right:${{}^2G_2[\xi^5]}$] (l0) {};
\draw (-1.5,1.866) node [draw,label=left:${{}^2G_2^\mathrm{II}[\I]}$] (l0) {};
\draw (-0.5,0.1339) node [draw,label=right:${{}^2G_2[\xi^7]}$] (l0) {};
\draw (-1.5,0.1339) node [draw,label=left:${{}^2G_2^\mathrm{II}[-\I]}$] (l0) {};

\draw (-0.5,2.03) node{$4$};
\draw (-1.5,2.03) node{$4$};
\draw (-0.5,-0.02) node{$4$};
\draw (-1.5,-0.02) node{$4$};
\end{tikzpicture}
\end{center}
The canonical ordering here is $\phi_{1,0},{}^2G_2[\xi^5],{}^2G_2^{\mathrm{II}}[\I],\phi_{1,6},{}^2G_2^{\mathrm{II}}[-\I],{}^2G_2[\xi^7]$, and all characters apart from the trivial $\phi_{1,0}$ have had $4$ added to them, so the new ordering should be $\phi_{1,0},{}^2G_2^{\mathrm{II}}[-\I],{}^2G_2[\xi^7],{}^2G_2[\xi^5],{}^2G_2^{\mathrm{II}}[\I],\phi_{1,6}$ according to Theorem \ref{thm:perversebijchange}. This is the ordering given in the table, and so the combinatorial form of Brou\'e's conjecture is verified in this case. The case $\kappa=7$ is similar, but with another application of Theorem \ref{thm:perversebijchange}.

\subsection{${}^2\!F_4(q)$, $d=24'$}

To be consistent with the previous section, set $\Phi_8''$ and $\Phi_{24}''$ to be the factors of $\Phi_8$ and $\Phi_{24}$ (which are reducible over $\Z[\sqrt2]$) which take zero on $\e^{2\pi\I/8}$ and $\e^{2\pi\I/24}$ respectively. Because there are misprints in the table of degrees in \cite{carterfinite}, we give the degrees of those characters for which $\Phi_{24}'$ does not divide their degree here. Let $\psi=\e^{2\pi\I/8}$.
\begin{center}\begin{tabular}{llcccccc}
\hline Name & Degree & $A(-)$ & $\omega_i q^{aA/e}$ & $\kappa=5$ & $\kappa=11$ & $\kappa=13$ & $\kappa=19$
\\ \hline $\phi_{1,0}$ & $1$ & $0$ & $1$ & $0$ & $0$ & $0$ & $0$
\\ ${}^2\!B_2[\psi^3],1$ & $q\Phi_1\Phi_2\Phi_4^2\Phi_{12}/\sqrt2$ & $11$ & $\psi^7q$ & $4$ & $10$ & $12$ & $18$
\\ ${}^2\!F_4^{\mathrm{II}}[-\I]$ & $q^4\Phi_1^2\Phi_2^2\Phi_4^2\Phi_{12}\Phi_{24}''/4$ & $20$ & $-\I q^2$ & $8$ & $18$ & $22$ & $32$
\\ ${}^2\!F_4[-\theta^2]$ & $q^4\Phi_1^2\Phi_2^2\Phi_4^2\Phi_8^2/3$ & $20$ & $-\theta q^2$ & $8$ & $18$ & $22$ & $32$
\\ ${}^2\!B_2[\psi^5],1$ & $q\Phi_1\Phi_2\Phi_4^2\Phi_{12}/\sqrt2$ & $11$ & $\psi q$ & $4$ & $10$ & $12$ & $18$
\\ $\phi_{2,1}$ & $q^4\Phi_4^2\Phi_8'^2\Phi_{12}\Phi_{24}''/4$ & $20$ & $q^2$ & $7$ & $17$ & $21$ & $31$
\\ ${}^2\!B_2[\psi^3],\ep$ & $q^{13}\Phi_1\Phi_2\Phi_4^2\Phi_{12}/\sqrt2$ & $23$ & $\psi^7q^3$ & $9$ & $21$ & $25$ & $37$
\\ ${}^2\!F_4[-\theta]$ & $q^4\Phi_1^2\Phi_2^2\Phi_4^2\Phi_8^2/3$ & $20$ & $-\theta^2q^2$ & $8$ & $18$ & $22$ & $32$
\\ ${}^2\!F_4^{\mathrm{II}}[\I]$ & $q^4\Phi_1^2\Phi_2^2\Phi_4^2\Phi_{12}\Phi_{24}''/4$ & $20$ & $\I q^2$ & $8$ & $18$ & $22$ & $32$
\\ ${}^2\!B_2[\psi^5],\ep$ & $q^{13}\Phi_1\Phi_2\Phi_4^2\Phi_{12}/\sqrt2$ & $23$ & $\psi q^3$ & $9$ & $21$ & $25$ & $37$
\\ $\phi_{1,8}$ & $q^{24}$ & $24$ & $q^4$ & $10$ & $22$ & $26$ & $38$
\\ ${}^2\!F_4^{\mathrm{II}}[-1]$& $q^4\Phi_1^2\Phi_2^2\Phi_8''^2\Phi_{12}\Phi_{24}''/12$ & $20$ & $-q^2$ & $10$ & $20$ & $24$ & $32$
\\\hline\end{tabular}\end{center}
When $\ell\mid\Phi_{24}'$ we get the following tree, determined in \cite{hiss1991}, with the $\pi_{\kappa/d}$-function in the case $\kappa=5$ attached.

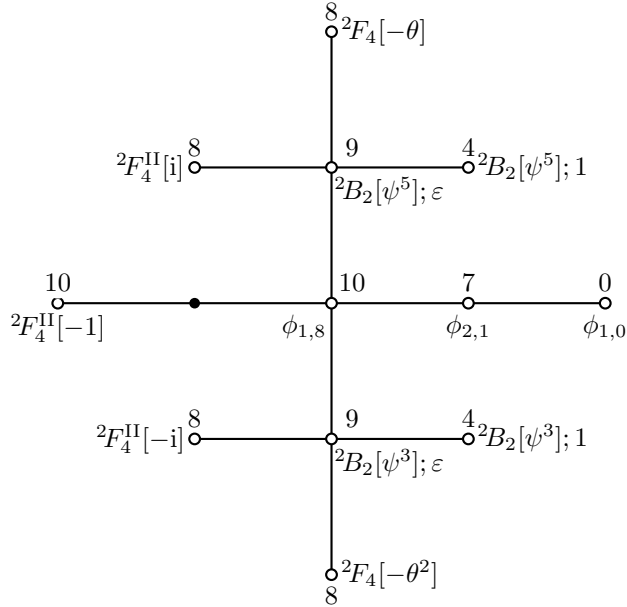
\begin{center}\begin{tikzpicture}[thick,scale=1.8]
\draw (2.422,2.82) node{${}^2\!B_2[\psi^5];\ep$};
\draw (2.422,0.82) node{${}^2\!B_2[\psi^3];\ep$};
\draw (1.8,1.82) node{$\phi_{1,8}$};
\draw (3,1.82) node{$\phi_{2,1}$};
\draw (4,1.82) node{$\phi_{1,0}$};
\draw (0,1.82) node{${}^2\!F_4^{\mathrm{II}}[-1]$};

\draw (0,2) -- (4,2);
\draw (1,1) -- (3,1);
\draw (1,3) -- (3,3);
\draw (2,0) -- (2,4);

\draw (4,2) node [draw] (l0) {};
\draw (3,2) node [draw] (l1) {};
\draw (2,2) node [draw] {};
\draw (1,2) node [fill=black!100] (ld) {};
\draw (0,2) node [draw] (l2) {};
\draw (3,1) node [draw,label=right:${{}^2\!B_2[\psi^3];1}$] {};
\draw (2,1) node [draw] (l4) {};
\draw (1,1) node [draw,label=left:${{}^2\!F_4^{\mathrm{II}}[-\I]}$] (l4) {};
\draw (2,0) node [draw,label=right:${{}^2\!F_4[-\theta^2]}$] (l4) {};
\draw (3,3) node [draw,label=right:${{}^2\!B_2[\psi^5];1}$] {};
\draw (2,3) node [draw] (l4) {};
\draw (1,3) node [draw,label=left:${{}^2\!F_4^{\mathrm{II}}[\I]}$] (l4) {};
\draw (2,4) node [draw,label=right:${{}^2\!F_4[-\theta]}$] (l4) {};

\draw (4,2.15) node{$0$};
\draw (3,2.15) node{$7$};
\draw (2.15,2.15) node{$10$};
\draw (0,2.15) node{$10$};
\draw (1,3.15) node{$8$};
\draw (1,1.15) node{$8$};
\draw (3,3.15) node{$4$};
\draw (3,1.15) node{$4$};
\draw (2,4.15) node{$8$};
\draw (2,-0.15) node{$8$};
\draw (2.15,3.15) node{$9$};
\draw (2.15,1.15) node{$9$};
\end{tikzpicture}\end{center}

The canonical ordering here is 
\[ \phi_{1,0},\phi_{2,1},{}^2\!B_2[\psi^5];1,{}^2\!F_4[-\theta],{}^2\!F_4^{\mathrm{II}}[\I],{}^2\!B_2[\psi^5];\ep,\phi_{1,8},{}^2\!F_4^{\mathrm{II}}[-1],{}^2\!F_4^{\mathrm{II}}[-\I],{}^2\!F_4[-\theta^2],{}^2\!B_2[\psi^3];1,{}^2\!B_2[\psi^3];\ep.\]
We add $4$ to each non-trivial character, and the ordering changes to 
\[ \phi_{1,0},{}^2\!B_2[\psi^3];1,{}^2\!B_2[\psi^3];\ep,\phi_{2,1},{}^2\!B_2[\psi^5];1,{}^2\!F_4[-\theta],{}^2\!F_4^{\mathrm{II}}[\I],{}^2\!B_2[\psi^5];\ep,\phi_{1,8},{}^2\!F_4^{\mathrm{II}}[-1],{}^2\!F_4^{\mathrm{II}}[-\I],{}^2\!F_4[-\theta^2].\]
At this point we have reached the correct $\pi_{\kappa/d}$-function for ${}^2\!B_2[\psi^i],1$, but adding another $2$ is needed for $\phi_{2,1}$ to be in place. This yields
\[ \phi_{1,0},{}^2\!B_2[\psi^3];1,{}^2\!F_4[-\theta^2],{}^2\!B_2[\psi^3];\ep,{}^2\!B_2[\psi^5];1,\phi_{2,1},{}^2\!F_4[-\theta],{}^2\!F_4^{\mathrm{II}}[\I],{}^2\!B_2[\psi^5];\ep,\phi_{1,8},{}^2\!F_4^{\mathrm{II}}[-1],{}^2\!F_4^{\mathrm{II}}[-\I].\]
We now fix $\phi_{1,0}$, $\phi_{2,1}$, ${}^2\!B_2[\psi^3],1$ and ${}^2\!B_2[\psi^5],1$, and adding $2$ to all remaining characters yields the correct bijection, which is
\[ \phi_{1,0},{}^2\!B_2[\psi^3];1,{}^2\!F_4^{\mathrm{II}}[-\I],{}^2\!F_4[-\theta^2],{}^2\!B_2[\psi^5];1,\phi_{2,1},{}^2\!B_2[\psi^3];\ep,{}^2\!F_4[-\theta],{}^2\!F_4^{\mathrm{II}}[\I],{}^2\!B_2[\psi^5];\ep,\phi_{1,8},{}^2\!F_4^{\mathrm{II}}[-1].\]

\subsection{$E_7(q)$, $d=14$}
\label{ssec:e7}

This example is included because it is one of the few blocks of exceptional groups for which there are parameters whose power of $q$ is a semi-integer, and so replacing $\kappa$ by $\kappa+d$ \emph{does} alter the bijection.

Here there is a single unipotent block, the principal block, and it has fourteen unipotent characters, so that the cyclotomic Weyl group is $Z_{14}$. Hence substituting $\zeta=\e^{2\kappa\pi\I/14}$ to the parameters should produce the set of $14$th roots of unity. We give the table below, ordered so that substitution $q\mapsto\e^{2\pi\I/14}$ (i.e., $\kappa=1$) gives the $14$th roots of unity in order.
\begin{center}\begin{tabular}{lcccccccc}
\hline Character & $A(-)$ & $\omega_i q^{aA/e}$ & $\kappa=1$ & $\kappa=3$ & $\kappa=5$ & $\kappa=9$ & $\kappa=11$ & $\kappa=13$
\\\hline $\phi_{1,0}$ & $0$ & $q^{0}$ & $0$ & $0$ & $0$ & $0$ & $0$ & $0$
\\ $E_7[-\I]$ & $52$ & $-\I q^{9/2}$ & $8$ & $22$ & $38$ & $66$ & $82$ & $96$
\\ $\phi_{27,2}$ & $26$ & $q^{2}$ & $3$ & $11$ & $19$ & $33$ & $41$ & $49$
\\ $\phi_{105,5}$ & $38$ & $q^{3}$ & $4$ & $16$ & $26$ & $50$ & $60$ & $72$
\\ $\phi_{189,10}$ & $48$ & $q^{4}$ & $5$ & $21$ & $35$ & $61$ & $75$ & $91$
\\ $\phi_{189,17}$ & $55$ & $q^{5}$ & $6$ & $24$ & $40$ & $70$ & $86$ & $104$
\\ $\phi_{105,26}$ & $59$ & $q^{6}$ & $7$ & $25$ & $41$ & $77$ & $93$ & $111$
\\ $\phi_{27,37}$ & $61$ & $q^{7}$ & $8$ & $26$ & $44$ & $78$ & $96$ & $114$
\\ $E_7[\I]$ & $52$ & $\I q^{9/2}$ & $8$ & $22$ & $38$ & $66$ & $82$ & $96$
\\ $\phi_{1,63}$ & $63$ & $q^{9}$ & $9$ & $27$ & $45$ & $81$ & $99$ & $117$
\\ $D_4;\ep_1$ & $38$ & $-q^{3}$ & $6$ & $16$ & $28$ & $48$ & $60$ & $70$
\\ $D_4;r\ep_1$ & $48$ & $-q^{4}$ & $7$ & $21$ & $35$ & $61$ & $75$ & $89$
\\ $D_4;r\ep_2$ & $55$ & $-q^{5}$ & $8$ & $24$ & $40$ & $70$ & $86$ & $102$
\\ $D_4;\ep_2$ & $59$ & $-q^{6}$ & $9$ & $25$ & $43$ & $75$ & $93$ & $109$
\\\hline\end{tabular}\end{center}
When $\ell\mid\Phi_{14}$ we get the following tree, determined in as-yet unpublished work of Dudas, Rouquier and the author, with the $\pi_{\kappa/d}$-function in the case $\kappa=1$ attached.
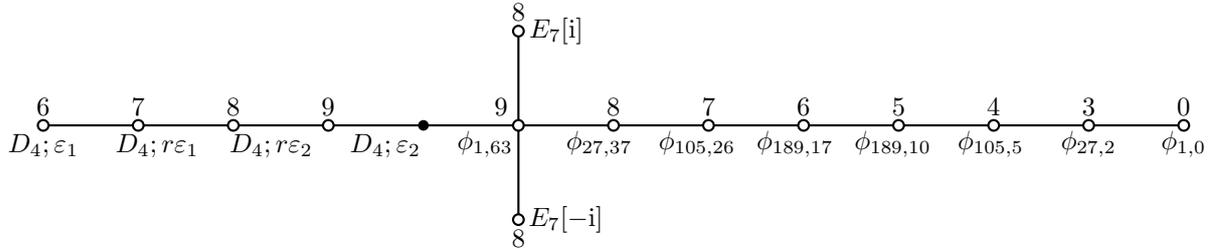
\begin{center}\begin{tikzpicture}[thick,scale=1.25]
\draw (0,-0.22) node{${D_4;\ep_1}$};
\draw (1.2,-0.22) node{$D_4;r\ep_1$};
\draw (2.4,-0.22) node{$D_4;r\ep_2$};
\draw (3.6,-0.22) node{$D_4;\ep_2$};
\draw (4.65,-0.2) node{$\phi_{1,63}$};
\draw (5.85,-0.2) node{$\phi_{27,37}$};
\draw (6.88,-0.2) node{$\phi_{105,26}$};
\draw (7.91,-0.2) node{$\phi_{189,17}$};
\draw (8.94,-0.2) node{$\phi_{189,10}$};
\draw (9.97,-0.2) node{$\phi_{105,5}$};
\draw (11,-0.2) node{$\phi_{27,2}$};
\draw (12,-0.2) node{$\phi_{1,0}$};

\draw (5.4,1) node {$E_7[\I]$};
\draw (5.5,-1) node {$E_7[-\I]$};

\draw (0,0.2) node{$6$};
\draw (1,0.2) node{$7$};
\draw (2,0.2) node{$8$};
\draw (3,0.2) node{$9$};
\draw (4.82,0.2) node{$9$};

\draw (5,-1.2) node{$8$};
\draw (5,1.2) node{$8$};

\draw (6,0.2) node{$8$};
\draw (7,0.2) node{$7$};
\draw (8,0.2) node{$6$};
\draw (9,0.2) node{$5$};
\draw (10,0.2) node{$4$};
\draw (11,0.2) node{$3$};
\draw (12,0.2) node{$0$};

\draw (0,0) -- (12,0);
\draw (5,-1) -- (5,1);

\draw (12,0) node [draw] (l0) {};
\draw (11,0) node [draw] (l0) {};
\draw (10,0) node [draw] (l0) {};
\draw (9,0) node [draw] (l0) {};
\draw (8,0) node [draw] (l0) {};
\draw (7,0) node [draw] (l1) {};
\draw (6,0) node [draw] (l2) {};
\draw (5,0) node [draw] (l2) {};
\draw (5,1) node [draw] (l2) {};
\draw (5,-1) node [draw] (l2) {};
\draw (4,0) node [fill=black!100] (ld) {};
\draw (3,0) node [draw] (l2) {};
\draw (2,0) node [draw] (l2) {};
\draw (1,0) node [draw] (l4) {};
\draw (0,0) node [draw] (l4) {};
\end{tikzpicture}\end{center}
The single application of Theorem \ref{thm:perversebijchange} -- as for the case of $G_2$ -- is easy, and omitted. What is of interest here is the change from $\kappa$ to $\kappa+d$. The sizes of the sets $X_{\geq i}$ are
\[ |X_{\geq 26}|=13,\;\; |X_{\geq 38}|=12,\;\; |X_{\geq 48}|=10,\;\; |X_{\geq 52}|=8,\;\; |X_{\geq 55}|=6,\;\; |X_{\geq 59}|=4,\;\; |X_{\geq 61}|=2,\;\; |X_{\geq 63}|=1.\]
Starting from $X_{\geq 63}=\{\phi_{1,63}\}$, we see that all three conditions at the start of this section are satisfied until we reach the jump between $X_{\geq 55}$ and $X_{\geq 52}$, and $X_{\geq 52}$ and $X_{\geq 48}$. We start with the list $X_{\geq 55}$, namely
\[ \phi_{189,17}, \;\phi_{105,26},\;\phi_{27,37},\;\phi_{1,63},\;D_4;r\ep_2,\;D_4;\ep_2,\]
then rotate by three iterations of a $6$-cycle (from $X_{\geq 55}$ to $X_{\geq 52}$) to get
\[ \phi_{1,63},\;D_4;r\ep_2,\;D_4;\ep_2,\;\phi_{189,17},\; \phi_{105,26},\;\phi_{27,37},\]
and then insert $E_7[\pm\I]$ in their appropriate places, to get
\[ E_7[-\I],\;\phi_{1,63},\;D_4;r\ep_2,\;D_4;\ep_2,\;E_7[\I],\;\phi_{189,17},\; \phi_{105,26},\;\phi_{27,37},\]
finally rotating by four iterations of an $8$-cycle (from $X_{\geq 52}$ to $X_{\geq 48}$) to get
\[ E_7[\I],\;\phi_{189,17},\; \phi_{105,26},\;\phi_{27,37},\;E_7[-\I],\;\phi_{1,63},\;D_4;r\ep_2,\;D_4;\ep_2.\]
We see that the relative positions of all characters other than $E_7[\pm\I]$ are the same, and that the $E_7[\pm\I]$ are swapped. All subsequent differences in the sizes of the $X_{\geq i}$ satisfy the second condition, and so the change in the bijection when moving from $\kappa$ to $\kappa+d$ is to swap $E_7[\I]$ and $E_7[-\I]$, consistent with the change upon substitution $q\mapsto \e^{2\pi\I\kappa/d}$ to $q\mapsto \e^{2\pi\I(\kappa+d)/d}$.

\medskip

This completes the proof of the combinatorial form of Brou\'e's conjecture for a representative sample of unipotent blocks with cyclic defect group in exceptional groups of Lie type.

\bigskip\bigskip

\noindent\textbf{Acknowledgments} I would like to thank Gunter Malle for suggesting that I consider the $d$-cuspidal pair in a preliminary version of Conjecture \ref{conj:DLcohom}. I thank Olivier Dudas, both for many useful conversations and for allowing me access to his unpublished work. Jean Michel also performed several calculations with Deligne--Lusztig varieties that backed up Conjecture \ref{conj:DLcohom}, and more importantly helped me to straighten out the problems I was having with the `very twisted' Ree and Suzuki groups. I would finally like to thank Rapha\"el Rouquier for letting me bounce ideas off him, and for introducing me to the problem of determining the `geometric perversity function'.

\appendix

\bigskip\bigskip

\noindent\section*{Appendix: Unipotent Blocks of Weight 1 in Exceptional Groups}
In this appendix we collate all information, both known and conjectural, about unipotent blocks of weight $1$ in exceptional groups of Lie type. For each block $B$ we produce the list of unipotent characters belonging to $B$, their degrees, their eigenvalues of a suitable root/power of the Frobenius (i.e., the parameters of the cyclotomic Hecke algebra), its $d$-cuspidal pair $(\bo L,\bo\lambda)$ (and the generic degree of $\bo\lambda$), and the values of the $\pi_{\kappa/d}$-function for all appropriate $\kappa$ (recalling that the difference between $\pi_{(\zeta+d)/d}(-)$ and $\pi_{\zeta/d}(-)$ is $2A(-)$, and noting that we also give $A(-)$), together with the Brauer tree, either known or conjectured.

In each case, we either state that it is conjectured, or give a reference for its proof. For groups of types $G_2$, ${}^2\!G_2$, ${}^3\!D_4$, $F_4$, ${}^2\!F_4$, $E_6$ and ${}^2\!E_6$, all Brauer trees are in the published literature, together with the planar embedding (except for ${}^2\!E_6(q^2)$, $\ell\mid\Phi_{12}(q^2)$, and $q\nequiv 1\bmod 3$), and we simply give a reference at the start of the respective section. For many of the others, from $E_7$ and $E_8$, we use four particular arguments.
\begin{enumerate}
\item Geck--Pfeiffer argument: in \cite[Table F]{geckpfeiffer}, the blocks of positive defect of the Hecke algebra of various Lie-type groups are given, and one uses Table C of \cite{geckpfeiffer} to reconstruct the labelling of the principal series characters (the second argument in $\phi_{n,b}$ is the $b$-function of that table). The Brauer trees of the Hecke algebra are lines with the characters appearing in order of increasing $a$-function; these are subtrees of the Brauer tree of the relevant block. This allows us to construct the subtree consisting of the principal characters.
\item Degree argument: the degree (as a polynomial in $q$) should increase towards the exceptional node. This, together with a Geck--Pfeiffer argument, normally allows us to construct the entire real stem.
\item Morita argument: Harish-Chandra induction from a particular Levi subgroup will provide a Morita equivalence. When we use this argument, we will provide the particular Levi subgroup. In general, if Harish-Chandra induction from the unipotent characters $\chi$ of a block $b$ of a Levi subgroup $\bo L$ of $\bo G$, when cut by a block $B$ of $\bo G$, result in irreducible unipotent characters of $B$, then this induction produces a Morita equivalence from $b$ to $B$. (This requires the normalizers of the appropriate tori to coincide.) Indeed, the bijection between the labels of the (planar embedded) Brauer trees of $b$ and $B$ matches that of Harish-Chandra induction.
\item Projective-Induction argument: if $\chi$ is a projective character of a Levi subgroup with the same power of $\Phi_d$ as in the group itself, then the Harish-Chandra (in general, the Deligne--Lusztig) induced character will also be a projective. This means that it is the sum of adjacent vertices of the Brauer tree. This is used to isolate some of the non-real characters for $d=9$ and $E_8(q)$.
\end{enumerate}

Occasionally more specialized arguments are needed, and each section will contain those, but we will appeal to these four arguments without extra comment when they apply.

Of the remaining trees for $E_7$ and $E_8$, all but two at the present time, have been solved by Dudas, Rouquier and the author in \cite{cdr2012un}, and we will just reference that article when it applies. The only remaining unknown trees are for $E_8(q)$, with the principal block for $d=15$ and a single non-principal block for $d=18$.
%
%

\tikzstyle{every node}=[circle, fill=black!0,
                        inner sep=0pt, minimum width=4pt]

\newpage

\section{${}^2\!B_2(q^2)$}

The group ${}^2\!B_2(q^2)$ has order $q^4\Phi_1\Phi_2\Phi_8$. Since $\sqrt 2$ lies inside our field of definition, $\Phi_8$ factorizes as
\[ \Phi_8=(q^2+\sqrt2q+1)(q^2-\sqrt2+1).\]
Denote by $\Phi_8'$ the first and by $\Phi_8''$ the second. If $\psi=\e^{\pi\I/4}$, then $\Phi_8'=(q-\psi^3)(q-\psi^5)$ and $\Phi_8''=(q-\psi)(q-\psi^7)$. This is the same way as Carter defines them in \cite{carterfinite}. This way is consistent with $\Phi_d''$ being the Coxeter case, and having $\e^{2\pi\I/d}$ as a zero.

The Brauer trees for this group are for the principal block and $\ell$ dividing any torus. It should be noted that $\ell\mid (q^2-1)$, as $\Phi_1=\sqrt2^n-1$ and $\Phi_2=\sqrt2^n+1$ with $n$ odd, so an integer cannot divide either. In \cite{burkhardt1979}, the Brauer trees for ${}^2\!B_2(q)$ are determined.

There are four unipotent characters of $G$, with the following degrees.

\begin{center}\begin{tabular}{ll}
\hline Character & Degree
\\ \hline $\phi_{1,0}$ & $1$
\\ $\phi_{1,2}$ & $q^4$
\\ ${}^2\!B_2[\psi^3]$ & $q\Phi_1\Phi_2/\sqrt2$
\\ ${}^2\!B_2[\psi^5]$ & $q\Phi_1\Phi_2/\sqrt2$
\\\hline\end{tabular}\end{center}

\newpage

\subsection{$d=1$}

For Suzuki groups, $d=1$ corresponds to $(q^2-1)$, and we can either let $d=2$ or $d=1$ and consider $q^2$ instead of $q$. For $d=1$ there is a single unipotent block of weight $1$, and two unipotent blocks of defect zero.
\\[1cm]
\noindent(i)\;\; \textbf{Block 1:} Cuspidal pair is $(\Phi_1\Phi_2,1)$, of degree $1$. There are two unipotent characters in the block, both real.

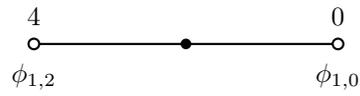
\begin{center}\begin{tabular}{lccc}
\hline Character & $A(-)$ & $\omega_i q^{aA/e}$ & $\kappa=1$
\\ \hline $\phi_{1,0}$ & $0$ & $1$ & $0$
\\ $\phi_{1,2}$ & $4$ & $-q^4$ & $4$
\\ \hline\end{tabular}\end{center}

\begin{center}\begin{tikzpicture}[thick,scale=2]
\draw (0,0) -- (2,0);
\draw (0,0) node [draw,label=below:$\phi_{1,2}$] (l0) {};
\draw (2,0) node [draw,label=below:$\phi_{1,0}$] (l2) {};
\draw (1,0) node [fill=black!100] (ld) {};
\draw (0,0.18) node{$4$};
\draw (2,0.18) node{$0$};
\end{tikzpicture}\end{center}

\newpage

\subsection{$d=8'$}

For ${}^2\!B_2(q^2)$ and $d=8'$ there is a single unipotent block of weight $1$, and no other unipotent blocks.
\\[1cm]
\noindent(i)\;\; \textbf{Block 1:} Cuspidal pair is $(\Phi_8',1)$, of degree $1$. There are four unipotent characters in the block, two of which are non-real.

\begin{center}\begin{tabular}{lcccc}
\hline Character & $A(-)$ & $\omega_i q^{aA/e}$ & $\kappa=3$ & $\kappa=5$
\\ \hline $\phi_{1,0}$ & $0$ & $1$ & $0$ & $0$
\\ ${}^2\!B_2[\psi^5]$ & $3$ & $\psi^7 q$ & $2$ & $4$
\\ ${}^2\!B_2[\psi^3]$ & $3$ & $\psi q$ & $2$ & $4$
\\ $\phi_{1,2}$ & $4$ & $q^2$ & $3$ & $5$
\\ \hline\end{tabular}\end{center}

\begin{center}\begin{tikzpicture}[thick,scale=2]
\draw (0,0.5) -- (-2,0.5);
\draw (-1,0) -- (-1,1);
\draw (0,0.5) node [draw,label=below:$\phi_{1,0}$] (l0) {};
\draw (-1,0.5) node [draw,label=below left:$\phi_{1,2}$] (l2) {};
\draw (-2,0.5) node [fill=black!100] (ld) {};
\draw (0,0.68) node{$0$};
\draw (-0.85,0.68) node{$3$};

\draw (-1,1) node [draw,label=left:${^2\!B_2[\psi^3]}$] (l0) {};
\draw (-1,0) node [draw,label=left:${^2\!B_2[\psi^5]}$] (l0) {};

\draw (-1,1.15) node{$2$};
\draw (-1,-0.15) node{$2$};
\end{tikzpicture}\end{center}

\newpage

\subsection{$d=8''$}

For ${}^2\!B_2(q^2)$ and $d=8''$ there is a single unipotent block of weight $1$, and no other unipotent blocks.
\\[1cm]
\noindent(i)\;\; \textbf{Block 1:} Cuspidal pair is $(\Phi_8'',1)$, of degree $1$. There are four unipotent characters in the block, two of which are non-real.

\begin{center}\begin{tabular}{lcccc}
\hline Character & $A(-)$ & $\omega_i q^{aA/e}$ & $\kappa=1$ & $\kappa=7$
\\ \hline $\phi_{1,0}$ & $0$ & $1$ & $0$ & $0$
\\ $\phi_{1,2}$ & $4$ & $q^2$ & $1$ & $7$
\\ ${}^2\!B_2[\psi^3]$ & $3$ & $\psi^3 q$ & $1$ & $5$
\\ ${}^2\!B_2[\psi^5]$ & $3$ & $\psi^5 q$ & $1$ & $5$
\\ \hline
\end{tabular}\end{center}

\begin{center}\begin{tikzpicture}[thick,scale=2]
\draw (0,0.5) -- (-2,0.5);
\draw (-2,0) -- (-2,1);
\draw (0,0.5) node [draw,label=below:$\phi_{1,0}$] (l0) {};
\draw (-1,0.5) node [draw,label=below:$\phi_{1,2}$] (l2) {};
\draw (-2,0.5) node [fill=black!100] (ld) {};
\draw (0,0.68) node{$0$};
\draw (-1,0.68) node{$1$};

\draw (-2,1) node [draw,label=left:${^2\!B_2[\psi^3]}$] (l0) {};
\draw (-2,0) node [draw,label=left:${^2\!B_2[\psi^5]}$] (l0) {};

\draw (-2,1.15) node{$1$};
\draw (-2,-0.15) node{$1$};
\end{tikzpicture}\end{center}

\newpage

\section{$G_2(q)$}
\label{sec:G2}

The group $G_2(q)$ has order $q^6\Phi_1^2\Phi_2^2\Phi_3\Phi_6$. The only Brauer trees for $G_2(q)$ are for the principal block and $d=3,6$. By \cite{shamash1989}, the Brauer trees for $G_2(q)$ are determined, along with the planar embedding.

There are ten unipotent characters of $G$, with the following degrees.

\begin{center}\begin{tabular}{ll}
\hline Character & Degree
\\ \hline $\phi_{1,0}$ & $1$
\\ $\phi_{2,1}$ & $q\Phi_2^2\Phi_3/6$
\\ $\phi_{2,2}$ & $q\Phi_2^2\Phi_6/2$
\\ $\phi_{1,3}'$ & $q\Phi_3\Phi_6/3$
\\ $\phi_{1,3}''$ & $q\Phi_3\Phi_6/3$
\\ $\phi_{1,6}$ & $q^6$
\\ $G_2[1]$ & $q\Phi_1^2\Phi_6/6$
\\ $G_2[-1]$ & $q\Phi_1^2\Phi_3/2$
\\ $G_2[\theta]$ & $q\Phi_1^2\Phi_2^2/3$
\\ $G_2[\theta^2]$ & $q\Phi_1^2\Phi_2^2/3$
\\\hline\end{tabular}\end{center}

\newpage

\subsection{$d=3$}

For $G_2(q)$ and $d=3$ there is a single unipotent block of weight $1$, together with four unipotent blocks of defect zero.
\\[1cm]
\noindent(i)\;\; \textbf{Block 1:} Cuspidal pair is $(\Phi_3,1)$, of degree $1$. There are six unipotent characters in the block, two of which are non-real.

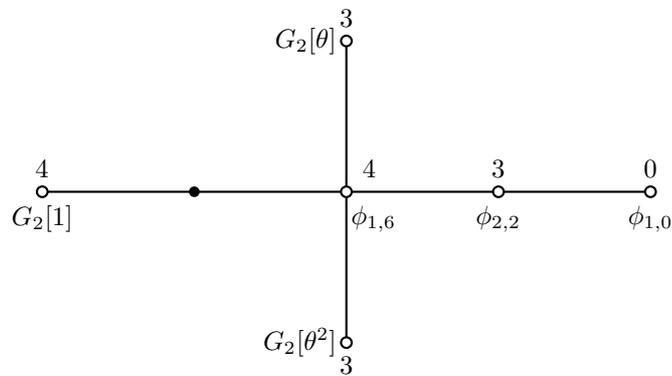
\begin{center}\begin{tabular}{lcccc}
\hline Character & $A(-)$ & $\omega_i q^{aA/e}$ & $\kappa=1$ & $\kappa=2$
\\\hline $\phi_{1,0}$ & $0$ & $q^{0}$ & $0$ & $0$
\\ $G_2[\theta^2]$ & $5$ & $-\theta q$ & $3$ & $7$
\\ $\phi_{2,2}$ & $5$ & $q$ & $3$ & $7$
\\ $G_2[\theta]$ & $5$ & $-\theta^2 q$ & $3$ & $7$
\\ $\phi_{1,6}$ & $6$ & $q^{2}$ & $4$ & $8$
\\ $G_2[1]$ & $5$ & $-q$ & $4$ & $6$
\\\hline\end{tabular}\end{center}

\begin{center}\begin{tikzpicture}[thick,scale=2]
\draw (0,0.32) node {$\phi_{1,0}$};
\draw (-1,0.32) node {$\phi_{2,2}$};
\draw (-1.82,0.32) node {$\phi_{1,6}$};
\draw (-4,0.32) node {$G_2[1]$};

\draw (0,0.5) -- (-4,0.5);
\draw (-2,-0.5) -- (-2,1.5);

\draw (0,0.5) node [draw] (l0) {};
\draw (-1,0.5) node [draw] (l1) {};
\draw (-4,0.5) node [draw] (l2) {};
\draw (-3,0.5) node [fill=black!100] (ld) {};
\draw (-2,-0.5) node [draw,label=left:${G_2[\theta^2]}$] (l4) {};
\draw (-2,0.5) node [draw] (l4) {};
\draw (-2,1.5) node [draw,label=left:${G_2[\theta]}$] (l4) {};
\draw (0,0.65) node{$0$};
\draw (-1,0.65) node{$3$};
\draw (-1.85,0.65) node{$4$};
\draw (-4,0.65) node{$4$};
\draw (-2,1.65) node{$3$};
\draw (-2,-0.65) node{$3$};
\end{tikzpicture}\end{center}
\newpage

\subsection{$d=6$}

For $G_2(q)$ and $d=6$ there is a single unipotent block of weight $1$, together with four unipotent blocks of defect zero.
\\[1cm]
\noindent(i)\;\; \textbf{Block 1:} Cuspidal pair is $(\Phi_6,1)$, of degree $1$. There are six unipotent characters in the block, two of which are non-real.

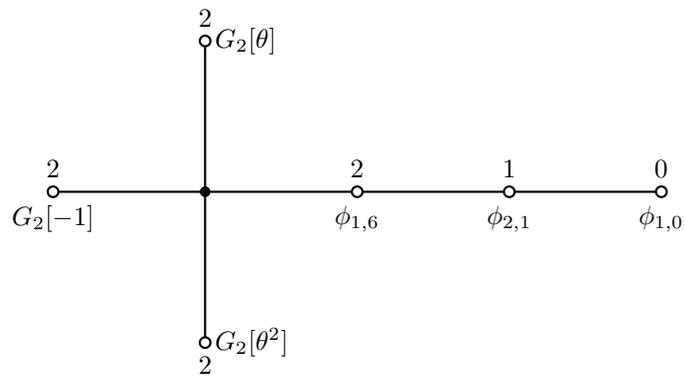
\begin{center}\begin{tabular}{lcccc}
\hline Character & $A(-)$ & $\omega_i q^{aA/e}$ & $\kappa=1$ & $\kappa=5$
\\\hline $\phi_{1,0}$ & $0$ & $q^{0}$ & $0$ & $0$
\\ $\phi_{2,1}$ & $5$ & $q$ & $1$ & $9$
\\ $\phi_{1,6}$ & $6$ & $q^{2}$ & $2$ & $10$
\\ $G_2[\theta]$ & $5$ & $\theta q$ & $2$ & $8$
\\ $G_2[-1]$ & $5$ & $-q$ & $2$ & $8$
\\ $G_2[\theta^2]$ & $5$ & $\theta^2q$ & $2$ & $8$
\\\hline\end{tabular}\end{center}

\begin{center}\begin{tikzpicture}[thick,scale=2]
\draw (0,0.32) node {$\phi_{1,0}$};
\draw (-1,0.32) node {$\phi_{2,1}$};
\draw (-2,0.32) node {$\phi_{1,6}$};
\draw (-4,0.32) node {$G_2[-1]$};

\draw (0,0.5) -- (-4,0.5);
\draw (-3,-0.5) -- (-3,1.5);

\draw (0,0.5) node [draw] (l0) {};
\draw (-1,0.5) node [draw] (l1) {};
\draw (-2,0.5) node [draw] (l1) {};
\draw (-4,0.5) node [draw] (l2) {};
\draw (-3,0.5) node [fill=black!100] (ld) {};
\draw (-3,-0.5) node [draw,label=right:${G_2[\theta^2]}$] (l4) {};
\draw (-3,1.5) node [draw,label=right:${G_2[\theta]}$] (l4) {};

\draw (-0,0.65) node{$0$};
\draw (-1,0.65) node{$1$};
\draw (-2,0.65) node{$2$};
\draw (-4,0.65) node{$2$};
\draw (-3,1.65) node{$2$};
\draw (-3,-0.65) node{$2$};
\end{tikzpicture}\end{center}

\newpage

\section{${}^2\!G_2(q^2)$}

The group ${}^2\!G_2(q^2)$, for $q$ an odd power of $\sqrt3$,  has order $q^6\Phi_1\Phi_2\Phi_4\Phi_{12}$. Since $\sqrt 3$ lies inside our field of definition, $\Phi_{12}$ factorizes as
\[ \Phi_{12}=(q^2+\sqrt3q+1)(q^2-\sqrt3+1).\]
Denote by $\Phi_{12}'$ the first and by $\Phi_{12}''$ the second. If $\xi=\e^{\pi\I/6}$, then $\Phi_{12}'=(q-\xi^5)(q-\xi^7)$ and $\Phi_{12}''=(q-\xi)(q-\xi^{11})$. \textbf{This is the opposite of how Carter defines them in \cite{carterfinite}.} This way is consistent with $\Phi_d''$ being the Coxeter case, and having $\e^{2\pi\I/d}$ as a zero.

The Brauer trees for this group are for the principal block and $\ell$ dividing any torus. It should be noted that $\ell\mid (q^2-1)$, as $\Phi_1=\sqrt3^n-1$ and $\Phi_2=\sqrt3^n+1$ with $n$ odd, so an integer cannot divide either. In \cite{hiss1991}, the Brauer trees for ${}^2\!G_2(q)$ are partially determined, but the planar embedding was not determined for $\ell\mid \Phi_{12}''$: this case was completed by Dudas in \cite{dudas2010un2}.

There are eight unipotent characters for ${}^2\!G_2(q^2)$, with degrees given below.

\begin{center}\begin{tabular}{ll}
\hline Name & Degree
\\ \hline $\phi_{1,0}$ & $1$
\\ $\phi_{1,2}$ & $q^6$
\\ ${}^2\!G_2[\xi^5]$ & $q\Phi_1\Phi_2\Phi_4/\sqrt3$
\\ ${}^2\!G_2[\xi^7]$ & $q\Phi_1\Phi_2\Phi_4/\sqrt3$
\\ ${}^2\!G_2^\mathrm{I}[\I]$ & $q\Phi_1\Phi_2\Phi_{12}'/2\sqrt3$
\\ ${}^2\!G_2^\mathrm{I}[-\I]$ & $q\Phi_1\Phi_2\Phi_{12}'/2\sqrt3$
\\ ${}^2\!G_2^\mathrm{II}[\I]$ & $q\Phi_1\Phi_2\Phi_{12}''/2\sqrt3$
\\ ${}^2\!G_2^\mathrm{II}[-\I]$ & $q\Phi_1\Phi_2\Phi_{12}''/2\sqrt3$
\\\hline\end{tabular}\end{center}

\newpage

\subsection{$d=1$}

For the small Ree groups, $d=1$ corresponds to $(q^2-1)$, and we can either let $d=2$ or $d=1$ and consider $q^2$ instead of $q$. For $d=1$ there is a single unipotent block of weight $1$, and six unipotent blocks of defect zero.
\\[1cm]
\noindent(i)\;\; \textbf{Block 1:} Cuspidal pair is $(\Phi_1\Phi_2,1)$, of degree $1$. There are two unipotent characters in the block, both real.

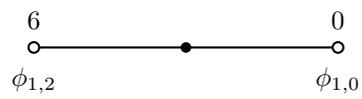
\begin{center}\begin{tabular}{lccc}
\hline Character & $A(-)$ & $\omega_i q^{aA/e}$ & $\kappa=1$
\\ \hline $\phi_{1,0}$ & $0$ & $1$ & $0$
\\ $\phi_{1,2}$ & $6$ & $-q^6$ & $6$
\\ \hline\end{tabular}\end{center}

\begin{center}\begin{tikzpicture}[thick,scale=2]
\draw (0,0) -- (2,0);
\draw (0,0) node [draw,label=below:$\phi_{1,2}$] (l0) {};
\draw (2,0) node [draw,label=below:$\phi_{1,0}$] (l2) {};
\draw (1,0) node [fill=black!100] (ld) {};
\draw (0,0.18) node{$6$};
\draw (2,0.18) node{$0$};
\end{tikzpicture}\end{center}

\newpage

\subsection{$d=4$}

For ${}^2\!G_2(q^2)$ and $d=4$ there is a single unipotent block of weight $1$, together with four unipotent blocks of defect zero.
\\[1cm]
\noindent(i)\;\; \textbf{Block 1:} Cuspidal pair is $(\Phi_4,1)$, of degree $1$. There are six unipotent characters in the block, four of which are non-real.

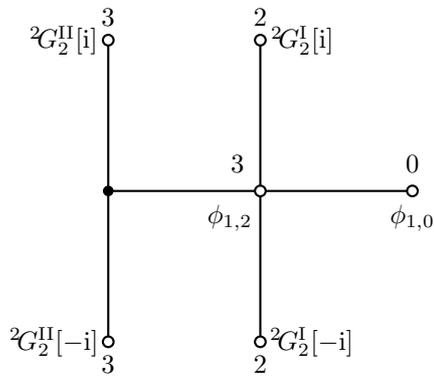
\begin{center}\begin{tabular}{lcccc}
\hline Character & $A(-)$ & $\omega_i q^{aA/e}$ & $\kappa=1$ & $\kappa=3$
\\ \hline $\phi_{1,0}$ & $0$ & $1$ & $0$ & $0$
\\ ${}^2\!G_2^{\mathrm{I}}[-\I]$ & $5$ & $\xi^{11}q$ & $2$ & $8$
\\ ${}^2\!G_2^{\mathrm{I}}[\I]$ & $5$ & $\xi q$ & $2$ & $8$
\\ $\phi_{1,2}$ & $6$ & $q^2$ & $3$ & $9$
\\ ${}^2\!G_2^{\mathrm{II}}[\I]$ & $5$ & $\xi^5 q$ & $3$ & $7$
\\ ${}^2\!G_2^{\mathrm{II}}[-\I]$ & $5$ & $\xi^7 q$ & $3$ & $7$
\\ \hline
\end{tabular}\end{center}

\begin{center}\begin{tikzpicture}[thick,scale=2]
\draw (0,0.32) node {$\phi_{1,0}$};
\draw (-1.2,0.32) node {$\phi_{1,2}$};

\draw (0,0.5) -- (-2,0.5);
\draw (-1,-0.5) -- (-1,1.5);
\draw (-2,-0.5) -- (-2,1.5);
\draw (0,0.5) node [draw] (l0) {};
\draw (-1,0.5) node [draw] (l2) {};
\draw (-2,0.5) node [fill=black!100] (ld) {};
\draw (0,0.68) node{$0$};
\draw (-1.15,0.68) node{$3$};

\draw (-1,1.5) node [draw,label=right:${{}^2\!G_2^\mathrm{I}[\I]}$] (l0) {};
\draw (-1,-0.5) node [draw,label=right:${{}^2\!G_2^\mathrm{I}[-\I]}$] (l0) {};

\draw (-2,1.5) node [draw,label=left:${{}^2\!G_2^\mathrm{II}[\I]}$] (l0) {};
\draw (-2,-0.5) node [draw,label=left:${{}^2\!G_2^\mathrm{II}[-\I]}$] (l0) {};

\draw (-1,1.65) node{$2$};
\draw (-1,-0.65) node{$2$};
\draw (-2,1.65) node{$3$};
\draw (-2,-0.65) node{$3$};
\end{tikzpicture}\end{center}

\newpage

\subsection{$d=12'$}

For ${}^2\!G_2(q^2)$ and $d=12'$ there is a single unipotent block of weight $1$, together with two unipotent blocks of defect zero.
\\[1cm]
\noindent(i)\;\; \textbf{Block 1:} Cuspidal pair is $(\Phi_{12}',1)$, of degree $1$. There are six unipotent characters in the block, four of which are non-real.

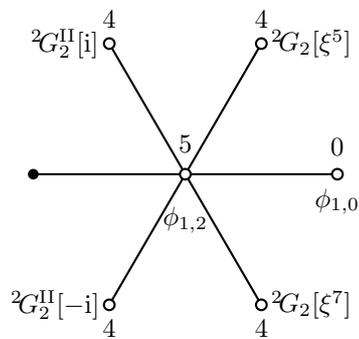
\begin{center}\begin{tabular}{lcccc}
\hline Character & $A(-)$ & $\omega_i q^{aA/e}$ & $\kappa=1$ & $\kappa=11$
\\ \hline $\phi_{1,0}$ & $0$ & $1$ & $0$ & $0$
\\ ${}^2\!G_2^{\mathrm{II}}[-\I]$ & $5$ & $-\I q$ & $4$ & $6$
\\ ${}^2\!G_2[\xi^7]$ & $5$ & $\xi^{11} q$ & $4$ & $6$
\\ ${}^2\!G_2[\xi^5]$ & $5$ & $\xi q$ & $4$ & $6$
\\ ${}^2\!G_2^{\mathrm{II}}[\I]$ & $5$ & $\I q$ & $4$ & $6$
\\ $\phi_{1,2}$ & $6$ & $q^2$ & $5$ & $7$
\\ \hline
\end{tabular}\end{center}

\begin{center}\begin{tikzpicture}[thick,scale=2]
\draw (0,0.82) node {$\phi_{1,0}$};
\draw (-1,0.7) node{$\phi_{1,2}$};

\draw (0,1) -- (-2,1);
\draw (-0.5,0.1339) -- (-1.5,1.866);
\draw (-0.5,1.866) -- (-1.5,0.1339);
\draw (0,1) node [draw] (l0) {};
\draw (-1,1) node [draw] (l2) {};
\draw (-2,1) node [fill=black!100] (ld) {};
\draw (0,1.18) node{$0$};
\draw (-1,1.2) node{$5$};

\draw (-0.5,1.866) node [draw,label=right:${{}^2\!G_2[\xi^5]}$] (l0) {};
\draw (-1.5,1.866) node [draw,label=left:${{}^2\!G_2^\mathrm{II}[\I]}$] (l0) {};
\draw (-0.5,0.1339) node [draw,label=right:${{}^2\!G_2[\xi^7]}$] (l0) {};
\draw (-1.5,0.1339) node [draw,label=left:${{}^2\!G_2^\mathrm{II}[-\I]}$] (l0) {};

\draw (-0.5,2.03) node{$4$};
\draw (-1.5,2.03) node{$4$};
\draw (-0.5,-0.02) node{$4$};
\draw (-1.5,-0.02) node{$4$};
\end{tikzpicture}
\end{center}

\newpage

\subsection{$d=12''$}

For ${}^2\!G_2(q^2)$ and $d=12''$ there is a single unipotent block of weight $1$, together with two unipotent blocks of defect zero.
\\[1cm]
\noindent(i)\;\; \textbf{Block 1:} Cuspidal pair is $(\Phi_{12}'',1)$, of degree $1$. There are six unipotent characters in the block, four of which are non-real.

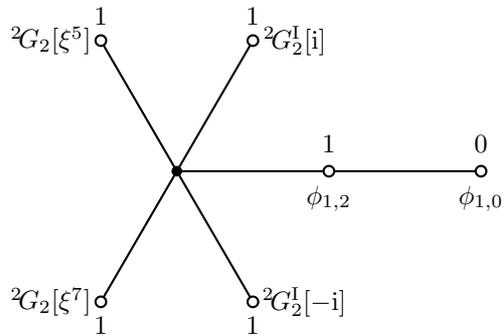
\begin{center}\begin{tabular}{lcccc}
\hline Character & $A(-)$ & $\omega_i q^{aA/e}$ & $\kappa=5$ & $\kappa=7$
\\ \hline $\phi_{1,0}$ & $0$ & $1$ & $0$ & $0$
\\ $\phi_{1,2}$ & $6$ & $q^2$ & $1$ & $11$
\\ ${}^2\!G_2^{\mathrm{I}}[\I]$ & $5$ & $\I q$ & $1$ & $9$
\\ ${}^2\!G_2[\xi^5]$ & $5$ & $\xi^5 q$ & $1$ & $9$
\\ ${}^2\!G_2[\xi^7]$ & $5$ & $\xi^7 q$ & $1$ & $9$
\\ ${}^2\!G_2^{\mathrm{I}}[-\I]$ & $5$ & $-\I q$ & $1$ & $9$
\\ \hline
\end{tabular}\end{center}
\begin{center}\begin{tikzpicture}[thick,scale=2]
\draw (1,0.82) node {$\phi_{1,0}$};
\draw (0,0.82) node{$\phi_{1,2}$};

\draw (-1,1) -- (1,1);
\draw (-0.5,0.1339) -- (-1.5,1.866);
\draw (-0.5,1.866) -- (-1.5,0.1339);
\draw (1,1) node [draw] (l0) {};
\draw (0,1) node [draw] (l0) {};
\draw (-1,1) node [fill=black!100] (ld) {};
\draw (1,1.18) node{$0$};
\draw (0,1.18) node{$1$};

\draw (-1.5,1.866) node [draw,label=left:${{}^2\!G_2[\xi^5]}$] (l0) {};
\draw (-0.5,1.866) node [draw,label=right:${{}^2\!G_2^\mathrm{I}[\I]}$] (l0) {};
\draw (-1.5,0.1339) node [draw,label=left:${{}^2\!G_2[\xi^7]}$] (l0) {};
\draw (-0.5,0.1339) node [draw,label=right:${{}^2\!G_2^\mathrm{I}[-\I]}$] (l0) {};

\draw (-0.5,2.03) node{$1$};
\draw (-1.5,2.03) node{$1$};
\draw (-0.5,-0.02) node{$1$};
\draw (-1.5,-0.02) node{$1$};
\end{tikzpicture}\end{center}

\newpage

\section{${}^3\!D_4(q^3)$}

The group ${}^3\!D_4(q^3)$ has order $q^{12}\Phi_1^2\Phi_2^2\Phi_3^2\Phi_6^2\Phi_{12}$.

The only Brauer tree for this group is for the principal block and $\ell$ dividing $\Phi_{12}$, given in \cite{geck1991}.

There are eight unipotent characters for ${}^3\!D_4(q^3)$, with degrees given below.

\begin{center}\begin{tabular}{ll}
\hline Name & Degree
\\ \hline $\phi_{1,0}$ & $1$
\\ $\phi_{1,6}$ & $q^{12}$
\\ $\phi_{2,1}$ & $q^3\Phi_2^2\Phi_6^2/2$
\\ $\phi_{2,2}$ & $q^3\Phi_2^2\Phi_{12}/2$
\\ $\phi_{1,3}'$ & $q\Phi_{12}$
\\ $\phi_{1,3}''$ & $q^7\Phi_{12}$
\\ ${}^3\!D_4[1]$ & $q^3\Phi_1^2\Phi_{12}/2$
\\ ${}^3\!D_4[-1]$ & $q^3\Phi_1^2\Phi_3^2/2$
\\\hline\end{tabular}\end{center}

\newpage

\subsection{$d=12$}

For ${}^3\!D_4(q^3)$ and $d=12$ there is a single unipotent block of weight $1$, together with two unipotent blocks of defect zero.
\\[1cm]
\noindent(i)\;\; \textbf{Block 1:} Cuspidal pair is $(\Phi_{12},1)$, of degree $1$. There are four unipotent characters in the block, all of which are real.

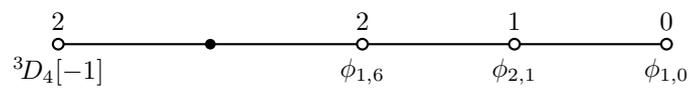
\begin{center}\begin{tabular}{lcccccc}
\hline Character & $A(-)$ & $\omega_i q^{aA/e}$ & $\kappa=1$ & $\kappa=5$ & $\kappa=7$ & $\kappa=11$
\\\hline $\phi_{1,0}$ & $0$ & $q^{0}$ & $0$ & $0$ & $0$ & $0$
\\ $\phi_{2,1}$ & $9$ & $q^3$ & $1$ & $7$ & $11$ & $17$
\\ $\phi_{1,6}$ & $12$ & $q^6$ & $2$ & $10$ & $14$ & $22$
\\ ${}^3\!D_4[-1]$ & $9$ & $-q^3$ & $2$ & $8$ & $10$ & $16$
\\ \hline
\end{tabular}\end{center}

\begin{center}\begin{tikzpicture}[thick,scale=2]
\draw (4,-0.18) node {$\phi_{1,0}$};
\draw (3,-0.18) node {$\phi_{2,1}$};
\draw (2,-0.18) node {$\phi_{1,6}$};
\draw (0,-0.18) node {${}^3\!D_4[-1]$};

\draw \foreach \x in {0,1,2,3}{
(\x,0) -- (\x+1,0)};
\draw (4,0) node [draw] (l0) {};
\draw (3,0) node [draw] (l0) {};
\draw (2,0) node [draw] (l0) {};
\draw (0,0) node [draw] (l0) {};
\draw (1,0) node [fill=black!100] (ld) {};
\draw (4,0.15) node{$0$};
\draw (3,0.15) node{$1$};
\draw (2,0.15) node{$2$};
\draw (0,0.15) node{$2$};
\end{tikzpicture}\end{center}

\newpage

\section{$F_4(q)$}
\vspace{-0.05cm}
The group $F_4(q)$ has order $q^{24}\Phi_1^4\Phi_2^4\Phi_3^3\Phi_4^2\Phi_6^2\Phi_8\Phi_{12}$. There are Brauer trees for $\Phi_8$ and $\Phi_{12}$ obviously, and also non-principal unipotent blocks of weight $1$ when $d=4$. The Brauer trees were determined in \cite{hisslubeck1998} in almost all cases: the planar embedding for $d=12$ was not given there, and was constructed by Dudas in \cite{dudas2010un2}.

There are thirty-seven unipotent characters of $F_4(q)$, with degrees given below.
\begin{center}\begin{tabular}{ll}
\hline Name & Degree
\\ \hline $\phi_{1,0}$ & $1$
\\ $\phi_{1,12}''$ & $q^{4}\Phi_{4}^{2}\Phi_{6}^{2}\Phi_{8}\Phi_{12}/8$
\\ $\phi_{1,12}'$ & $q^{4}\Phi_{4}^{2}\Phi_{6}^{2}\Phi_{8}\Phi_{12}/8$
\\ $\phi_{1,24}$ & $q^{24}$
\\ $\phi_{2,4}''$ & $q\Phi_{4}\Phi_{8}\Phi_{12}/2$
\\ $\phi_{2,16}'$ & $q^{13}\Phi_{4}\Phi_{8}\Phi_{12}/2$
\\ $\phi_{2,4}'$ & $q\Phi_{4}\Phi_{8}\Phi_{12}/2$
\\ $\phi_{2,16}''$ & $q^{13}\Phi_{4}\Phi_{8}\Phi_{12}/2$
\\ $\phi_{4,8}$ & $q^{4}\Phi_{2}^{4}\Phi_{6}^{2}\Phi_{8}\Phi_{12}/8$
\\ $\phi_{9,2}$ & $q^{2}\Phi_{3}^{2}\Phi_{6}^{2}\Phi_{12}$
\\ $\phi_{9,6}''$ & $q^{4}\Phi_{3}^{2}\Phi_{4}^{2}\Phi_{8}\Phi_{12}/8$
\\ $\phi_{9,6}'$ & $q^{4}\Phi_{3}^{2}\Phi_{4}^{2}\Phi_{8}\Phi_{12}/8$
\\ $\phi_{9,10}$ & $q^{10}\Phi_{3}^{2}\Phi_{6}^{2}\Phi_{12}$
\\ $\phi_{6,6}'$ & $q^{4}\Phi_{3}^{2}\Phi_{6}^{2}\Phi_{8}\Phi_{12}/3$
\\ $\phi_{6,6}''$ & $q^{4}\Phi_{3}^{2}\Phi_{4}^{2}\Phi_{6}^{2}\Phi_{8}/12$
\\ $\phi_{12,4}$ & $q^{4}\Phi_{2}^{4}\Phi_{3}^{2}\Phi_{8}\Phi_{12}/24$
\\ $\phi_{4,1}$ & $q\Phi_{2}^{2}\Phi_{6}^{2}\Phi_{8}/2$
\\ $\phi_{4,7}''$ & $q^{4}\Phi_{2}^{2}\Phi_{4}\Phi_{6}^{2}\Phi_{8}\Phi_{12}/4$
\\ $\phi_{4,7}'$ & $q^{4}\Phi_{2}^{2}\Phi_{4}\Phi_{6}^{2}\Phi_{8}\Phi_{12}/4$
\\ $\phi_{4,13}$ & $q^{13}\Phi_{2}^{2}\Phi_{6}^{2}\Phi_{8}/2$
\\ $\phi_{8,3}''$ & $q^{3}\Phi_{4}^{2}\Phi_{8}\Phi_{12}$
\\ $\phi_{8,9}'$ & $q^{9}\Phi_{4}^{2}\Phi_{8}\Phi_{12}$
\\ $\phi_{8,3}'$ & $q^{3}\Phi_{4}^{2}\Phi_{8}\Phi_{12}$
\\ $\phi_{8,9}''$ & $q^{9}\Phi_{4}^{2}\Phi_{8}\Phi_{12}$
\\ $\phi_{16,5}$ & $q^{4}\Phi_{2}^{4}\Phi_{4}^{2}\Phi_{6}^{2}\Phi_{12}/4$
\\ $B_2;\ep$ & $q^{13}\Phi_{1}^{2}\Phi_{3}^{2}\Phi_{8}/2$
\\ $B_2;r$ & $q^{4}\Phi_{1}^{2}\Phi_{2}^{2}\Phi_{3}^{2}\Phi_{6}^{2}\Phi_{8}/4$
\\ $B_2;\ep''$ & $q^{4}\Phi_{1}^{2}\Phi_{3}^{2}\Phi_{4}\Phi_{8}\Phi_{12}/4$
\\ $B_2;\ep'$ & $q^{4}\Phi_{1}^{2}\Phi_{3}^{2}\Phi_{4}\Phi_{8}\Phi_{12}/4$
\\ $B_2;1$ & $q\Phi_{1}^{2}\Phi_{3}^{2}\Phi_{8}/2$
\\ $F_4[-1]$ & $q^{4}\Phi_{1}^{4}\Phi_{3}^{2}\Phi_{4}^{2}\Phi_{12}/4$
\\ $F_4[-\I]$ & $q^{4}\Phi_{1}^{4}\Phi_{2}^{4}\Phi_{3}^{2}\Phi_{6}^{2}/4$
\\ $F_4[\I]$ & $q^{4}\Phi_{1}^{4}\Phi_{2}^{4}\Phi_{3}^{2}\Phi_{6}^{2}/4$
\\ $F_4[\theta]$ & $q^{4}\Phi_{1}^{4}\Phi_{2}^{4}\Phi_{4}^{2}\Phi_{8}/3$
\\ $F_4[\theta^2]$ & $q^{4}\Phi_{1}^{4}\Phi_{2}^{4}\Phi_{4}^{2}\Phi_{8}/3$
\\ $F_4^\mathrm{I}[1]$ & $q^{4}\Phi_{1}^{4}\Phi_{3}^{2}\Phi_{8}\Phi_{12}/8$
\\ $F_4^\mathrm{II}[1]$ & $q^{4}\Phi_{1}^{4}\Phi_{6}^{2}\Phi_{8}\Phi_{12}/24$
\\\hline\end{tabular}\end{center}

\newpage
\subsection{$d=4$}

For $F_4(q)$ and $d=4$ there are two unipotent blocks of weight $1$, together with the principal block (cyclotomic Weyl group $G_8$) and thirteen unipotent blocks of defect zero.
\\[1cm]
\noindent(i)\;\; \textbf{Block 1:} Cuspidal pair is $(\Phi_4.B_2(q),\phi_{11,-})$, of degree $q\Phi_4/2$. There are four unipotent characters in the block, all of which are real.

\begin{center}\begin{tabular}{lcccc}
\hline Character & $A(-)$ & $\omega_i q^{aA/e}$ & $\kappa=1$ & $\kappa=3$
\\\hline $\phi_{2,4}'$ & $8$ & $q^{2}$ & $4$ & $12$
\\ $B_2;\ep'$ & $17$ & $-q^{5}$ & $8$ & $26$
\\ $\phi_{2,16}'$ & $20$ & $q^{8}$ & $10$ & $30$
\\ $\phi_{4,7}'$ & $17$ & $q^{5}$ & $9$ & $25$
\\\hline\end{tabular}\end{center}

\begin{center}\begin{tikzpicture}[thick,scale=2]
\draw (0,0) -- (-4,0);
\draw (0,0) node [draw,label=below:$\phi_{2,4}'$] (l0) {};
\draw (-1,0) node [draw,label=below:$\phi_{4,7}'$] (l1) {};
\draw (-2,0) node [draw,label=below:$\phi_{2,16}'$] (l2) {};
\draw (-3,0) node [fill=black!100] (ld) {};
\draw (-4,0) node [draw,label=below:${B_2;\ep'}$] (l4) {};
\draw (0,0.18) node{$4$};
\draw (-1,0.18) node{$9$};
\draw (-2,0.18) node{$10$};
\draw (-4,0.18) node{$8$};
\end{tikzpicture}\end{center}

\noindent(ii)\;\; \textbf{Block 2:} Cuspidal pair is $(\Phi_4.B_2(q),\phi_{-,2})$, of degree $q\Phi_4/2$. There are four unipotent characters in the block, all of which are real.

\begin{center}\begin{tabular}{lcccc}
\hline Character & $A(-)$ & $\omega_i q^{aA/e}$ & $\kappa=1$ & $\kappa=3$
\\\hline $\phi_{2,4}''$ & $8$ & $q^{2}$ & $4$ & $12$
\\ $B_2;\ep''$ & $17$ & $-q^{5}$ & $8$ & $26$
\\ $\phi_{2,16}''$ & $20$ & $q^{8}$ & $10$ & $30$
\\ $\phi_{4,7}''$ & $17$ & $q^{5}$ & $9$ & $25$
\\\hline\end{tabular}\end{center}

\begin{center}\begin{tikzpicture}[thick,scale=2]
\draw (0,0) -- (-4,0);
\draw (0,0) node [draw,label=below:$\phi_{2,4}''$] (l0) {};
\draw (-1,0) node [draw,label=below:$\phi_{4,7}''$] (l1) {};
\draw (-2,0) node [draw,label=below:$\phi_{2,16}''$] (l2) {};
\draw (-3,0) node [fill=black!100] (ld) {};
\draw (-4,0) node [draw,label=below:${B_2;\ep''}$] (l4) {};
\draw (0,0.18) node{$4$};
\draw (-1,0.18) node{$9$};
\draw (-2,0.18) node{$10$};
\draw (-4,0.18) node{$8$};
\end{tikzpicture}\end{center}

\newpage
\subsection{$d=8$}

For $F_4(q)$ and $d=8$ there is a single unipotent block of weight $1$, and twenty-nine unipotent blocks of defect zero.
\\[1cm]
\noindent\textbf{Block 1:} Cuspidal pair is $(\Phi_8,1)$, of degree $1$. There are eight unipotent characters in the block, two of which -- $F_4[\pm\I]$ -- are non-real.

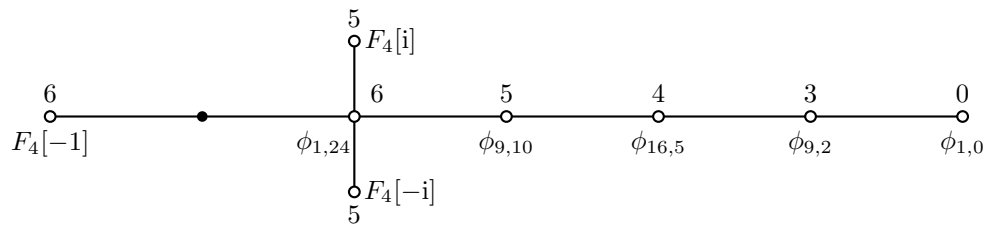
\begin{center}\begin{tabular}{lcccccc}
\hline Character & $A(-)$ & $\omega_i q^{aA/e}$ & $\kappa=1$ & $\kappa=3$ & $\kappa=5$ & $\kappa=7$
\\ \hline$\phi_{1,0}$ & $0$ & $q^{0}$ & $0$ & $0$ & $0$ & $0$
\\ $F_4[-\I]$ & $20$ & $-\I q^{3}$ & $5$ & $15$ & $25$ & $35$
\\ $\phi_{9,2}$ & $14$ & $q^{2}$ & $3$ & $11$ & $17$ & $25$
\\ $\phi_{16,5}$ & $20$ & $q^{3}$ & $4$ & $14$ & $26$ & $36$
\\ $\phi_{9,10}$ & $22$ & $q^{4}$ & $5$ & $17$ & $27$ & $39$
\\ $F_4[\I]$ & $20$ & $\I q^{3}$ & $5$ & $15$ & $25$ & $35$
\\ $\phi_{1,24}$ & $24$ & $q^{6}$ & $6$ & $18$ & $30$ & $42$
\\ $F_4[-1]$ & $20$ & $-q^{3}$ & $6$ & $16$ & $24$ & $34$
\\\hline\end{tabular}\end{center}

\begin{center}\begin{tikzpicture}[thick,scale=2]
\draw (0,0.32) node {$\phi_{1,0}$};
\draw (-1,0.32) node {$\phi_{9,2}$};
\draw (-2,0.32) node {$\phi_{16,5}$};
\draw (-3,0.32) node {$\phi_{9,10}$};
\draw (-4.2,0.32) node {$\phi_{1,24}$};

\draw (-6,0.32) node {$F_4[-1]$};

\draw (0,0.5) -- (-6,0.5);
\draw (-4,0) -- (-4,1);

\draw (0,0.5) node [draw] (l0) {};
\draw (-1,0.5) node [draw] (l1) {};
\draw (-2,0.5) node [draw] (l2) {};
\draw (-3,0.5) node [draw] (l2) {};
\draw (-4,0) node [draw,label=right:${F_4[-\I]}$] (l4) {};
\draw (-4,0.5) node [draw] (l4) {};
\draw (-4,1) node [draw,label=right:${F_4[\I]}$] (l4) {};
\draw (-5,0.5) node [fill=black!100] (ld) {};
\draw (-6,0.5) node [draw] (l2) {};
\draw (-0,0.65) node{$0$};
\draw (-1,0.65) node{$3$};
\draw (-2,0.65) node{$4$};
\draw (-3,0.65) node{$5$};
\draw (-3.85,0.65) node{$6$};
\draw (-4,1.15) node{$5$};
\draw (-4,-0.15) node{$5$};
\draw (-6,0.65) node{$6$};
\end{tikzpicture}\end{center}

\newpage

\subsection{$d=12$}

For $F_4(q)$ and $d=12$ there is a single unipotent block of weight $1$, and twenty-five unipotent blocks of defect zero.
\\[1cm]
\noindent\textbf{Block 1:} Cuspidal pair is $(\Phi_{12},1)$, of degree $1$. There are twelve unipotent characters in the block, four of which -- $F_4[\pm\I]$ and $F_4[\theta^i]$ -- are non-real.

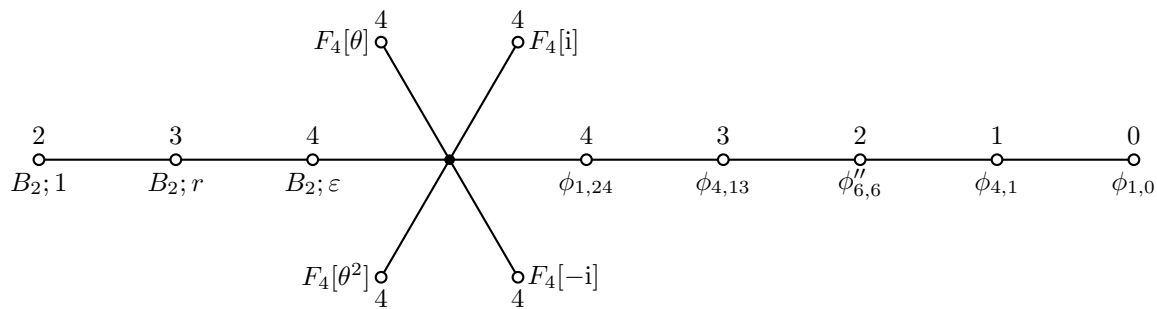
\begin{center}\begin{tabular}{lcccccc}
\hline Character & $A(-)$ & $\omega_i q^{aA/e}$ & $\kappa=1$ & $\kappa=5$ & $\kappa=7$ & $\kappa=11$
\\\hline $\phi_{1,0}$ & $0$ & $q^{0}$ & $0$ & $0$ & $0$ & $0$
\\ $\phi_{4,1}$ & $11$ & $q$ & $1$ & $9$ & $13$ & $21$
\\ $\phi_{6,6}''$ & $20$ & $q^{2}$ & $2$ & $18$ & $22$ & $38$
\\ $\phi_{4,13}$ & $23$ & $q^{3}$ & $3$ & $19$ & $27$ & $43$
\\ $\phi_{1,24}$ & $24$ & $q^{4}$ & $4$ & $20$ & $28$ & $44$
\\ $F_4[\I]$ & $20$ & $\I q^{2}$ & $4$ & $16$ & $24$ & $36$
\\ $F_4[\theta]$ & $20$ & $\theta q^{2}$ & $4$ & $16$ & $24$ & $36$
\\ $B_2;1$ & $11$ & $-q$ & $2$ & $10$ & $12$ & $20$
\\ $B_2;r$ & $20$ & $-q^{2}$ & $3$ & $17$ & $23$ & $37$
\\ $B_2;\ep$ & $23$ & $-q^{3}$ & $4$ & $20$ & $26$ & $42$
\\ $F_4[\theta^2]$ & $20$ & $\theta^2q^{2}$ & $4$ & $16$ & $24$ & $36$
\\ $F_4[-\I]$ & $20$ & $-\I q^{2}$ & $4$ & $16$ & $24$ & $36$
\\\hline\end{tabular}\end{center}

\begin{center}\begin{tikzpicture}[thick,scale=1.8]

\draw (0,0.82) node {$\phi_{1,0}$};
\draw (-1,0.82) node {$\phi_{4,1}$};
\draw (-2,0.82) node {$\phi_{6,6}''$};
\draw (-3,0.82) node {$\phi_{4,13}$};
\draw (-4,0.82) node {$\phi_{1,24}$};

\draw (-6,0.82) node {$B_2;\ep$};
\draw (-7,0.82) node {$B_2;r$};
\draw (-8,0.82) node {$B_2;1$};

\draw (0,1) -- (-8,1);
\draw (-4.5,0.1339) -- (-5.5,1.866);
\draw (-4.5,1.866) -- (-5.5,0.1339);
\draw (0,1) node [draw] (l0) {};
\draw (-1,1) node [draw] (l0) {};
\draw (-2,1) node [draw] (l0) {};
\draw (-3,1) node [draw] (l0) {};
\draw (-4,1) node [draw] (l0) {};
\draw (-5,1) node [fill=black!100] (ld) {};
\draw (-6,1) node [draw] (l0) {};
\draw (-7,1) node [draw] (l0) {};
\draw (-8,1) node [draw] (l0) {};
\draw (0,1.18) node{$0$};
\draw (-1,1.18) node{$1$};
\draw (-2,1.18) node{$2$};
\draw (-3,1.18) node{$3$};
\draw (-4,1.18) node{$4$};
\draw (-6,1.18) node{$4$};
\draw (-7,1.18) node{$3$};
\draw (-8,1.18) node{$2$};

\draw (-4.5,1.866) node [draw,label=right:${F_4[\I]}$] (l0) {};
\draw (-5.5,1.866) node [draw,label=left:${F_4[\theta]}$] (l0) {};
\draw (-4.5,0.1339) node [draw,label=right:${F_4[-\I]}$] (l0) {};
\draw (-5.5,0.1339) node [draw,label=left:${F_4[\theta^2]}$] (l0) {};

\draw (-4.5,2.03) node{$4$};
\draw (-5.5,2.03) node{$4$};
\draw (-4.5,-0.02) node{$4$};
\draw (-5.5,-0.02) node{$4$};
\end{tikzpicture}\end{center}

\newpage

\section{${}^2\!F_4(q^2)$}

The group ${}^2\!F_4(q^2)$, for $q$ an odd power of $\sqrt2$,  has order $q^{24}\Phi_1^2\Phi_2^2\Phi_4^2\Phi_8^2\Phi_{12}\Phi_{24}$. Since $\sqrt 2$ lies inside our field of definition, $\Phi_{8}$ factorizes as in the Suzuki group case, and $\Phi_{24}$ factorizes as
\[ \Phi_{12}=(q^4+\sqrt2q^3+q^2+\sqrt2q+1)(q^4-\sqrt2q^3+q^2-\sqrt2q+1).\]
Denote by $\Phi_{24}'$ the first and by $\Phi_{24}''$ the second. If $\lambda=\e^{\pi\I/12}$, then $\Phi_{12}'=(q-\lambda^5)(q-\lambda^{11})(q-\lambda^{13})(q-\lambda^{19})$ and $\Phi_{24}''=(q-\lambda)(q-\lambda^7)(q-\lambda^{17})(q-\lambda^{23})$. This is the same as how Carter defines them in \cite{carterfinite}. This way is consistent with $\Phi_d''$ being the Coxeter case, and having $\e^{2\pi\I/d}$ as a zero.

The Brauer trees for this group are for the principal block and $\ell$ dividing $\Phi_{12}$, $\Phi_{24}'$ and $\Phi_{24}''$, and non-principal blocks for $\Phi_1\Phi_2$: as with the other Suzuki and Ree groups, that if $\ell$ is an integer dividing $q^2-1$, as $\Phi_1=\sqrt2^n-1$ and $\Phi_2=\sqrt2^n+1$ with $n$ odd, so an integer cannot divide either of these. In \cite{hiss1991}, the Brauer trees for ${}^2\!F_4(q)$ are mostly determined, up to questions of choosing your field of definitions, but the planar embedding was not determined for $\ell\mid \Phi_{24}''$: this case was completed by Dudas in \cite{dudas2010un2}. Here we make consistent choices of the various cuspidal characters.

There are misprints in the table given in \cite{carterfinite}: the correct degrees for the twenty-one unipotent characters are as follows.

\begin{center}
\begin{tabular}{ll}
\hline Name & Degree
\\ \hline $\phi_{1,0}$ & $1$
\\ $\phi_{1,4}''$ & $q^2\Phi_{12}\Phi_{24}$
\\ $\phi_{1,4}'$ & $q^{10}\Phi_{12}\Phi_{24}$
\\ $\phi_{1,8}$ & $q^{24}$
\\ $\phi_{2,3}$ & $q^4\Phi_4^2\Phi_8''^2\Phi_{12}\Phi_{24}'/4$
\\ $\phi_{2,1}$ & $q^4\Phi_4^2\Phi_8'^2\Phi_{12}\Phi_{24}''/4$
\\ $\phi_{2,2}$ & $q^4\Phi_8^2\Phi_{24}/2$
\\ ${}^2\!B_2[\psi^3];1$ & $q\Phi_1\Phi_2\Phi_4^2\Phi_{12}/\sqrt2$
\\ ${}^2\!B_2[\psi^5];1$ & $q\Phi_1\Phi_2\Phi_4^2\Phi_{12}/\sqrt2$
\\ ${}^2\!B_2[\psi^3];\ep$ & $q^{13}\Phi_1\Phi_2\Phi_4^2\Phi_{12}/\sqrt2$
\\ ${}^2\!B_2[\psi^5];\ep$ & $q^{13}\Phi_1\Phi_2\Phi_4^2\Phi_{12}/\sqrt2$
\\ ${}^2F_4^{\mathrm{I}}[-1]$ & $q^4\Phi_1^2\Phi_2^2\Phi_4^2\Phi_{24}/6$
\\ ${}^2\!F_4^{\mathrm{II}}[-1]$& $q^4\Phi_1^2\Phi_2^2\Phi_8''^2\Phi_{12}\Phi_{24}''/12$
\\ ${}^2\!F_4^{\mathrm{III}}[-1]$ & $q^4\Phi_1^2\Phi_2^2\Phi_8'^2\Phi_{12}\Phi_{24}'/12$
\\ ${}^2\!F_4^{\mathrm{IV}}[-1]$&$q^4\Phi_1^2\Phi_2^2\Phi_{12}\Phi_{24}/3$
\\ ${}^2\!F_4^{\mathrm{I}}[\I]$ & $q^4\Phi_1^2\Phi_2^2\Phi_4^2\Phi_{12}\Phi_{24}'/4$
\\ ${}^2\!F_4^{\mathrm{II}}[\I]$ & $q^4\Phi_1^2\Phi_2^2\Phi_4^2\Phi_{12}\Phi_{24}''/4$
\\ ${}^2\!F_4^{\mathrm{I}}[-\I]$ & $q^4\Phi_1^2\Phi_2^2\Phi_4^2\Phi_{12}\Phi_{24}'/4$
\\ ${}^2\!F_4^{\mathrm{II}}[-\I]$ & $q^4\Phi_1^2\Phi_2^2\Phi_4^2\Phi_{12}\Phi_{24}''/4$
\\ ${}^2\!F_4[-\theta]$ & $q^4\Phi_1^2\Phi_2^2\Phi_4^2\Phi_8^2/3$
\\ ${}^2\!F_4[-\theta^2]$ & $q^4\Phi_1^2\Phi_2^2\Phi_4^2\Phi_8^2/3$
\\\hline\end{tabular}\end{center}

\newpage

\subsection{$d=1$}

For ${}^2\!F_4(q)$ and $\ell\mid (q^2-1)$ there are two unipotent blocks of weight $1$, together with seventeen unipotent blocks of defect zero.
\\[1cm]
\noindent(i)\;\; \textbf{Block 1:} Cuspidal pair is $(\Phi_1\Phi_2.{}^2\!B_2(q),{}^2\!B_2[\psi^3])$, of degree $q\Phi_1\Phi_2/\sqrt2$. There are two unipotent characters in the block, both of which are real.

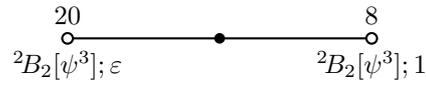
\begin{center}\begin{tabular}{lccc}
\hline Character & $A(-)$ & $\omega_i q^{aA/e}$ & $\kappa=1$
\\ \hline ${}^2\!B_2[\psi^3];1$ & $11$ & $q^4$ & $8$
\\ ${}^2\!B_2[\psi^3];\ep$ & $23$ & $-q^{16}$ & $20$
\\ \hline\end{tabular}%
\end{center}

\noindent(ii)\;\; \textbf{Block 2:}  Cuspidal pair is $(\Phi_1\Phi_2.{}^2\!B_2(q),{}^2\!B_2[\psi^5])$, of degree $q\Phi_1\Phi_2/\sqrt2$. There are two unipotent characters in the block, both of which are real.
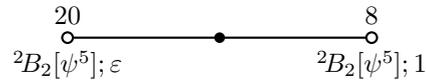
\begin{center}\begin{tabular}{lccc}
\hline Character & $A(-)$ & $\omega_i q^{aA/e}$ & $\kappa=1$
\\ \hline ${}^2\!B_2[\psi^5];1$ & $11$ & $q^4$ & $8$
\\ ${}^2\!B_2[\psi^5];\ep$ & $23$ & $-q^{16}$ & $20$
\\ \hline\end{tabular}%
\end{center}

\newpage

\subsection{$d=12$}

For ${}^2\!F_4(q)$ and $d=12$ there is a single unipotent block of weight $1$, together with fifteen unipotent blocks of defect zero.
\\[1cm]
\noindent(i)\;\;\textbf{Block 1:} Cuspidal pair is $(\Phi_{12},1)$, of degree $1$. There are six unipotent characters in the block, two of which -- ${}^2\!F_4[-\theta^i]$ -- are non-real.

\begin{center}\begin{tabular}{lcccccc}
\hline Character & $A(-)$ & $\omega_i q^{aA/e}$ & $\kappa=1$ & $\kappa=5$ & $\kappa=7$ & $\kappa=11$
\\ \hline $\phi_{1,0}$ & $0$ & $1$ & $0$ & $0$ & $0$ & $0$
\\ ${}^2\!F_4[-\theta]$ & $20$ & $-\theta q^4$ & $3$ & $17$ & $23$ & $35$
\\ $\phi_{2,2}$ & $20$ & $q^4$ & $3$ & $17$ & $23$ & $37$
\\ ${}^2\!F_4[-\theta^2]$ & $20$ & $-\theta^2 q^4$ & $3$ & $17$ & $23$ & $35$
\\ $\phi_{1,8}$ & $24$ & $q^8$ & $4$ & $20$ & $28$ & $44$
\\ ${}^2\!F_4^{\mathrm{I}}[-1]$ & $20$ & $-q^4$ & $4$ & $16$ & $24$ & $36$
\\ \hline
\end{tabular}%
\end{center}

\newpage

\subsection{$d=24'$}

For ${}^2\!F_4(q)$ and $d=24'$ there is a single unipotent block of weight $1$, together with nine unipotent blocks of defect zero.
\\[1cm]
\noindent(i)\;\;\textbf{Block 1:} Cuspidal pair is $(\Phi_{24}',1)$, of degree $1$. There are twelve unipotent characters in the block, eight of which are non-real.

\begin{center}\begin{tabular}{lcccccc}
\hline Character & $A(-)$ & $\omega_i q^{aA/e}$ & $\kappa=5$ & $\kappa=11$ & $\kappa=13$ & $\kappa=19$
\\ \hline $\phi_{1,0}$ & $0$ & $1$ & $0$ & $0$ & $0$ & $0$
\\ ${}^2\!B_2[\psi^3],1$ & $11$ & $\psi^7q$ & $4$ & $10$ & $12$ & $18$
\\ ${}^2\!F_4^{\mathrm{II}}[-\I]$ & $20$ & $-\I q^2$ & $8$ & $18$ & $22$ & $32$
\\ ${}^2\!F_4[-\theta^2]$ & $20$ & $-\theta q^2$ & $8$ & $18$ & $22$ & $32$
\\ ${}^2\!B_2[\psi^5],1$ &  $11$ & $\psi q$ & $4$ & $10$ & $12$ & $18$
\\ $\phi_{2,1}$ & $20$ & $q^2$ & $7$ & $17$ & $21$ & $31$
\\ ${}^2\!B_2[\psi^3],\ep$ & $23$ & $\psi^7q^3$ & $9$ & $21$ & $25$ & $37$
\\ ${}^2\!F_4[-\theta]$ & $20$ & $-\theta^2q^2$ & $8$ & $18$ & $22$ & $32$
\\ ${}^2\!F_4^{\mathrm{II}}[\I]$ & $20$ & $\I q^2$ & $8$ & $18$ & $22$ & $32$
\\ ${}^2\!B_2[\psi^5],\ep$ & $23$ & $\psi q^3$ & $9$ & $21$ & $25$ & $37$
\\ $\phi_{1,8}$ & $24$ & $q^4$ & $10$ & $22$ & $26$ & $38$
\\ ${}^2\!F_4^{\mathrm{II}}[-1]$ & $20$ & $-q^2$ & $10$ & $20$ & $24$ & $32$
\\\hline\end{tabular}%
\end{center}

\newpage

\subsection{$d=24''$}

For ${}^2\!F_4(q)$ and $d=24''$ there is a single unipotent block of weight $1$, together with nine unipotent blocks of defect zero.
\\[1cm]
\noindent(i)\;\;\textbf{Block 1:} Cuspidal pair is $(\Phi_{24}'',1)$, of degree $1$. There are twelve unipotent characters in the block, eight of which are non-real.

\begin{center}\begin{tabular}{lcccccc}
\hline Character & $A(-)$ & $\omega_i q^{aA/e}$ & $\kappa=1$ & $\kappa=7$ & $\kappa=17$ & $\kappa=23$
\\ \hline $\phi_{1,0}$ & $0$ & $1$ & $0$ & $0$ & $0$ & $0$
\\ $\phi_{2,3}$ & $20$ & $q^2$ & $1$ & $13$ & $27$ & $39$
\\ $\phi_{1,8}$ & $24$ & $q^4$ & $2$ & $14$ & $34$ & $46$
\\ ${}^2\!F_4[-\theta^2]$ & $20$ & $-\theta^2 q^2$ & $2$ & $12$ & $28$ & $38$ 
\\ ${}^2\!F_4^{\mathrm{I}}[\I]$ & $20$ & $\I q^2$ & $2$ & $12$ & $28$ & $38$
\\ ${}^2\!B_2[\psi^3];1$ & $11$ & $\psi^3 q$ & $1$ & $7$ & $15$ & $21$
\\ ${}^2\!B_2[\psi^3];\ep$ & $23$ & $\psi^3 q^3$ & $2$ & $14$ & $32$ & $44$
\\ ${}^2\!F_4^{\mathrm{III}}[-1]$ & $20$ & $-q^2$ & $2$ &  $10$ & $30$ & $38$
\\ ${}^2\!B_2[\psi^5];1$ & $11$ & $\psi^5 q$ & $1$ & $7$ & $15$ & $21$
\\ ${}^2\!B_2[\psi^5];\ep$ & $23$ & $\psi^5 q^3$ & $2$ & $14$ & $32$ & $44$
\\ ${}^2\!F_4^{\mathrm{I}}[-\I]$ & $20$ & $-\I q^2$ & $2$ & $12$ & $28$ & $38$
\\ ${}^2\!F_4[-\theta]$ & $20$ & $-\theta q^2$ & $2$ & $12$ & $28$ & $38$
\\\hline\end{tabular}%
\end{center}
\newpage

\section{$E_6(q)$}

The group $E_6(q)$ has order $q^{36}\Phi_1^6\Phi_2^4\Phi_3^3\Phi_4^2\Phi_5\Phi_6^2\Phi_8\Phi_9\Phi_{12}$. There are Brauer trees for $\Phi_5$, $\Phi_8$, $\Phi_9$ and $\Phi_{12}$ obviously, and also non-principal unipotent blocks of weight $1$ when $d=2$, $3$ and $4$. The Brauer trees for all unipotent blocks of weight $1$, together with their planar embeddings, were constructed in \cite{hlm1995}.

There are thirty unipotent characters of $E_6(q)$, with degrees given below.

\begin{center}\begin{tabular}{ll}
\hline Name & Degree
\\ \hline $\phi_{1,0}$ & $1$
\\ $\phi_{1,36}$ & $q^{36}$
\\ $\phi_{10,9}$ & $q^{7}\Phi_{5}\Phi_{6}^{2}\Phi_{8}\Phi_{9}\Phi_{12}/3$
\\ $\phi_{6,1}$ & $q\Phi_{8}\Phi_{9}$
\\ $\phi_{6,25}$ & $q^{25}\Phi_{8}\Phi_{9}$
\\ $\phi_{20,10}$ & $q^{7}\Phi_{4}^{2}\Phi_{5}\Phi_{6}^{2}\Phi_{8}\Phi_{9}/6$
\\ $\phi_{15,5}$ & $q^{3}\Phi_{5}\Phi_{6}^{2}\Phi_{8}\Phi_{9}/2$
\\
\hline\end{tabular}\end{center}

\newpage

\subsection{$d=2$}

For $E_6(q)$ and $d=2$ there is a single unipotent block of weight $1$, together with the principal block, with cyclotomic Weyl group $F_4$, and three unipotent blocks of defect zero.
\\[1cm]
\noindent(i)\;\; \textbf{Block 1:} Cuspidal pair is $(\Phi_2.A_5(q),\phi_{321})$, of degree $q^4\Phi_2^3\Phi_4\Phi_6$. There are two unipotent characters in the block, both of which are real.

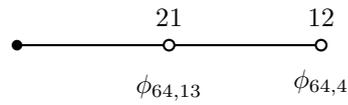
\begin{center}\begin{tabular}{lccc}
\hline Character & $A(-)$ & $\omega_i q^{aA/e}$ & $\kappa=1$
\\\hline $\phi_{64,4}$ & $12$ & $q^{6}$ & $12$
\\ $\phi_{64,13}$ & $21$ & $q^{15}$ & $21$
\\\hline\end{tabular}%
\end{center}

\newpage

\subsection{$d=3$}

For $E_6(q)$ and $d=3$ there is a single unipotent block of weight $1$, together with the principal block, with cyclotomic Weyl group $F_4$, and three unipotent blocks of defect zero.
\\[1cm]
\noindent(i)\;\; \textbf{Block 1:} Cuspidal pair is $(\Phi_3.{}^3\!D_4(q),{}^3\!D_4[-1])$, of degree $q^3\Phi_1^2\Phi_3^2/2$. There are three unipotent characters in the block, all of which are real.

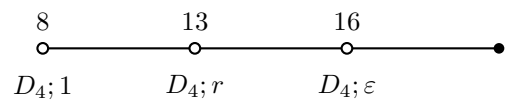
\begin{center}\begin{tabular}{lcccc}
\hline Character & $A(-)$ & $\omega_i q^{aA/e}$ & $\kappa=1$ & $\kappa=2$
\\\hline $D_4;1$ & $12$ & $-q^{4}$ & $8$ & $16$
\\ $D_4;r$ & $20$ &$-q^{8}$ & $13$ & $27$
\\ $D_4;\ep$ & $24$ &$-q^{12}$ & $16$ & $32$
\\\hline\end{tabular}%
\end{center}

\newpage

\subsection{$d=4$}

For $E_6(q)$ and $d=4$ there is a single unipotent block of weight $1$, together with the principal block, with cyclotomic Weyl group $G_8$, and ten unipotent blocks of defect zero.
\\[1cm]
\noindent(i)\;\; \textbf{Block 1:} Cuspidal pair is $(\Phi_1\Phi_4.{}^2\!A_3(q),\phi_{22})$, of degree $q^2\Phi_4$. There are four unipotent characters in the block, all of which are real.

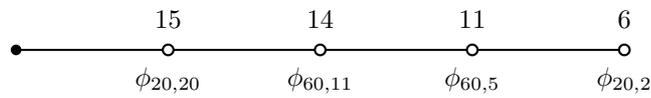
\begin{center}\begin{tabular}{lcccc}
\hline Character & $A(-)$ & $\omega_i q^{aA/e}$ & $\kappa=1$ & $\kappa=3$
\\\hline $\phi_{20,2}$ & $12$ & $q^{3}$ & $6$ & $18$
\\ $\phi_{20,20}$ & $30$ & $q^{12}$ & $15$ & $45$
\\ $\phi_{60,11}$ & $27$ & $q^{9}$ & $14$ & $40$
\\ $\phi_{60,5}$ & $21$ & $q^{6}$ & $11$ & $31$
\\\hline\end{tabular}%
\end{center}

\newpage

\subsection{$d=5$}

For $E_6(q)$ and $d=5$ there are two unipotent blocks of weight $1$, together with twenty unipotent blocks of defect zero.
\\[1cm]
\noindent(i)\;\; \textbf{Block 1:} Cuspidal pair is $(\Phi_1\Phi_5.A_1(q),\phi_2)$, of degree $1$. There are five unipotent characters in the block, all of which are real.

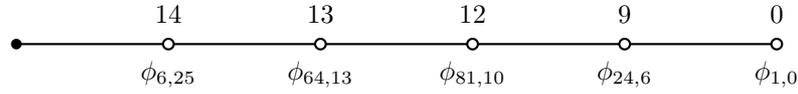
\begin{center}\begin{tabular}{lcccccc}
\hline Character & $A(-)$ & $\omega_i q^{aA/e}$ & $\kappa=1$ & $\kappa=2$ & $\kappa=3$ & $\kappa=4$
\\\hline $\phi_{1,0}$ & $0$ & $q^{0}$ & $0$ & $0$ & $0$ & $0$
\\ $\phi_{24,6}$ & $24$ & $q^{6}$ & $9$ & $19$ & $29$ & $39$
\\ $\phi_{6,25}$ & $35$ & $q^{12}$ & $14$ & $28$ & $42$ & $56$
\\ $\phi_{81,10}$ & $30$ & $q^{8}$ & $12$ & $24$ & $36$ & $48$
\\ $\phi_{64,13}$ & $32$ & $q^{9}$ & $13$ & $25$ & $39$ & $51$
\\\hline\end{tabular}%
\end{center}

\noindent(ii)\;\; \textbf{Block 2:} Cuspidal pair is $(\Phi_1\Phi_5.A_1(q),\phi_{11})$, of degree $q$. There are five unipotent characters in the block, all of which are real.

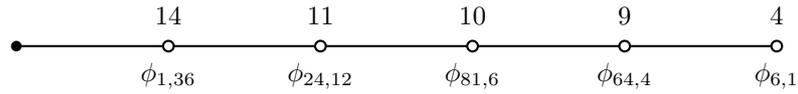
\begin{center}\begin{tabular}{lcccccc}
\hline Character & $A(-)$ & $\omega_i q^{aA/e}$ & $\kappa=1$ & $\kappa=2$ & $\kappa=3$ & $\kappa=4$
\\\hline $\phi_{6,1}$ & $10$ & $q^{2}$ & $4$ & $8$ & $12$ & $16$
\\ $\phi_{24,12}$ & $29$ & $q^{8}$ & $11$ & $23$ & $35$ & $47$
\\ $\phi_{1,36}$ & $35$ & $q^{14}$ & $14$ & $28$ & $42$ & $56$
\\ $\phi_{64,4}$ & $22$ & $q^{5}$ & $9$ & $17$ & $27$ & $35$
\\ $\phi_{81,6}$ & $25$ & $q^{6}$ & $10$ & $20$ & $30$ & $40$
\\\hline\end{tabular}%
\end{center}

\newpage

\subsection{$d=8$}

For $E_6(q)$ and $d=8$ there is a single unipotent block of weight $1$, together with twenty-two unipotent blocks of defect zero.
\\[1cm]
\noindent(i)\;\;\textbf{Block 1:} Cuspidal pair is $(\Phi_1\Phi_2\Phi_8,1)$, of degree $1$. There are eight unipotent characters in the block, all of which are real.

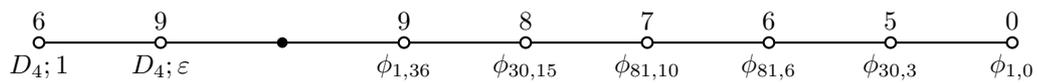
\begin{center}\begin{tabular}{lcccccc}
\hline Character & $A(-)$ & $\omega_i q^{aA/e}$ & $\kappa=1$ & $\kappa=3$ & $\kappa=5$ & $\kappa=7$
\\\hline $\phi_{1,0}$ & $0$ & $q^{0}$ & $0$ & $0$ & $0$ & $0$
\\ $\phi_{1,36}$ & $36$ & $q^{9}$ & $9$ & $27$ & $45$ & $63$
\\ $D_4;\ep$ & $33$ & $-q^{6}$ & $9$ & $25$ & $41$ & $57$
\\ $\phi_{30,3}$ & $21$ & $q^{3}$ & $5$ & $15$ & $27$ & $37$
\\ $\phi_{81,6}$ & $26$ & $q^{4}$ & $6$ & $20$ & $32$ & $46$
\\ $\phi_{81,10}$ & $30$ & $q^{5}$ & $7$ & $23$ & $37$ & $53$
\\ $\phi_{30,15}$ & $33$ & $q^{6}$ & $8$ & $24$ & $42$ & $58$
\\ $D_4;1$ & $21$ & $-q^{3}$ & $6$ & $16$ & $26$ & $36$
\\\hline\end{tabular}%
\end{center}

\newpage

\subsection{$d=9$}

For $E_6(q)$ and $d=9$ there is a single unipotent block of weight $1$, together with twenty-one unipotent blocks of defect zero.
\\[1cm]
\noindent(i)\;\;\textbf{Block 1:} Cuspidal pair is $(\Phi_9,1)$, of degree $1$. There are nine characters in the block, two of which -- $E_6[\theta^i]$ -- are non-real.

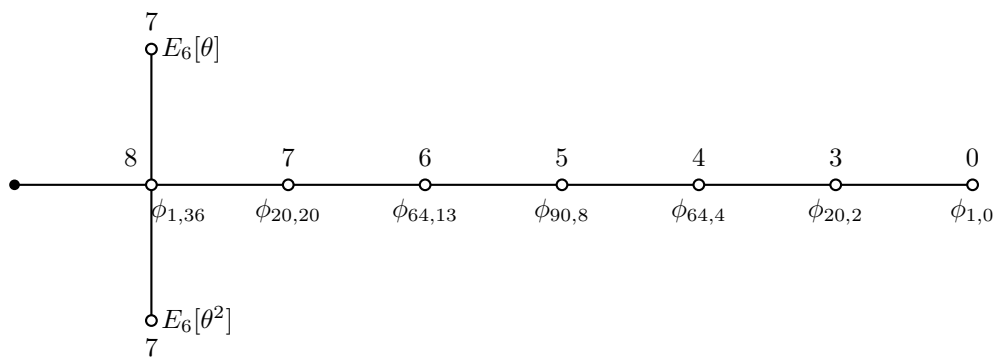
\begin{center}\begin{tabular}{lcccccccc}
\hline Character & $A(-)$ & $\omega_i q^{aA/e}$ & $\kappa=1$ & $\kappa=2$ & $\kappa=4$ & $\kappa=5$ & $\kappa=7$ & $\kappa=8$
\\\hline $\phi_{1,0}$ & $0$ & $q^{0}$ & $0$ & $0$ & $0$ & $0$ & $0$ & $0$
\\ $E_6[\theta^2]$ & $29$ & $\theta^2 q^{4}$ & $7$ & $13$ & $25$ & $33$ & $45$ & $51$
\\ $\phi_{20,2}$ & $16$ & $q^{2}$ & $3$ & $7$ & $15$ & $17$ & $25$ & $29$
\\ $\phi_{64,4}$ & $23$ & $q^{3}$ & $4$ & $10$ & $20$ & $26$ & $36$ & $42$
\\ $\phi_{90,8}$ & $29$ & $q^{4}$ & $5$ & $13$ & $27$ & $31$ & $45$ & $53$
\\ $\phi_{64,13}$ & $32$ & $q^{5}$ & $6$ & $14$ & $28$ & $36$ & $50$ & $58$
\\ $\phi_{20,20}$ & $34$ & $q^{6}$ & $7$ & $15$ & $31$ & $37$ & $53$ & $61$
\\ $E_6[\theta]$ & $29$ & $\theta q^{4}$ & $7$ & $13$ & $25$ & $33$ & $45$ & $51$
\\ $\phi_{1,36}$ & $36$ & $q^{8}$ & $8$ & $16$ & $32$ & $40$ & $56$ & $64$
\\\hline\end{tabular}%
\end{center}

\newpage

\subsection{$d=12$}

For $E_6(q)$ and $d=12$ there is a single unipotent block of weight $1$, together with eighteen unipotent blocks of defect zero.
\\[1cm]
\noindent(i)\;\;\textbf{Block 1:} Cuspidal pair is $(\Phi_3\Phi_{12},1)$, of degree $1$. There are twelve characters in the block, two of which -- $E_6[\theta^i]$ -- are non-real.

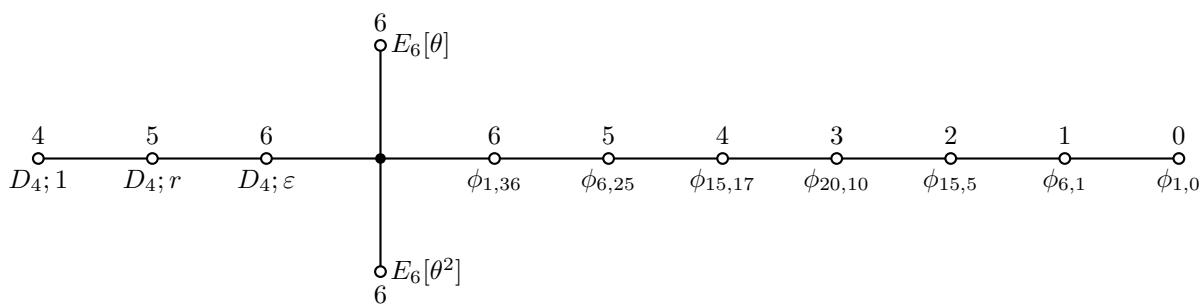
\begin{center}\begin{tabular}{lcccccc}
\hline Character & $A(-)$ & $\omega_i q^{aA/e}$ & $\kappa=1$ & $\kappa=5$ & $\kappa=7$ & $\kappa=11$
\\\hline $\phi_{1,0}$ & $0$ & $q^{0}$ & $0$ & $0$ & $0$ & $0$
\\ $\phi_{6,1}$ & $11$ & $q$ & $1$ & $9$ & $13$ & $21$
\\ $\phi_{15,5}$ & $21$ & $q^{2}$ & $2$ & $18$ & $24$ & $40$
\\ $\phi_{20,10}$ & $29$ & $q^{3}$ & $3$ & $25$ & $33$ & $55$
\\ $\phi_{15,17}$ & $33$ & $q^{4}$ & $4$ & $28$ & $38$ & $62$
\\ $\phi_{6,25}$ & $35$ & $q^{5}$ & $5$ & $29$ & $41$ & $65$
\\ $\phi_{1,36}$ & $36$ & $q^{6}$ & $6$ & $30$ & $42$ & $66$
\\ $E_6[\theta]$ & $29$ & $\theta q^{3}$ & $6$ & $24$ & $34$ & $52$
\\ $D_4;1$ & $21$ & $-q^{2}$ & $4$ & $18$ & $24$ & $38$
\\ $D_4;r$ & $29$ & $-q^{3}$ & $5$ & $25$ & $33$ & $53$
\\ $D_4;\ep$ & $33$ & $-q^{4}$ & $6$ & $28$ & $38$ & $60$
\\ $E_6[\theta^2]$ & $29$ & $\theta^2 q^{3}$ & $6$ & $24$ & $34$ & $52$
\\\hline\end{tabular}%
\end{center}

\newpage

\section{${}^2\!E_6(q^2)$}

The group ${}^2\!E_6(q^2)$ has order $q^{36}\Phi_1^4\Phi_2^6\Phi_3^2\Phi_4^2\Phi_6^3\Phi_8\Phi_{10}\Phi_{12}\Phi_{18}$. There are Brauer trees for $\Phi_8$, $\Phi_{10}$, $\Phi_{12}$ and $\Phi_{18}$ obviously, and also non-principal unipotent blocks of weight $1$ when $d=1$, $4$ and $6$. In \cite{hisslubeck1998} almost all Brauer trees were constructed: firstly, in the case where $d=1$ two blocks are suggested there, but this is incorrect and there is a single block; for $d=2$, the result below was only obtained for $q\equiv -1\bmod 3$, and the tree was not determined otherwise. This last case appears amenable to a proof via Deligne--Lusztig theory, and Dudas, Rouquier and the author are currently investigating this case.

There are thirty unipotent characters of ${}^2\!E_6(q)$, with degrees given below.

\begin{center}\begin{tabular}{ll}
\hline Name & Degree
\\\hline $\phi_{1,0}$ & $1$
\\ $\phi_{1,24}$ & $q^{36}$
\\ $\phi_{6,6}'$ & $q^{7}\Phi_{3}^{2}\Phi_{8}\Phi_{10}\Phi_{12}\Phi_{18}/3$
\\ $\phi_{2,4}'$ & $q\Phi_{8}\Phi_{18}$
\\ $\phi_{2,16}''$ & $q^{25}\Phi_{8}\Phi_{18}$
\\ $\phi_{12,4}$ & $q^{7}\Phi_{3}^{2}\Phi_{4}^{2}\Phi_{8}\Phi_{10}\Phi_{18}/6$
\\ $\phi_{9,2}$ & $q^{3}\Phi_{3}^{2}\Phi_{8}\Phi_{10}\Phi_{18}/2$
\\ $\phi_{9,10}$ & $q^{15}\Phi_{3}^{2}\Phi_{8}\Phi_{10}\Phi_{18}/2$
\\ $\phi_{1,12}''$ & $q^{3}\Phi_{8}\Phi_{10}\Phi_{12}\Phi_{18}/2$
\\ $\phi_{1,12}'$ & $q^{15}\Phi_{8}\Phi_{10}\Phi_{12}\Phi_{18}/2$
\\ $\phi_{4,1}$ & $q^{2}\Phi_{4}\Phi_{8}\Phi_{10}\Phi_{12}$
\\ $\phi_{4,13}$ & $q^{20}\Phi_{4}\Phi_{8}\Phi_{10}\Phi_{12}$
\\ $\phi_{8,3}''$ & $q^{6}\Phi_{4}^{2}\Phi_{8}\Phi_{12}\Phi_{18}$
\\ $\phi_{8,9}'$ & $q^{12}\Phi_{4}^{2}\Phi_{8}\Phi_{12}\Phi_{18}$
\\ $\phi_{2,4}''$ & $q^{3}\Phi_{4}^{2}\Phi_{10}\Phi_{12}\Phi_{18}/2$
\\ $\phi_{2,16}'$ & $q^{15}\Phi_{4}^{2}\Phi_{10}\Phi_{12}\Phi_{18}/2$
\\ $\phi_{4,8}$ & $q^{7}\Phi_{4}^{2}\Phi_{8}\Phi_{10}\Phi_{12}\Phi_{18}/2$
\\ ${}^2\!E_6[1]$ & $q^{7}\Phi_{1}^{4}\Phi_{8}\Phi_{10}\Phi_{12}\Phi_{18}/6$
\\ $\phi_{6,6}''$ & $q^{7}\Phi_{3}^{2}\Phi_{6}^{3}\Phi_{8}\Phi_{10}\Phi_{12}/3$
\\ $\phi_{4,7}'$ & $q^{5}\Phi_{4}\Phi_{8}\Phi_{10}\Phi_{12}\Phi_{18}$
\\ $\phi_{4,7}''$ & $q^{11}\Phi_{4}\Phi_{8}\Phi_{10}\Phi_{12}\Phi_{18}$
\\ ${}^2\!A_5;1$ & $q^{4}\Phi_{1}^{3}\Phi_{3}^{2}\Phi_{4}^{2}\Phi_{8}\Phi_{12}$
\\ ${}^2\!A_5;\ep$ & $q^{13}\Phi_{1}^{3}\Phi_{3}^{2}\Phi_{4}^{2}\Phi_{8}\Phi_{12}$
\\ $\phi_{9,6}'$ & $q^{6}\Phi_{3}^{2}\Phi_{6}^{3}\Phi_{12}\Phi_{18}$
\\ $\phi_{9,6}''$ & $q^{10}\Phi_{3}^{2}\Phi_{6}^{3}\Phi_{12}\Phi_{18}$
\\ $\phi_{8,3}'$ & $q^{3}\Phi_{2}^{4}\Phi_{6}^{2}\Phi_{10}\Phi_{18}/2$
\\ $\phi_{8,9}''$ & $q^{15}\Phi_{2}^{4}\Phi_{6}^{2}\Phi_{10}\Phi_{18}/2$
\\ $\phi_{16,5}$ & $q^{7}\Phi_{2}^{4}\Phi_{6}^{2}\Phi_{8}\Phi_{10}\Phi_{18}/2$
\\ ${}^2\!E_6[\theta]$ & $q^{7}\Phi_{1}^{4}\Phi_{2}^{6}\Phi_{4}^{2}\Phi_{8}\Phi_{10}/3$
\\ ${}^2\!E_6[\theta^2]$ & $q^{7}\Phi_{1}^{4}\Phi_{2}^{6}\Phi_{4}^{2}\Phi_{8}\Phi_{10}/3$
\\\hline\end{tabular}\end{center}

\newpage

\subsection{$d=1$}

For ${}^2\!E_6(q)$ and $d=1$ there is a single unipotent block of weight $1$, together with the principal block, with cyclotomic Weyl group $F_4$, and three unipotent blocks of defect zero.
\\[1cm]
\noindent(i)\;\;\textbf{Block 1:} Cuspidal pair is $(\Phi_1.{}^2\!A_5(q),\phi_{321})$, of degree $q^4\Phi_1^3\Phi_3\Phi_4$. There are two characters in the block, both of which are real.

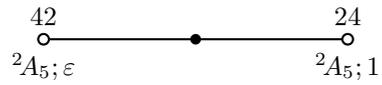
\begin{center}\begin{tabular}{lccc}
\hline Character & $A(-)$ & $\omega_i q^{aA/e}$ & $\kappa=1$
\\\hline ${}^2\!A_5;1$ & $12$ & $q^{6}$ & $24$
\\ ${}^2\!A_5;\ep$ & $21$ & $-q^{15}$ & $42$
\\\hline\end{tabular}\end{center}

\begin{center}\begin{tikzpicture}[thick,scale=2]
\draw (2,-0.18) node {${}^2\!A_5;1$};
\draw (0,-0.18) node {${}^2\!A_5;\ep$};
\draw (0,0) -- (2,0);
\draw (2,0) node [draw] (l0) {};
\draw (0,0) node [draw] (l1) {};
\draw (1,0) node [fill=black!100] (ld) {};
\draw (2,0.15) node{$24$};
\draw (0,0.15) node{$42$};
\end{tikzpicture}\end{center}

\newpage

\subsection{$d=4$}

For ${}^2\!E_6(q)$ and $d=4$ there is a single unipotent block of weight $1$, together with the principal block, with cyclotomic Weyl group $G_8$, and ten unipotent blocks of defect zero.
\\[1cm]
\noindent(i)\;\;\textbf{Block 1:} Cuspidal pair is $(\Phi_2\Phi_4.A_3(q),\phi_{22})$, of degree $q^2\Phi_4$. There are four characters in the block, all of which are real.

\begin{center}\begin{tabular}{lcccc}
\hline Character & $A(-)$ & $\omega_i q^{aA/e}$ & $\kappa=1$ & $\kappa=3$
\\\hline $\phi_{4,1}$ & $12$ & $q^{3}$ & $6$ & $18$
\\ $\phi_{4,13}$ & $30$ & $-q^{12}$ & $15$ & $45$
\\ $\phi_{4,7}'$ & $21$ & $-q^{6}$ & $10$ & $32$
\\ $\phi_{4,7}''$ & $27$ & $q^{9}$ & $13$ & $41$
\\\hline\end{tabular}\end{center}

\begin{center}\begin{tikzpicture}[thick,scale=2]
\draw (0,0) -- (-4,0);
\draw (0,0) node [draw,label=below:$\phi_{4,1}$] (l0) {};
\draw (-1,0) node [draw,label=below:$\phi_{4,7}''$] (l1) {};
\draw (-3,0) node [draw,label=below:$\phi_{4,13}$] (l2) {};
\draw (-4,0) node [draw,label=below:$\phi_{4,7}'$] (l2) {};
\draw (-2,0) node [fill=black!100] (ld) {};
\draw (0,0.18) node{$6$};
\draw (-1,0.18) node{$13$};
\draw (-3,0.18) node{$15$};
\draw (-4,0.18) node{$10$};
\end{tikzpicture}\end{center}

\newpage

\subsection{$d=6$}

For ${}^2\!E_6(q)$ and $d=6$ there is a single unipotent block of weight $1$, together with the principal block, with cyclotomic Weyl group $G_{25}$, and three unipotent blocks of defect zero.
\\[1cm]
\noindent(i)\;\; \textbf{Block 1:} Cuspidal pair is $(\Phi_6.{}^3\!D_4(q),\phi_{2,1})$, of degree $q^3\Phi_2^2\Phi_6^2/2$. There are five unipotent characters in the block, all of which are real.

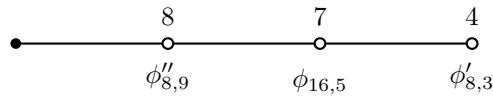
\begin{center}\begin{tabular}{lcccc}
\hline Character & $A(-)$ & $\omega_i q^{aA/e}$ & $\kappa=1$ & $\kappa=5$
\\\hline $\phi_{8,3}'$ & $12$ & $q^{4}$ & $4$ & $20$
\\ $\phi_{16,5}$ & $20$ & $q^{8}$ & $7$ & $33$
\\ $\phi_{8,9}''$ & $24$ & $q^{12}$ & $8$ & $40$
\\\hline\end{tabular}\end{center}
\begin{center}\begin{tikzpicture}[thick,scale=2]
\draw (0,0) -- (-3,0);
\draw (0,0) node [draw,label=below:$\phi_{8,3}'$] (l0) {};
\draw (-1,0) node [draw,label=below:$\phi_{16,5}$] (l1) {};
\draw (-2,0) node [draw,label=below:$\phi_{8,9}''$] (l2) {};
\draw (-3,0) node [fill=black!100] (ld) {};
\draw (-0,0.18) node{$4$};
\draw (-1,0.18) node{$7$};
\draw (-2,0.18) node{$8$};
\end{tikzpicture}\end{center}

\newpage

\subsection{$d=8$}

For ${}^2\!E_6(q)$ and $d=8$ there is a single unipotent block of weight $1$, together with twenty-two unipotent blocks of defect zero.
\\[1cm]
\noindent(i)\;\;\textbf{Block 1:} Cuspidal pair is $(\Phi_1\Phi_2\Phi_8,1)$, of degree $1$. There are eight unipotent characters in the block, all of which are real.

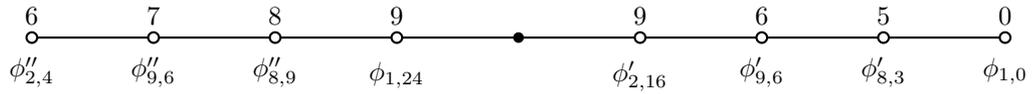
\begin{center}\begin{tabular}{lcccccc}
\hline Character & $A(-)$ & $\omega_i q^{aA/e}$ & $\kappa=1$ & $\kappa=3$ & $\kappa=5$ & $\kappa=7$
\\\hline $\phi_{1,0}$ & $0$ & $q^{0}$ & $0$ & $0$ & $0$ & $0$
\\ $\phi_{9,6}''$ & $30$ & $-q^{5}$ & $7$ & $23$ & $37$ & $53$
\\ $\phi_{8,9}''$ & $33$ & $-q^{6}$ & $8$ & $24$ & $42$ & $58$
\\ $\phi_{8,3}'$ & $21$ & $q^{3}$ & $5$ & $15$ & $27$ & $37$
\\ $\phi_{9,6}'$ & $26$ & $q^{4}$ & $6$ & $20$ & $32$ & $46$
\\ $\phi_{1,24}$ & $36$ & $-q^{9}$ & $9$ & $27$ & $45$ & $63$
\\ $\phi_{2,16}'$ & $33$ & $q^{6}$ & $9$ & $25$ & $41$ & $57$
\\ $\phi_{2,4}''$ & $21$ & $-q^{3}$ & $6$ & $16$ & $26$ & $36$
\\\hline\end{tabular}\end{center}
\begin{center}\begin{tikzpicture}[thick,scale=1.6]
\draw (0,0) -- (-8,0);
\draw (0,0) node [draw,label=below:$\phi_{1,0}$] (l0) {};
\draw (-1,0) node [draw,label=below:$\phi_{8,3}'$] (l1) {};
\draw (-2,0) node [draw,label=below:$\phi_{9,6}'$] (l2) {};
\draw (-3,0) node [draw,label=below:$\phi_{2,16}'$] (l2) {};
\draw (-4,0) node [fill=black!100] (ld) {};
\draw (-5,0) node [draw,label=below:$\phi_{1,24}$] (l2) {};
\draw (-6,0) node [draw,label=below:$\phi_{8,9}''$] (l2) {};
\draw (-7,0) node [draw,label=below:$\phi_{9,6}''$] (l2) {};
\draw (-8,0) node [draw,label=below:$\phi_{2,4}''$] (l2) {};
\draw (0,0.18) node{$0$};
\draw (-1,0.18) node{$5$};
\draw (-2,0.18) node{$6$};
\draw (-3,0.18) node{$9$};
\draw (-5,0.18) node{$9$};
\draw (-6,0.18) node{$8$};
\draw (-7,0.18) node{$7$};
\draw (-8,0.18) node{$6$};
\end{tikzpicture}\end{center}

\newpage

\subsection{$d=10$}

For ${}^2\!E_6(q)$ and $d=10$ there are two unipotent blocks of weight $1$, together with twenty unipotent blocks of defect zero.
\\[1cm]
\noindent(i)\;\; \textbf{Block 1:} Cuspidal pair is $(\Phi_2\Phi_{10}.A_1(q),\phi_2)$, of degree $1$. There are five unipotent characters in the block, all of which are real.

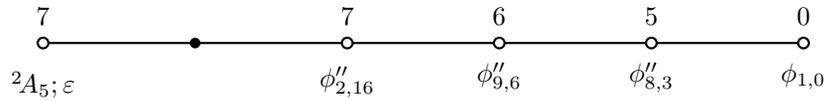
\begin{center}\begin{tabular}{lcccccc}
\hline Character & $A(-)$ & $\omega_i q^{aA/e}$ & $\kappa=1$ & $\kappa=3$ & $\kappa=7$ & $\kappa=9$
\\\hline $\phi_{1,0}$ & $0$ & $q^{0}$ & $0$ & $0$ & $0$ & $0$
\\ $\phi_{2,16}''$ & $35$ & $q^{12}$ & $7$ & $21$ & $49$ & $63$
\\ ${}^2\!A_5:\ep$ & $32$ & $-q^{9}$ & $7$ & $19$ & $45$ & $57$
\\ $\phi_{8,3}''$ & $24$ & $q^{6}$ & $5$ & $15$ & $33$ & $43$
\\ $\phi_{9,6}''$ & $30$ & $q^{8}$ & $6$ & $18$ & $42$ & $54$
\\\hline\end{tabular}\end{center}
\begin{center}\begin{tikzpicture}[thick,scale=2]
\draw (0,1.1) -- (-5,1.1);
\draw (0,1.1) node [draw,label=below:$\phi_{1,0}$] (l0) {};
\draw (-1,1.1) node [draw,label=below:$\phi_{8,3}''$] (l1) {};
\draw (-2,1.1) node [draw,label=below:$\phi_{9,6}''$] (l2) {};
\draw (-3,1.1) node [draw,label=below:$\phi_{2,16}''$] (l2) {};
\draw (-5,1.1) node [draw,label=below:${^2\!A_5;\ep}$] (l2) {};
\draw (-4,1.1) node [fill=black!100] (ld) {};
\draw (0,1.28) node{$0$};
\draw (-1,1.28) node{$5$};
\draw (-2,1.28) node{$6$};
\draw (-3,1.28) node{$7$};
\draw (-5,1.28) node{$7$};
\end{tikzpicture}\end{center}
\noindent(ii)\;\; \textbf{Block 2:} Cuspidal pair is $(\Phi_2\Phi_{10}.A_1(q),\phi_{11})$, of degree $q$. There are five unipotent characters in the block, all of which are real.

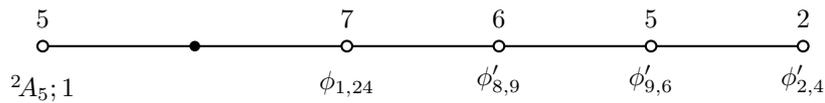
\begin{center}\begin{tabular}{lcccccc}
\hline Character & $A(-)$ & $\omega_i q^{aA/e}$ & $\kappa=1$ & $\kappa=3$ & $\kappa=7$ & $\kappa=9$
\\\hline $\phi_{2,4}'$ & $10$ & $q^{2}$ & $2$ & $6$ & $14$ & $18$
\\ $\phi_{1,24}$ & $35$ & $q^{14}$ & $7$ & $21$ & $49$ & $63$
\\ $\phi_{9,6}'$ & $25$ & $q^{6}$ & $5$ & $15$ & $35$ & $45$
\\ $\phi_{8,9}'$ & $29$ & $q^{8}$ & $6$ & $18$ & $40$ & $52$
\\ ${}^2\!A_5:1$ & $22$ & $-q^{5}$ & $5$ & $13$ & $31$ & $39$
\\\hline\end{tabular}\end{center}
\begin{center}\begin{tikzpicture}[thick,scale=2]

\draw (0,0) -- (-5,0);
\draw (0,0) node [draw,label=below:$\phi_{2,4}'$] (l0) {};
\draw (-1,0) node [draw,label=below:$\phi_{9,6}'$] (l1) {};
\draw (-2,0) node [draw,label=below:$\phi_{8,9}'$] (l2) {};
\draw (-3,0) node [draw,label=below:$\phi_{1,24}$] (l2) {};
\draw (-5,0) node [draw,label=below:${^2\!A_5;1}$] (l2) {};
\draw (-4,0) node [fill=black!100] (ld) {};
\draw (0,0.18) node{$2$};
\draw (-1,0.18) node{$5$};
\draw (-2,0.18) node{$6$};
\draw (-3,0.18) node{$7$};
\draw (-5,0.18) node{$5$};
\end{tikzpicture}\end{center}

\newpage

\subsection{$d=12$}

For ${}^2\!E_6(q)$ and $d=12$ there is a single unipotent block of weight $1$, together with eighteen unipotent blocks of defect zero.
\\[1cm]
\noindent(i)\;\; \textbf{Block 1:} Cuspidal pair is $(\Phi_6\Phi_{12},1)$, of degree $1$. There are twelve unipotent characters in the block, two of which  -- ${}^2\!E_6[\theta^i]$ -- are non-real. This Brauer tree was conjectural for $q\not\equiv 1\bmod 3$, and is proved in full generality in \cite{cdr2012un}.

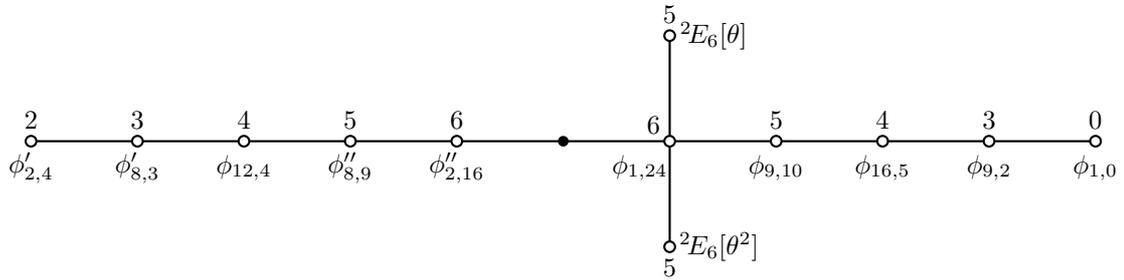
\begin{center}\begin{tabular}{lcccccc}
\hline Character & $A(-)$ & $\omega_i q^{aA/e}$ & $\kappa=1$ & $\kappa=5$ & $\kappa=7$ & $\kappa=11$
\\\hline $\phi_{1,0}$ & $0$ & $q^{0}$ & $0$ & $0$ & $0$ & $0$
\\ ${}^2\!E_6[\theta^2]$ & $29$ & $-\theta q^{3}$ & $5$ & $23$ & $35$ & $53$
\\ $\phi_{9,2}$ & $21$ & $q^{2}$ & $3$ & $19$ & $23$ & $39$
\\ $\phi_{16,5}$ & $29$ & $q^{3}$ & $4$ & $24$ & $34$ & $54$
\\ $\phi_{9,10}$ & $33$ & $q^{4}$ & $5$ & $29$ & $37$ & $61$
\\ ${}^2\!E_6[\theta]$ & $29$ & $-\theta^2 q^{3}$ & $5$ & $23$ & $35$ & $53$
\\ $\phi_{1,24}$ & $36$ & $q^{6}$ & $6$ & $30$ & $42$ & $66$
\\ $\phi_{2,4}'$ & $11$ & $-q$ & $2$ & $10$ & $12$ & $20$
\\ $\phi_{8,3}'$ & $21$ & $-q^{2}$ & $3$ & $17$ & $25$ & $39$
\\ $\phi_{12,4}$ & $29$ & $-q^{3}$ & $4$ & $26$ & $32$ & $54$
\\ $\phi_{8,9}''$ & $33$ & $-q^{4}$ & $5$ & $27$ & $39$ & $61$
\\ $\phi_{2,16}''$ & $35$ & $-q^{5}$ & $6$ & $30$ & $40$ & $64$
\\\hline\end{tabular}\end{center}
\begin{center}\begin{tikzpicture}[thick,scale=1.4]
\draw (0,0.25) node {$\phi_{1,0}$};
\draw (-1,0.25) node {$\phi_{9,2}$};
\draw (-2,0.25) node {$\phi_{16,5}$};
\draw (-3,0.25) node {$\phi_{9,10}$};
\draw (-4.28,0.25) node {$\phi_{1,24}$};
\draw (-6,0.25) node {$\phi_{2,16}''$};
\draw (-7,0.25) node {$\phi_{8,9}''$};
\draw (-8,0.25) node {$\phi_{12,4}$};
\draw (-9,0.25) node {$\phi_{8,3}'$};
\draw (-10,0.25) node {$\phi_{2,4}'$};

\draw (0,0.5) -- (-10,0.5);
\draw (-4,-0.5) -- (-4,1.5);

\draw (0,0.5) node [draw] (l0) {};
\draw (-1,0.5) node [draw] (l0) {};
\draw (-2,0.5) node [draw] (l0) {};
\draw (-3,0.5) node [draw] (l0) {};
\draw (-4,0.5) node [draw] (l0) {};
\draw (-4,-0.5) node [draw,label=right:${{}^2\!E_6[\theta^2]}$] (l4) {};
\draw (-4,1.5) node [draw,label=right:${{}^2\!E_6[\theta]}$] (l4) {};
\draw (-6,0.5) node [draw] (l0) {};
\draw (-7,0.5) node [draw] (l0) {};
\draw (-8,0.5) node [draw] (l0) {};
\draw (-9,0.5) node [draw] (l0) {};
\draw (-10,0.5) node [draw] (l0) {};

\draw (-5,0.5) node [fill=black!100] (ld) {};
\draw (0,0.7) node{$0$};
\draw (-1,0.7) node{$3$};
\draw (-2,0.7) node{$4$};
\draw (-3,0.7) node{$5$};
\draw (-4.15,0.65) node{$6$};
\draw (-4,1.7) node{$5$};
\draw (-4,-0.7) node{$5$};
\draw (-6,0.7) node{$6$};
\draw (-7,0.7) node{$5$};
\draw (-8,0.7) node{$4$};
\draw (-9,0.7) node{$3$};
\draw (-10,0.7) node{$2$};
\end{tikzpicture}\end{center}

\newpage

\subsection{$d=18$}
For ${}^2\!E_6(q)$ and $d=18$ there is a single unipotent block of weight $1$, together with twenty-one unipotent blocks of defect zero.
\\[1cm]
\noindent(i)\;\;\textbf{Block 1:} Cuspidal pair is $(\Phi_{18},1)$, of degree $1$. There are nine characters in the block, two of which -- ${}^2\!E_6[\theta^i]$ -- are non-real.

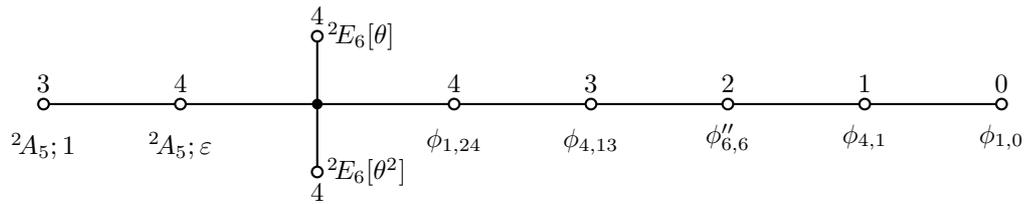
\begin{center}\begin{tabular}{lcccccccc}
\hline Character & $A(-)$ & $\omega_i q^{aA/e}$ & $\kappa=1$ & $\kappa=5$ & $\kappa=7$ & $\kappa=11$ & $\kappa=13$ & $\kappa=17$
\\\hline $\phi_{1,0}$ & $0$ & $q^{0}$ & $0$ & $0$ & $0$ & $0$ & $0$ & $0$
\\ $\phi_{4,1}$ & $16$ & $q^{2}$ & $1$ & $9$ & $13$ & $19$ & $23$ & $31$
\\ $\phi_{6,6}''$ & $29$ & $q^{4}$ & $2$ & $16$ & $24$ & $34$ & $42$ & $56$
\\ $\phi_{4,13}$ & $34$ & $q^{6}$ & $3$ & $19$ & $27$ & $41$ & $49$ & $65$
\\ $\phi_{1,24}$ & $36$ & $q^{8}$ & $4$ & $20$ & $28$ & $44$ & $52$ & $68$
\\ ${}^2\!E_6[\theta]$ & $29$ & $\theta q^{4}$ & $4$ & $16$ & $22$ & $36$ & $42$ & $54$
\\ ${}^2\!A_5:1$ & $23$ & $-q^{3}$ & $3$ & $13$ & $19$ & $27$ & $33$ & $43$
\\ ${}^2\!A_5:\ep$ & $32$ & $-q^{5}$ & $4$ & $18$ & $26$ & $38$ & $46$ & $60$
\\ ${}^2\!E_6[\theta^2]$ & $29$ & $\theta^2 q^{4}$ & $4$ & $16$ & $22$ & $36$ & $42$ & $54$
\\\hline\end{tabular}\end{center}
\begin{center}\begin{tikzpicture}[thick,scale=1.8]
\draw (0,0.5) -- (-7,0.5);
\draw (-5,0) -- (-5,1);

\draw (0,0.5) node [draw,label=below:$\phi_{1,0}$] (l0) {};
\draw (-1,0.5) node [draw,label=below:$\phi_{4,1}$] (l1) {};
\draw (-2,0.5) node [draw,label=below:$\phi_{6,6}''$] (l2) {};
\draw (-3,0.5) node [draw,label=below:$\phi_{4,13}$] (l2) {};
\draw (-4,0.5) node [draw,label=below:$\phi_{1,24}$] (l2) {};
\draw (-5,0) node [draw,label=right:${^2\!E_6[\theta^2]}$] (l4) {};
\draw (-5,1) node [draw,label=right:${^2\!E_6[\theta]}$] (l4) {};
\draw (-5,0.5) node [fill=black!100] (ld) {};
\draw (-6,0.5) node [draw,label=below:${^2\!A_5;\ep}$] (l4) {};
\draw (-7,0.5) node [draw,label=below:${^2\!A_5;1}$] (l2) {};
\draw (-0,0.65) node{$0$};
\draw (-1,0.65) node{$1$};
\draw (-2,0.65) node{$2$};
\draw (-3,0.65) node{$3$};
\draw (-4,0.65) node{$4$};
\draw (-7,0.65) node{$3$};
\draw (-6,0.65) node{$4$};
\draw (-5,1.15) node{$4$};
\draw (-5,-0.15) node{$4$};
\end{tikzpicture}\end{center}

\newpage

\section{$E_7(q)$}

The group $E_7(q)$ has order $q^{63}\Phi_1^7\Phi_2^7\Phi_3^3\Phi_4^2\Phi_5\Phi_6^3\Phi_7\Phi_8\Phi_9\Phi_{10}\Phi_{12}\Phi_{14}\Phi_{18}$. For this group, there are Brauer trees of unipotent blocks for all possible $d$ except $d=4$. Apart from the Brauer trees for $\Phi_{18}$ that have been constructed by Dudas and Rouquier \cite{dudasrouquier2012un}, no other Brauer trees have appeared in the literature. However, most of them can be deduced simply from the structure of the blocks of the Hecke algebra, which we have referred to as the Geck--Pfeiffer argument.

The remaining trees are more problematic: for $d=9$ there is an unfolding of the corresponding tree for $E_6(q)$ which can be used, but the details are still being worked out. For $d=10$ and $d=14$ there is another approach that Dudas, Rouquier and the author are taking, which should generate results. (This method might well be applicable to other trees in $E_8$.)

There are seventy-six unipotent characters of $E_7(q)$, whose degrees are given overleaf.
\begin{center}\begin{tabular}{llll}
\hline Name & Degree & Name & Degree
\\\hline$\phi_{1,0}$ & $1$
&$\phi_{210,6}$ & $q^{6}\Phi_{5}\Phi_{7}\Phi_{8}\Phi_{9}\Phi_{10}\Phi_{14}\Phi_{18}$
\\$\phi_{1,63}$ & $q^{63}$
&$\phi_{210,21}$ & $q^{21}\Phi_{5}\Phi_{7}\Phi_{8}\Phi_{9}\Phi_{10}\Phi_{14}\Phi_{18}$
\\$\phi_{7,46}$ & $q^{46}\Phi_{7}\Phi_{12}\Phi_{14}$
&$\phi_{210,10}$ & $q^{10}\Phi_{5}\Phi_{7}\Phi_{8}\Phi_{9}\Phi_{10}\Phi_{12}\Phi_{14}\Phi_{18}$
\\$\phi_{7,1}$ & $q\Phi_{7}\Phi_{12}\Phi_{14}$
&$\phi_{210,13}$ & $q^{13}\Phi_{5}\Phi_{7}\Phi_{8}\Phi_{9}\Phi_{10}\Phi_{12}\Phi_{14}\Phi_{18}$
\\$\phi_{15,28}$ & $q^{25}\Phi_{5}\Phi_{8}\Phi_{9}\Phi_{10}\Phi_{12}\Phi_{14}\Phi_{18}/2$
&$\phi_{216,16}$ & $q^{15}\Phi_{2}^{4}\Phi_{3}^{2}\Phi_{6}^{3}\Phi_{9}\Phi_{10}\Phi_{12}\Phi_{14}\Phi_{18}/2$
\\$\phi_{15,7}$ & $q^{4}\Phi_{5}\Phi_{8}\Phi_{9}\Phi_{10}\Phi_{12}\Phi_{14}\Phi_{18}/2$
&$\phi_{216,9}$ & $q^{8}\Phi_{2}^{4}\Phi_{3}^{2}\Phi_{6}^{3}\Phi_{9}\Phi_{10}\Phi_{12}\Phi_{14}\Phi_{18}/2$
\\$\phi_{21,6}$ & $q^{3}\Phi_{7}\Phi_{8}\Phi_{9}\Phi_{10}\Phi_{12}\Phi_{14}/2$
&$\phi_{280,18}$ & $q^{16}\Phi_{4}^{2}\Phi_{5}\Phi_{7}\Phi_{8}\Phi_{9}\Phi_{10}\Phi_{14}\Phi_{18}/3$
\\$\phi_{21,33}$ & $q^{30}\Phi_{7}\Phi_{8}\Phi_{9}\Phi_{10}\Phi_{12}\Phi_{14}/2$
&$\phi_{280,9}$ & $q^{7}\Phi_{4}^{2}\Phi_{5}\Phi_{7}\Phi_{8}\Phi_{9}\Phi_{10}\Phi_{14}\Phi_{18}/3$
\\$\phi_{21,36}$ & $q^{36}\Phi_{7}\Phi_{9}\Phi_{14}\Phi_{18}$
&$\phi_{280,8}$ & $q^{7}\Phi_{2}^{4}\Phi_{5}\Phi_{6}^{3}\Phi_{7}\Phi_{10}\Phi_{12}\Phi_{14}\Phi_{18}/2$
\\$\phi_{21,3}$ & $q^{3}\Phi_{7}\Phi_{9}\Phi_{14}\Phi_{18}$
&$\phi_{280,17}$ & $q^{16}\Phi_{2}^{4}\Phi_{5}\Phi_{6}^{3}\Phi_{7}\Phi_{10}\Phi_{12}\Phi_{14}\Phi_{18}/2$
\\$\phi_{27,2}$ & $q^{2}\Phi_{3}^{2}\Phi_{6}^{2}\Phi_{9}\Phi_{12}\Phi_{18}$
&$\phi_{315,16}$ & $q^{16}\Phi_{3}^{3}\Phi_{5}\Phi_{7}\Phi_{8}\Phi_{10}\Phi_{12}\Phi_{14}\Phi_{18}/6$
\\$\phi_{27,37}$ & $q^{37}\Phi_{3}^{2}\Phi_{6}^{2}\Phi_{9}\Phi_{12}\Phi_{18}$
&$\phi_{315,7}$ & $q^{7}\Phi_{3}^{3}\Phi_{5}\Phi_{7}\Phi_{8}\Phi_{10}\Phi_{12}\Phi_{14}\Phi_{18}/6$
\\$\phi_{35,22}$ & $q^{16}\Phi_{5}\Phi_{6}^{3}\Phi_{7}\Phi_{8}\Phi_{9}\Phi_{10}\Phi_{12}\Phi_{14}/6$
&$\phi_{336,14}$ & $q^{13}\Phi_{2}^{4}\Phi_{6}^{2}\Phi_{7}\Phi_{8}\Phi_{9}\Phi_{10}\Phi_{14}\Phi_{18}/2$
\\$\phi_{35,13}$ & $q^{7}\Phi_{5}\Phi_{6}^{3}\Phi_{7}\Phi_{8}\Phi_{9}\Phi_{10}\Phi_{12}\Phi_{14}/6$
&$\phi_{336,11}$ & $q^{10}\Phi_{2}^{4}\Phi_{6}^{2}\Phi_{7}\Phi_{8}\Phi_{9}\Phi_{10}\Phi_{14}\Phi_{18}/2$
\\$\phi_{35,4}$ & $q^{3}\Phi_{5}\Phi_{7}\Phi_{8}\Phi_{12}\Phi_{14}\Phi_{18}/2$
&$\phi_{378,14}$ & $q^{14}\Phi_{3}^{2}\Phi_{6}^{2}\Phi_{7}\Phi_{8}\Phi_{9}\Phi_{12}\Phi_{14}\Phi_{18}$
\\$\phi_{35,31}$ & $q^{30}\Phi_{5}\Phi_{7}\Phi_{8}\Phi_{12}\Phi_{14}\Phi_{18}/2$
&$\phi_{378,9}$ & $q^{9}\Phi_{3}^{2}\Phi_{6}^{2}\Phi_{7}\Phi_{8}\Phi_{9}\Phi_{12}\Phi_{14}\Phi_{18}$
\\$\phi_{56,30}$ & $q^{30}\Phi_{2}^{4}\Phi_{6}^{2}\Phi_{7}\Phi_{10}\Phi_{14}\Phi_{18}/2$
&$\phi_{405,8}$ & $q^{8}\Phi_{3}^{3}\Phi_{5}\Phi_{6}^{2}\Phi_{8}\Phi_{9}\Phi_{12}\Phi_{14}\Phi_{18}/2$
\\$\phi_{56,3}$ & $q^{3}\Phi_{2}^{4}\Phi_{6}^{2}\Phi_{7}\Phi_{10}\Phi_{14}\Phi_{18}/2$
&$\phi_{405,15}$ & $q^{15}\Phi_{3}^{3}\Phi_{5}\Phi_{6}^{2}\Phi_{8}\Phi_{9}\Phi_{12}\Phi_{14}\Phi_{18}/2$
\\$\phi_{70,18}$ & $q^{16}\Phi_{5}\Phi_{7}\Phi_{8}\Phi_{9}\Phi_{10}\Phi_{12}\Phi_{14}\Phi_{18}/3$
&$\phi_{420,10}$ & $q^{10}\Phi_{4}^{2}\Phi_{5}\Phi_{7}\Phi_{8}\Phi_{9}\Phi_{12}\Phi_{14}\Phi_{18}/2$
\\$\phi_{70,9}$ & $q^{7}\Phi_{5}\Phi_{7}\Phi_{8}\Phi_{9}\Phi_{10}\Phi_{12}\Phi_{14}\Phi_{18}/3$
&$\phi_{420,13}$ & $q^{13}\Phi_{4}^{2}\Phi_{5}\Phi_{7}\Phi_{8}\Phi_{9}\Phi_{12}\Phi_{14}\Phi_{18}/2$
\\$\phi_{84,12}$ & $q^{10}\Phi_{4}^{2}\Phi_{7}\Phi_{8}\Phi_{9}\Phi_{10}\Phi_{12}\Phi_{14}\Phi_{18}/2$
&$\phi_{512,12}$ & $q^{11}\Phi_{2}^{7}\Phi_{4}^{2}\Phi_{6}^{3}\Phi_{8}\Phi_{10}\Phi_{12}\Phi_{14}\Phi_{18}/2$
\\$\phi_{84,15}$ & $q^{13}\Phi_{4}^{2}\Phi_{7}\Phi_{8}\Phi_{9}\Phi_{10}\Phi_{12}\Phi_{14}\Phi_{18}/2$
&$\phi_{512,11}$ & $q^{11}\Phi_{2}^{7}\Phi_{4}^{2}\Phi_{6}^{3}\Phi_{8}\Phi_{10}\Phi_{12}\Phi_{14}\Phi_{18}/2$
\\$\phi_{105,26}$ & $q^{25}\Phi_{5}\Phi_{7}\Phi_{8}\Phi_{9}\Phi_{10}\Phi_{12}\Phi_{18}/2$
&$D_4;\ep$ & $q^{30}\Phi_{1}^{4}\Phi_{3}^{2}\Phi_{5}\Phi_{7}\Phi_{9}\Phi_{14}/2$
\\$\phi_{105,5}$ & $q^{4}\Phi_{5}\Phi_{7}\Phi_{8}\Phi_{9}\Phi_{10}\Phi_{12}\Phi_{18}/2$
&$D_4;\sigma_2'$ & $q^{13}\Phi_{1}^{4}\Phi_{3}^{2}\Phi_{5}\Phi_{7}\Phi_{8}\Phi_{9}\Phi_{14}\Phi_{18}/2$
\\$\phi_{105,6}$ & $q^{6}\Phi_{5}\Phi_{7}\Phi_{9}\Phi_{10}\Phi_{12}\Phi_{14}\Phi_{18}$
&$D_4;\ep_1$ & $q^{4}\Phi_{1}^{4}\Phi_{3}^{2}\Phi_{5}\Phi_{7}\Phi_{9}\Phi_{10}\Phi_{18}/2$
\\$\phi_{105,21}$ & $q^{21}\Phi_{5}\Phi_{7}\Phi_{9}\Phi_{10}\Phi_{12}\Phi_{14}\Phi_{18}$
&$D_4;r\ep$ & $q^{16}\Phi_{1}^{4}\Phi_{3}^{3}\Phi_{5}\Phi_{7}\Phi_{9}\Phi_{10}\Phi_{12}\Phi_{14}/2$
\\$\phi_{105,12}$ & $q^{12}\Phi_{5}\Phi_{7}\Phi_{9}\Phi_{10}\Phi_{12}\Phi_{14}\Phi_{18}$
&$D_4;r\ep_1$ & $q^{8}\Phi_{1}^{4}\Phi_{3}^{3}\Phi_{5}\Phi_{6}^{2}\Phi_{7}\Phi_{9}\Phi_{12}\Phi_{18}/2$
\\$\phi_{105,15}$ & $q^{15}\Phi_{5}\Phi_{7}\Phi_{9}\Phi_{10}\Phi_{12}\Phi_{14}\Phi_{18}$
&$D_4;r\ep_2$ & $q^{15}\Phi_{1}^{4}\Phi_{3}^{3}\Phi_{5}\Phi_{6}^{2}\Phi_{7}\Phi_{9}\Phi_{12}\Phi_{18}/2$
\\$\phi_{120,4}$ & $q^{4}\Phi_{2}^{4}\Phi_{5}\Phi_{6}^{2}\Phi_{9}\Phi_{10}\Phi_{14}\Phi_{18}/2$
&$D_4;\ep_2$ & $q^{25}\Phi_{1}^{4}\Phi_{3}^{2}\Phi_{5}\Phi_{7}\Phi_{9}\Phi_{10}\Phi_{18}/2$
\\$\phi_{120,25}$ & $q^{25}\Phi_{2}^{4}\Phi_{5}\Phi_{6}^{2}\Phi_{9}\Phi_{10}\Phi_{14}\Phi_{18}/2$
&$D_4;r$ & $q^{7}\Phi_{1}^{4}\Phi_{3}^{3}\Phi_{5}\Phi_{7}\Phi_{9}\Phi_{10}\Phi_{12}\Phi_{14}/2$
\\$\phi_{168,6}$ & $q^{6}\Phi_{4}^{2}\Phi_{7}\Phi_{8}\Phi_{9}\Phi_{12}\Phi_{14}\Phi_{18}$
&$D_4;\sigma_2$ & $q^{10}\Phi_{1}^{4}\Phi_{3}^{2}\Phi_{5}\Phi_{7}\Phi_{8}\Phi_{9}\Phi_{14}\Phi_{18}/2$
\\$\phi_{168,21}$ & $q^{21}\Phi_{4}^{2}\Phi_{7}\Phi_{8}\Phi_{9}\Phi_{12}\Phi_{14}\Phi_{18}$
&$D_4;1$ & $q^{3}\Phi_{1}^{4}\Phi_{3}^{2}\Phi_{5}\Phi_{7}\Phi_{9}\Phi_{14}/2$
\\$\phi_{189,10}$ & $q^{8}\Phi_{3}^{2}\Phi_{6}^{3}\Phi_{7}\Phi_{8}\Phi_{9}\Phi_{10}\Phi_{12}\Phi_{18}/2$
&$E_6[\theta];\ep$ & $q^{16}\Phi_{1}^{6}\Phi_{2}^{6}\Phi_{4}^{2}\Phi_{5}\Phi_{7}\Phi_{8}\Phi_{10}\Phi_{14}/3$
\\$\phi_{189,17}$ & $q^{15}\Phi_{3}^{2}\Phi_{6}^{3}\Phi_{7}\Phi_{8}\Phi_{9}\Phi_{10}\Phi_{12}\Phi_{18}/2$
&$E_6[\theta];1$ & $q^{7}\Phi_{1}^{6}\Phi_{2}^{6}\Phi_{4}^{2}\Phi_{5}\Phi_{7}\Phi_{8}\Phi_{10}\Phi_{14}/3$
\\$\phi_{189,22}$ & $q^{22}\Phi_{3}^{2}\Phi_{6}^{2}\Phi_{7}\Phi_{9}\Phi_{12}\Phi_{14}\Phi_{18}$
&$E_6[\theta^2];\ep$ & $q^{16}\Phi_{1}^{6}\Phi_{2}^{6}\Phi_{4}^{2}\Phi_{5}\Phi_{7}\Phi_{8}\Phi_{10}\Phi_{14}/3$
\\$\phi_{189,5}$ & $q^{5}\Phi_{3}^{2}\Phi_{6}^{2}\Phi_{7}\Phi_{9}\Phi_{12}\Phi_{14}\Phi_{18}$
&$E_6[\theta^2];1$ & $q^{7}\Phi_{1}^{6}\Phi_{2}^{6}\Phi_{4}^{2}\Phi_{5}\Phi_{7}\Phi_{8}\Phi_{10}\Phi_{14}/3$
\\$\phi_{189,20}$ & $q^{20}\Phi_{3}^{2}\Phi_{6}^{2}\Phi_{7}\Phi_{9}\Phi_{12}\Phi_{14}\Phi_{18}$
&$E_7[-\I]$ & $q^{11}\Phi_{1}^{7}\Phi_{3}^{3}\Phi_{4}^{2}\Phi_{5}\Phi_{7}\Phi_{8}\Phi_{9}\Phi_{12}/2$
\\$\phi_{189,7}$ & $q^{7}\Phi_{3}^{2}\Phi_{6}^{2}\Phi_{7}\Phi_{9}\Phi_{12}\Phi_{14}\Phi_{18}$
&$E_7[\I]$ & $q^{11}\Phi_{1}^{7}\Phi_{3}^{3}\Phi_{4}^{2}\Phi_{5}\Phi_{7}\Phi_{8}\Phi_{9}\Phi_{12}/2$
\\\hline\end{tabular}\end{center}

\newpage

\subsection{$d=1$}

For $E_7(q)$ and $d=1$ there are two unipotent blocks of weight $1$, together with the principal block with cyclotomic Weyl group $E_7$, a block of weight $3$ with cyclotomic Weyl group $B_3$, and two unipotent blocks of defect zero.
\\[1cm]
\noindent(i)\;\;\textbf{Block 1/2:} Cuspidal pair is $(\Phi_1.E_6(q),E_6[\theta^i])$, of degree $q^7\Phi_1^6\Phi_2^4\Phi_4^2\Phi_5\Phi_8/3$. There are two unipotent characters in each block, neither of which is real, and the blocks are conjugate.
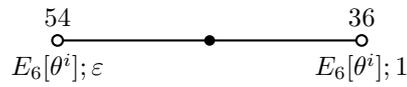
\begin{center}\begin{tabular}{lccc}
\hline Character & $A(-)$ & $\omega_i q^{aA/e}$ & $\kappa=1$
\\\hline $E_6[\theta^i];1$ & $18$ & $q^9$ & $36$
\\ $E_6[\theta^i];\ep$ & $27$ & $-q^{18}$ & $54$
\\\hline\end{tabular}\end{center}

\begin{center}\begin{tikzpicture}[thick,scale=2]
\draw (2,-0.18) node {$E_6[\theta^i];1$};
\draw (0,-0.18) node {$E_6[\theta^i];\ep$};
\draw (0,0) -- (2,0);
\draw (2,0) node [draw] (l0) {};
\draw (0,0) node [draw] (l1) {};
\draw (1,0) node [fill=black!100] (ld) {};
\draw (2,0.15) node{$36$};
\draw (0,0.15) node{$54$};
\end{tikzpicture}\end{center}
\noindent\textbf{Proof of Brauer trees}: Since both characters have the same parity, they cannot be connected; this completes the proof.

\newpage

\subsection{$d=2$}
For $E_7(q)$ and $d=2$ there are two unipotent blocks of weight $1$, together with the principal block with cyclotomic Weyl group $E_7$, a block of weight $3$ with cyclotomic Weyl group $B_3$, and two unipotent blocks of defect zero.
\\[1cm]
\noindent(i)\;\;\textbf{Block 1/2:} Cuspidal pair is $(\Phi_2.{}^2\!E_6(q^2),{}^2\!E_6[\theta^i])$, of degree $q^7\Phi_1^4\Phi_2^6\Phi_4^2\Phi_8\Phi_{10}/3$. There are two unipotent characters in each block, neither of which is real, and the blocks are conjugate.
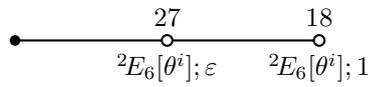
\begin{center}\begin{tabular}{lccc}
\hline Character & $A(-)$ & $\omega_i q^{aA/e}$ & $\kappa=1$
\\\hline $E_6[\theta^i];1$ & $18$ & $q^9$ & $18$
\\ $E_6[\theta^i];\ep$ & $27$ & $-q^{18}$ & $27$
\\\hline\end{tabular}\end{center}

\begin{center}\begin{tikzpicture}[thick,scale=2]
\draw (2,-0.18) node {${}^2\!E_6[\theta^i];1$};
\draw (1,-0.18) node {${}^2\!E_6[\theta^i];\ep$};
\draw (0,0) -- (2,0);
\draw (1,0) node [draw] (l0) {};
\draw (2,0) node [draw] (l1) {};
\draw (0,0) node [fill=black!100] (ld) {};
\draw (2,0.15) node{$18$};
\draw (1,0.15) node{$27$};
\end{tikzpicture}\end{center}

\noindent\textbf{Proof of Brauer trees}: Since the characters have different parities, they cannot both be connected to the exceptional node; this completes the proof.

\newpage

\subsection{$d=3$}

For $E_7(q)$ and $d=3$ there are three unipotent blocks of weight $1$, together with the principal block with cyclotomic Weyl group $G_{26}$ and ten unipotent blocks of defect zero.
\\[1cm]
\noindent(i)\;\;\textbf{Block 1:} Cuspidal pair is $(\Phi_1\Phi_3.{}^3\!D_4(q^3),{}^3\!D_4[-1])$, of degree $q^3\Phi_1^2\Phi_3^2/2$. There are six unipotent characters in the block, all of which are real.
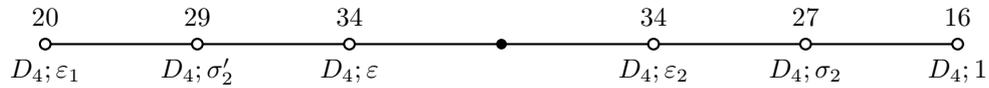
\begin{center}\begin{tabular}{lcccc}
\hline Character & $A(-)$ & $\omega_i q^{aA/e}$ & $\kappa=1$ & $\kappa=2$
\\\hline $D_4;1$ & $24$ & $q^{4}$ & $16$ & $32$
\\ $D_4;\sigma_2'$ & $44$ & $-q^{9}$ & $29$ & $59$
\\ $D_4;\sigma_2$ & $41$ & $q^{8}$ & $27$ & $55$
\\ $D_4;\ep$ & $51$ & $-q^{13}$ & $34$ & $68$
\\ $D_4;\ep_2$ & $50$ & $q^{12}$ & $34$ & $66$
\\ $D_4;\ep_1$ & $29$ & $-q^{5}$ & $20$ & $38$
\\\hline\end{tabular}\end{center}
\begin{center}\begin{tikzpicture}[thick,scale=2]
\draw (0,-0.18) node{$D_4;\ep_1$};
\draw (1,-0.18) node{$D_4;\sigma_2'$};
\draw (2,-0.18) node{$D_4;\ep$};

\draw (4,-0.18) node{$D_4;\ep_2$};
\draw (5,-0.18) node{$D_4;\sigma_2$};
\draw (6,-0.18) node{$D_4;1$};
\draw (0,0) -- (6,0);
\draw (3,0) node [fill=black!100] (ld) {};
\draw (0,0) node [draw] (l4) {};
\draw (1,0) node [draw] (l4) {};
\draw (2,0) node [draw] (l4) {};

\draw (4,0) node [draw] (l4) {};
\draw (5,0) node [draw] (l4) {};
\draw (6,0) node [draw] (l4) {};
\draw (0,0.18) node{$20$};
\draw (1,0.18) node{$29$};
\draw (2,0.18) node{$34$};
\draw (4,0.18) node{$34$};
\draw (5,0.18) node{$27$};
\draw (6,0.18) node{$16$};
\end{tikzpicture}\end{center}

\noindent(ii)\;\;\textbf{Block 2:} Cuspidal pair is $(\Phi_3.A_5(q),\phi_{42})$, of degree $q^2\Phi_3^2\Phi_6$. There are six unipotent characters in the block, all of which are real.
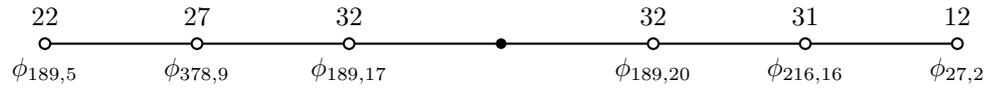
\begin{center}\begin{tabular}{lcccc}
\hline Character & $A(-)$ & $\omega_i q^{aA/e}$ & $\kappa=1$ & $\kappa=2$
\\\hline $\phi_{27,2}$ & $18$ & $q^{3}$ & $12$ & $24$
\\ $\phi_{378,9}$ & $41$ & $-q^{8}$ & $27$ & $55$
\\ $\phi_{216,16}$ & $47$ & $q^{10}$ & $31$ & $63$
\\ $\phi_{189,5}$ & $33$ & $-q^{6}$ & $22$ & $44$
\\ $\phi_{189,20}$ & $48$ & $q^{11}$ & $32$ & $64$
\\ $\phi_{189,17}$ & $47$ & $-q^{10}$ & $32$ & $62$
\\\hline\end{tabular}\end{center}
\begin{center}\begin{tikzpicture}[thick,scale=2]
\draw (0,-0.18) node{$\phi_{189,5}$};
\draw (1,-0.18) node{$\phi_{378,9}$};
\draw (2,-0.18) node{$\phi_{189,17}$};

\draw (4,-0.18) node{$\phi_{189,20}$};
\draw (5,-0.18) node{$\phi_{216,16}$};
\draw (6,-0.18) node{$\phi_{27,2}$};
\draw (0,0) -- (6,0);
\draw (3,0) node [fill=black!100] (ld) {};
\draw (0,0) node [draw] (l4) {};
\draw (1,0) node [draw] (l4) {};
\draw (2,0) node [draw] (l4) {};

\draw (4,0) node [draw] (l4) {};
\draw (5,0) node [draw] (l4) {};
\draw (6,0) node [draw] (l4) {};
\draw (0,0.18) node{$22$};
\draw (1,0.18) node{$27$};
\draw (2,0.18) node{$32$};
\draw (4,0.18) node{$32$};
\draw (5,0.18) node{$31$};
\draw (6,0.18) node{$12$};
\end{tikzpicture}\end{center}

\noindent(iii)\;\;\textbf{Block 3:} Cuspidal pair is $(\Phi_3.A_5(q),\phi_{2211})$, of degree $q^7\Phi_3^2\Phi_6$. There are six unipotent characters in the block, all of which are real.
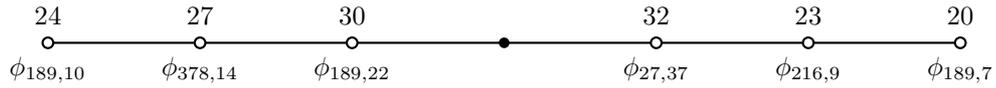
\begin{center}\begin{tabular}{lcccc}
\hline Character & $A(-)$ & $\omega_i q^{aA/e}$ & $\kappa=1$ & $\kappa=2$
\\\hline $\phi_{189,7}$ & $30$ & $q^{5}$ & $20$ & $40$
\\ $\phi_{189,22}$ & $45$ & $-q^{10}$ & $30$ & $60$
\\ $\phi_{216,9}$ & $35$ & $q^{6}$ & $23$ & $47$
\\ $\phi_{378,14}$ & $41$ & $-q^{8}$ & $27$ & $55$
\\ $\phi_{27,37}$ & $48$ & $q^{13}$ & $32$ & $64$
\\ $\phi_{189,10}$ & $35$ & $-q^{6}$ & $24$ & $46$
\\\hline\end{tabular}\end{center}
\begin{center}\begin{tikzpicture}[thick,scale=2]
\draw (0,-0.18) node{$\phi_{189,10}$};
\draw (1,-0.18) node{$\phi_{378,14}$};
\draw (2,-0.18) node{$\phi_{189,22}$};

\draw (4,-0.18) node{$\phi_{27,37}$};
\draw (5,-0.18) node{$\phi_{216,9}$};
\draw (6,-0.18) node{$\phi_{189,7}$};
\draw (0,0) -- (6,0);
\draw (3,0) node [fill=black!100] (ld) {};
\draw (0,0) node [draw] (l4) {};
\draw (1,0) node [draw] (l4) {};
\draw (2,0) node [draw] (l4) {};

\draw (4,0) node [draw] (l4) {};
\draw (5,0) node [draw] (l4) {};
\draw (6,0) node [draw] (l4) {};
\draw (0,0.18) node{$24$};
\draw (1,0.18) node{$27$};
\draw (2,0.18) node{$30$};
\draw (4,0.18) node{$32$};
\draw (5,0.18) node{$23$};
\draw (6,0.18) node{$20$};
\end{tikzpicture}\end{center}

\noindent\textbf{Proof of Brauer trees}: Since all characters are real, the Brauer trees are lines. For the second and third blocks, a Geck--Pfeiffer argument then a degree argument is enough.

For the first tree, we use a projective-induction argument, using the $D_6$ Levi subgroup. The $D_4$-series for $D_6$ consists solely of projective characters, and so any of these, Harish-Chandra induced to $E_7$, will be a projective character. Using GAP (and keeping GAP notation for the characters of $D_6$), we calculate the Harish-Chandra induced characters of certain $D_4$-series characters:
\begin{align*}
D_4;2.&\mapsto [D_4;1+D_4;\sigma_2]+(D_4;r)
\\ D_4;.2&\mapsto [D_4;\sigma_2+D_4;\ep_2]+(D_4;r\ep_2)
\\ D_4;11.&\mapsto [D_4;\ep_1+D_4;\sigma_2']+(D_4;r\ep_1)
\\ D_4;.11&\mapsto [D_4;\sigma_2'+D_4;\ep]+(D_4;r\ep)
\end{align*}
These computations prove that $D_4;1$, $D_4;\sigma_2$ and $D_4;\ep_2$ form a line in the Brauer tree, as do $D_4;\ep_1$, $D_4;\sigma_2'$ and $D_4;\ep$. A degree argument now completes the proof, by determining to which characters the exceptional node is connected.

\newpage
\subsection{$d=5$}

For $E_7(q)$ and $d=5$ there are three unipotent blocks of weight $1$, together with 46 unipotent blocks of defect zero.
\\[1cm]
\noindent(i)\;\;\textbf{Block 1:} Cuspidal pair is $(\Phi_5.A_2(q),\phi_{3})$, of degree $1$. There are ten unipotent characters in the block, all of which are real.

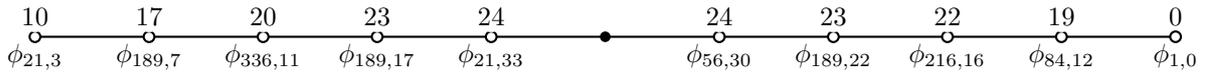
\begin{center}\begin{tabular}{lcccccc}
\hline Character & $A(-)$ & $\omega_i q^{aA/e}$ & $\kappa=1$ & $\kappa=2$ & $\kappa=3$ & $\kappa=4$
\\\hline $\phi_{1,0}$ & $0$ & $q^{0}$ & $0$ & $0$ & $0$ & $0$
\\ $\phi_{21,3}$ & $27$ & $-q^{3}$ & $10$ & $22$ & $32$ & $44$
\\ $\phi_{84,12}$ & $50$ & $q^{6}$ & $19$ & $41$ & $59$ & $81$
\\ $\phi_{21,33}$ & $60$ & $-q^{9}$ & $24$ & $48$ & $72$ & $96$
\\ $\phi_{216,16}$ & $55$ & $q^{7}$ & $22$ & $44$ & $66$ & $88$
\\ $\phi_{189,7}$ & $43$ & $-q^{5}$ & $17$ & $35$ & $51$ & $69$
\\ $\phi_{189,22}$ & $58$ & $q^{8}$ & $23$ & $47$ & $69$ & $93$
\\ $\phi_{336,11}$ & $50$ & $-q^{6}$ & $20$ & $40$ & $60$ & $80$
\\ $\phi_{56,30}$ & $60$ & $q^{9}$ & $24$ & $48$ & $72$ & $96$
\\ $\phi_{189,17}$ & $55$ & $-q^{7}$ & $23$ & $45$ & $65$ & $87$
\\\hline\end{tabular}\end{center}
\begin{center}\begin{tikzpicture}[thick,scale=1.5]
\draw (0,-0.18) node{$\phi_{21,3}$};
\draw (1,-0.18) node{$\phi_{189,7}$};
\draw (2,-0.18) node{$\phi_{336,11}$};
\draw (3,-0.18) node{$\phi_{189,17}$};
\draw (4,-0.18) node{$\phi_{21,33}$};

\draw (6,-0.18) node{$\phi_{56,30}$};
\draw (7,-0.18) node{$\phi_{189,22}$};
\draw (8,-0.18) node{$\phi_{216,16}$};
\draw (9,-0.18) node{$\phi_{84,12}$};
\draw (10,-0.18) node{$\phi_{1,0}$};
\draw (0,0) -- (10,0);
\draw (5,0) node [fill=black!100] (ld) {};
\draw (0,0) node [draw] (l4) {};
\draw (1,0) node [draw] (l4) {};
\draw (2,0) node [draw] (l4) {};
\draw (3,0) node [draw] (l4) {};
\draw (4,0) node [draw] (l4) {};
\draw (6,0) node [draw] (l4) {};
\draw (7,0) node [draw] (l4) {};
\draw (8,0) node [draw] (l4) {};
\draw (9,0) node [draw] (l4) {};
\draw (10,0) node [draw] (l4) {};

\draw (0,0.18) node{$10$};
\draw (1,0.18) node{$17$};
\draw (2,0.18) node{$20$};
\draw (3,0.18) node{$23$};
\draw (4,0.18) node{$24$};
\draw (6,0.18) node{$24$};
\draw (7,0.18) node{$23$};
\draw (8,0.18) node{$22$};
\draw (9,0.18) node{$19$};
\draw (10,0.18) node{$0$};
\end{tikzpicture}\end{center}

\noindent(ii)\;\;\textbf{Block 2:} Cuspidal pair is $(\Phi_5.A_2(q),\phi_{21})$, of degree $q\Phi_2$. There are ten unipotent characters in the block, all of which are real.
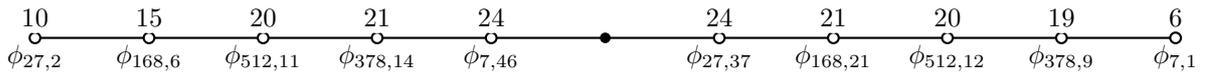
\begin{center}\begin{tabular}{lcccccc}
\hline Character & $A(-)$ & $\omega_i q^{aA/e}$ & $\kappa=1$ & $\kappa=2$ & $\kappa=3$ & $\kappa=4$
\\\hline $\phi_{7,1}$ & $15$ & $q^{3/2}$ & $6$ & $12$ & $18$ & $24$
\\ $\phi_{168,6}$ & $40$ & $-q^{9/2}$ & $15$ & $33$ & $47$ & $65$
\\ $\phi_{168,21}$ & $55$ & $q^{15/2}$ & $21$ & $45$ & $65$ & $89$
\\ $\phi_{7,46}$ & $60$ & $-q^{21/2}$ & $24$ & $48$ & $72$ & $96$
\\ $\phi_{512,11}$ & $50$ & $-q^{6}$ & $20$ & $40$ & $60$ & $80$
\\ $\phi_{378,14}$ & $52$ & $-q^{13/2}$ & $21$ & $43$ & $61$ & $83$
\\ $\phi_{27,37}$ & $59$ & $q^{19/2}$ & $24$ & $48$ & $70$ & $94$
\\ $\phi_{27,2}$ & $24$ & $-q^{5/2}$ & $10$ & $20$ & $28$ & $38$
\\ $\phi_{378,9}$ & $47$ & $q^{11/2}$ & $19$ & $39$ & $55$ & $75$
\\ $\phi_{512,12}$ & $50$ & $q^{6}$ & $20$ & $40$ & $60$ & $80$
\\\hline\end{tabular}\end{center}
\begin{center}\begin{tikzpicture}[thick,scale=1.5]
\draw (0,-0.18) node{$\phi_{27,2}$};
\draw (1,-0.18) node{$\phi_{168,6}$};
\draw (2,-0.18) node{$\phi_{512,11}$};
\draw (3,-0.18) node{$\phi_{378,14}$};
\draw (4,-0.18) node{$\phi_{7,46}$};

\draw (6,-0.18) node{$\phi_{27,37}$};
\draw (7,-0.18) node{$\phi_{168,21}$};
\draw (8,-0.18) node{$\phi_{512,12}$};
\draw (9,-0.18) node{$\phi_{378,9}$};
\draw (10,-0.18) node{$\phi_{7,1}$};
\draw (0,0) -- (10,0);
\draw (5,0) node [fill=black!100] (ld) {};
\draw (0,0) node [draw] (l4) {};
\draw (1,0) node [draw] (l4) {};
\draw (2,0) node [draw] (l4) {};
\draw (3,0) node [draw] (l4) {};
\draw (4,0) node [draw] (l4) {};
\draw (6,0) node [draw] (l4) {};
\draw (7,0) node [draw] (l4) {};
\draw (8,0) node [draw] (l4) {};
\draw (9,0) node [draw] (l4) {};
\draw (10,0) node [draw] (l4) {};

\draw (0,0.18) node{$10$};
\draw (1,0.18) node{$15$};
\draw (2,0.18) node{$20$};
\draw (3,0.18) node{$21$};
\draw (4,0.18) node{$24$};
\draw (6,0.18) node{$24$};
\draw (7,0.18) node{$21$};
\draw (8,0.18) node{$20$};
\draw (9,0.18) node{$19$};
\draw (10,0.18) node{$6$};
\end{tikzpicture}\end{center}

\noindent(iii)\;\;\textbf{Block 3:} Cuspidal pair is $(\Phi_5.A_2(q),\phi_{111})$, of degree $q^3$. There are ten unipotent characters in the block, all of which are real.
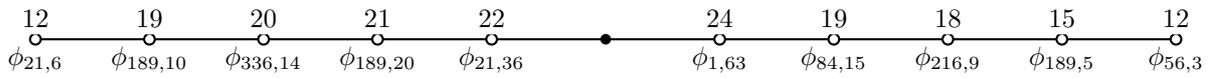
\begin{center}\begin{tabular}{lcccccc}
\hline Character & $A(-)$ & $\omega_i q^{aA/e}$ & $\kappa=1$ & $\kappa=2$ & $\kappa=3$ & $\kappa=4$
\\\hline $\phi_{56,3}$ & $30$ & $q^{3}$ & $12$ & $24$ & $36$ & $48$
\\ $\phi_{336,14}$ & $50$ & $-q^{6}$ & $20$ & $40$ & $60$ & $80$
\\ $\phi_{189,5}$ & $38$ & $q^{4}$ & $15$ & $31$ & $45$ & $61$
\\ $\phi_{189,20}$ & $53$ & $-q^{7}$ & $21$ & $43$ & $63$ & $85$
\\ $\phi_{216,9}$ & $45$ & $q^{5}$ & $18$ & $36$ & $54$ & $72$
\\ $\phi_{21,6}$ & $30$ & $-q^{3}$ & $12$ & $24$ & $36$ & $48$
\\ $\phi_{84,15}$ & $50$ & $q^{6}$ & $19$ & $41$ & $59$ & $81$
\\ $\phi_{21,36}$ & $57$ & $-q^{9}$ & $22$ & $46$ & $68$ & $92$
\\ $\phi_{1,63}$ & $60$ & $q^{12}$ & $24$ & $48$ & $72$ & $96$
\\ $\phi_{189,10}$ & $45$ & $-q^{5}$ & $19$ & $37$ & $53$ & $71$
\\\hline\end{tabular}\end{center}
\begin{center}\begin{tikzpicture}[thick,scale=1.5]
\draw (0,-0.18) node{$\phi_{21,6}$};
\draw (1,-0.18) node{$\phi_{189,10}$};
\draw (2,-0.18) node{$\phi_{336,14}$};
\draw (3,-0.18) node{$\phi_{189,20}$};
\draw (4,-0.18) node{$\phi_{21,36}$};

\draw (6,-0.18) node{$\phi_{1,63}$};
\draw (7,-0.18) node{$\phi_{84,15}$};
\draw (8,-0.18) node{$\phi_{216,9}$};
\draw (9,-0.18) node{$\phi_{189,5}$};
\draw (10,-0.18) node{$\phi_{56,3}$};
\draw (0,0) -- (10,0);
\draw (5,0) node [fill=black!100] (ld) {};
\draw (0,0) node [draw] (l4) {};
\draw (1,0) node [draw] (l4) {};
\draw (2,0) node [draw] (l4) {};
\draw (3,0) node [draw] (l4) {};
\draw (4,0) node [draw] (l4) {};
\draw (6,0) node [draw] (l4) {};
\draw (7,0) node [draw] (l4) {};
\draw (8,0) node [draw] (l4) {};
\draw (9,0) node [draw] (l4) {};
\draw (10,0) node [draw] (l4) {};

\draw (0,0.18) node{$12$};
\draw (1,0.18) node{$19$};
\draw (2,0.18) node{$20$};
\draw (3,0.18) node{$21$};
\draw (4,0.18) node{$22$};
\draw (6,0.18) node{$24$};
\draw (7,0.18) node{$19$};
\draw (8,0.18) node{$18$};
\draw (9,0.18) node{$15$};
\draw (10,0.18) node{$12$};
\end{tikzpicture}\end{center}

\noindent\textbf{Proof of Brauer trees}: Since all characters are real, the Brauer trees are lines. A Geck--Pfeiffer argument then a degree argument is enough.

\newpage
\subsection{$d=6$}

For $E_7(q)$ and $d=6$ there are three unipotent blocks of weight $1$, together with the principal block, with cyclotomic Weyl group $G_{26}$, and ten unipotent blocks of defect zero.
\\[1cm]
\noindent(i)\;\;\textbf{Block 1:} Cuspidal pair is $(\Phi_2\Phi_6.{}^3\!D_4(q^3),\phi_{2,1})$, of degree $q^3\Phi_2^2\Phi_6^2/2$. There are six unipotent characters in the block, all of which are real.
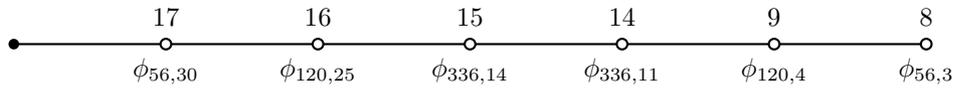
\begin{center}
\begin{tabular}{lcccc}
\hline Character & $A(-)$ & $\omega_i q^{aA/e}$ & $\kappa=1$ & $\kappa=5$
\\\hline $\phi_{56,3}$ & $24$ & $q^{4}$ & $8$ & $40$
\\ $\phi_{120,4}$ & $29$ & $q^{5}$ & $9$ & $49$
\\ $\phi_{120,25}$ & $50$ & $q^{12}$ & $16$ & $84$
\\ $\phi_{56,30}$ & $51$ & $q^{13}$ & $17$ & $85$
\\ $\phi_{336,11}$ & $41$ & $q^{8}$ & $14$ & $68$
\\ $\phi_{336,14}$ & $44$ & $q^{9}$ & $15$ & $73$
\\\hline\end{tabular}%
\end{center}

\noindent(ii)\;\;\textbf{Block 2:} Cuspidal pair is $(\Phi_6.{}^2\!A_5(q^2),\phi_{42})$, of degree $q^2\Phi_3\Phi_6^2$. There are six unipotent characters in the block, all of which are real.
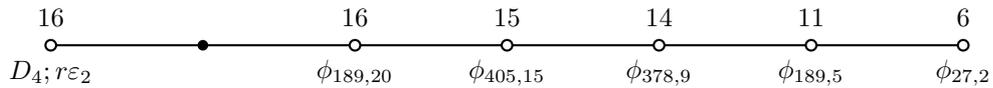
\begin{center}\begin{tabular}{lcccc}
\hline Character & $A(-)$ & $\omega_i q^{aA/e}$ & $\kappa=1$ & $\kappa=5$
\\\hline $\phi_{27,2}$ & $18$ & $q^{3}$ & $6$ & $30$
\\ $\phi_{405,15}$ & $47$ & $q^{10}$ & $15$ & $79$
\\ $\phi_{189,20}$ & $48$ & $q^{11}$ & $16$ & $80$
\\ $\phi_{189,5}$ & $33$ & $q^{6}$ & $11$ & $55$
\\ $D_4;r\ep_2$ & $47$ & $-q^{10}$ & $16$ & $78$
\\ $\phi_{378,9}$ & $41$ & $q^{8}$ & $14$ & $68$
\\\hline\end{tabular}%
\end{center}

\noindent(iii)\;\;\textbf{Block 3:} Cuspidal pair is $(\Phi_6.{}^2\!A_5(q^2),\phi_{2211})$, of degree $q^7\Phi_3\Phi_6^2$. There are six unipotent characters in the block, all of which are real.
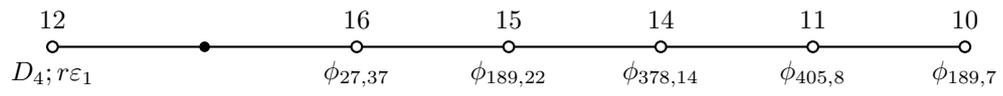
\begin{center}\begin{tabular}{lcccc}
\hline Character & $A(-)$ & $\omega_i q^{aA/e}$ & $\kappa=1$ & $\kappa=5$
\\\hline $\phi_{189,7}$ & $30$ & $q^{5}$ & $10$ & $50$
\\ $\phi_{405,8}$ & $35$ & $q^{6}$ & $11$ & $59$
\\ $\phi_{27,37}$ & $48$ & $q^{13}$ & $16$ & $80$
\\ $\phi_{378,14}$ & $41$ & $q^{8}$ & $14$ & $68$
\\ $D_4;r\ep_1$ & $35$ & $-q^{6}$ & $12$ & $58$
\\ $\phi_{189,22}$ & $45$ & $q^{10}$ & $15$ & $75$
\\\hline\end{tabular}%
\end{center}

\noindent\textbf{Proof of Brauer trees}: Since all characters are real, the Brauer trees are lines. A Geck--Pfeiffer argument then a degree argument is enough.

\newpage

\subsection{$d=7$}

For $E_7(q)$ and $d=7$ there is a single unipotent block of weight $1$, together with 62 unipotent blocks of defect zero.
\\[1cm]
\noindent(i)\;\;\textbf{Block 1:} Cuspidal pair is $(\Phi_1\Phi_6,1)$, of degree $1$. There are fourteen unipotent characters in the block, all of which are real.
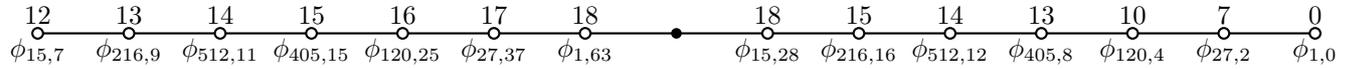
\begin{center}\begin{tabular}{lcccccccc}
\hline Character & $A(-)$ & $\omega_i q^{aA/e}$ & $\kappa=1$ & $\kappa=2$ & $\kappa=3$ & $\kappa=4$ & $\kappa=5$ & $\kappa=6$
\\\hline $\phi_{1,0}$ & $0$ & $q^{0}$ & $0$ & $0$ & $0$ & $0$ & $0$ & $0$
\\ $\phi_{216,9}$ & $48$ & $-q^{4}$ & $13$ & $27$ & $41$ & $55$ & $69$ & $83$
\\ $\phi_{512,11}$ & $52$ & $-q^{9/2}$ & $14$ & $30$ & $44$ & $60$ & $74$ & $90$
\\ $\phi_{405,15}$ & $55$ & $-q^{5}$ & $15$ & $31$ & $49$ & $61$ & $79$ & $95$
\\ $\phi_{27,2}$ & $26$ & $q^{2}$ & $7$ & $15$ & $23$ & $29$ & $37$ & $45$
\\ $\phi_{120,25}$ & $59$ & $-q^{6}$ & $16$ & $34$ & $50$ & $68$ & $84$ & $102$
\\ $\phi_{120,4}$ & $38$ & $q^{3}$ & $10$ & $22$ & $32$ & $44$ & $54$ & $66$
\\ $\phi_{27,37}$ & $61$ & $-q^{7}$ & $17$ & $35$ & $53$ & $69$ & $87$ & $105$
\\ $\phi_{405,8}$ & $48$ & $q^{4}$ & $13$ & $27$ & $43$ & $53$ & $69$ & $83$
\\ $\phi_{512,12}$ & $52$ & $q^{9/2}$ & $14$ & $30$ & $44$ & $60$ & $74$ & $90$
\\ $\phi_{216,16}$ & $55$ & $q^{5}$ & $15$ & $31$ & $47$ & $63$ & $79$ & $95$
\\ $\phi_{1,63}$ & $63$ & $-q^{9}$ & $18$ & $36$ & $54$ & $72$ & $90$ & $108$
\\ $\phi_{15,28}$ & $59$ & $q^{6}$ & $18$ & $34$ & $52$ & $66$ & $84$ & $100$
\\ $\phi_{15,7}$ & $38$ & $-q^{3}$ & $12$ & $22$ & $34$ & $42$ & $54$ & $64$
\\\hline\end{tabular}%
\end{center}

\noindent\textbf{Proof of Brauer tree}: Since all characters are real, the Brauer tree is a line. A Geck--Pfeiffer argument then a degree argument is enough.

\newpage

\subsection{$d=8$}

For $E_7(q)$ and $d=8$ there are four unipotent blocks of weight $1$, together with 44 unipotent blocks of defect zero.
\\[1cm]
\noindent(i)\;\;\textbf{Block 1:} Cuspidal pair is $(\Phi_8.A_1(q^2)A_1(q),\phi_2^2)$, of degree $1$. There are eight unipotent characters in the block, all of which are real.
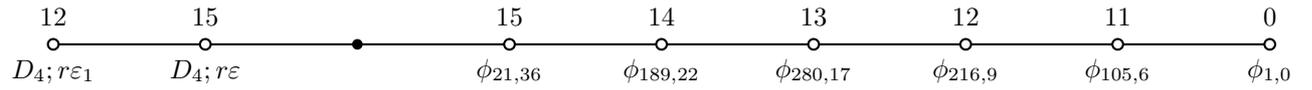
\begin{center}\begin{tabular}{lcccccc}
\hline Character & $A(-)$ & $\omega_i q^{aA/e}$ & $\kappa=1$ & $\kappa=3$ & $\kappa=5$ & $\kappa=7$
\\\hline $\phi_{1,0}$ & $0$ & $q^{0}$ & $0$ & $0$ & $0$ & $0$
\\ $\phi_{280,17}$ & $56$ & $q^{9}$ & $13$ & $41$ & $71$ & $99$
\\ $\phi_{189,22}$ & $58$ & $q^{10}$ & $14$ & $44$ & $72$ & $102$
\\ $D_4;r\ep_1$ & $48$ & $-q^{7}$ & $12$ & $36$ & $60$ & $84$
\\ $\phi_{21,36}$ & $60$ & $q^{12}$ & $15$ & $45$ & $75$ & $105$
\\ $D_4;r\ep$ & $56$ & $-q^{9}$ & $15$ & $43$ & $69$ & $97$
\\ $\phi_{105,6}$ & $42$ & $q^{6}$ & $11$ & $31$ & $53$ & $73$
\\ $\phi_{216,9}$ & $48$ & $q^{7}$ & $12$ & $36$ & $60$ & $84$
\\\hline\end{tabular}%
\end{center}

\noindent(ii)\;\;\textbf{Block 2:} Cuspidal pair is $(\Phi_8.A_1(q^2)A_1(q),\phi_2\phi_{11})$, of degree $q$. There are eight unipotent characters in the block, all of which are real.
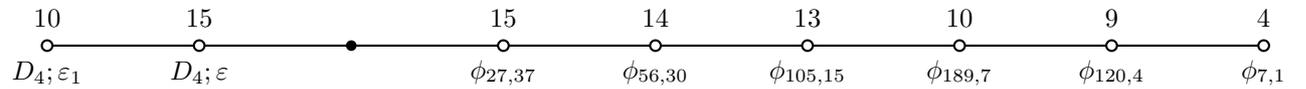
\begin{center}\begin{tabular}{lcccccc}
\hline Character & $A(-)$ & $\omega_i q^{aA/e}$ & $\kappa=1$ & $\kappa=3$ & $\kappa=5$ & $\kappa=7$
\\\hline $\phi_{7,1}$ & $16$ & $q^{2}$ & $4$ & $12$ & $20$ & $28$
\\ $\phi_{56,30}$ & $59$ & $q^{11}$ & $14$ & $44$ & $74$ & $104$
\\ $\phi_{27,37}$ & $60$ & $q^{12}$ & $15$ & $45$ & $75$ & $105$
\\ $\phi_{120,4}$ & $37$ & $q^{5}$ & $9$ & $27$ & $47$ & $65$
\\ $\phi_{189,7}$ & $42$ & $q^{6}$ & $10$ & $32$ & $52$ & $74$
\\ $D_4;\ep$ & $59$ & $-q^{11}$ & $15$ & $45$ & $73$ & $103$
\\ $\phi_{105,15}$ & $50$ & $q^{8}$ & $13$ & $37$ & $63$ & $87$
\\ $D_4;\ep_1$ & $37$ & $-q^{5}$ & $10$ & $28$ & $46$ & $64$
\\\hline\end{tabular}%
\end{center}

\noindent(iii)\;\;\textbf{Block 3:} Cuspidal pair is $(\Phi_8.A_1(q^2)A_1(q),\phi_{11}\phi_2)$, of degree $q^2$. There are eight unipotent characters in the block, all of which are real.
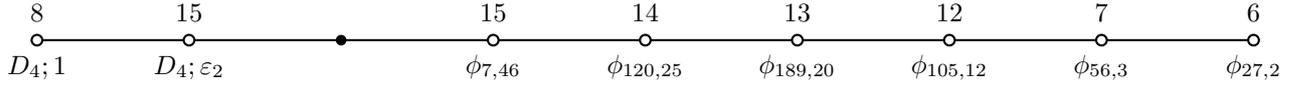
\begin{center}\begin{tabular}{lcccccc}
\hline Character & $A(-)$ & $\omega_i q^{aA/e}$ & $\kappa=1$ & $\kappa=3$ & $\kappa=5$ & $\kappa=7$
\\\hline $\phi_{27,2}$ & $24$ & $q^{3}$ & $6$ & $18$ & $30$ & $42$
\\ $\phi_{56,3}$ & $31$ & $q^{4}$ & $7$ & $23$ & $39$ & $55$
\\ $\phi_{7,46}$ & $60$ & $q^{13}$ & $15$ & $45$ & $75$ & $105$
\\ $D_4;\ep_2$ & $57$ & $-q^{10}$ & $15$ & $43$ & $71$ & $99$
\\ $\phi_{105,12}$ & $46$ & $q^{7}$ & $12$ & $34$ & $58$ & $80$
\\ $D_4;1$ & $31$ & $-q^{4}$ & $8$ & $24$ & $38$ & $54$
\\ $\phi_{189,20}$ & $54$ & $q^{9}$ & $13$ & $41$ & $67$ & $95$
\\ $\phi_{120,25}$ & $57$ & $q^{10}$ & $14$ & $42$ & $72$ & $100$
\\\hline\end{tabular}%
\end{center}

\noindent(iv)\;\;\textbf{Block 4:} Cuspidal pair is $(\Phi_8.A_1(q^2)A_1(q),\phi_{11}^2)$, of degree $q^3$. There are eight unipotent characters in the block, all of which are real.
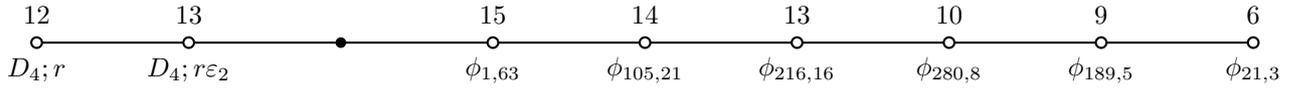
\begin{center}\begin{tabular}{lcccccc}
\hline Character & $A(-)$ & $\omega_i q^{aA/e}$ & $\kappa=1$ & $\kappa=3$ & $\kappa=5$ & $\kappa=7$
\\ $\phi_{21,3}$ & $24$ & $q^{3}$ & $6$ & $18$ & $30$ & $42$
\\ $D_4;r\ep_2$ & $52$ & $-q^{8}$ & $13$ & $39$ & $65$ & $91$
\\ $\phi_{189,5}$ & $38$ & $q^{5}$ & $9$ & $29$ & $47$ & $67$
\\ $\phi_{280,8}$ & $44$ & $q^{6}$ & $10$ & $32$ & $56$ & $78$
\\ $\phi_{1,63}$ & $60$ & $q^{15}$ & $15$ & $45$ & $75$ & $105$
\\ $\phi_{216,16}$ & $52$ & $q^{8}$ & $13$ & $39$ & $65$ & $91$
\\ $\phi_{105,21}$ & $54$ & $q^{9}$ & $14$ & $40$ & $68$ & $94$
\\ $D_4;r$ & $44$ & $-q^{6}$ & $12$ & $34$ & $54$ & $76$
\\\hline\end{tabular}%
\end{center}

\noindent\textbf{Proof of Brauer trees}: Since all characters are real, the Brauer tree is a line. A Geck--Pfeiffer argument then a degree argument is enough.

\newpage

\subsection{$d=9$}

For $E_7(q)$ and $d=9$ there is a single unipotent block of weight $1$, together with 58 unipotent blocks of defect zero.
\\[1cm]
\noindent(i)\;\;\textbf{Block 1:} Cuspidal pair is $(\Phi_1\Phi_9,1)$, of degree $1$. There are eighteen unipotent characters in the block, four of which -- $E_6[\theta^i];1$ and $E_6[\theta^i];\ep$ -- are non-real.
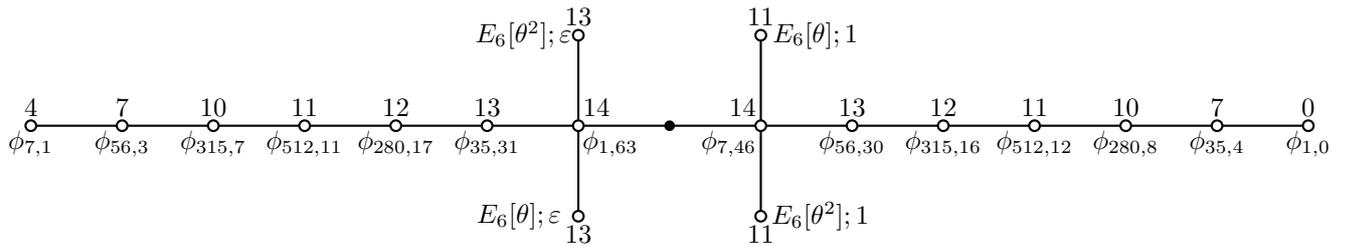
\begin{center}\begin{tabular}{lcccccccc}
\hline Character & $A(-)$ & $\omega_i q^{aA/e}$ & $\kappa=1$ & $\kappa=2$ & $\kappa=4$ & $\kappa=5$ & $\kappa=7$ & $\kappa=8$
\\\hline $\phi_{1,0}$ & $0$ & $q^{0}$ & $0$ & $0$ & $0$ & $0$ & $0$ & $0$
\\ $\phi_{35,31}$ & $60$ & $-q^{5}$ & $13$ & $27$ & $55$ & $65$ & $93$ & $107$
\\ $E_6[\theta];\ep$ & $56$ & $\theta^2 q^{4}$ & $13$ & $25$ & $49$ & $63$ & $87$ & $99$
\\ $E_6[\theta^2];1$ & $47$ & $-\theta q^{3}$ & $11$ & $21$ & $41$ & $53$ & $73$ & $83$
\\ $\phi_{35,4}$ & $33$ & $q^{2}$ & $7$ & $15$ & $31$ & $35$ & $51$ & $59$
\\ $\phi_{1,63}$ & $63$ & $-q^{7}$ & $14$ & $28$ & $56$ & $70$ & $98$ & $112$
\\ $\phi_{280,8}$ & $47$ & $q^{3}$ & $10$ & $22$ & $42$ & $52$ & $72$ & $84$
\\ $\phi_{512,12}$ & $52$ & $q^{7/2}$ & $11$ & $23$ & $45$ & $59$ & $81$ & $93$
\\ $\phi_{315,16}$ & $56$ & $q^{4}$ & $12$ & $24$ & $52$ & $60$ & $88$ & $100$
\\ $E_6[\theta];1$ & $47$ & $-\theta^2 q^{3}$ & $11$ & $21$ & $41$ & $53$ & $73$ & $83$
\\ $\phi_{56,30}$ & $60$ & $q^{5}$ & $13$ & $27$ & $53$ & $67$ & $93$ & $107$
\\ $\phi_{7,1}$ & $17$ & $-q$ & $4$ & $8$ & $16$ & $18$ & $26$ & $30$
\\ $\phi_{7,46}$ & $62$ & $q^{6}$ & $14$ & $28$ & $56$ & $68$ & $96$ & $110$
\\ $\phi_{56,3}$ & $33$ & $-q^{2}$ & $7$ & $15$ & $29$ & $37$ & $51$ & $59$
\\ $E_6[\theta^2];\ep$ & $56$ & $\theta q^{4}$ & $13$ & $25$ & $49$ & $63$ & $87$ & $99$
\\ $\phi_{315,7}$ & $47$ & $-q^{3}$ & $10$ & $20$ & $44$ & $50$ & $74$ & $84$
\\ $\phi_{512,11}$ & $52$ & $-q^{7/2}$ & $11$ & $23$ & $45$ & $59$ & $81$ & $93$
\\ $\phi_{280,17}$ & $56$ & $-q^{4}$ & $12$ & $26$ & $50$ & $62$ & $86$ & $100$
\\\hline\end{tabular}%
\end{center}

\noindent\textbf{Proof of Brauer tree}: This is given in \cite{cdr2012un}.

\newpage

\subsection{$d=10$}

For $E_7(q)$ and $d=10$ there are three unipotent blocks of weight $1$, together with $46$ unipotent blocks of defect zero.
\\[1cm]
\noindent(i)\;\; \textbf{Block 1:} Cuspidal pair is $(\Phi_{10}.{}^2\!A_2(q),\phi_{3})$, of degree $1$. There are ten unipotent characters in the block, all of which are real.

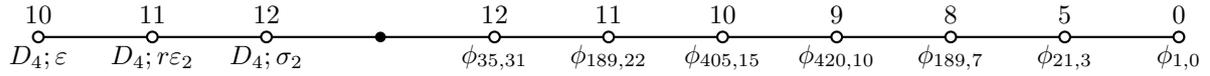
\begin{center}\begin{tabular}{lcccccc}
\hline Character & $A(-)$ & $\omega_i q^{aA/e}$ & $\kappa=1$ & $\kappa=3$ & $\kappa=7$ & $\kappa=9$
\\\hline $\phi_{1,0}$ & $0$ & $q^{0}$ & $0$ & $0$ & $0$ & $0$
\\ $D_4;\sigma_2$ & $50$ & $-q^{6}$ & $10$ & $30$ & $70$ & $90$
\\ $D_4;r\ep_2$ & $55$ & $-q^{7}$ & $11$ & $33$ & $77$ & $99$
\\ $\phi_{21,3}$ & $27$ & $q^{3}$ & $5$ & $17$ & $37$ & $49$
\\ $D_4;\ep$ & $60$ & $-q^{9}$ & $12$ & $36$ & $84$ & $108$
\\ $\phi_{189,7}$ & $43$ & $q^{5}$ & $8$ & $26$ & $60$ & $78$
\\ $\phi_{420,10}$ & $50$ & $q^{6}$ & $9$ & $31$ & $69$ & $91$
\\ $\phi_{405,15}$ & $55$ & $q^{7}$ & $10$ & $32$ & $78$ & $100$
\\ $\phi_{189,22}$ & $58$ & $q^{8}$ & $11$ & $35$ & $81$ & $105$
\\ $\phi_{35,31}$ & $60$ & $q^{9}$ & $12$ & $36$ & $84$ & $108$
\\\hline\end{tabular}%
\end{center}

\noindent(ii)\;\; \textbf{Block 2:} Cuspidal pair is $(\Phi_{10}.{}^2\!A_2(q),\phi_{21})$, of degree $q\Phi_1$. There are ten unipotent characters in the block, two of which -- $E_7[\pm\I]$ --  are non-real.

\begin{center}\begin{tabular}{lcccccc}
\hline Character & $A(-)$ & $\omega_i q^{aA/e}$ & $\kappa=1$ & $\kappa=3$ & $\kappa=7$ & $\kappa=9$
\\\hline $\phi_{7,1}$ & $15$ & $q^{3/2}$ & $3$ & $9$ & $21$ & $27$
\\ $\phi_{27,2}$ & $24$ & $q^{5/2}$ & $4$ & $14$ & $34$ & $44$
\\ $E_7[-\I]$ & $50$ & $-\I q^{6}$ & $10$ & $30$ & $70$ & $90$
\\ $\phi_{168,6}$ & $40$ & $q^{9/2}$ & $7$ & $25$ & $55$ & $73$
\\ $\phi_{378,9}$ & $47$ & $q^{11/2}$ & $8$ & $28$ & $66$ & $86$
\\ $\phi_{378,14}$ & $52$ & $q^{13/2}$ & $9$ & $31$ & $73$ & $95$
\\ $\phi_{168,21}$ & $55$ & $q^{15/2}$ & $10$ & $34$ & $76$ & $100$
\\ $E_7[\I]$ & $50$ & $\I q^{6}$ & $10$ & $30$ & $70$ & $90$
\\ $\phi_{27,37}$ & $59$ & $q^{19/2}$ & $11$ & $35$ & $83$ & $107$
\\ $\phi_{7,46}$ & $60$ & $q^{21/2}$ & $12$ & $36$ & $84$ & $108$
\\\hline\end{tabular}\end{center}
\textbf{WARNING:} this tree is conjectural. The real stem is definitely correct, but the precise location of the pair of non-real characters is not currently known.

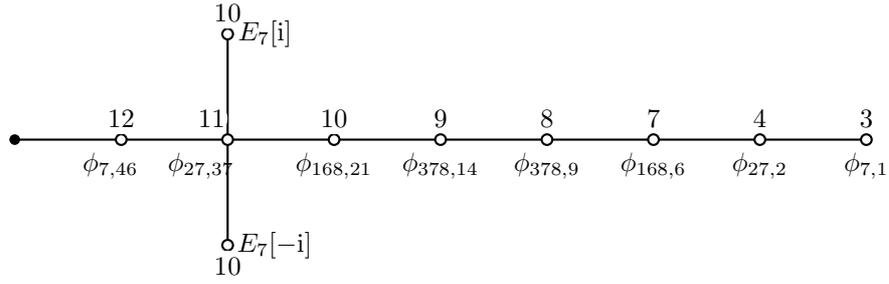
\begin{center}%
\end{center}

\noindent(iii)\;\; \textbf{Block 3:} Cuspidal pair is $(\Phi_{10}.{}^2\!A_2(q),\phi_{111})$, of degree $q^3$. There are ten unipotent characters in the block, all of which are real.

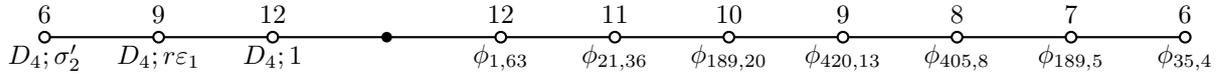
\begin{center}\begin{tabular}{lcccccc}
\hline Character & $A(-)$ & $\omega_i q^{aA/e}$ & $\kappa=1$ & $\kappa=3$ & $\kappa=7$ & $\kappa=9$
\\\hline $\phi_{35,4}$ & $30$ & $q^{3}$ & $6$ & $18$ & $42$ & $54$
\\ $\phi_{189,5}$ & $38$ & $q^{4}$ & $7$ & $23$ & $53$ & $69$
\\ $\phi_{405,8}$ & $45$ & $q^{5}$ & $8$ & $26$ & $64$ & $82$
\\ $\phi_{420,13}$ & $50$ & $q^{6}$ & $9$ & $31$ & $69$ & $91$
\\ $\phi_{189,20}$ & $53$ & $q^{7}$ & $10$ & $32$ & $74$ & $96$
\\ $D_4;1$ & $30$ & $-q^{3}$ & $6$ & $18$ & $42$ & $54$
\\ $\phi_{21,36}$ & $57$ & $q^{9}$ & $11$ & $35$ & $79$ & $103$
\\ $D_4;r\ep_1$ & $45$ & $-q^{5}$ & $9$ & $27$ & $63$ & $81$
\\ $D_4;\sigma_2'$ & $50$ & $-q^{6}$ & $10$ & $30$ & $70$ & $90$
\\ $\phi_{1,63}$ & $60$ & $q^{12}$ & $12$ & $36$ & $84$ & $108$
\\\hline\end{tabular}%
\end{center}

\noindent\textbf{Proof of Brauer trees}: For the first and third trees, since all characters are real, the Brauer tree is a line. A Geck--Pfeiffer argument then a degree argument is enough. For the second tree, see \cite{cdr2012un}.

\newpage
\subsection{$d=12$}

For $E_7(q)$ and $d=12$ there are two unipotent blocks of weight $1$, together with $52$ unipotent blocks of defect zero.
\\[1cm]
\noindent(i)\;\; \textbf{Block 1:} Cuspidal pair is $(\Phi_{12}.A_1(q^3),\phi_{2})$, of degree $1$. There are twelve unipotent characters in the block, two of which -- $E_6[\theta^i];\ep$ -- are non-real.

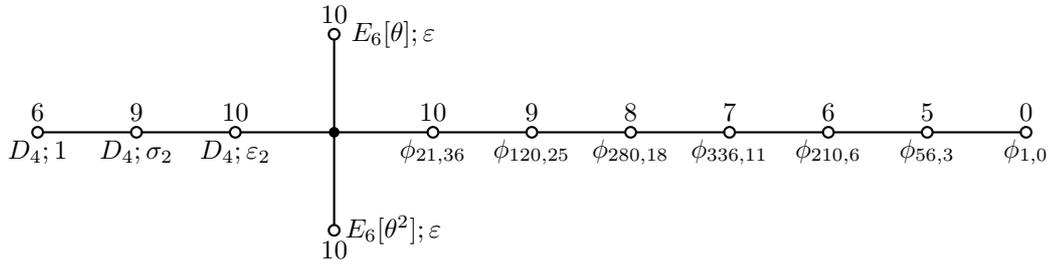
\begin{center}\begin{tabular}{lcccccc}
\hline Character & $A(-)$ & $\omega_i q^{aA/e}$ & $\kappa=1$ & $\kappa=5$ & $\kappa=7$ & $\kappa=11$
\\\hline $\phi_{1,0}$ & $0$ & $q^{0}$ & $0$ & $0$ & $0$ & $0$
\\ $D_4;\ep_2$ & $59$ & $-q^{7}$ & $10$ & $50$ & $68$ & $108$
\\ $E_6[\theta^2];\ep$ & $56$ & $\theta^2 q^{6}$ & $10$ & $46$ & $66$ & $102$
\\ $\phi_{56,3}$ & $33$ & $q^{3}$ & $5$ & $27$ & $39$ & $61$
\\ $\phi_{210,6}$ & $42$ & $q^{4}$ & $6$ & $36$ & $48$ & $78$
\\ $\phi_{336,11}$ & $50$ & $q^{5}$ & $7$ & $41$ & $59$ & $93$
\\ $\phi_{280,18}$ & $56$ & $q^{6}$ & $8$ & $48$ & $64$ & $104$
\\ $\phi_{120,25}$ & $59$ & $q^{7}$ & $9$ & $49$ & $69$ & $109$
\\ $\phi_{21,36}$ & $60$ & $q^{8}$ & $10$ & $50$ & $70$ & $110$
\\ $D_4;1$ & $33$ & $-q^{3}$ & $6$ & $28$ & $38$ & $60$
\\ $E_6[\theta];\ep$ & $56$ & $\theta q^{6}$ & $10$ & $46$ & $66$ & $102$
\\ $D_4;\sigma_2$ & $50$ & $-q^{5}$ & $9$ & $43$ & $57$ & $91$
\\\hline\end{tabular}\end{center}

\begin{center}\begin{tikzpicture}[thick,scale=1.3]
\draw (0,-0.2) node{$D_4;1$};
\draw (1,-0.2) node{$D_4;\sigma_2$};
\draw (2,-0.2) node{$D_4;\ep_2$};
\draw (4,-0.2) node{$\phi_{21,36}$};
\draw (5,-0.2) node{$\phi_{120,25}$};
\draw (6,-0.2) node{$\phi_{280,18}$};
\draw (7,-0.2) node{$\phi_{336,11}$};
\draw (8,-0.2) node{$\phi_{210,6}$};
\draw (9,-0.2) node{$\phi_{56,3}$};
\draw (10,-0.2) node{$\phi_{1,0}$};

\draw (3.6,1) node {$E_6[\theta];\ep$};
\draw (3.6,-1) node {$E_6[\theta^2];\ep$};

\draw (0,0.2) node{$6$};
\draw (1,0.2) node{$9$};
\draw (2,0.2) node{$10$};

\draw (3,-1.2) node{$10$};
\draw (3,1.2) node{$10$};

\draw (4,0.2) node{$10$};
\draw (5,0.2) node{$9$};
\draw (6,0.2) node{$8$};
\draw (7,0.2) node{$7$};
\draw (8,0.2) node{$6$};
\draw (9,0.2) node{$5$};
\draw (10,0.2) node{$0$};

\draw (0,0) -- (10,0);
\draw (3,-1) -- (3,1);

\draw (10,0) node [draw] (l0) {};
\draw (9,0) node [draw] (l0) {};
\draw (8,0) node [draw] (l0) {};
\draw (7,0) node [draw] (l1) {};
\draw (6,0) node [draw] (l2) {};
\draw (5,0) node [draw] (l2) {};
\draw (4,0) node [draw] (l2) {};
\draw (3,1) node [draw] (l2) {};
\draw (3,-1) node [draw] (l2) {};
\draw (3,0) node [fill=black!100] (ld) {};
\draw (2,0) node [draw] (l2) {};
\draw (1,0) node [draw] (l4) {};
\draw (0,0) node [draw] (l4) {};
\end{tikzpicture}\end{center}

\noindent(ii)\;\; \textbf{Block 2:} Cuspidal pair is $(\Phi_{12}.A_1(q^3),\phi_{11})$, of degree $q^3$. There are twelve unipotent characters in the block, two of which -- $E_6[\theta^i];1$ -- are non-real.

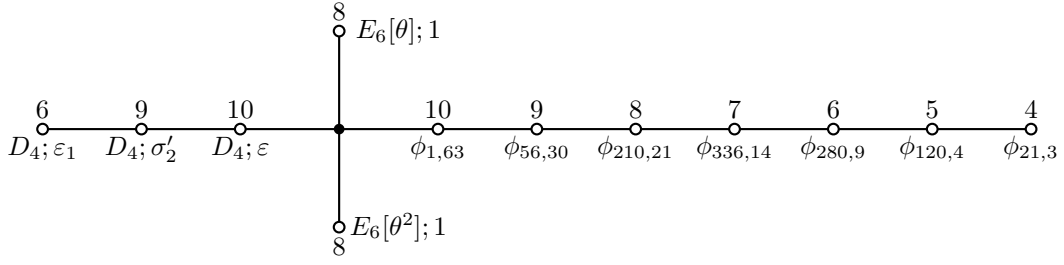
\begin{center}\begin{tabular}{lcccccc}
\hline Character & $A(-)$ & $\omega_i q^{aA/e}$ & $\kappa=1$ & $\kappa=5$ & $\kappa=7$ & $\kappa=11$
\\\hline $\phi_{21,3}$ & $24$ & $q^{2}$ & $4$ & $20$ & $28$ & $44$
\\ $\phi_{120,4}$ & $35$ & $q^{3}$ & $5$ & $29$ & $41$ & $65$
\\ $\phi_{280,9}$ & $44$ & $q^{4}$ & $6$ & $38$ & $50$ & $82$
\\ $\phi_{336,14}$ & $50$ & $q^{5}$ & $7$ & $41$ & $59$ & $93$
\\ $\phi_{210,21}$ & $54$ & $q^{6}$ & $8$ & $46$ & $62$ & $100$
\\ $\phi_{56,30}$ & $57$ & $q^{7}$ & $9$ & $47$ & $67$ & $105$
\\ $E_6[\theta];1$ & $44$ & $\theta q^{4}$ & $8$ & $36$ & $52$ & $80$
\\ $D_4;\ep_1$ & $35$ & $-q^{3}$ & $6$ & $30$ & $40$ & $64$
\\ $\phi_{1,63}$ & $60$ & $q^{10}$ & $10$ & $50$ & $70$ & $110$
\\ $D_4;\sigma_2'$ & $50$ & $-q^{5}$ & $9$ & $43$ & $57$ & $91$
\\ $E_6[\theta^2];1$ & $44$ & $\theta^2 q^{4}$ & $8$ & $36$ & $52$ & $80$
\\ $D_4;\ep$ & $57$ & $-q^{7}$ & $10$ & $48$ & $66$ & $104$
\\\hline\end{tabular}\end{center}

\begin{center}\begin{tikzpicture}[thick,scale=1.3]
\draw (0,-0.2) node{$D_4;\ep_1$};
\draw (1,-0.2) node{$D_4;\sigma_2'$};
\draw (2,-0.2) node{$D_4;\ep$};
\draw (4,-0.2) node{$\phi_{1,63}$};
\draw (5,-0.2) node{$\phi_{56,30}$};
\draw (6,-0.2) node{$\phi_{210,21}$};
\draw (7,-0.2) node{$\phi_{336,14}$};
\draw (8,-0.2) node{$\phi_{280,9}$};
\draw (9,-0.2) node{$\phi_{120,4}$};
\draw (10,-0.2) node{$\phi_{21,3}$};

\draw (3.6,1) node {$E_6[\theta];1$};
\draw (3.6,-1) node {$E_6[\theta^2];1$};

\draw (0,0.2) node{$6$};
\draw (1,0.2) node{$9$};
\draw (2,0.2) node{$10$};

\draw (3,-1.2) node{$8$};
\draw (3,1.2) node{$8$};

\draw (4,0.2) node{$10$};
\draw (5,0.2) node{$9$};
\draw (6,0.2) node{$8$};
\draw (7,0.2) node{$7$};
\draw (8,0.2) node{$6$};
\draw (9,0.2) node{$5$};
\draw (10,0.2) node{$4$};

\draw (0,0) -- (10,0);
\draw (3,-1) -- (3,1);

\draw (10,0) node [draw] (l0) {};
\draw (9,0) node [draw] (l0) {};
\draw (8,0) node [draw] (l0) {};
\draw (7,0) node [draw] (l1) {};
\draw (6,0) node [draw] (l2) {};
\draw (5,0) node [draw] (l2) {};
\draw (4,0) node [draw] (l2) {};
\draw (3,1) node [draw] (l2) {};
\draw (3,-1) node [draw] (l2) {};
\draw (3,0) node [fill=black!100] (ld) {};
\draw (2,0) node [draw] (l2) {};
\draw (1,0) node [draw] (l4) {};
\draw (0,0) node [draw] (l4) {};
\end{tikzpicture}\end{center}

\noindent\textbf{Proof of Brauer trees}: For both blocks we need a Morita argument. We claim that Harish-Chandra induction from the $E_6(q)$ Levi subgroup, together with cutting by the relevant block, induces a Morita equivalence. To prove this, we use GAP to find that induction from the principal block of $E_6(q)$ of a unipotent character to $E_7(q)$ has only a single constituent in each unipotent block. The precise calculations, together with the other constituents that lie in blocks of defect zero, are given below, starting with the principal series, then the $D_4$-series, then the $E_6$-series:

\begin{align*} \phi_{1,0} &\mapsto \phi_{1,0}+\phi_{21,3}(+\phi_{7,1}+\phi_{27,2})
\\ \phi_{6,1}&\mapsto \phi_{56,3}+\phi_{120,4}(+\phi_{7,1}+\phi_{21,6}+\phi_{27,2}+\phi_{105,5})
\\ \phi_{15,5}&\mapsto\phi_{210,6}+\phi_{280,9}(+\phi_{21,6}+\phi_{35,13}+\phi_{105,5}+\phi_{189,10})
\\ \phi_{20,10}&\mapsto \phi_{336,11}+\phi_{336,14}(+\phi_{35,22}+\phi_{35,13}+\phi_{189,10}+\phi_{189,17})
\\ \phi_{15,17}&\mapsto \phi_{280,18}+\phi_{210,21}(+\phi_{21,33}+\phi_{35,22}+\phi_{105,26}+\phi_{189,17})
\\ \phi_{6,25}&\mapsto \phi_{120,25}+\phi_{56,30}(+\phi_{7,46}+\phi_{21,33}+\phi_{27,37}+\phi_{105,26})
\\ \phi_{1,36}&\mapsto \phi_{21,36}+\phi_{1,63}(+\phi_{7,46}+\phi_{27,37})
\\ D_4;1&\mapsto D_4;1+D_4;\ep_1(+D_4;r\ep_1+D_4;r)
\\ D_4;r&\mapsto D_4;\sigma_2+D_4;\sigma_2'(+D_4;r\ep+D_4;r\ep_1+D_4;r\ep_2+D_4;r)
\\ D_4;\ep&\mapsto D_4;\ep_2+D_4;\ep(+D_4;r\ep+D_4;r\ep_2)
\\ E_6[\theta]&\mapsto E_6[\theta];1+E_6[\theta];\ep
\\ E_6[\theta]&\mapsto E_6[\theta^2];1+E_6[\theta^2];\ep
\end{align*}

\newpage
\subsection{$d=14$}

For $E_7(q)$ and $d=14$ there is a single unipotent block of weight $1$, together with 62 unipotent blocks of defect zero.
\\[1cm]
\noindent(i)\;\;\textbf{Block 1:} Cuspidal pair is $(\Phi_2\Phi_{14},1)$, of degree $1$. There are fourteen unipotent characters in the block, two of which -- $E_7[\pm\I]$ -- are non-real.

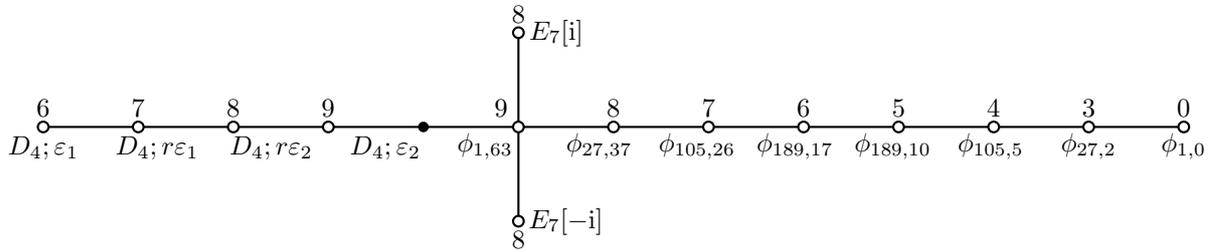
\begin{center}\begin{tabular}{lcccccccc}
\hline Character & $A(-)$ & $\omega_i q^{aA/e}$ & $\kappa=1$ & $\kappa=3$ & $\kappa=5$ & $\kappa=9$ & $\kappa=11$ & $\kappa=13$
\\\hline $\phi_{1,0}$ & $0$ & $q^{0}$ & $0$ & $0$ & $0$ & $0$ & $0$ & $0$
\\ $E_7[-\I]$ & $52$ & $-\I q^{9/2}$ & $8$ & $22$ & $38$ & $66$ & $82$ & $96$
\\ $\phi_{27,2}$ & $26$ & $q^{2}$ & $3$ & $11$ & $19$ & $33$ & $41$ & $49$
\\ $\phi_{105,5}$ & $38$ & $q^{3}$ & $4$ & $16$ & $26$ & $50$ & $60$ & $72$
\\ $\phi_{189,10}$ & $48$ & $q^{4}$ & $5$ & $21$ & $35$ & $61$ & $75$ & $91$
\\ $\phi_{189,17}$ & $55$ & $q^{5}$ & $6$ & $24$ & $40$ & $70$ & $86$ & $104$
\\ $\phi_{105,26}$ & $59$ & $q^{6}$ & $7$ & $25$ & $41$ & $77$ & $93$ & $111$
\\ $\phi_{27,37}$ & $61$ & $q^{7}$ & $8$ & $26$ & $44$ & $78$ & $96$ & $114$
\\ $E_7[\I]$ & $52$ & $\I q^{9/2}$ & $8$ & $22$ & $38$ & $66$ & $82$ & $96$
\\ $\phi_{1,63}$ & $63$ & $q^{9}$ & $9$ & $27$ & $45$ & $81$ & $99$ & $117$
\\ $D_4,\ep_1$ & $38$ & $-q^{3}$ & $6$ & $16$ & $28$ & $48$ & $60$ & $70$
\\ $D_4,r\ep_1$ & $48$ & $-q^{4}$ & $7$ & $21$ & $35$ & $61$ & $75$ & $89$
\\ $D_4,r\ep_2$ & $55$ & $-q^{5}$ & $8$ & $24$ & $40$ & $70$ & $86$ & $102$
\\ $D_4,\ep_2$ & $59$ & $-q^{6}$ & $9$ & $25$ & $43$ & $75$ & $93$ & $109$
\\\hline\end{tabular}\end{center}
\begin{center}\begin{tikzpicture}[thick,scale=1.25]
\draw (0,-0.22) node{${D_4;\ep_1}$};
\draw (1.2,-0.22) node{$D_4;r\ep_1$};
\draw (2.4,-0.22) node{$D_4;r\ep_2$};
\draw (3.6,-0.22) node{$D_4;\ep_2$};
\draw (4.65,-0.2) node{$\phi_{1,63}$};
\draw (5.85,-0.2) node{$\phi_{27,37}$};
\draw (6.88,-0.2) node{$\phi_{105,26}$};
\draw (7.91,-0.2) node{$\phi_{189,17}$};
\draw (8.94,-0.2) node{$\phi_{189,10}$};
\draw (9.97,-0.2) node{$\phi_{105,5}$};
\draw (11,-0.2) node{$\phi_{27,2}$};
\draw (12,-0.2) node{$\phi_{1,0}$};

\draw (5.4,1) node {$E_7[\I]$};
\draw (5.5,-1) node {$E_7[-\I]$};

\draw (0,0.2) node{$6$};
\draw (1,0.2) node{$7$};
\draw (2,0.2) node{$8$};
\draw (3,0.2) node{$9$};
\draw (4.82,0.2) node{$9$};

\draw (5,-1.2) node{$8$};
\draw (5,1.2) node{$8$};

\draw (6,0.2) node{$8$};
\draw (7,0.2) node{$7$};
\draw (8,0.2) node{$6$};
\draw (9,0.2) node{$5$};
\draw (10,0.2) node{$4$};
\draw (11,0.2) node{$3$};
\draw (12,0.2) node{$0$};

\draw (0,0) -- (12,0);
\draw (5,-1) -- (5,1);

\draw (12,0) node [draw] (l0) {};
\draw (11,0) node [draw] (l0) {};
\draw (10,0) node [draw] (l0) {};
\draw (9,0) node [draw] (l0) {};
\draw (8,0) node [draw] (l0) {};
\draw (7,0) node [draw] (l1) {};
\draw (6,0) node [draw] (l2) {};
\draw (5,0) node [draw] (l2) {};
\draw (5,1) node [draw] (l2) {};
\draw (5,-1) node [draw] (l2) {};
\draw (4,0) node [fill=black!100] (ld) {};
\draw (3,0) node [draw] (l2) {};
\draw (2,0) node [draw] (l2) {};
\draw (1,0) node [draw] (l4) {};
\draw (0,0) node [draw] (l4) {};
\end{tikzpicture}\end{center}
\noindent\textbf{Proof of Brauer tree}: This is given in \cite{cdr2012un}.

\newpage

\subsection{$d=18$}

For $E_7(q)$ and $d=18$ there is a single unipotent block of weight $1$, together with 58 unipotent blocks of defect zero.
\\[1cm]
\noindent(i)\;\;\textbf{Block 1:} Cuspidal pair is $(\Phi_2\Phi_{18},1)$, of degree $1$. There are eighteen unipotent characters in the block, six of which -- $E_6[\theta^i];1$, $E_6[\theta^i];\ep$, and $E_7[\pm\I]$ -- are non-real.

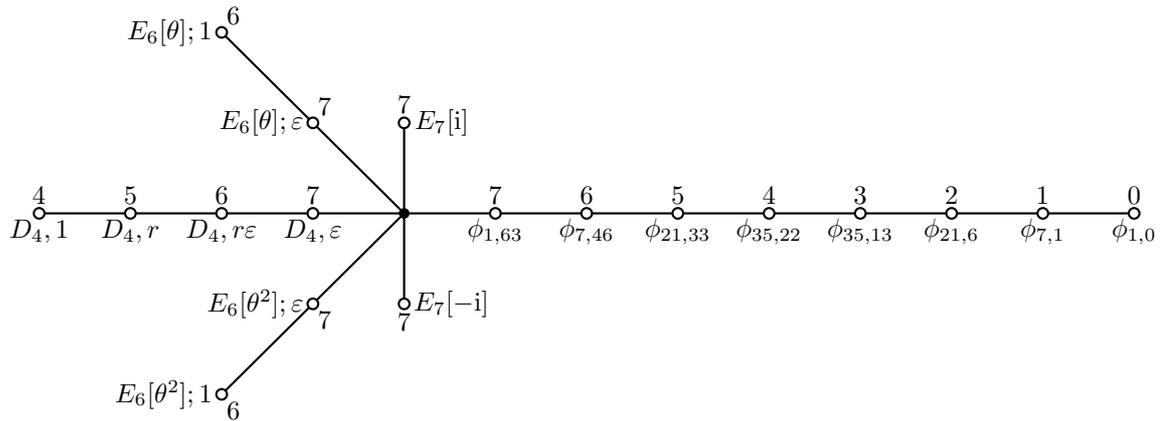
\begin{center}\begin{tabular}{lcccccccc}
\hline Character & $A(-)$ & $\omega_i q^{aA/e}$ & $\kappa=1$ & $\kappa=5$ & $\kappa=7$ & $\kappa=11$ & $\kappa=13$ & $\kappa=17$
\\\hline $\phi_{1,0}$ & $0$ & $q^{0}$ & $0$ & $0$ & $0$ & $0$ & $0$ & $0$
\\ $\phi_{7,1}$ & $17$ & $q$ & $1$ & $9$ & $13$ & $21$ & $25$ & $33$
\\ $\phi_{21,6}$ & $33$ & $q^{2}$ & $2$ & $18$ & $26$ & $40$ & $48$ & $64$
\\ $\phi_{35,13}$ & $47$ & $q^{3}$ & $3$ & $27$ & $37$ & $57$ & $67$ & $91$
\\ $\phi_{35,22}$ & $56$ & $q^{4}$ & $4$ & $32$ & $44$ & $68$ & $80$ & $108$
\\ $\phi_{21,33}$ & $60$ & $q^{5}$ & $5$ & $33$ & $47$ & $73$ & $87$ & $115$
\\ $\phi_{7,46}$ & $62$ & $q^{6}$ & $6$ & $34$ & $48$ & $76$ & $90$ & $118$
\\ $\phi_{1,63}$ & $63$ & $q^{7}$ & $7$ & $35$ & $49$ & $77$ & $91$ & $119$
\\ $E_7[\I]$ & $52$ & $\I q^{7/2}$ & $7$ & $29$ & $41$ & $63$ & $75$ & $97$
\\ $E_6[\theta];1$ & $47$ & $\theta q^{3}$ & $6$ & $26$ & $36$ & $58$ & $68$ & $88$
\\ $E_6[\theta];\ep$ & $56$ & $\theta q^{4}$ & $7$ & $31$ & $43$ & $69$ & $81$ & $105$
\\ $D_4,1$ & $33$ & $-q^{2}$ & $4$ & $18$ & $26$ & $40$ & $48$ & $62$
\\ $D_4,r$ & $47$ & $-q^{3}$ & $5$ & $25$ & $37$ & $57$ & $69$ & $89$
\\ $D_4,r\ep$ & $56$ & $-q^{4}$ & $6$ & $30$ & $44$ & $68$ & $82$ & $106$
\\ $D_4,\ep$ & $60$ & $-q^{5}$ & $7$ & $33$ & $47$ & $73$ & $87$ & $113$
\\ $E_6[\theta^2];1$ & $47$ & $\theta^2 q^{3}$ & $6$ & $26$ & $36$ & $58$ & $68$ & $88$
\\ $E_6[\theta^2];\ep$ & $56$ & $\theta^2 q^{4}$ & $7$ & $31$ & $43$ & $69$ & $81$ & $105$
\\ $E_7[-\I]$ & $52$ & $-\I q^{7/2}$ & $7$ & $29$ & $41$ & $63$ & $75$ & $97$
\\\hline\end{tabular}\end{center}

\begin{center}\begin{tikzpicture}[thick,scale=1.2]

\draw (12,0.8) node {$\phi_{1,0}$};
\draw (11,0.8) node {$\phi_{7,1}$};
\draw (10,0.8) node {$\phi_{21,6}$};
\draw (9,0.8) node {$\phi_{35,13}$};
\draw (8,0.8) node {$\phi_{35,22}$};
\draw (7,0.8) node {$\phi_{21,33}$};
\draw (6,0.8) node {$\phi_{7,46}$};
\draw (5,0.8) node {$\phi_{1,63}$};

\draw (3,0.8) node {$D_4,\ep$};
\draw (2,0.8) node {$D_4,r\ep$};
\draw (1,0.8) node {$D_4,r$};
\draw (0,0.8) node {$D_4,1$};

\draw (0,1.2) node {$4$};
\draw (1,1.2) node {$5$};
\draw (2,1.2) node {$6$};
\draw (3,1.2) node {$7$};

\draw (5,1.2) node {$7$};
\draw (6,1.2) node {$6$};
\draw (7,1.2) node {$5$};
\draw (8,1.2) node {$4$};
\draw (9,1.2) node {$3$};
\draw (10,1.2) node {$2$};
\draw (11,1.2) node {$1$};
\draw (12,1.2) node {$0$};

\draw (4,2.2) node {$7$};
\draw (4,-0.2) node {$7$};

\draw (3.13,2.18) node {$7$};
\draw (3.13,-0.18) node {$7$};
\draw (2.13,3.18) node {$6$};
\draw (2.13,-1.18) node {$6$};

\draw (0,1) -- (12,1);
\draw (2,-1) -- (4,1);
\draw (2,3) -- (4,1);
\draw (4,0) -- (4,2);

\draw (0,1) node [draw] (l1) {};
\draw (1,1) node [draw] (l0) {};
\draw (2,1) node [draw] (l1) {};
\draw (3,1) node [draw] (l1) {};
\draw (4,1) node [fill=black!100] (ld) {};
\draw (5,1) node [draw] (l1) {};
\draw (6,1) node [draw] (l1) {};
\draw (7,1) node [draw] (l1) {};
\draw (8,1) node [draw] (l1) {};
\draw (9,1) node [draw] (l1) {};
\draw (10,1) node [draw] (l1) {};
\draw (11,1) node [draw] (l1) {};
\draw (12,1) node [draw] (l1) {};

\draw (4,2) node [draw,label=right:${E_7[\I]}$] (l1) {};
\draw (4,0) node [draw,label=right:${E_7[-\I]}$] (l1) {};

\draw (3,2) node [draw,label=left:${E_6[\theta];\ep}$] (l1) {};
\draw (3,0) node [draw,label=left:${E_6[\theta^2];\ep}$] (l1) {};
\draw (2,3) node [draw,label=left:${E_6[\theta];1}$] (l1) {};
\draw (2,-1) node [draw,label=left:${E_6[\theta^2];1}$] (l1) {};
\end{tikzpicture}\end{center}

\newpage
\section{$E_8(q)$}

The group $E_8(q)$ has order $q^{120}\Phi_1^8\Phi_2^8\Phi_3^4\Phi_4^4\Phi_5^2\Phi_6^4\Phi_7\Phi_8^2\Phi_9\Phi_{10}^2\Phi_{12}^2\Phi_{14}\Phi_{15}\Phi_{18}\Phi_{20}\Phi_{24}\Phi_{30}$. There are Brauer trees for $\Phi_7$, $\Phi_9$, $\Phi_{14}$, $\Phi_{15}$, $\Phi_{18}$, $\Phi_{20}$, $\Phi_{24}$ and $\Phi_{30}$ obviously, and also non-principal unipotent blocks of weight $1$ when $d=2$, $3$, $5$, $6$, $8$, $10$ and $12$. The only one to appear in the literature so far is for $d=30$, which was constructed by Dudas and Rouquier \cite{dudasrouquier2012un}.

There are 166 unipotent characters of $E_8(q)$.

\newpage

\subsection{$d=1$}

For $E_8(q)$ and $d=1$ there are two unipotent blocks of weight $1$, together with the principal block with cyclotomic Weyl group $E_8$,  a unipotent block of weight $4$ with cyclotomic Weyl group $F_4$, two unipotent blocks of weight 2 with cyclotomic Weyl group $G(6,6,2)$, and 13 unipotent blocks of defect zero
\\[1cm]
\noindent(i)\;\;\textbf{Block 1/2:} Cuspidal pair is $(\Phi_1.E_7(q),E_7[\pm\I])$, of degree $q^{11}\Phi_1^7\Phi_3^3\Phi_4^2\Phi_5\Phi_7\Phi_8\Phi_9\Phi_{12}/2$. There are two unipotent characters in each block, neither of which is real, and the blocks are conjugate.
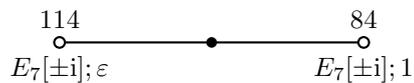
\begin{center}\begin{tabular}{lccc}
\hline Character & $A(-)$ & $\omega_i q^{aA/e}$ & $\kappa=1$
\\\hline $E_7[\pm\I],1$ & $42$ & $q^{21}$ & $84$
\\ $E_7[\pm\I],\ep$ & $57$ & $q^{36}$ & $114$
\\\hline\end{tabular}\end{center}

\begin{center}\begin{tikzpicture}[thick,scale=2]
\draw (2,0.15) node[rectangle]{$84$};
\draw (0,0.15) node[rectangle]{$114$};
\draw (2,-0.18) node[rectangle] {$E_7[\pm\I];1$};
\draw (0,-0.18) node[rectangle] {$E_7[\pm\I];\ep$};
\draw (0,0) -- (2,0);
\draw (2,0) node [draw] (l0) {};
\draw (0,0) node [draw] (l1) {};
\draw (1,0) node [fill=black!100] (ld) {};
\end{tikzpicture}\end{center}
\noindent\textbf{Proof of Brauer trees}: Since both characters have the same parity, they cannot be connected; this completes the proof.

\newpage
\subsection{$d=2$}
For $E_8(q)$ and $d=2$ there are two unipotent blocks of weight $1$, together with the principal block with cyclotomic Weyl group $E_8$, a unipotent block of weight $4$ with cyclotomic Weyl group $F_4$, two unipotent blocks of weight 2 with cyclotomic Weyl group $G(6,6,2)$, and 13 unipotent blocks of defect zero.
\\[1cm]
\noindent(i)\;\; \textbf{Block 1:} Cuspidal pair is $(\Phi_2.E_7(q),\phi_{512,11})$, of degree $q^{11}\Phi_2^7\Phi_4^2\Phi_6^3\Phi_8\Phi_{10}\Phi_{12}\Phi_{14}\Phi_{18}/2$. There are two unipotent characters in the block, both of which are real.

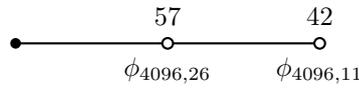
\begin{center}\begin{tabular}{lccc}
\hline Character & $A(-)$ & $\omega_i q^{aA/e}$ & $\kappa=1$
\\\hline $\phi_{4096,11}$ & $42$ & $q^{21}$ & $42$
\\ $\phi_{4096,26}$ & $57$ & $q^{36}$ & $57$
\\\hline\end{tabular}\end{center}
\begin{center}\begin{tikzpicture}[thick,scale=2]
\draw (1,-0.18) node{$\phi_{4096,26}$};
\draw (2,-0.18) node{$\phi_{4096,11}$};

\draw (0,0) -- (2,0);
\draw (0,0) node [fill=black!100] (ld) {};
\draw (1,0) node [draw] (l4) {};
\draw (2,0) node [draw] (l4) {};
\draw (2,0.18) node{$42$};
\draw (1,0.18) node{$57$};
\end{tikzpicture}\end{center}

\noindent(ii)\;\; \textbf{Block 2:} Cuspidal pair is $(\Phi_2.E_7(q),\phi_{512,12})$, of degree $q^{11}\Phi_2^7\Phi_4^2\Phi_6^3\Phi_8\Phi_{10}\Phi_{12}\Phi_{14}\Phi_{18}/2$. There are two unipotent characters in the block, both of which are real.

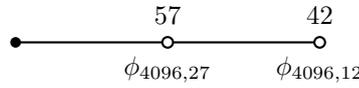
\begin{center}\begin{tabular}{lccc}
\hline Character & $A(-)$ & $\omega_i q^{aA/e}$ & $\kappa=1$
\\\hline $\phi_{4096,12}$ & $42$ & $q^{21}$ & $42$
\\ $\phi_{4096,27}$ & $57$ & $q^{36}$ & $57$
\\\hline\end{tabular}\end{center}
\begin{center}\begin{tikzpicture}[thick,scale=2]
\draw (1,-0.18) node{$\phi_{4096,27}$};
\draw (2,-0.18) node{$\phi_{4096,12}$};

\draw (0,0) -- (2,0);
\draw (0,0) node [fill=black!100] (ld) {};
\draw (1,0) node [draw] (l4) {};
\draw (2,0) node [draw] (l4) {};
\draw (2,0.18) node{$42$};
\draw (1,0.18) node{$57$};
\end{tikzpicture}\end{center}

\noindent\textbf{Proof of Brauer trees}: Since all characters are real, the Brauer tree is a line. A Geck--Pfeiffer argument then a degree argument is enough.

\newpage
\subsection{$d=3$}
For $E_8(q)$ and $d=3$ there are three unipotent blocks of weight $1$, together with the principal block with cyclotomic Weyl group $G_{32}$, a unipotent block of weight $2$ with cyclotomic Weyl group $G_5$, and 25 unipotent blocks of defect zero.
\\[1cm]
\noindent(i)\;\; \textbf{Block 1:} Cuspidal pair is $(\Phi_3.E_6(q)$,$\phi_{81,6})$, of degree $q^6\Phi_3^3\Phi_6^2\Phi_9\Phi_{12}$. There are six unipotent characters in the block, all of which are real.

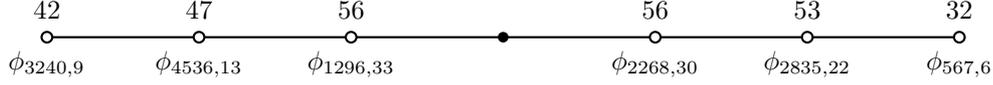
\begin{center}\begin{tabular}{lcccc}
\hline Character & $A(-)$ & $\omega_i q^{aA/e}$ & $\kappa=1$ & $\kappa=2$
\\\hline $\phi_{567,6}$ & $48$ & $q^{8}$ & $32$ & $64$
\\ $\phi_{4536,13}$ & $71$ & $-q^{13}$ & $47$ & $95$
\\ $\phi_{2268,30}$ & $84$ & $q^{18}$ & $56$ & $112$
\\ $\phi_{3240,9}$ & $$63 & $-q^{11}$ & $42$ & $84$
\\ $\phi_{2835,22}$ & $80$ & $q^{16}$ & $53$ & $107$
\\ $\phi_{1296,33}$ & $84$ & $-q^{18}$ & $56$ & $112$
\\\hline\end{tabular}\end{center}
\begin{center}\begin{tikzpicture}[thick,scale=2]
\draw (0,-0.18) node{$\phi_{3240,9}$};
\draw (1,-0.18) node{$\phi_{4536,13}$};
\draw (2,-0.18) node{$\phi_{1296,33}$};

\draw (4,-0.18) node{$\phi_{2268,30}$};
\draw (5,-0.18) node{$\phi_{2835,22}$};
\draw (6,-0.18) node{$\phi_{567,6}$};
\draw (0,0) -- (6,0);
\draw (3,0) node [fill=black!100] (ld) {};
\draw (0,0) node [draw] (l4) {};
\draw (1,0) node [draw] (l4) {};
\draw (2,0) node [draw] (l4) {};

\draw (4,0) node [draw] (l4) {};
\draw (5,0) node [draw] (l4) {};
\draw (6,0) node [draw] (l4) {};
\draw (0,0.18) node{$42$};
\draw (1,0.18) node{$47$};
\draw (2,0.18) node{$56$};
\draw (4,0.18) node{$56$};
\draw (5,0.18) node{$53$};
\draw (6,0.18) node{$32$};
\end{tikzpicture}\end{center}

\noindent(ii)\;\; \textbf{Block 2:} Cuspidal pair is $(\Phi_3.E_6(q)$,$\phi_{81,10})$, of degree $q^{10}\Phi_3^3\Phi_6^2\Phi_9\Phi_{12}$. There are six unipotent characters in the block, all of which are real.

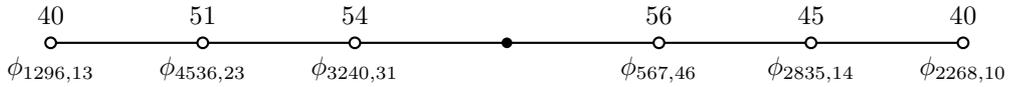
\begin{center}\begin{tabular}{lcccc}
\hline Character & $A(-)$ & $\omega_i q^{aA/e}$ & $\kappa=1$ & $\kappa=2$
\\\hline $\phi_{2268,10}$ & $60$ & $q^{10}$ & $40$ & $80$
\\ $\phi_{4536,23}$ & $77$ & $-q^{15}$ & $51$ & $103$
\\ $\phi_{567,46}$ & $84$ & $q^{20}$ & $56$ & $112$
\\ $\phi_{1296,13}$ & $60$ & $-q^{10}$ & $40$ & $80$
\\ $\phi_{2835,14}$ & $68$ & $q^{12}$ & $45$ & $91$
\\ $\phi_{3240,31}$ & $81$ & $-q^{17}$ & $54$ & $108$
\\\hline\end{tabular}\end{center}
\begin{center}\begin{tikzpicture}[thick,scale=2]
\draw (0,-0.18) node{$\phi_{1296,13}$};
\draw (1,-0.18) node{$\phi_{4536,23}$};
\draw (2,-0.18) node{$\phi_{3240,31}$};

\draw (4,-0.18) node{$\phi_{567,46}$};
\draw (5,-0.18) node{$\phi_{2835,14}$};
\draw (6,-0.18) node{$\phi_{2268,10}$};
\draw (0,0) -- (6,0);
\draw (3,0) node [fill=black!100] (ld) {};
\draw (0,0) node [draw] (l4) {};
\draw (1,0) node [draw] (l4) {};
\draw (2,0) node [draw] (l4) {};

\draw (4,0) node [draw] (l4) {};
\draw (5,0) node [draw] (l4) {};
\draw (6,0) node [draw] (l4) {};
\draw (0,0.18) node{$40$};
\draw (1,0.18) node{$51$};
\draw (2,0.18) node{$54$};
\draw (4,0.18) node{$56$};
\draw (5,0.18) node{$45$};
\draw (6,0.18) node{$40$};
\end{tikzpicture}\end{center}

\noindent(iii)\;\; \textbf{Block 3:} Cuspidal pair is $(\Phi_3.E_6(q)$,$\phi_{90,8})$, of degree $q^7\Phi_3^3\Phi_5\Phi_6^2\Phi_8\Phi_{12}/3$. There are six unipotent characters in the block, all of which are real.

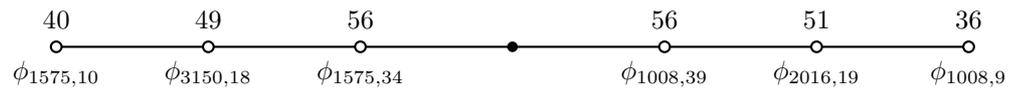
\begin{center}\begin{tabular}{lcccc}
\hline Character & $A(-)$ & $\omega_i q^{aA/e}$ & $\kappa=1$ & $\kappa=2$
\\\hline $\phi_{1008,9}$ & $54$ & $q^9$ & $36$ & $72$
\\ $\phi_{3150,18}$ & $75$ & $-q^{14}$ & $49$ & $101$
\\ $\phi_{1008,39}$ & $84$ & $q^{19}$ & $56$ & $112$
\\ $\phi_{1575,34}$ & $83$ & $-q^{18}$ & $56$ & $110$
\\ $\phi_{2016,19}$ & $75$ & $q^{14}$ & $51$ & $99$
\\ $\phi_{1575,10}$ & $59$ & $-q^{10}$ & $40$ & $78$
\\\hline\end{tabular}\end{center}
\begin{center}\begin{tikzpicture}[thick,scale=2]
\draw (0,-0.18) node{$\phi_{1575,10}$};
\draw (1,-0.18) node{$\phi_{3150,18}$};
\draw (2,-0.18) node{$\phi_{1575,34}$};

\draw (4,-0.18) node{$\phi_{1008,39}$};
\draw (5,-0.18) node{$\phi_{2016,19}$};
\draw (6,-0.18) node{$\phi_{1008,9}$};
\draw (0,0) -- (6,0);
\draw (3,0) node [fill=black!100] (ld) {};
\draw (0,0) node [draw] (l4) {};
\draw (1,0) node [draw] (l4) {};
\draw (2,0) node [draw] (l4) {};

\draw (4,0) node [draw] (l4) {};
\draw (5,0) node [draw] (l4) {};
\draw (6,0) node [draw] (l4) {};
\draw (0,0.18) node{$40$};
\draw (1,0.18) node{$49$};
\draw (2,0.18) node{$56$};
\draw (4,0.18) node{$56$};
\draw (5,0.18) node{$51$};
\draw (6,0.18) node{$36$};
\end{tikzpicture}\end{center}

\bigskip

\noindent\textbf{Proof of Brauer trees}: Since all characters are real, the Brauer tree is a line. A Geck--Pfeiffer argument then a degree argument is enough.

\newpage
\subsection{$d=5$}

For $E_8(q)$ and $d=5$ there are two unipotent blocks of weight $1$, together with the principal block with cyclotomic Weyl group $G_{16}$, and 101 unipotent blocks of defect zero.
\\[1cm]
\noindent(i)\;\; \textbf{Block 1:} Cuspidal pair is $(\Phi_4.A_4(q)$,$\phi_{32})$, of degree $q^2\Phi_5$. There are ten unipotent characters in the block, all of which are real.

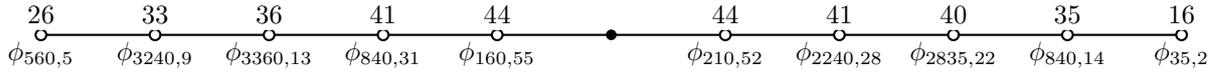
\begin{center}\begin{tabular}{lcccccc}
\hline Character & $A(-)$ & $\omega_i q^{aA/e}$ & $\kappa=1$ & $\kappa=2$ & $\kappa=3$ & $\kappa=4$
\\\hline $\phi_{35,2}$ & $40$ & $q^{4}$ & $16$ & $32$ & $48$ & $64$
\\ $\phi_{560,5}$ & $67$ & $-q^{7}$ & $26$ & $54$ & $80$ & $108$
\\ $\phi_{840,14}$ & $90$ & $q^{10}$ & $35$ & $73$ & $107$ & $145$
\\ $\phi_{840,31}$ & $104$ & $-q^{13}$ & $41$ & $83$ & $125$ & $167$
\\ $\phi_{210,52}$ & $110$ & $q^{16}$ & $44$ & $88$ & $132$ & $176$
\\ $\phi_{3240,9}$ & $83$ & $-q^{9}$ & $33$ & $67$ & $99$ & $133$
\\ $\phi_{2835,22}$ & $100$ & $q^{12}$ & $40$ & $80$ & $120$ & $160$
\\ $\phi_{3360,13}$ & $90$ & $-q^{10}$ & $36$ & $72$ & $108$ & $144$
\\ $\phi_{2240,28}$ & $104$ & $q^{13}$ & $41$ & $83$ & $125$ & $167$
\\ $\phi_{160,55}$ & $110$ & $-q^{16}$ & $44$ & $88$ & $132$ & $176$
\\\hline\end{tabular}\end{center}

\begin{center}\begin{tikzpicture}[thick,scale=1.5]
\draw (0,-0.18) node{$\phi_{560,5}$};
\draw (1,-0.18) node{$\phi_{3240,9}$};
\draw (2,-0.18) node{$\phi_{3360,13}$};
\draw (3,-0.18) node{$\phi_{840,31}$};
\draw (4,-0.18) node{$\phi_{160,55}$};

\draw (6,-0.18) node{$\phi_{210,52}$};
\draw (7,-0.18) node{$\phi_{2240,28}$};
\draw (8,-0.18) node{$\phi_{2835,22}$};
\draw (9,-0.18) node{$\phi_{840,14}$};
\draw (10,-0.18) node{$\phi_{35,2}$};
\draw (0,0) -- (10,0);
\draw (5,0) node [fill=black!100] (ld) {};
\draw (0,0) node [draw] (l4) {};
\draw (1,0) node [draw] (l4) {};
\draw (2,0) node [draw] (l4) {};
\draw (3,0) node [draw] (l4) {};
\draw (4,0) node [draw] (l4) {};
\draw (6,0) node [draw] (l4) {};
\draw (7,0) node [draw] (l4) {};
\draw (8,0) node [draw] (l4) {};
\draw (9,0) node [draw] (l4) {};
\draw (10,0) node [draw] (l4) {};

\draw (0,0.18) node{$26$};
\draw (1,0.18) node{$33$};
\draw (2,0.18) node{$36$};
\draw (3,0.18) node{$41$};
\draw (4,0.18) node{$44$};
\draw (6,0.18) node{$44$};
\draw (7,0.18) node{$41$};
\draw (8,0.18) node{$40$};
\draw (9,0.18) node{$35$};
\draw (10,0.18) node{$16$};
\end{tikzpicture}\end{center}

\noindent(ii)\;\; \textbf{Block 2:} Cuspidal pair is $(\Phi_4.A_4(q)$,$\phi_{221})$, of degree $q^4\Phi_5$. There are ten unipotent characters in the block, all of which are real.

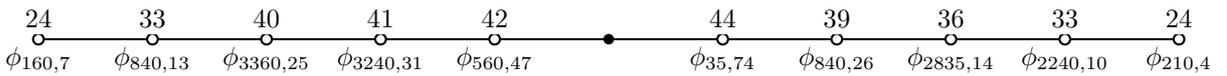
\begin{center}\begin{tabular}{lcccccc}
\hline Character & $A(-)$ & $\omega_i q^{aA/e}$ & $\kappa=1$ & $\kappa=2$ & $\kappa=3$ & $\kappa=4$
\\\hline $\phi_{210,4}$ & $60$ & $q^{6}$ & $24$ & $48$ & $72$ & $96$
\\ $\phi_{840,13}$ & $84$ & $-q^{9}$ & $33$ & $67$ & $101$ & $135$
\\ $\phi_{840,26}$ & $100$ & $q^{12}$ & $39$ & $81$ & $119$ & $161$
\\ $\phi_{560,47}$ & $107$ & $-q^{15}$ & $42$ & $86$ & $128$ & $172$
\\ $\phi_{35,74}$ & $110$ & $q^{18}$ & $44$ & $88$ & $132$ & $176$
\\ $\phi_{160,7}$ & $60$ & $-q^{6}$ & $24$ & $48$ & $72$ & $96$
\\ $\phi_{2240,10}$ & $84$ & $q^{9}$ & $33$ & $67$ & $101$ & $135$
\\ $\phi_{3360,25}$ & $100$ & $-q^{12}$ & $40$ & $80$ & $120$ & $160$
\\ $\phi_{2835,14}$ & $90$ & $q^{10}$ & $36$ & $72$ & $108$ & $144$
\\ $\phi_{3240,31}$ & $103$ & $-q^{13}$ & $41$ & $83$ & $123$ & $165$
\\\hline\end{tabular}\end{center}

\begin{center}\begin{tikzpicture}[thick,scale=1.5]
\draw (0,-0.18) node{$\phi_{160,7}$};
\draw (1,-0.18) node{$\phi_{840,13}$};
\draw (2,-0.18) node{$\phi_{3360,25}$};
\draw (3,-0.18) node{$\phi_{3240,31}$};
\draw (4,-0.18) node{$\phi_{560,47}$};

\draw (6,-0.18) node{$\phi_{35,74}$};
\draw (7,-0.18) node{$\phi_{840,26}$};
\draw (8,-0.18) node{$\phi_{2835,14}$};
\draw (9,-0.18) node{$\phi_{2240,10}$};
\draw (10,-0.18) node{$\phi_{210,4}$};
\draw (0,0) -- (10,0);
\draw (5,0) node [fill=black!100] (ld) {};
\draw (0,0) node [draw] (l4) {};
\draw (1,0) node [draw] (l4) {};
\draw (2,0) node [draw] (l4) {};
\draw (3,0) node [draw] (l4) {};
\draw (4,0) node [draw] (l4) {};
\draw (6,0) node [draw] (l4) {};
\draw (7,0) node [draw] (l4) {};
\draw (8,0) node [draw] (l4) {};
\draw (9,0) node [draw] (l4) {};
\draw (10,0) node [draw] (l4) {};

\draw (0,0.18) node{$24$};
\draw (1,0.18) node{$33$};
\draw (2,0.18) node{$40$};
\draw (3,0.18) node{$41$};
\draw (4,0.18) node{$42$};
\draw (6,0.18) node{$44$};
\draw (7,0.18) node{$39$};
\draw (8,0.18) node{$36$};
\draw (9,0.18) node{$33$};
\draw (10,0.18) node{$24$};
\end{tikzpicture}\end{center}

\bigskip

\noindent\textbf{Proof of Brauer trees}: Since all characters are real, the Brauer tree is a line. A Geck--Pfeiffer argument then a degree argument is enough.

\newpage
\subsection{$d=6$}

For $E_8(q)$ and $d=6$ there are three unipotent blocks of weight $1$, together with the principal block with cyclotomic Weyl group $G_{32}$, a unipotent block of weight $2$ with cyclotomic Weyl group $G_5$, and 25 unipotent blocks of defect zero.
\\[1cm]
\noindent(i)\;\; \textbf{Block 1:} Cuspidal pair is $(\Phi_6.{}^2\!E_6(q),\phi_{9,6}')$, of degree $q^6\Phi_3^2\Phi_6^3\Phi_{12}\Phi_{18}$. There are six unipotent characters in the block, all of which are real.

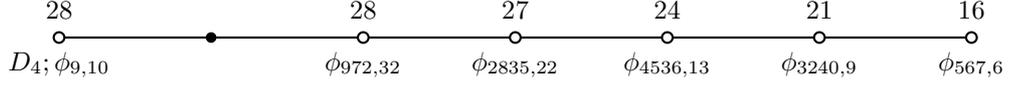
\begin{center}\begin{tabular}{lcccc}
\hline Character & $A(-)$ & $\omega_i q^{aA/e}$ & $\kappa=1$ & $\kappa=5$
\\\hline $\phi_{567,6}$ & $48$ & $q^{8}$ & $16$ & $80$
\\ $D_4;\phi_{9,10}$ & $84$ & $-q^{18}$ & $28$ & $140$
\\ $\phi_{2835,22}$ & $80$ & $q^{16}$ & $27$ & $133$
\\ $\phi_{3240,9}$ & $63$ & $q^{11}$ & $21$ & $105$
\\ $\phi_{972,32}$ & $84$ & $q^{18}$ & $28$ & $140$
\\ $\phi_{4536,13}$ & $71$ & $q^{13}$ & $24$ & $118$
\\\hline\end{tabular}\end{center}
\begin{center}\begin{tikzpicture}[thick,scale=2]
\draw (0,-0.18) node{$D_4;\phi_{9,10}$};
\draw (2,-0.18) node{$\phi_{972,32}$};
\draw (3,-0.18) node{$\phi_{2835,22}$};
\draw (4,-0.18) node{$\phi_{4536,13}$};
\draw (5,-0.18) node{$\phi_{3240,9}$};
\draw (6,-0.18) node{$\phi_{567,6}$};

\draw (0,0) -- (6,0);
\draw (6,0) node [draw] (l2) {};
\draw (5,0) node [draw] (l2) {};
\draw (4,0) node [draw] (l2) {};
\draw (3,0) node [draw] (l2) {};
\draw (1,0) node [fill=black!100] (ld) {};
\draw (2,0) node [draw] (l4) {};
\draw (0,0) node [draw] (l4) {};
\draw (0,0.18) node{$28$};
\draw (2,0.18) node{$28$};
\draw (3,0.18) node{$27$};
\draw (4,0.18) node{$24$};
\draw (5,0.18) node{$21$};
\draw (6,0.18) node{$16$};
\end{tikzpicture}\end{center}

\noindent(ii)\;\; \textbf{Block 2:} Cuspidal pair is $(\Phi_6.{}^2\!E_6(q),\phi_{9,6}'')$, of degree $q^{10}\Phi_3^2\Phi_6^3\Phi_{12}\Phi_{18}$. There are six unipotent characters in the block, all of which are real.
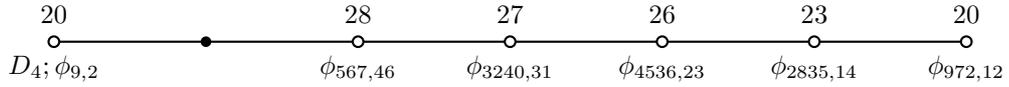
\begin{center}\begin{tabular}{lcccc}
\hline Character & $A(-)$ & $\omega_i q^{aA/e}$ & $\kappa=1$ & $\kappa=5$
\\\hline $\phi_{972,12}$ & $60$ & $q^{10}$ & $20$ & $100$
\\ $\phi_{3240,31}$ & $81$ & $q^{17}$ & $27$ & $135$
\\ $\phi_{2835,14}$ & $68$ & $q^{12}$ & $23$ & $113$
\\ $D_4;\phi_{9,2}$ & $60$ & $-q^{10}$ & $20$ & $100$
\\ $\phi_{567,46}$ & $84$ & $q^{20}$ & $28$ & $140$
\\ $\phi_{4536,23}$ & $77$ & $q^{15}$ & $26$ & $128$
\\\hline\end{tabular}\end{center}
\begin{center}\begin{tikzpicture}[thick,scale=2]
\draw (0,-0.18) node{$D_4;\phi_{9,2}$};
\draw (2,-0.18) node{$\phi_{567,46}$};
\draw (3,-0.18) node{$\phi_{3240,31}$};
\draw (4,-0.18) node{$\phi_{4536,23}$};
\draw (5,-0.18) node{$\phi_{2835,14}$};
\draw (6,-0.18) node{$\phi_{972,12}$};

\draw (0,0) -- (6,0);
\draw (6,0) node [draw] (l2) {};
\draw (5,0) node [draw] (l2) {};
\draw (4,0) node [draw] (l2) {};
\draw (3,0) node [draw] (l2) {};
\draw (1,0) node [fill=black!100] (ld) {};
\draw (2,0) node [draw] (l4) {};
\draw (0,0) node [draw] (l4) {};
\draw (0,0.18) node{$20$};
\draw (2,0.18) node{$28$};
\draw (3,0.18) node{$27$};
\draw (4,0.18) node{$26$};
\draw (5,0.18) node{$23$};
\draw (6,0.18) node{$20$};
\end{tikzpicture}\end{center}

\noindent(iii)\;\; \textbf{Block 3:} Cuspidal pair is $(\Phi_6.{}^2\!E_6(q),\phi_{6,6}'')$, of degree $q^7\Phi_3^2\Phi_6^3\Phi_8\Phi_{10}\Phi_{12}/3$. There are six unipotent characters in the block, all of which are real.
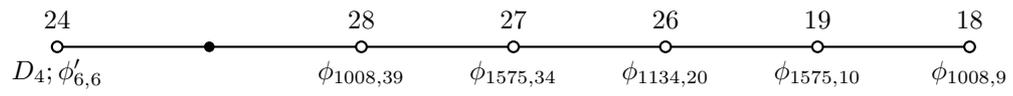
\begin{center}\begin{tabular}{lcccc}
\hline Character & $A(-)$ & $\omega_i q^{aA/e}$ & $\kappa=1$ & $\kappa=5$
\\\hline $\phi_{1008,9}$ & $54$ & $q^{9}$ & $18$ & $90$
\\ $\phi_{1575,10}$ & $59$ & $q^{10}$ & $19$ & $99$
\\ $D_4;\phi_{6,6}'$ & $75$ & $-q^{14}$ & $24$ & $126$
\\ $\phi_{1575,34}$ & $83$ & $q^{18}$ & $27$ & $139$
\\ $\phi_{1008,39}$ & $84$ & $q^{19}$ & $28$ & $140$
\\ $\phi_{1134,20}$ & $75$ & $q^{14}$ & $26$ & $124$
\\\hline\end{tabular}\end{center}
\begin{center}\begin{tikzpicture}[thick,scale=2]
\draw (0,-0.18) node{$D_4;\phi_{6,6}'$};
\draw (2,-0.18) node{$\phi_{1008,39}$};
\draw (3,-0.18) node{$\phi_{1575,34}$};
\draw (4,-0.18) node{$\phi_{1134,20}$};
\draw (5,-0.18) node{$\phi_{1575,10}$};
\draw (6,-0.18) node{$\phi_{1008,9}$};

\draw (0,0) -- (6,0);
\draw (6,0) node [draw] (l2) {};
\draw (5,0) node [draw] (l2) {};
\draw (4,0) node [draw] (l2) {};
\draw (3,0) node [draw] (l2) {};
\draw (1,0) node [fill=black!100] (ld) {};
\draw (2,0) node [draw] (l4) {};
\draw (0,0) node [draw] (l4) {};
\draw (0,0.18) node{$24$};
\draw (2,0.18) node{$28$};
\draw (3,0.18) node{$27$};
\draw (4,0.18) node{$26$};
\draw (5,0.18) node{$19$};
\draw (6,0.18) node{$18$};
\end{tikzpicture}\end{center}

\bigskip

\noindent\textbf{Proof of Brauer trees}: Since all characters are real, the Brauer tree is a line. A Geck--Pfeiffer argument then a degree argument is enough.

\newpage
\subsection{$d=7$}

For $E_8(q)$ and $d=7$ there are two unipotent blocks of weight $1$, together with $138$ unipotent blocks of defect zero.
\\[1cm]
\noindent(i)\;\; \textbf{Block 1:} Cuspidal pair is $(\Phi_1\Phi_7.A_1(q),\phi_{2})$, of degree $1$. There are fourteen unipotent characters in the block, all of which are real.

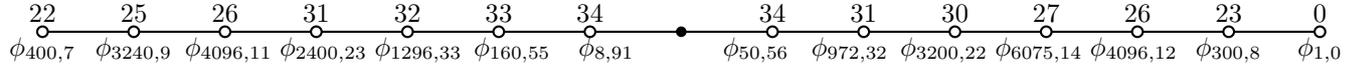
\begin{center}\begin{tabular}{lcccccccc}
\hline Character & $A(-)$ & $\omega_i q^{aA/e}$ & $\kappa=1$ & $\kappa=2$ & $\kappa=3$ & $\kappa=4$ & $\kappa=5$ & $\kappa=6$
\\\hline $\phi_{1,0}$ & $0$ & $q^{0}$ & $0$ & $0$ & $0$ & $0$ & $0$ & $0$
\\ $\phi_{4096,12}$ & $94$ & $q^{15/2}$ & $26$ & $54$ & $80$ & $108$ & $134$ & $162$
\\ $\phi_{6075,14}$ & $98$ & $q^{8}$ & $27$ & $55$ & $85$ & $111$ & $141$ & $169$
\\ $\phi_{160,55}$ & $116$ & $-q^{12}$ & $33$ & $67$ & $99$ & $133$ & $165$ & $199$
\\ $\phi_{3200,22}$ & $105$ & $q^{9}$ & $30$ & $60$ & $90$ & $120$ & $150$ & $180$
\\ $\phi_{400,7}$ & $78$ & $-q^{6}$ & $22$ & $44$ & $66$ & $90$ & $112$ & $134$
\\ $\phi_{972,32}$ & $110$ & $q^{10}$ & $31$ & $63$ & $95$ & $125$ & $157$ & $189$
\\ $\phi_{3240,9}$ & $89$ & $-q^{7}$ & $25$ & $51$ & $77$ & $101$ & $127$ & $153$
\\ $\phi_{4096,11}$ & $94$ & $-q^{15/2}$ & $26$ & $54$ & $80$ & $108$ & $134$ & $162$
\\ $\phi_{8,91}$ & $119$ & $-q^{15}$ & $34$ & $68$ & $102$ & $136$ & $170$ & $204$
\\ $\phi_{50,56}$ & $116$ & $q^{12}$ & $34$ & $66$ & $100$ & $132$ & $166$ & $198$
\\ $\phi_{2400,23}$ & $105$ & $-q^{9}$ & $31$ & $61$ & $91$ & $119$ & $149$ & $179$
\\ $\phi_{300,8}$ & $78$ & $q^{6}$ & $23$ & $45$ & $67$ & $89$ & $111$ & $133$
\\ $\phi_{1296,33}$ & $110$ & $-q^{10}$ & $32$ & $62$ & $94$ & $126$ & $158$ & $188$
\\\hline\end{tabular}\end{center}
\begin{center}\begin{tikzpicture}[thick,scale=1.2]
\draw (0,-0.2) node{$\phi_{400,7}$};
\draw (1,-0.2) node{$\phi_{3240,9}$};
\draw (2.04,-0.2) node{$\phi_{4096,11}$};
\draw (3.08,-0.2) node{$\phi_{2400,23}$};
\draw (4.12,-0.2) node{$\phi_{1296,33}$};
\draw (5.15,-0.2) node{$\phi_{160,55}$};
\draw (6.18,-0.2) node{$\phi_{8,91}$};

\draw (7.82,-0.2) node{$\phi_{50,56}$};
\draw (8.85,-0.2) node{$\phi_{972,32}$};
\draw (9.88,-0.2) node{$\phi_{3200,22}$};
\draw (10.92,-0.2) node{$\phi_{6075,14}$};
\draw (11.96,-0.2) node{$\phi_{4096,12}$};
\draw (13,-0.2) node{$\phi_{300,8}$};
\draw (14,-0.2) node{$\phi_{1,0}$};

\draw (0,0.2) node{$22$};
\draw (1,0.2) node{$25$};
\draw (2,0.2) node{$26$};
\draw (3,0.2) node{$31$};
\draw (4,0.2) node{$32$};
\draw (5,0.2) node{$33$};
\draw (6,0.2) node{$34$};

\draw (14,0.2) node{$0$};
\draw (13,0.2) node{$23$};
\draw (12,0.2) node{$26$};
\draw (11,0.2) node{$27$};
\draw (10,0.2) node{$30$};
\draw (9,0.2) node{$31$};
\draw (8,0.2) node{$34$};

\draw (0,0) -- (14,0);
\draw (7,0) node [fill=black!100] (ld) {};

\draw (14,0) node [draw] (l2) {};
\draw (13,0) node [draw] (l2) {};
\draw (12,0) node [draw] (l2) {};
\draw (11,0) node [draw] (l2) {};
\draw (10,0) node [draw] (l2) {};
\draw (9,0) node [draw] (l2) {};
\draw (8,0) node [draw] (l2) {};
\draw (6,0) node [draw] (l2) {};
\draw (5,0) node [draw] (l2) {};
\draw (4,0) node [draw] (l2) {};
\draw (3,0) node [draw] (l2) {};
\draw (2,0) node [draw] (l2) {};
\draw (1,0) node [draw] (l2) {};
\draw (0,0) node [draw] (l2) {};

\end{tikzpicture}\end{center}

\noindent(ii)\;\; \textbf{Block 2:} Cuspidal pair is $(\Phi_1\Phi_7.A_1(q)$,$\phi_{11})$, of degree $q$. There are fourteen unipotent characters in the block, all of which are real.

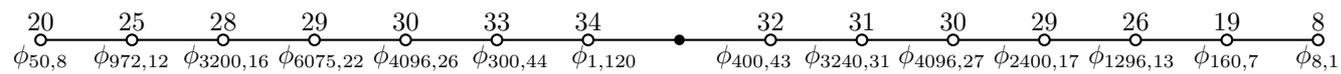
\begin{center}\begin{tabular}{lcccccccc}
\hline Character & $A(-)$ & $\omega_i q^{aA/e}$ & $\kappa=1$ & $\kappa=2$ & $\kappa=3$ & $\kappa=4$ & $\kappa=5$ & $\kappa=6$
\\\hline $\phi_{8,1}$ & $28$ & $q^{2}$ & $8$ & $16$ & $24$ & $32$ & $40$ & $48$
\\ $\phi_{4096,27}$ & $108$ & $q^{19/2}$ & $30$ & $62$ & $92$ & $124$ & $154$ & $186$
\\ $\phi_{3240,31}$ & $110$ & $q^{10}$ & $31$ & $63$ & $95$ & $125$ & $157$ & $189$
\\ $\phi_{972,12}$ & $89$ & $-q^{7}$ & $25$ & $51$ & $77$ & $101$ & $127$ & $153$
\\ $\phi_{400,43}$ & $113$ & $q^{11}$ & $32$ & $64$ & $96$ & $130$ & $162$ & $194$
\\ $\phi_{3200,16}$ & $98$ & $-q^{8}$ & $28$ & $56$ & $84$ & $112$ & $140$ & $168$
\\ $\phi_{160,7}$ & $67$ & $q^{5}$ & $19$ & $39$ & $57$ & $77$ & $95$ & $115$
\\ $\phi_{6075,22}$ & $105$ & $-q^{9}$ & $29$ & $59$ & $91$ & $119$ & $151$ & $181$
\\ $\phi_{4096,26}$ & $108$ & $-q^{19/2}$ & $30$ & $62$ & $92$ & $124$ & $154$ & $186$
\\ $\phi_{1,120}$ & $120$ & $-q^{17}$ & $34$ & $68$ & $102$ & $136$ & $170$ & $204$
\\ $\phi_{1296,13}$ & $89$ & $q^{7}$ & $26$ & $50$ & $76$ & $102$ & $128$ & $152$
\\ $\phi_{300,44}$ & $113$ & $-q^{11}$ & $33$ & $65$ & $97$ & $129$ & $161$ & $193$
\\ $\phi_{2400,17}$ & $98$ & $q^{8}$ & $29$ & $57$ & $85$ & $111$ & $139$ & $167$
\\ $\phi_{50,8}$ & $67$ & $-q^{5}$ & $20$ & $38$ & $58$ & $76$ & $96$ & $114$
\\\hline\end{tabular}\end{center}
\begin{center}\begin{tikzpicture}[thick,scale=1.2]
\draw (0,-0.2) node{$\phi_{50,8}$};
\draw (1,-0.2) node{$\phi_{972,12}$};
\draw (2.04,-0.2) node{$\phi_{3200,16}$};
\draw (3.08,-0.2) node{$\phi_{6075,22}$};
\draw (4.12,-0.2) node{$\phi_{4096,26}$};
\draw (5.15,-0.2) node{$\phi_{300,44}$};
\draw (6.18,-0.2) node{$\phi_{1,120}$};

\draw (7.82,-0.2) node{$\phi_{400,43}$};
\draw (8.85,-0.2) node{$\phi_{3240,31}$};
\draw (9.88,-0.2) node{$\phi_{4096,27}$};
\draw (10.92,-0.2) node{$\phi_{2400,17}$};
\draw (11.96,-0.2) node{$\phi_{1296,13}$};
\draw (13,-0.2) node{$\phi_{160,7}$};
\draw (14,-0.2) node{$\phi_{8,1}$};

\draw (0,0.2) node{$20$};
\draw (1,0.2) node{$25$};
\draw (2,0.2) node{$28$};
\draw (3,0.2) node{$29$};
\draw (4,0.2) node{$30$};
\draw (5,0.2) node{$33$};
\draw (6,0.2) node{$34$};

\draw (14,0.2) node{$8$};
\draw (13,0.2) node{$19$};
\draw (12,0.2) node{$26$};
\draw (11,0.2) node{$29$};
\draw (10,0.2) node{$30$};
\draw (9,0.2) node{$31$};
\draw (8,0.2) node{$32$};

\draw (0,0) -- (14,0);
\draw (7,0) node [fill=black!100] (ld) {};

\draw (14,0) node [draw] (l2) {};
\draw (13,0) node [draw] (l2) {};
\draw (12,0) node [draw] (l2) {};
\draw (11,0) node [draw] (l2) {};
\draw (10,0) node [draw] (l2) {};
\draw (9,0) node [draw] (l2) {};
\draw (8,0) node [draw] (l2) {};
\draw (6,0) node [draw] (l2) {};
\draw (5,0) node [draw] (l2) {};
\draw (4,0) node [draw] (l2) {};
\draw (3,0) node [draw] (l2) {};
\draw (2,0) node [draw] (l2) {};
\draw (1,0) node [draw] (l2) {};
\draw (0,0) node [draw] (l2) {};

\end{tikzpicture}\end{center}

\bigskip

\noindent\textbf{Proof of Brauer trees}: Since all characters are real, the Brauer tree is a line. A Geck--Pfeiffer argument then a degree argument is enough.

\newpage

\subsection{$d=8$}

For $E_8(q)$ and $d=8$ there are six unipotent blocks of weight $1$, together with the principal block with cyclotomic Weyl group $G_9$, and $86$ unipotent blocks of defect zero.
\\[1cm]
\noindent(i)\;\; \textbf{Block 1:} Cuspidal pair is $(\Phi_8.{}^2\!D_4(q)$,$\phi_{13,-})$, of degree $q\Phi_8$. There are eight unipotent characters in the block, all of which are real.

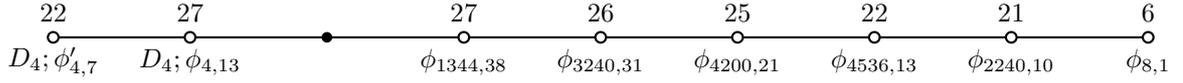
\begin{center}\begin{tabular}{lcccccc}
\hline Character & $A(-)$ & $\omega_i q^{aA/e}$ & $\kappa=1$ & $\kappa=3$ & $\kappa=5$ & $\kappa=7$
\\\hline $\phi_{8,1}$ & $24$ & $q^{3}$ & $6$ & $18$ & $30$ & $42$
\\ $\phi_{2240,10}$ & $87$ & $q^{12}$ & $21$ & $65$ & $109$ & $153$
\\ $\phi_{4536,13}$ & $92$ & $q^{13}$ & $22$ & $70$ & $114$ & $162$
\\ $D_4;\phi_{4,13}$ & $108$ & $-q^{18}$ & $27$ & $81$ & $135$ & $189$
\\ $\phi_{4200,21}$ & $100$ & $q^{15}$ & $25$ & $75$ & $125$ & $175$
\\ $D_4;\phi_{4,7}'$ & $87$ & $-q^{12}$ & $22$ & $66$ & $108$ & $152$
\\ $\phi_{3240,31}$ & $106$ & $q^{17}$ & $26$ & $80$ & $132$ & $186$
\\ $\phi_{1344,38}$ & $108$ & $q^{18}$ & $27$ & $81$ & $135$ & $189$
\\\hline\end{tabular}\end{center}
\begin{center}\begin{tikzpicture}[thick,scale=1.8]
\draw (0,-0.18) node{${D_4;\phi_{4,7}'}$};
\draw (1,-0.18) node{$D_4;\phi_{4,13}$};
\draw (3,-0.18) node{$\phi_{1344,38}$};
\draw (4,-0.18) node{$\phi_{3240,31}$};
\draw (5,-0.18) node{$\phi_{4200,21}$};
\draw (6,-0.18) node{$\phi_{4536,13}$};
\draw (7,-0.18) node{$\phi_{2240,10}$};
\draw (8,-0.18) node{$\phi_{8,1}$};

\draw (0,0) -- (8,0);
\draw (8,0) node [draw] (l0) {};
\draw (7,0) node [draw] (l1) {};
\draw (6,0) node [draw] (l2) {};
\draw (5,0) node [draw] (l2) {};
\draw (4,0) node [draw] (l2) {};
\draw (3,0) node [draw] (l2) {};
\draw (2,0) node [fill=black!100] (ld) {};
\draw (1,0) node [draw] (l4) {};
\draw (0,0) node [draw] (l4) {};
\draw (0,0.18) node{$22$};
\draw (1,0.18) node{$27$};
\draw (3,0.18) node{$27$};
\draw (4,0.18) node{$26$};
\draw (5,0.18) node{$25$};
\draw (6,0.18) node{$22$};
\draw (7,0.18) node{$21$};
\draw (8,0.18) node{$6$};
\end{tikzpicture}\end{center}

\noindent(ii)\;\; \textbf{Block 2:} Cuspidal pair $(\Phi_8.{}^2\!D_4(q)$,$\phi_{0123,13})$, of degree $q^7\Phi_8$. There are eight unipotent characters in the block, all of which are real.

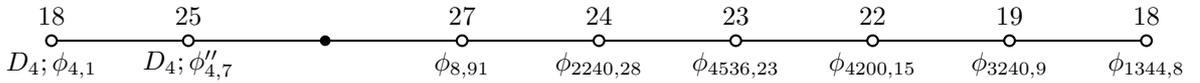
\begin{center}\begin{tabular}{lcccccc}
\hline Character & $A(-)$ & $\omega_i q^{aA/e}$ & $\kappa=1$ & $\kappa=3$ & $\kappa=5$ & $\kappa=7$
\\\hline $\phi_{1344,8}$ & $72$ & $q^{9}$ & $18$ & $54$ & $90$ & $126$
\\ $\phi_{3240,9}$ & $78$ & $q^{10}$ & $19$ & $59$ & $97$ & $137$
\\ $D_4;\phi_{4,7}''$ & $99$ & $-q^{15}$ & $25$ & $75$ & $123$ & $173$
\\ $\phi_{4200,15}$ & $88$ & $q^{12}$ & $22$ & $66$ & $110$ & $154$
\\ $D_4;\phi_{4,1}$ & $72$ & $-q^{9}$ & $18$ & $54$ & $90$ & $126$
\\ $\phi_{4536,23}$ & $96$ & $q^{14}$ & $23$ & $73$ & $119$ & $169$
\\ $\phi_{2240,28}$ & $99$ & $q^{15}$ & $24$ & $74$ & $124$ & $174$
\\ $\phi_{8,91}$ & $108$ & $q^{24}$ & $27$ & $81$ & $135$ & $189$
\\\hline\end{tabular}\end{center}

\begin{center}\begin{tikzpicture}[thick,scale=1.8]
\draw (0,-0.18) node{${D_4;\phi_{4,1}}$};
\draw (1,-0.18) node{$D_4;\phi_{4,7}''$};
\draw (3,-0.18) node{$\phi_{8,91}$};
\draw (4,-0.18) node{$\phi_{2240,28}$};
\draw (5,-0.18) node{$\phi_{4536,23}$};
\draw (6,-0.18) node{$\phi_{4200,15}$};
\draw (7,-0.18) node{$\phi_{3240,9}$};
\draw (8,-0.18) node{$\phi_{1344,8}$};

\draw (0,0) -- (8,0);
\draw (8,0) node [draw] (l0) {};
\draw (7,0) node [draw] (l1) {};
\draw (6,0) node [draw] (l2) {};
\draw (5,0) node [draw] (l2) {};
\draw (4,0) node [draw] (l2) {};
\draw (3,0) node [draw] (l2) {};
\draw (2,0) node [fill=black!100] (ld) {};
\draw (1,0) node [draw] (l4) {};
\draw (0,0) node [draw] (l4) {};
\draw (0,0.18) node{$18$};
\draw (1,0.18) node{$25$};
\draw (3,0.18) node{$27$};
\draw (4,0.18) node{$24$};
\draw (5,0.18) node{$23$};
\draw (6,0.18) node{$22$};
\draw (7,0.18) node{$19$};
\draw (8,0.18) node{$18$};
\end{tikzpicture}\end{center}

\noindent(iii)\;\; \textbf{Block 3:} Cuspidal pair $(\Phi_8.{}^2\!D_4(q)$,$\phi_{023,1})$, of degree $q^3\Phi_3\Phi_8$. There are eight unipotent characters in the block, all of which are real.

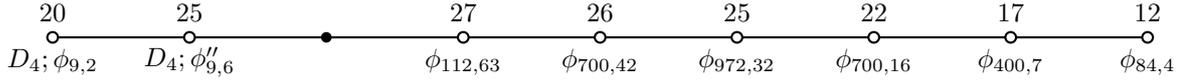
\begin{center}\begin{tabular}{lcccccc}
\hline Character & $A(-)$ & $\omega_i q^{aA/e}$ & $\kappa=1$ & $\kappa=3$ & $\kappa=5$ & $\kappa=7$
\\\hline $\phi_{84,4}$ & $48$ & $q^{6}$ & $12$ & $36$ & $60$ & $84$
\\ $D_4;\phi_{9,2}$ & $81$ & $-q^{11}$ & $20$ & $60$ & $102$ & $142$
\\ $\phi_{972,32}$ & $101$ & $q^{16}$ & $25$ & $77$ & $125$ & $177$
\\ $\phi_{400,7}$ & $69$ & $q^{9}$ & $17$ & $51$ & $87$ & $121$
\\ $\phi_{700,42}$ & $105$ & $q^{18}$ & $26$ & $78$ & $132$ & $184$
\\ $D_4;\phi_{9,6}''$ & $98$ & $-q^{15}$ & $25$ & $75$ & $121$ & $171$
\\ $\phi_{700,16}$ & $86$ & $q^{12}$ & $22$ & $64$ & $108$ & $150$
\\ $\phi_{112,63}$ & $108$ & $q^{21}$ & $27$ & $81$ & $135$ & $189$
\\\hline\end{tabular}\end{center}

\begin{center}\begin{tikzpicture}[thick,scale=1.8]
\draw (0,-0.18) node{${D_4;\phi_{9,2}}$};
\draw (1,-0.18) node{$D_4;\phi_{9,6}''$};
\draw (3,-0.18) node{$\phi_{112,63}$};
\draw (4,-0.18) node{$\phi_{700,42}$};
\draw (5,-0.18) node{$\phi_{972,32}$};
\draw (6,-0.18) node{$\phi_{700,16}$};
\draw (7,-0.18) node{$\phi_{400,7}$};
\draw (8,-0.18) node{$\phi_{84,4}$};

\draw (0,0) -- (8,0);
\draw (8,0) node [draw] (l0) {};
\draw (7,0) node [draw] (l1) {};
\draw (6,0) node [draw] (l2) {};
\draw (5,0) node [draw] (l2) {};
\draw (4,0) node [draw] (l2) {};
\draw (3,0) node [draw] (l2) {};
\draw (2,0) node [fill=black!100] (ld) {};
\draw (1,0) node [draw] (l4) {};
\draw (0,0) node [draw] (l4) {};
\draw (0,0.18) node{$20$};
\draw (1,0.18) node{$25$};
\draw (3,0.18) node{$27$};
\draw (4,0.18) node{$26$};
\draw (5,0.18) node{$25$};
\draw (6,0.18) node{$22$};
\draw (7,0.18) node{$17$};
\draw (8,0.18) node{$12$};
\end{tikzpicture}\end{center}

\noindent(iv)\;\; \textbf{Block 4:} Cuspidal pair $(\Phi_8.{}^2\!D_4(q)$,$\phi_{123,0})$, of degree $q^3\Phi_6\Phi_8$. There are eight unipotent characters in the block, all of which are real.

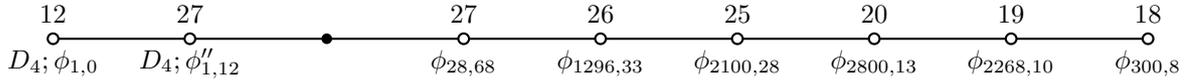
\begin{center}\begin{tabular}{lcccccc}
\hline Character & $A(-)$ & $\omega_i q^{aA/e}$ & $\kappa=1$ & $\kappa=3$ & $\kappa=5$ & $\kappa=7$
\\\hline $D_4;\phi_{1,0}$ & $48$ & $-q^{6}$ & $12$ & $36$ & $60$ & $84$
\\ $\phi_{2268,10}$ & $81$ & $q^{11}$ & $19$ & $61$ & $101$ & $143$
\\ $\phi_{2800,13}$ & $86$ & $q^{12}$ & $20$ & $64$ & $108$ & $152$
\\ $\phi_{28,68}$ & $108$ & $q^{21}$ & $27$ & $81$ & $135$ & $189$
\\ $D_4;\phi_{1,12}''$ & $105$ & $-q^{18}$ & $27$ & $79$ & $131$ & $183$
\\ $\phi_{2100,28}$ & $98$ & $q^{15}$ & $25$ & $73$ & $123$ & $171$
\\ $\phi_{1296,33}$ & $101$ & $q^{16}$ & $26$ & $76$ & $126$ & $176$
\\ $\phi_{300,8}$ & $69$ & $q^{9}$ & $18$ & $52$ & $86$ & $120$
\\\hline\end{tabular}\end{center}

\begin{center}\begin{tikzpicture}[thick,scale=1.8]
\draw (0,-0.18) node{${D_4;\phi_{1,0}}$};
\draw (1,-0.18) node{$D_4;\phi_{1,12}''$};
\draw (3,-0.18) node{$\phi_{28,68}$};
\draw (4,-0.18) node{$\phi_{1296,33}$};
\draw (5,-0.18) node{$\phi_{2100,28}$};
\draw (6,-0.18) node{$\phi_{2800,13}$};
\draw (7,-0.18) node{$\phi_{2268,10}$};
\draw (8,-0.18) node{$\phi_{300,8}$};

\draw (0,0) -- (8,0);
\draw (8,0) node [draw] (l0) {};
\draw (7,0) node [draw] (l1) {};
\draw (6,0) node [draw] (l2) {};
\draw (5,0) node [draw] (l2) {};
\draw (4,0) node [draw] (l2) {};
\draw (3,0) node [draw] (l2) {};
\draw (2,0) node [fill=black!100] (ld) {};
\draw (1,0) node [draw] (l4) {};
\draw (0,0) node [draw] (l4) {};
\draw (0,0.18) node{$12$};
\draw (1,0.18) node{$27$};
\draw (3,0.18) node{$27$};
\draw (4,0.18) node{$26$};
\draw (5,0.18) node{$25$};
\draw (6,0.18) node{$20$};
\draw (7,0.18) node{$19$};
\draw (8,0.18) node{$18$};
\end{tikzpicture}\end{center}

\noindent(v)\;\; \textbf{Block 5:} Cuspidal pair $(\Phi_8.{}^2\!D_4(q)$,$\phi_{013,2})$, of degree $q^3\Phi_3\Phi_8$. There are eight unipotent characters in the block, all of which are real.

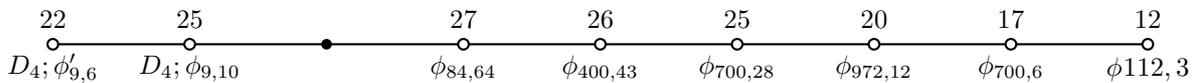
\begin{center}\begin{tabular}{lcccccc}
\hline Character & $A(-)$ & $\omega_i q^{aA/e}$ & $\kappa=1$ & $\kappa=3$ & $\kappa=5$ & $\kappa=7$
\\\hline $\phi_{112,3}$ & $48$ & $q^{6}$ & $12$ & $36$ & $60$ & $84$
\\ $\phi_{700,28}$ & $98$ & $q^{15}$ & $25$ & $73$ & $123$ & $171$
\\ $D_4;\phi_{9,6}'$ & $86$ & $-q^{12}$ & $22$ & $66$ & $106$ & $150$
\\ $\phi_{700,6}$ & $69$ & $q^{9}$ & $17$ & $51$ & $87$ & $121$
\\ $\phi_{400,43}$ & $105$ & $q^{18}$ & $26$ & $78$ & $132$ & $184$
\\ $\phi_{972,12}$ & $81$ & $q^{11}$ & $20$ & $62$ & $100$ & $142$
\\ $D_4;\phi_{9,10}$ & $101$ & $-q^{16}$ & $25$ & $75$ & $127$ & $177$
\\ $\phi_{84,64}$ & $108$ & $q^{21}$ & $27$ & $81$ & $135$ & $189$
\\\hline\end{tabular}\end{center}

\begin{center}\begin{tikzpicture}[thick,scale=1.8]
\draw (0,-0.18) node{${D_4;\phi_{9,6}'}$};
\draw (1,-0.18) node{$D_4;\phi_{9,10}$};
\draw (3,-0.18) node{$\phi_{84,64}$};
\draw (4,-0.18) node{$\phi_{400,43}$};
\draw (5,-0.18) node{$\phi_{700,28}$};
\draw (6,-0.18) node{$\phi_{972,12}$};
\draw (7,-0.18) node{$\phi_{700,6}$};
\draw (8,-0.18) node{$\phi_{}112,3$};

\draw (0,0) -- (8,0);
\draw (8,0) node [draw] (l0) {};
\draw (7,0) node [draw] (l1) {};
\draw (6,0) node [draw] (l2) {};
\draw (5,0) node [draw] (l2) {};
\draw (4,0) node [draw] (l2) {};
\draw (3,0) node [draw] (l2) {};
\draw (2,0) node [fill=black!100] (ld) {};
\draw (1,0) node [draw] (l4) {};
\draw (0,0) node [draw] (l4) {};
\draw (0,0.18) node{$22$};
\draw (1,0.18) node{$25$};
\draw (3,0.18) node{$27$};
\draw (4,0.18) node{$26$};
\draw (5,0.18) node{$25$};
\draw (6,0.18) node{$20$};
\draw (7,0.18) node{$17$};
\draw (8,0.18) node{$12$};
\end{tikzpicture}\end{center}

\noindent(vi)\;\; \textbf{Block 6:} Cuspidal pair $(\Phi_8.{}^2\!D_4(q)$,$\phi_{012,3})$, of degree $q^3\Phi_6\Phi_8$. There are eight unipotent characters in the block, all of which are real.

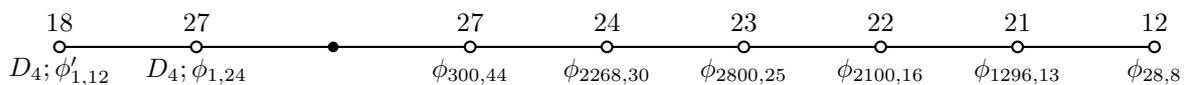
\begin{center}\begin{tabular}{lcccccc}
\hline Character & $A(-)$ & $\omega_i q^{aA/e}$ & $\kappa=1$ & $\kappa=3$ & $\kappa=5$ & $\kappa=7$
\\\hline $\phi_{28,8}$ & $48$ & $q^{6}$ & $12$ & $36$ & $60$ & $84$
\\ $\phi_{2800,25}$ & $98$ & $q^{15}$ & $23$ & $73$ & $123$ & $173$
\\ $\phi_{2268,30}$ & $$101 & $q^{16}$ & $24$ & $76$ & $126$ & $178$
\\ $D_4;\phi_{1,24}$ & $108$ & $-q^{21}$ & $27$ & $81$ & $135$ & $189$
\\ $\phi_{300,44}$ & $105$ & $q^{18}$ & $27$ & $79$ & $131$ & $183$
\\ $\phi_{1296,13}$ & $81$ & $q^{11}$ & $21$ & $61$ & $101$ & $141$
\\ $\phi_{2100,16}$ & $86$ & $q^{12}$ & $22$ & $64$ & $108$ & $150$
\\ $D_4;\phi_{1,12}'$ & $69$ & $-q^{9}$ & $18$ & $52$ & $86$ & $120$
\\\hline\end{tabular}\end{center}

\begin{center}\begin{tikzpicture}[thick,scale=1.8]
\draw (0,-0.18) node{${D_4;\phi_{1,12}'}$};
\draw (1,-0.18) node{$D_4;\phi_{1,24}$};
\draw (3,-0.18) node{$\phi_{300,44}$};
\draw (4,-0.18) node{$\phi_{2268,30}$};
\draw (5,-0.18) node{$\phi_{2800,25}$};
\draw (6,-0.18) node{$\phi_{2100,16}$};
\draw (7,-0.18) node{$\phi_{1296,13}$};
\draw (8,-0.18) node{$\phi_{28,8}$};

\draw (0,0) -- (8,0);
\draw (8,0) node [draw] (l0) {};
\draw (7,0) node [draw] (l1) {};
\draw (6,0) node [draw] (l2) {};
\draw (5,0) node [draw] (l2) {};
\draw (4,0) node [draw] (l2) {};
\draw (3,0) node [draw] (l2) {};
\draw (2,0) node [fill=black!100] (ld) {};
\draw (1,0) node [draw] (l4) {};
\draw (0,0) node [draw] (l4) {};
\draw (0,0.18) node{$18$};
\draw (1,0.18) node{$27$};
\draw (3,0.18) node{$27$};
\draw (4,0.18) node{$24$};
\draw (5,0.18) node{$23$};
\draw (6,0.18) node{$22$};
\draw (7,0.18) node{$21$};
\draw (8,0.18) node{$12$};
\end{tikzpicture}\end{center}

\bigskip

\noindent\textbf{Proof of Brauer trees}: Since all characters are real, the Brauer tree is a line. A Geck--Pfeiffer argument then a degree argument is enough.

\newpage
\subsection{$d=9$}

For $E_8(q)$ and $d=9$ there are three unipotent blocks of weight $1$, together with $112$ unipotent blocks of defect zero. Note that two of the Brauer trees in this case depend on the corresponding Brauer tree for $E_7$, $d=9$ being known. Since the shape is known for that tree, the shapes of each of these trees is correct, and it is the labelling of the cuspidal characters that is needed.
\\[1cm]
\noindent(i)\;\; \textbf{Block 1:} Cuspidal pair is $(\Phi_{9}.A_2(q)$,$\phi_{3})$, of degree $1$. There are eighteen unipotent characters in the block, four of which are non-real.

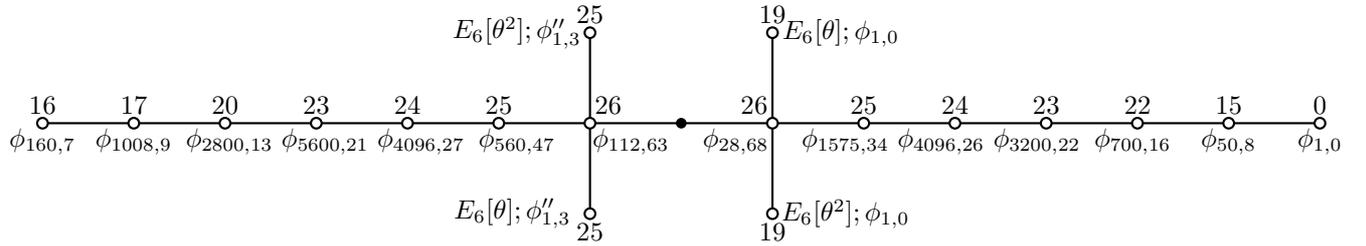
\begin{center}\begin{tabular}{lcccccccc}
\hline Character & $A(-)$ & $\omega_i q^{aA/e}$ & $\kappa=1$ & $\kappa=2$ & $\kappa=4$ & $\kappa=5$ & $\kappa=7$ & $\kappa=8$
\\\hline $\phi_{1,0}$ & $0$ & $q^{0}$ & $0$ & $0$ & $0$ & $0$ & $0$ & $0$
\\ $\phi_{1008,9}$ & $83$ & $-q^{5}$ & $17$ & $37$ & $75$ & $91$ & $129$ & $149$
\\ $\phi_{28,68}$ & $117$ & $q^{10}$ & $26$ & $52$ & $104$ & $130$ & $182$ & $208$
\\ $\phi_{2800,13}$ & $95$ & $-q^{6}$ & $20$ & $44$ & $84$ & $106$ & $146$ & $170$
\\ $E_6[\theta^2];\phi_{1,3}''$ & $112$ & $\theta q^{8}$ & $25$ & $49$ & $99$ & $125$ & $175$ & $199$
\\ $\phi_{5600,21}$ & $105$ & $-q^{7}$ & $23$ & $45$ & $95$ & $115$ & $165$ & $187$
\\ $\phi_{4096,27}$ & $109$ & $-q^{15/2}$ & $24$ & $48$ & $96$ & $122$ & $170$ & $194$
\\ $E_6[\theta^2];\phi_{1,0}$ & $83$ & $-\theta q^{5}$ & $19$ & $37$ & $73$ & $93$ & $129$ & $147$
\\ $\phi_{50,8}$ & $68$ & $q^{4}$ & $15$ & $31$ & $61$ & $75$ & $105$ & $121$
\\ $\phi_{560,47}$ & $115$ & $-q^{9}$ & $25$ & $51$ & $103$ & $127$ & $179$ & $205$
\\ $E_6[\theta];\phi_{1,3}''$ & $112$ & $\theta^2 q^{8}$ & $25$ & $49$ & $99$ & $125$ & $175$ & $199$
\\ $\phi_{112,63}$ & $117$ & $-q^{10}$ & $26$ & $52$ & $104$ & $130$ & $182$ & $208$
\\ $\phi_{700,16}$ & $95$ & $q^{6}$ & $22$ & $42$ & $86$ & $104$ & $148$ & $168$
\\ $E_6[\theta];\phi_{1,0}$ & $83$ & $-\theta^2q^{5}$ & $19$ & $37$ & $73$ & $93$ & $129$ & $147$
\\ $\phi_{3200,22}$ & $105$ & $q^{7}$ & $23$ & $47$ & $93$ & $117$ & $163$ & $187$
\\ $\phi_{4096,26}$ & $109$ & $q^{15/2}$ & $24$ & $48$ & $96$ & $122$ & $170$ & $194$
\\ $\phi_{1575,34}$ & $112$ & $q^{8}$ & $25$ & $51$ & $101$ & $123$ & $173$ & $199$
\\ $\phi_{160,7}$ & $68$ & $-q^{4}$ & $16$ & $30$ & $60$ & $76$ & $106$ & $120$
\\\hline\end{tabular}\end{center}
\begin{center}\begin{tikzpicture}[thick,scale=1.2]

\draw (8.75,1) node{$E_6[\theta];\phi_{1,0}$};
\draw (8.8,-1) node{$E_6[\theta^2];\phi_{1,0}$};

\draw (5.2,1) node{$E_6[\theta^2];\phi_{1,3}''$};
\draw (5.15,-1) node{$E_6[\theta];\phi_{1,3}''$};

\draw (0,-0.2) node{$\phi_{160,7}$};
\draw (1,-0.2) node{$\phi_{1008,9}$};
\draw (2.05,-0.2) node{$\phi_{2800,13}$};
\draw (3.1,-0.2) node{$\phi_{5600,21}$};
\draw (4.15,-0.2) node{$\phi_{4096,27}$};
\draw (5.2,-0.2) node{$\phi_{560,47}$};
\draw (6.45,-0.2) node{$\phi_{112,63}$};

\draw (7.6,-0.2) node{$\phi_{28,68}$};
\draw (8.8,-0.2) node{$\phi_{1575,34}$};
\draw (9.85,-0.2) node{$\phi_{4096,26}$};
\draw (10.9,-0.2) node{$\phi_{3200,22}$};
\draw (11.95,-0.2) node{$\phi_{700,16}$};
\draw (13,-0.2) node{$\phi_{50,8}$};
\draw (14,-0.2) node{$\phi_{1,0}$};

\draw (0,0.2) node{$16$};
\draw (1,0.2) node{$17$};
\draw (2,0.2) node{$20$};
\draw (3,0.2) node{$23$};
\draw (4,0.2) node{$24$};
\draw (5,0.2) node{$25$};
\draw (6.2,0.2) node{$26$};

\draw (6,1.2) node{$25$};
\draw (6,-1.2) node{$25$};
\draw (8,1.2) node{$19$};
\draw (8,-1.2) node{$19$};

\draw (14,0.2) node{$0$};
\draw (13,0.2) node{$15$};
\draw (12,0.2) node{$22$};
\draw (11,0.2) node{$23$};
\draw (10,0.2) node{$24$};
\draw (9,0.2) node{$25$};
\draw (7.8,0.2) node{$26$};

\draw (0,0) -- (14,0);
\draw (6,-1) -- (6,1);
\draw (8,-1) -- (8,1);

\draw (7,0) node [fill=black!100] (ld) {};

\draw (14,0) node [draw] (l2) {};
\draw (13,0) node [draw] (l2) {};
\draw (12,0) node [draw] (l2) {};
\draw (11,0) node [draw] (l2) {};
\draw (10,0) node [draw] (l2) {};
\draw (9,0) node [draw] (l2) {};
\draw (8,0) node [draw] (l2) {};
\draw (6,0) node [draw] (l2) {};
\draw (5,0) node [draw] (l2) {};
\draw (4,0) node [draw] (l2) {};
\draw (3,0) node [draw] (l2) {};
\draw (2,0) node [draw] (l2) {};
\draw (1,0) node [draw] (l2) {};
\draw (0,0) node [draw] (l2) {};
\draw (6,1) node [draw] (l2) {};
\draw (6,-1) node [draw] (l2) {};
\draw (8,1) node [draw] (l2) {};
\draw (8,-1) node [draw] (l2) {};
\end{tikzpicture}\end{center}

\noindent(ii)\;\; \textbf{Block 2:} Cuspidal pair is $(\Phi_{9}.A_2(q)$,$\phi_{21})$, of degree $q\Phi_2$. There are eighteen unipotent characters in the block, four of which are non-real.

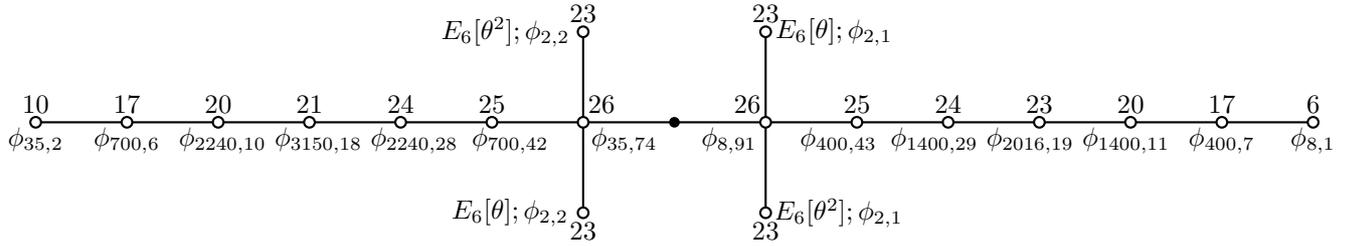
\begin{center}\begin{tabular}{lcccccccc}
\hline Character & $A(-)$ & $\omega_i q^{aA/e}$ & $\kappa=1$ & $\kappa=2$ & $\kappa=4$ & $\kappa=5$ & $\kappa=7$ & $\kappa=8$
\\\hline $\phi_{8,1}$ & $27$ & $q^{3/2}$ & $6$ & $12$ & $24$ & $30$ & $42$ & $48$
\\ $\phi_{3150,18}$ & $102$ & $-q^{13/2}$ & $21$ & $47$ & $93$ & $111$ & $157$ & $183$
\\ $\phi_{8,91}$ & $117$ & $q^{23/2}$ & $26$ & $52$ & $104$ & $130$ & $182$ & $208$
\\ $\phi_{2240,28}$ & $108$ & $-q^{15/2}$ & $24$ & $48$ & $96$ & $120$ & $168$ & $192$
\\ $E_6[\theta];\phi_{2,2}$ & $102$ & $\theta^2 q^{13/2}$ & $23$ & $45$ & $91$ & $113$ & $159$ & $181$
\\ $\phi_{700,42}$ & $112$ & $-q^{17/2}$ & $25$ & $49$ & $101$ & $123$ & $175$ & $199$
\\ $\phi_{400,7}$ & $76$ & $q^{9/2}$ & $17$ & $35$ & $67$ & $85$ & $117$ & $135$
\\ $E_6[\theta^2];\phi_{2,1}$ & $102$ & $-\theta q^{13/2}$ & $23$ & $45$ & $89$ & $115$ & $159$ & $181$
\\ $\phi_{1400,11}$ & $90$ & $q^{11/2}$ & $20$ & $40$ & $82$ & $98$ & $140$ & $160$
\\ $\phi_{35,74}$ & $116$ & $-q^{21/2}$ & $26$ & $52$ & $104$ & $128$ & $180$ & $206$
\\ $\phi_{2016,19}$ & $102$ & $q^{13/2}$ & $23$ & $45$ & $91$ & $113$ & $159$ & $181$
\\ $\phi_{35,2}$ & $44$ & $-q^{5/2}$ & $10$ & $20$ & $40$ & $48$ & $68$ & $78$
\\ $\phi_{1400,29}$ & $108$ & $q^{15/2}$ & $24$ & $48$ & $98$ & $118$ & $168$ & $192$
\\ $E_6[\theta];\phi_{2,1}$ & $102$ & $-\theta^2 q^{13/2}$ & $23$ & $45$ & $89$ & $115$ & $159$ & $181$
\\ $\phi_{400,43}$ & $112$ & $q^{17/2}$ & $25$ & $51$ & $99$ & $125$ & $173$ & $199$
\\ $\phi_{700,6}$ & $76$ & $-q^{9/2}$ & $17$ & $33$ & $69$ & $83$ & $119$ & $135$
\\ $E_6[\theta^2];\phi_{2,2}$ & $102$ & $\theta q^{13/2}$ & $23$ & $45$ & $91$ & $113$ & $159$ & $181$
\\ $\phi_{2240,10}$ & $90$ & $-q^{11/2}$ & $20$ & $40$ & $80$ & $100$ & $140$ & $160$
\\\hline\end{tabular}\end{center}
\begin{center}\begin{tikzpicture}[thick,scale=1.2]
\draw (8.75,1) node{$E_6[\theta];\phi_{2,1}$};
\draw (8.8,-1) node{$E_6[\theta^2];\phi_{2,1}$};

\draw (5.2,-1) node{$E_6[\theta];\phi_{2,2}$};
\draw (5.15,1) node{$E_6[\theta^2];\phi_{2,2}$};

\draw (0,-0.2) node{$\phi_{35,2}$};
\draw (1,-0.2) node{$\phi_{700,6}$};
\draw (2.05,-0.2) node{$\phi_{2240,10}$};
\draw (3.1,-0.2) node{$\phi_{3150,18}$};
\draw (4.15,-0.2) node{$\phi_{2240,28}$};
\draw (5.2,-0.2) node{$\phi_{700,42}$};
\draw (6.45,-0.2) node{$\phi_{35,74}$};

\draw (7.6,-0.2) node{$\phi_{8,91}$};
\draw (8.8,-0.2) node{$\phi_{400,43}$};
\draw (9.85,-0.2) node{$\phi_{1400,29}$};
\draw (10.9,-0.2) node{$\phi_{2016,19}$};
\draw (11.95,-0.2) node{$\phi_{1400,11}$};
\draw (13,-0.2) node{$\phi_{400,7}$};
\draw (14,-0.2) node{$\phi_{8,1}$};

\draw (0,0.2) node{$10$};
\draw (1,0.2) node{$17$};
\draw (2,0.2) node{$20$};
\draw (3,0.2) node{$21$};
\draw (4,0.2) node{$24$};
\draw (5,0.2) node{$25$};
\draw (6.2,0.2) node{$26$};

\draw (6,1.2) node{$23$};
\draw (6,-1.2) node{$23$};
\draw (8,1.2) node{$23$};
\draw (8,-1.2) node{$23$};

\draw (14,0.2) node{$6$};
\draw (13,0.2) node{$17$};
\draw (12,0.2) node{$20$};
\draw (11,0.2) node{$23$};
\draw (10,0.2) node{$24$};
\draw (9,0.2) node{$25$};
\draw (7.8,0.2) node{$26$};

\draw (0,0) -- (14,0);
\draw (6,-1) -- (6,1);
\draw (8,-1) -- (8,1);

\draw (7,0) node [fill=black!100] (ld) {};

\draw (14,0) node [draw] (l2) {};
\draw (13,0) node [draw] (l2) {};
\draw (12,0) node [draw] (l2) {};
\draw (11,0) node [draw] (l2) {};
\draw (10,0) node [draw] (l2) {};
\draw (9,0) node [draw] (l2) {};
\draw (8,0) node [draw] (l2) {};
\draw (6,0) node [draw] (l2) {};
\draw (5,0) node [draw] (l2) {};
\draw (4,0) node [draw] (l2) {};
\draw (3,0) node [draw] (l2) {};
\draw (2,0) node [draw] (l2) {};
\draw (1,0) node [draw] (l2) {};
\draw (0,0) node [draw] (l2) {};
\draw (6,1) node [draw] (l2) {};
\draw (6,-1) node [draw] (l2) {};
\draw (8,1) node [draw] (l2) {};
\draw (8,-1) node [draw] (l2) {};
\end{tikzpicture}\end{center}

\noindent(iii)\;\; \textbf{Block 3:} Cuspidal pair is $(\Phi_{9}.A_2(q)$,$\phi_{111})$, of degree $q^3$. There are eighteen unipotent characters in the block, four of which are non-real.
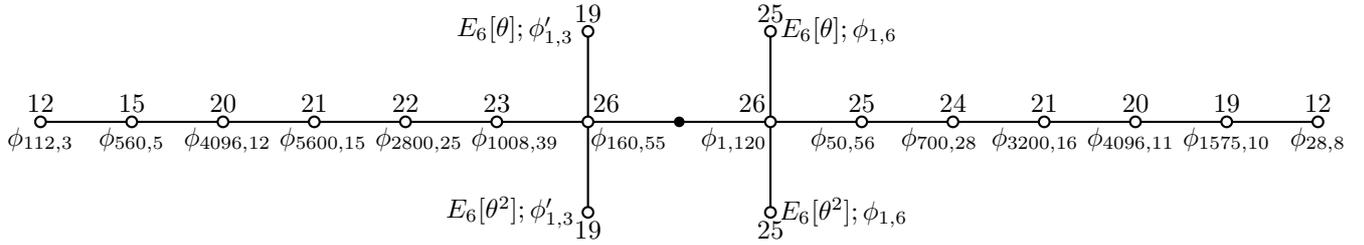
\begin{center}\begin{tabular}{lcccccccc}
\hline Character & $A(-)$ & $\omega_i q^{aA/e}$ & $\kappa=1$ & $\kappa=2$ & $\kappa=4$ & $\kappa=5$ & $\kappa=7$ & $\kappa=8$
\\\hline $\phi_{28,8}$ & $54$ & $q^{3}$ & $12$ & $24$ & $48$ & $60$ & $84$ & $96$
\\ $\phi_{1008,39}$ & $110$ & $-q^{8}$ & $23$ & $49$ & $99$ & $121$ & $171$ & $197$
\\ $\phi_{1,120}$ & $117$ & $q^{13}$ & $26$ & $52$ & $104$ & $130$ & $182$ & $208$
\\ $\phi_{160,55}$ & $113$ & $-q^{9}$ & $26$ & $50$ & $100$ & $126$ & $176$ & $200$
\\ $\phi_{1575,10}$ & $85$ & $q^{5}$ & $19$ & $39$ & $77$ & $93$ & $131$ & $151$
\\ $\phi_{4096,11}$ & $91$ & $q^{11/2}$ & $20$ & $40$ & $80$ & $102$ & $142$ & $162$
\\ $\phi_{3200,16}$ & $96$ & $q^{6}$ & $21$ & $43$ & $85$ & $107$ & $149$ & $171$
\\ $E_6[\theta^2];\phi_{1,6}$ & $110$ & $-\theta q^{8}$ & $25$ & $49$ & $97$ & $123$ & $171$ & $195$
\\ $\phi_{700,28}$ & $104$ & $q^{7}$ & $24$ & $46$ & $94$ & $114$ & $162$ & $184$
\\ $\phi_{112,3}$ & $54$ & $-q^{3}$ & $12$ & $24$ & $48$ & $60$ & $84$ & $96$
\\ $E_6[\theta^2];\phi_{1,3}'$ & $85$ & $\theta q^{5}$ & $19$ & $37$ & $75$ & $95$ & $133$ & $151$
\\ $\phi_{560,5}$ & $70$ & $-q^{4}$ & $15$ & $31$ & $63$ & $77$ & $109$ & $125$
\\ $\phi_{50,56}$ & $113$ & $q^{9}$ & $25$ & $51$ & $101$ & $125$ & $175$ & $201$
\\ $E_6[\theta];\phi_{1,6}$ & $110$ & $-\theta^2 q^{8}$ & $25$ & $49$ & $97$ & $123$ & $171$ & $195$
\\ $\phi_{4096,12}$ & $91$ & $-q^{11/2}$ & $20$ & $40$ & $80$ & $102$ & $142$ & $162$
\\ $\phi_{5600,15}$ & $96$ & $-q^{6}$ & $21$ & $41$ & $87$ & $105$ & $151$ & $171$
\\ $E_6[\theta];\phi_{1,3}'$ & $85$ & $\theta^2 q^{5}$ & $19$ & $37$ & $75$ & $95$ & $133$ & $151$
\\ $\phi_{2800,25}$ & $104$ & $-q^{7}$ & $22$ & $48$ & $92$ & $116$ & $160$ & $186$
\\\hline\end{tabular}\end{center}
\begin{center}\begin{tikzpicture}[thick,scale=1.2]
\draw (8.75,1) node{$E_6[\theta];\phi_{1,6}$};
\draw (8.8,-1) node{$E_6[\theta^2];\phi_{1,6}$};

\draw (5.2,1) node{$E_6[\theta];\phi_{1,3}'$};
\draw (5.15,-1) node{$E_6[\theta^2];\phi_{1,3}'$};

\draw (0,-0.2) node{$\phi_{112,3}$};
\draw (1,-0.2) node{$\phi_{560,5}$};
\draw (2.05,-0.2) node{$\phi_{4096,12}$};
\draw (3.1,-0.2) node{$\phi_{5600,15}$};
\draw (4.15,-0.2) node{$\phi_{2800,25}$};
\draw (5.2,-0.2) node{$\phi_{1008,39}$};
\draw (6.45,-0.2) node{$\phi_{160,55}$};

\draw (7.6,-0.2) node{$\phi_{1,120}$};
\draw (8.8,-0.2) node{$\phi_{50,56}$};
\draw (9.85,-0.2) node{$\phi_{700,28}$};
\draw (10.9,-0.2) node{$\phi_{3200,16}$};
\draw (11.95,-0.2) node{$\phi_{4096,11}$};
\draw (13,-0.2) node{$\phi_{1575,10}$};
\draw (14,-0.2) node{$\phi_{28,8}$};

\draw (0,0.2) node{$12$};
\draw (1,0.2) node{$15$};
\draw (2,0.2) node{$20$};
\draw (3,0.2) node{$21$};
\draw (4,0.2) node{$22$};
\draw (5,0.2) node{$23$};
\draw (6.2,0.2) node{$26$};

\draw (6,1.2) node{$19$};
\draw (6,-1.2) node{$19$};
\draw (8,1.2) node{$25$};
\draw (8,-1.2) node{$25$};

\draw (14,0.2) node{$12$};
\draw (13,0.2) node{$19$};
\draw (12,0.2) node{$20$};
\draw (11,0.2) node{$21$};
\draw (10,0.2) node{$24$};
\draw (9,0.2) node{$25$};
\draw (7.8,0.2) node{$26$};

\draw (0,0) -- (14,0);
\draw (6,-1) -- (6,1);
\draw (8,-1) -- (8,1);

\draw (7,0) node [fill=black!100] (ld) {};

\draw (14,0) node [draw] (l2) {};
\draw (13,0) node [draw] (l2) {};
\draw (12,0) node [draw] (l2) {};
\draw (11,0) node [draw] (l2) {};
\draw (10,0) node [draw] (l2) {};
\draw (9,0) node [draw] (l2) {};
\draw (8,0) node [draw] (l2) {};
\draw (6,0) node [draw] (l2) {};
\draw (5,0) node [draw] (l2) {};
\draw (4,0) node [draw] (l2) {};
\draw (3,0) node [draw] (l2) {};
\draw (2,0) node [draw] (l2) {};
\draw (1,0) node [draw] (l2) {};
\draw (0,0) node [draw] (l2) {};
\draw (6,1) node [draw] (l2) {};
\draw (6,-1) node [draw] (l2) {};
\draw (8,1) node [draw] (l2) {};
\draw (8,-1) node [draw] (l2) {};
\end{tikzpicture}\end{center}

\noindent\textbf{Proof of Brauer trees}: This is given in \cite{cdr2012un}.

\newpage
\subsection{$d=10$}

For $E_8(q)$ and $d=10$ there are two unipotent blocks of weight $1$, together with the principal block with cyclotomic Weyl group $G_{16}$, and $101$ unipotent blocks of defect zero.
\\[1cm]
\noindent(i)\;\; \textbf{Block 1:} Cuspidal pair is $(\Phi_{10}.{}^2\!A_4(q)$,$\phi_{32})$, of degree $q^2\Phi_{10}$. There are ten unipotent characters in the block, all of which are real.

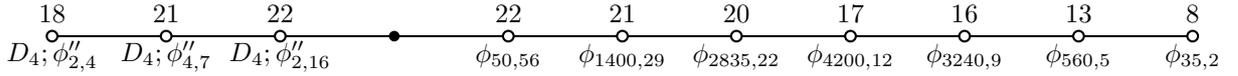
\begin{center}\begin{tabular}{lcccccc}
\hline Character & $A(-)$ & $\omega_i q^{aA/e}$ & $\kappa=1$ & $\kappa=3$ & $\kappa=7$ & $\kappa=9$
\\\hline $\phi_{35,2}$ & $40$ & $q^{4}$ & $8$ & $24$ & $56$ & $72$
\\ $D_4;\phi_{2,4}''$ & $90$ & $-q^{10}$ & $18$ & $54$ & $126$ & $162$
\\ $\phi_{50,56}$ & $110$ & $q^{16}$ & $22$ & $66$ & $154$ & $198$
\\ $\phi_{560,5}$ & $67$ & $q^{7}$ & $13$ & $41$ & $93$ & $121$
\\ $D_4;\phi_{4,7}''$ & $104$ & $-q^{13}$ & $21$ & $63$ & $145$ & $187$
\\ $\phi_{3240,9}$ & $83$ & $q^{9}$ & $16$ & $50$ & $116$ & $150$
\\ $\phi_{4200,12}$ & $90$ & $q^{10}$ & $17$ & $55$ & $125$ & $163$
\\ $D_4;\phi_{2,16}''$ & $110$ & $-q^{16}$ & $22$ & $66$ & $154$ & $198$
\\ $\phi_{2835,22}$ & $100$ & $q^{12}$ & $20$ & $60$ & $140$ & $180$
\\ $\phi_{1400,29}$ & $104$ & $q^{13}$ & $21$ & $63$ & $145$ & $187$
\\\hline\end{tabular}\end{center}
\begin{center}\begin{tikzpicture}[thick,scale=1.5]
\draw (0,-0.18) node{${D_4;\phi_{2,4}''}$};
\draw (1,-0.18) node{$D_4;\phi_{4,7}''$};
\draw (2,-0.18) node{$D_4;\phi_{2,16}''$};
\draw (4,-0.18) node{$\phi_{50,56}$};
\draw (5,-0.18) node{$\phi_{1400,29}$};
\draw (6,-0.18) node{$\phi_{2835,22}$};
\draw (7,-0.18) node{$\phi_{4200,12}$};
\draw (8,-0.18) node{$\phi_{3240,9}$};
\draw (9,-0.18) node{$\phi_{560,5}$};
\draw (10,-0.18) node{$\phi_{35,2}$};

\draw (0,0) -- (10,0);
\draw (10,0) node [draw] (l0) {};
\draw (9,0) node [draw] (l0) {};
\draw (8,0) node [draw] (l0) {};
\draw (7,0) node [draw] (l1) {};
\draw (6,0) node [draw] (l2) {};
\draw (5,0) node [draw] (l2) {};
\draw (4,0) node [draw] (l2) {};
\draw (2,0) node [draw] (l2) {};
\draw (3,0) node [fill=black!100] (ld) {};
\draw (1,0) node [draw] (l4) {};
\draw (0,0) node [draw] (l4) {};
\draw (0,0.2) node{$18$};
\draw (1,0.2) node{$21$};
\draw (2,0.2) node{$22$};
\draw (4,0.2) node{$22$};
\draw (5,0.2) node{$21$};
\draw (6,0.2) node{$20$};
\draw (7,0.2) node{$17$};
\draw (8,0.2) node{$16$};
\draw (9,0.2) node{$13$};
\draw (10,0.2) node{$8$};
\end{tikzpicture}\end{center}

\noindent(ii)\;\; \textbf{Block 2:} Cuspidal pair is $(\Phi_{10}.{}^2\!A_4(q)$,$\phi_{221})$, of degree $q^4\Phi_{10}$. There are ten unipotent characters in the block, all of which are real.

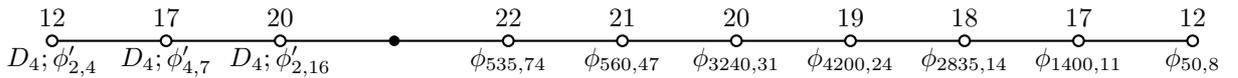
\begin{center}\begin{tabular}{lcccccc}
\hline Character & $A(-)$ & $\omega_i q^{aA/e}$ & $\kappa=1$ & $\kappa=3$ & $\kappa=7$ & $\kappa=9$
\\\hline $\phi_{50,8}$ & $60$ & $q^{6}$ & $12$ & $36$ & $84$ & $108$
\\ $D_4;\phi_{2,16}'$ & $100$ & $-q^{12}$ & $20$ & $60$ & $140$ & $180$
\\ $\phi_{35,74}$ & $110$ & $q^{18}$ & $22$ & $66$ & $154$ & $198$
\\ $\phi_{1400,11}$ & $84$ & $q^{9}$ & $17$ & $51$ & $117$ & $151$
\\ $\phi_{2835,14}$ & $90$ & $q^{10}$ & $18$ & $54$ & $126$ & $162$
\\ $D_4;\phi_{2,4}'$ & $60$ & $-q^{6}$ & $12$ & $36$ & $84$ & $108$
\\ $\phi_{4200,24}$ & $100$ & $q^{12}$ & $19$ & $61$ & $139$ & $181$
\\ $\phi_{3240,31}$ & $103$ & $q^{13}$ & $20$ & $62$ & $144$ & $186$
\\ $D_4;\phi_{4,7}'$ & $84$ & $-q^{9}$ & $17$ & $51$ & $117$ & $151$
\\ $\phi_{560,47}$ & $107$ & $q^{15}$ & $21$ & $65$ & $149$ & $193$
\\\hline\end{tabular}\end{center}
\begin{center}\begin{tikzpicture}[thick,scale=1.5]
\draw (0,-0.18) node{${D_4;\phi_{2,4}'}$};
\draw (1,-0.18) node{$D_4;\phi_{4,7}'$};
\draw (2,-0.18) node{$D_4;\phi_{2,16}'$};
\draw (4,-0.18) node{$\phi_{535,74}$};
\draw (5,-0.18) node{$\phi_{560,47}$};
\draw (6,-0.18) node{$\phi_{3240,31}$};
\draw (7,-0.18) node{$\phi_{4200,24}$};
\draw (8,-0.18) node{$\phi_{2835,14}$};
\draw (9,-0.18) node{$\phi_{1400,11}$};
\draw (10,-0.18) node{$\phi_{50,8}$};

\draw (0,0) -- (10,0);
\draw (10,0) node [draw] (l0) {};
\draw (9,0) node [draw] (l0) {};
\draw (8,0) node [draw] (l0) {};
\draw (7,0) node [draw] (l1) {};
\draw (6,0) node [draw] (l2) {};
\draw (5,0) node [draw] (l2) {};
\draw (4,0) node [draw] (l2) {};
\draw (2,0) node [draw] (l2) {};
\draw (3,0) node [fill=black!100] (ld) {};
\draw (1,0) node [draw] (l4) {};
\draw (0,0) node [draw] (l4) {};
\draw (0,0.2) node{$12$};
\draw (1,0.2) node{$17$};
\draw (2,0.2) node{$20$};
\draw (4,0.2) node{$22$};
\draw (5,0.2) node{$21$};
\draw (6,0.2) node{$20$};
\draw (7,0.2) node{$19$};
\draw (8,0.2) node{$18$};
\draw (9,0.2) node{$17$};
\draw (10,0.2) node{$12$};
\end{tikzpicture}\end{center}

\bigskip

\noindent\textbf{Proof of Brauer trees}: Since all characters are real, the Brauer tree is a line. A Geck--Pfeiffer argument then a degree argument is enough.

\newpage
\subsection{$d=12$}

For $E_8(q)$ and $d=12$ there are four unipotent blocks of weight $1$, together with the principal block, with cyclotomic Weyl group $G_{10}$, and $70$ unipotent blocks of defect zero.
\\[1cm]
\noindent(i)\;\;\textbf{Block 1:} Cuspidal pair is $(\Phi_{12}.{}^3\!D_4(q),\phi_{1,3}')$, of degree $q\Phi_{12}$. There are twelve unipotent characters in the block, two of which -- $E_6[\theta^i];\phi_{1,6}$ -- are non-real.

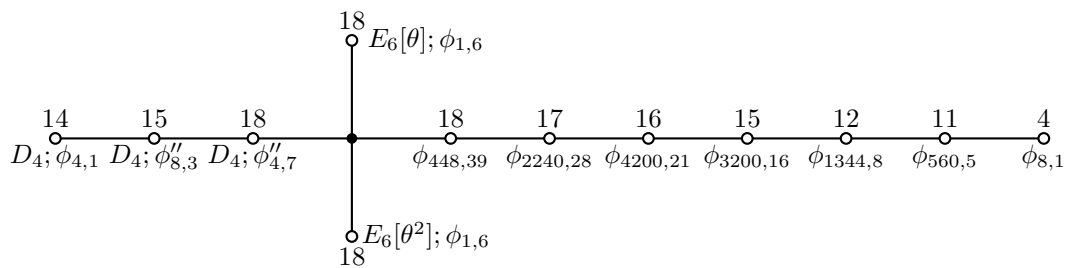
\begin{center}\begin{tabular}{lcccccc}
\hline Character & $A(-)$ & $\omega_i q^{aA/e}$ & $\kappa=1$ & $\kappa=5$ & $\kappa=7$ & $\kappa=11$
\\\hline $\phi_{8,1}$ & $24$ & $q^{2}$ & $4$ & $20$ & $28$ & $44$
\\ $D_4;\phi_{8,3}''$ & $94$ & $-q^{9}$ & $15$ & $79$ & $109$ & $173$
\\ $E_6[\theta];\phi_{1,6}$ & $108$ & $\theta q^{12}$ & $18$ & $90$ & $126$ & $198$
\\ $D_4;\phi_{4,7}''$ & $105$ & $-q^{11}$ & $18$ & $88$ & $122$ & $192$
\\ $\phi_{560,5}$ & $68$ & $q^{6}$ & $11$ & $57$ & $79$ & $125$
\\ $\phi_{1344,8}$ & $78$ & $q^{7}$ & $12$ & $64$ & $92$ & $144$
\\ $E_6[\theta^2];\phi_{1,6}$ & $108$ & $\theta^2 q^{12}$ & $18$ & $90$ & $126$ & $198$
\\ $\phi_{3200,16}$ & $94$ & $q^{9}$ & $15$ & $79$ & $109$ & $173$
\\ $\phi_{4200,21}$ & $100$ & $q^{10}$ & $16$ & $84$ & $116$ & $184$
\\ $\phi_{2240,28}$ & $105$ & $q^{11}$ & $17$ & $87$ & $123$ & $193$
\\ $\phi_{448,39}$ & $108$ & $q^{12}$ & $18$ & $90$ & $126$ & $198$
\\ $D_4;\phi_{4,1}$ & $78$ & $-q^{7}$ & $14$ & $66$ & $90$ & $142$
\\\hline\end{tabular}\end{center}
\begin{center}\begin{tikzpicture}[thick,scale=1.3]
\draw (0,-0.2) node{$D_4;\phi_{4,1}$};
\draw (1,-0.2) node{$D_4;\phi_{8,3}''$};
\draw (2,-0.2) node{$D_4;\phi_{4,7}''$};
\draw (4,-0.2) node{$\phi_{448,39}$};
\draw (5,-0.2) node{$\phi_{2240,28}$};
\draw (6,-0.2) node{$\phi_{4200,21}$};
\draw (7,-0.2) node{$\phi_{3200,16}$};
\draw (8,-0.2) node{$\phi_{1344,8}$};
\draw (9,-0.2) node{$\phi_{560,5}$};
\draw (10,-0.2) node{$\phi_{8,1}$};

\draw (3.75,1) node {$E_6[\theta];\phi_{1,6}$};
\draw (3.75,-1) node {$E_6[\theta^2];\phi_{1,6}$};

\draw (0,0.2) node{$14$};
\draw (1,0.2) node{$15$};
\draw (2,0.2) node{$18$};

\draw (3,-1.2) node{$18$};
\draw (3,1.2) node{$18$};

\draw (4,0.2) node{$18$};
\draw (5,0.2) node{$17$};
\draw (6,0.2) node{$16$};
\draw (7,0.2) node{$15$};
\draw (8,0.2) node{$12$};
\draw (9,0.2) node{$11$};
\draw (10,0.2) node{$4$};

\draw (0,0) -- (10,0);
\draw (3,-1) -- (3,1);

\draw (10,0) node [draw] (l0) {};
\draw (9,0) node [draw] (l0) {};
\draw (8,0) node [draw] (l0) {};
\draw (7,0) node [draw] (l1) {};
\draw (6,0) node [draw] (l2) {};
\draw (5,0) node [draw] (l2) {};
\draw (4,0) node [draw] (l2) {};
\draw (3,1) node [draw] (l2) {};
\draw (3,-1) node [draw] (l2) {};
\draw (3,0) node [fill=black!100] (ld) {};
\draw (2,0) node [draw] (l2) {};
\draw (1,0) node [draw] (l4) {};
\draw (0,0) node [draw] (l4) {};
\end{tikzpicture}\end{center}
\noindent(ii)\;\;\textbf{Block 2:} Cuspidal pair is $(\Phi_{12}.{}^3\!D_4(q),\phi_{2,2})$, of degree $q^3\Phi_2^2\Phi_{12}/2$. There are twelve unipotent characters in the block, two of which -- $E_6[\theta^i];\phi_{2,2}$ -- are non-real.

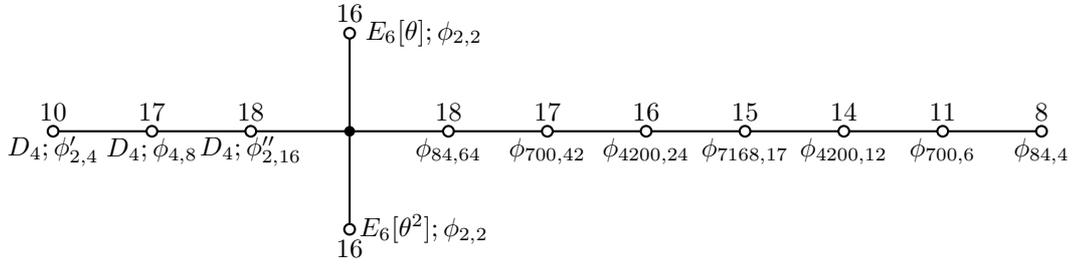
\begin{center}\begin{tabular}{lcccccc}
\hline Character & $A(-)$ & $\omega_i q^{aA/e}$ & $\kappa=1$ & $\kappa=5$ & $\kappa=7$ & $\kappa=11$
\\\hline $\phi_{84,4}$ & $48$ & $q^{4}$ & $8$ & $40$ & $56$ & $88$
\\ $E_6[\theta^2];\phi_{2,2}$ & $95$ & $\theta^2 q^{9}$ & $16$ & $80$ & $110$ & $174$
\\ $\phi_{700,6}$ & $69$ & $q^{6}$ & $11$ & $59$ & $79$ & $127$
\\ $D_4;\phi_{2,16}''$ & $107$ & $-q^{13}$ & $18$ & $90$ & $124$ & $196$
\\ $\phi_{4200,12}$ & $87$ & $q^{8}$ & $14$ & $74$ & $100$ & $160$
\\ $\phi_{7168,17}$ & $95$ & $q^{9}$ & $15$ & $77$ & $113$ & $175$
\\ $\phi_{4200,24}$ & $99$ & $q^{10}$ & $16$ & $84$ & $114$ & $182$
\\ $D_4;\phi_{2,4}'$ & $59$ & $-q^{5}$ & $10$ & $50$ & $68$ & $108$
\\ $\phi_{700,42}$ & $105$ & $q^{12}$ & $17$ & $89$ & $121$ & $193$
\\ $E_6[\theta];\phi_{2,2}$ & $95$ & $\theta q^{9}$ & $16$ & $80$ & $110$ & $174$
\\ $\phi_{84,64}$ & $108$ & $q^{14}$ & $18$ & $90$ & $126$ & $198$
\\ $D_4;\phi_{4,8}$ & $95$ & $-q^{9}$ & $17$ & $81$ & $109$ & $173$
\\\hline\end{tabular}\end{center}
\begin{center}\begin{tikzpicture}[thick,scale=1.3]
\draw (0,-0.2) node{$D_4;\phi_{2,4}'$};
\draw (1,-0.2) node{$D_4;\phi_{4,8}$};
\draw (2,-0.2) node{$D_4;\phi_{2,16}''$};
\draw (4,-0.2) node{$\phi_{84,64}$};
\draw (5,-0.2) node{$\phi_{700,42}$};
\draw (6,-0.2) node{$\phi_{4200,24}$};
\draw (7,-0.2) node{$\phi_{7168,17}$};
\draw (8,-0.2) node{$\phi_{4200,12}$};
\draw (9,-0.2) node{$\phi_{700,6}$};
\draw (10,-0.2) node{$\phi_{84,4}$};

\draw (3.75,1) node {$E_6[\theta];\phi_{2,2}$};
\draw (3.75,-1) node {$E_6[\theta^2];\phi_{2,2}$};

\draw (0,0.2) node{$10$};
\draw (1,0.2) node{$17$};
\draw (2,0.2) node{$18$};

\draw (3,-1.2) node{$16$};
\draw (3,1.2) node{$16$};

\draw (4,0.2) node{$18$};
\draw (5,0.2) node{$17$};
\draw (6,0.2) node{$16$};
\draw (7,0.2) node{$15$};
\draw (8,0.2) node{$14$};
\draw (9,0.2) node{$11$};
\draw (10,0.2) node{$8$};

\draw (0,0) -- (10,0);
\draw (3,-1) -- (3,1);

\draw (10,0) node [draw] (l0) {};
\draw (9,0) node [draw] (l0) {};
\draw (8,0) node [draw] (l0) {};
\draw (7,0) node [draw] (l1) {};
\draw (6,0) node [draw] (l2) {};
\draw (5,0) node [draw] (l2) {};
\draw (4,0) node [draw] (l2) {};
\draw (3,1) node [draw] (l2) {};
\draw (3,-1) node [draw] (l2) {};
\draw (3,0) node [fill=black!100] (ld) {};
\draw (2,0) node [draw] (l2) {};
\draw (1,0) node [draw] (l4) {};
\draw (0,0) node [draw] (l4) {};
\end{tikzpicture}\end{center}

\noindent(iii)\;\;\textbf{Block 3:} Cuspidal pair is $(\Phi_{12}.{}^3\!D_4(q),\phi_{1,3}'')$, of degree $q^7\Phi_{12}$. There are twelve unipotent characters in the block, two of which -- $E_6[\theta^i];\phi_{1,0}$ -- are non-real.

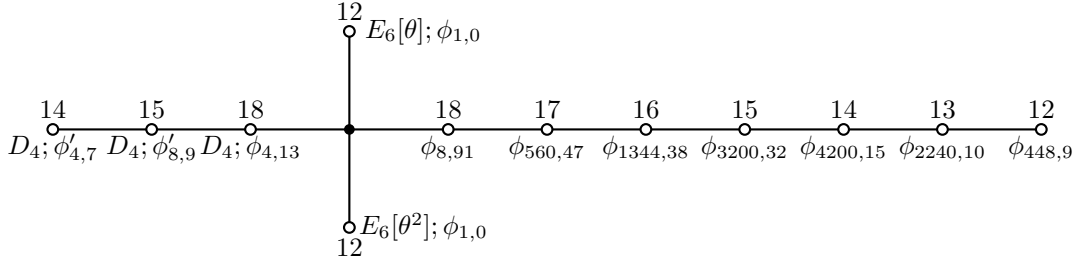
\begin{center}\begin{tabular}{lcccccc}
\hline Character & $A(-)$ & $\omega_i q^{aA/e}$ & $\kappa=1$ & $\kappa=5$ & $\kappa=7$ & $\kappa=11$
\\\hline $\phi_{448,9}$ & $72$ & $q^{6}$ & $12$ & $60$ & $84$ & $132$
\\ $\phi_{2240,10}$ & $81$ & $q^{7}$ & $13$ & $67$ & $95$ & $149$
\\ $\phi_{4200,15}$ & $88$ & $q^{8}$ & $14$ & $74$ & $102$ & $162$
\\ $\phi_{3200,22}$ & $94$ & $q^{9}$ & $15$ & $79$ & $109$ & $173$
\\ $E_6[\theta];\phi_{1,0}$ & $72$ & $\theta q^{6}$ & $12$ & $60$ & $84$ & $132$
\\ $\phi_{1344,38}$ & $102$ & $q^{11}$ & $16$ & $84$ & $120$ & $188$
\\ $\phi_{560,47}$ & $104$ & $q^{12}$ & $17$ & $87$ & $121$ & $191$
\\ $D_4;\phi_{4,7}'$ & $81$ & $-q^{7}$ & $14$ & $68$ & $94$ & $148$
\\ $E_6[\theta^2];\phi_{1,0}$ & $72$ & $\theta^2 q^{6}$ & $12$ & $60$ & $84$ & $132$
\\ $D_4;\phi_{8,9}'$ & $94$ & $-q^{9}$ & $15$ & $79$ & $109$ & $173$
\\ $\phi_{8,91}$ & $108$ & $q^{16}$ & $18$ & $90$ & $126$ & $198$
\\ $D_4;\phi_{4,13}$ & $102$ & $-q^{11}$ & $18$ & $86$ & $118$ & $186$
\\\hline\end{tabular}\end{center}
\begin{center}\begin{tikzpicture}[thick,scale=1.3]
\draw (0,-0.2) node{$D_4;\phi_{4,7}'$};
\draw (1,-0.2) node{$D_4;\phi_{8,9}'$};
\draw (2,-0.2) node{$D_4;\phi_{4,13}$};
\draw (4,-0.2) node{$\phi_{8,91}$};
\draw (5,-0.2) node{$\phi_{560,47}$};
\draw (6,-0.2) node{$\phi_{1344,38}$};
\draw (7,-0.2) node{$\phi_{3200,32}$};
\draw (8,-0.2) node{$\phi_{4200,15}$};
\draw (9,-0.2) node{$\phi_{2240,10}$};
\draw (10,-0.2) node{$\phi_{448,9}$};

\draw (3.75,1) node {$E_6[\theta];\phi_{1,0}$};
\draw (3.75,-1) node {$E_6[\theta^2];\phi_{1,0}$};

\draw (0,0.2) node{$14$};
\draw (1,0.2) node{$15$};
\draw (2,0.2) node{$18$};

\draw (3,-1.2) node{$12$};
\draw (3,1.2) node{$12$};

\draw (4,0.2) node{$18$};
\draw (5,0.2) node{$17$};
\draw (6,0.2) node{$16$};
\draw (7,0.2) node{$15$};
\draw (8,0.2) node{$14$};
\draw (9,0.2) node{$13$};
\draw (10,0.2) node{$12$};

\draw (0,0) -- (10,0);
\draw (3,-1) -- (3,1);

\draw (10,0) node [draw] (l0) {};
\draw (9,0) node [draw] (l0) {};
\draw (8,0) node [draw] (l0) {};
\draw (7,0) node [draw] (l1) {};
\draw (6,0) node [draw] (l2) {};
\draw (5,0) node [draw] (l2) {};
\draw (4,0) node [draw] (l2) {};
\draw (3,1) node [draw] (l2) {};
\draw (3,-1) node [draw] (l2) {};
\draw (3,0) node [fill=black!100] (ld) {};
\draw (2,0) node [draw] (l2) {};
\draw (1,0) node [draw] (l4) {};
\draw (0,0) node [draw] (l4) {};
\end{tikzpicture}\end{center}

\noindent(iv)\;\;\textbf{Block 4:} Cuspidal pair is $(\Phi_{12}.{}^3\!D_4(q),{}^3\!D_4[1])$, of degree $q^3\Phi_1^2\Phi_{12}/2$. There are twelve unipotent characters in the block, two of which -- $E_8[-\theta^i]$ -- are non-real.

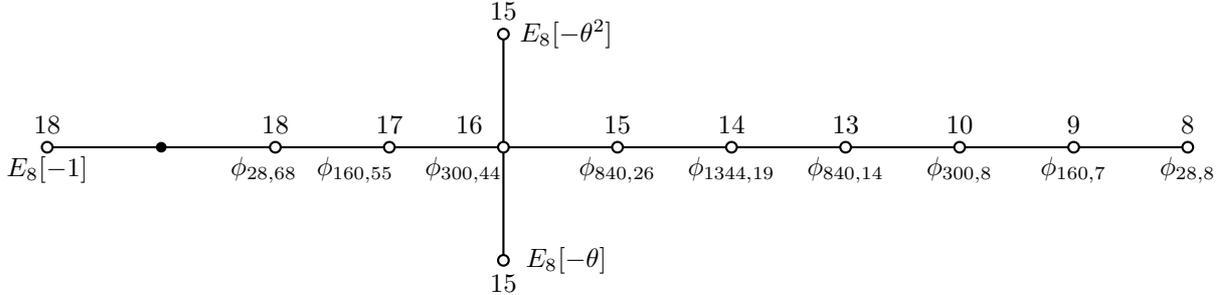
\begin{center}\begin{tabular}{lcccccc}
\hline Character & $A(-)$ & $\omega_i q^{aA/e}$ & $\kappa=1$ & $\kappa=5$ & $\kappa=7$ & $\kappa=11$
\\\hline $\phi_{28,8}$ & $48$ & $q^{4}$ & $8$ & $40$ & $56$ & $88$
\\ $\phi_{160,7}$ & $59$ & $q^{5}$ & $9$ & $49$ & $69$ & $109$
\\ $\phi_{300,8}$ & $69$ & $q^{6}$ & $10$ & $58$ & $80$ & $128$
\\ $E_8[-\theta]$ & $95$ & $-\theta q^{9}$ & $15$ & $79$ & $111$ & $175$
\\ $\phi_{840,14}$ & $87$ & $q^{8}$ & $13$ & $73$ & $101$ & $161$
\\ $\phi_{1344,19}$ & $95$ & $q^{9}$ & $14$ & $78$ & $112$ & $176$
\\ $\phi_{840,26}$ & $99$ & $q^{10}$ & $15$ & $83$ & $115$ & $183$
\\ $E_8[-\theta^2]$ & $95$ & $-\theta^2 q^{9}$ & $15$ & $79$ & $111$ & $175$
\\ $\phi_{300,44}$ & $105$ & $q^{12}$ & $16$ & $88$ & $122$ & $194$
\\ $\phi_{160,55}$ & $107$ & $q^{13}$ & $17$ & $89$ & $125$ & $197$
\\ $\phi_{28,68}$ & $108$ & $q^{14}$ & $18$ & $90$ & $126$ & $198$
\\ $E_8[-1]$ & $95$ & $-q^{9}$ & $18$ & $80$ & $110$ & $172$
\\\hline\end{tabular}\end{center}

\begin{center}\begin{tikzpicture}[thick,scale=1.5]
\draw (0,-0.2) node{$E_8[-1]$};
\draw (1.9,-0.2) node{$\phi_{28,68}$};
\draw (2.7,-0.2) node{$\phi_{160,55}$};
\draw (3.65,-0.2) node{$\phi_{300,44}$};
\draw (5,-0.2) node{$\phi_{840,26}$};
\draw (6,-0.2) node{$\phi_{1344,19}$};
\draw (7,-0.2) node{$\phi_{840,14}$};
\draw (8,-0.2) node{$\phi_{300,8}$};
\draw (9,-0.2) node{$\phi_{160,7}$};
\draw (10,-0.2) node{$\phi_{28,8}$};

\draw (4.55,-1) node {$E_8[-\theta]$};
\draw (4.55,1) node {$E_8[-\theta^2]$};

\draw (0,0.2) node{$18$};

\draw (4,-1.2) node{$15$};
\draw (4,1.2) node{$15$};

\draw (2,0.2) node{$18$};
\draw (3,0.2) node{$17$};
\draw (3.7,0.2) node{$16$};
\draw (5,0.2) node{$15$};
\draw (6,0.2) node{$14$};
\draw (7,0.2) node{$13$};
\draw (8,0.2) node{$10$};
\draw (9,0.2) node{$9$};
\draw (10,0.2) node{$8$};

\draw (0,0) -- (10,0);
\draw (4,-1) -- (4,1);

\draw (10,0) node [draw] (l0) {};
\draw (9,0) node [draw] (l0) {};
\draw (8,0) node [draw] (l0) {};
\draw (7,0) node [draw] (l1) {};
\draw (6,0) node [draw] (l2) {};
\draw (5,0) node [draw] (l2) {};
\draw (4,0) node [draw] (l2) {};
\draw (4,1) node [draw] (l2) {};
\draw (4,-1) node [draw] (l2) {};
\draw (1,0) node [fill=black!100] (ld) {};
\draw (2,0) node [draw] (l2) {};
\draw (3,0) node [draw] (l4) {};
\draw (0,0) node [draw] (l4) {};
\end{tikzpicture}\end{center}

\noindent\textbf{Proof of Brauer trees}: We use a Morita argument for the first three blocks. The fourth block's Brauer tree is given in \cite{cdr2012un}.

Each of the blocks for $d=12$ and $E_7(q)$, under Harish-Chandra induction, induces a Morita equivalence between it and two of the three blocks of $E_8(q)$, which will prove the shape of the tree (and planar embedding) for blocks 1, 2 and 3.

The principal unipotent block for $E_7(q)$, under Harish-Chandra induction, maps onto blocks 1 and 2, and yields a Morita equivalence. (It also has image in the other blocks, but this is not irreducible.) We collate the first and second blocks' terms in square brackets, and the remainder is in round brackets.

\begin{align*}
   \phi_{1,0}&\mapsto [\phi_{8,1}]+[\phi_{84,4}]+(\phi_{1,0}+\phi_{35,2}+\phi_{112,3})
\\ \phi_{56,3}&\mapsto [\phi_{560,5}]+[\phi_{700,6}]+(\phi_{210,4}+\phi_{567,6}+\phi_{1400,8}+\phi_{1575,10}+\phi_{2268,10}+\phi_{112,3}+\phi_{400,7}
\\ &\qquad+\phi_{1008,9}+\phi_{1400,7}+\phi_{3240,9})
\\ \phi_{210,6}&\mapsto [\phi_{1344,8}]+[\phi_{4200,12}]+(\phi_{525,12}+\phi_{567,6}+\phi_{1050,10}+\phi_{1400,8}+\phi_{1575,10}+\phi_{2100,16}+\phi_{2268,10}
\\ &\qquad+\phi_{4096,12}+\phi_{6075,14}+\phi_{1008,9}+\phi_{1296,13}+\phi_{1400,7}+\phi_{2400,17}+\phi_{2800,13}+\phi_{3240,9}
\\ &\qquad+\phi_{3360,13}+\phi_{4096,11}+\phi_{5600,15})
\\ \phi_{336,11}&\mapsto [\phi_{3200,16}]+[\phi_{7168,17}]+(\phi_{1680,22}+\phi_{2100,20}+\phi_{2688,20}+\phi_{2100,16}+\phi_{2100,28}+\phi_{4536,18}
\\ &\qquad+\phi_{5670,18}+\phi_{4096,12}+\phi_{6075,14}+\phi_{6075,22}+\phi_{1296,13}+\phi_{2400,17}+\phi_{2400,23}+\phi_{2800,13}
\\ &\qquad+\phi_{5600,19}+\phi_{3360,13}+\phi_{4096,11}+\phi_{5600,15}+\phi_{5600,21})
\\ \phi_{280,18}&\mapsto [\phi_{4200,21}]+[\phi_{4200,24}]+(\phi_{1400,20}+\phi_{1680,22}+\phi_{1575,34}+\phi_{2100,28}+\phi_{2268,30}+\phi_{4536,18}
\\ &\qquad+\phi_{5670,18}+\phi_{4096,26}+\phi_{6075,22}+\phi_{1008,39}+\phi_{1296,33}+\phi_{2400,23}+\phi_{2800,25}+\phi_{5600,19}
\\ &\qquad+\phi_{3240,31}+\phi_{3360,25}+\phi_{4096,27}+\phi_{5600,21})
\\ \phi_{120,25}&\mapsto [\phi_{2240,28}]+[\phi_{700,42}]+(\phi_{210,52}+\phi_{567,46}+\phi_{1344,38}+\phi_{1400,32}+\phi_{1575,34}+\phi_{2268,30}+\phi_{4096,26}
\\ &\qquad+\phi_{560,47}+\phi_{1008,39}+\phi_{1296,33}+\phi_{1400,37}+\phi_{2800,25}+\phi_{3240,31}+\phi_{4096,27})
\\ \phi_{21,36}&\mapsto [\phi_{448,39}]+[\phi_{84,64}]+(\phi_{525,36}+\phi_{567,46}+\phi_{1344,38}+\phi_{112,63}+\phi_{560,47}+\phi_{1400,37})
\\ D_4;1&\mapsto [D_4;\phi_{4,1}]+[D_4;\phi_{2,4}']+(D_4;\phi_{1,0}+D_4;\phi_{8,3}'+D_4;\phi_{9,2})
\\ D_4;\sigma_2&\mapsto [D_4;\phi_{8,3}'']+[D_4;\phi_{4,8}]+(D_4;\phi_{16,5}+D_4;\phi_{2,4}''+D_4;\phi_{9,2}+D_4;\phi_{9,6}'')
\\ D_4;\ep_2&\mapsto [D_4;\phi_{4,7}'']+[D_4;\phi_{2,16}'']+(D_4;\phi_{1,12}''+D_4;\phi_{8,9}''+D_4;\phi_{9,6}'')
\\ E_6[\theta];\ep&\mapsto [E_6[\theta];\phi_{1,6}]+[E_6[\theta];\phi_{2,2}]+(E_6[\theta];\phi_{2,1}+E_6[\theta];\phi_{1,3}'')
\\ E_6[\theta^2];\ep&\mapsto [E_6[\theta^2];\phi_{1,6}]+[E_6[\theta^2];\phi_{2,2}]+(E_6[\theta^2];\phi_{2,1}+E_6[\theta^2];\phi_{1,3}'')
\end{align*}

The non-principal unipotent block for $E_7(q)$, under Harish-Chandra induction, maps onto blocks 2 and 3, and yields a Morita equivalence. (It also has image in the other blocks, but this is not irreducible.) We collate the second and third blocks' terms in square brackets, and the remainder is in round brackets.
\begin{align*}
\phi_{21,3}&\mapsto [\phi_{84,4}]+[\phi_{448,9}]+(\phi_{525,12}+\phi_{567,6}+\phi_{1344,8}+\phi_{112,3}+\phi_{560,5}+\phi_{1400,7})
\\ \phi_{120,4}&\mapsto [\phi_{700,6}]+[\phi_{2240,10}]+(\phi_{210,4}+\phi_{567,6}+\phi_{1344,8}+\phi_{1400,8}+\phi_{1575,10}+\phi_{2268,10}+\phi_{4096,12}
\\ &\qquad+\phi_{560,5}+\phi_{1008,9}+\phi_{1296,13}+\phi_{1400,7}+\phi_{2800,13}+\phi_{3240,9}+\phi_{4096,11})
\\ \phi_{280,9}&\mapsto [\phi_{4200,12}]+[\phi_{4200,15}]+(\phi_{1400,20}+\phi_{1680,22}+\phi_{1575,10}+\phi_{2100,16}+\phi_{2268,10}+\phi_{4536,18}
\\ &\qquad+\phi_{5670,18}+\phi_{4096,12}+\phi_{6075,14}+\phi_{1008,9}+\phi_{1296,13}+\phi_{2400,17}+\phi_{2800,13}+\phi_{5600,19}
\\ &\qquad+\phi_{3240,9}+\phi_{3360,13}+\phi_{4096,11}+\phi_{5600,15})
\\ \phi_{336,14}&\mapsto [\phi_{7168,17}]+[\phi_{3200,32}]+(\phi_{1680,22}+\phi_{2100,20}+\phi_{2688,20}+\phi_{2100,16}+\phi_{2100,28}+\phi_{4536,18}
\\ &\qquad+\phi_{5670,18}+\phi_{4096,26}+\phi_{6075,14}+\phi_{6075,22}+\phi_{1296,33}+\phi_{2400,17}+\phi_{2400,23}+\phi_{2800,25}
\\ &\qquad+\phi_{5600,19}+\phi_{3360,25}+\phi_{4096,27}+\phi_{5600,15}+\phi_{5600,21})
\\ \phi_{210,21}&\mapsto [\phi_{4200,24}]+[\phi_{1344,38}]+(\phi_{525,36}+\phi_{567,46}+\phi_{1050,34}+\phi_{1400,32}+\phi_{1575,34}+\phi_{2100,28}+\phi_{2268,30}
\\ &\qquad+\phi_{4096,26}+\phi_{6075,22}+\phi_{1008,39}+\phi_{1296,33}+\phi_{1400,37}+\phi_{2400,23}+\phi_{2800,25}+\phi_{3240,31}
\\ &\qquad+\phi_{3360,25}+\phi_{4096,27}+\phi_{5600,21})
\\ \phi_{56,30}&\mapsto [\phi_{700,42}]+[\phi_{560,47}]+(\phi_{210,52}+\phi_{567,46}+\phi_{1400,32}+\phi_{1575,34}+\phi_{2268,30}+\phi_{112,63}
\\ &\qquad+\phi_{400,43}+\phi_{1008,39}+\phi_{1400,37}+\phi_{3240,31})
\\ \phi_{1,63}&\mapsto [\phi_{84,64}]+[\phi_{8,91}]+(\phi_{35,74}+\phi_{112,63}+\phi_{1,120})
\\ D_4;\ep_1&\mapsto [D_4;\phi_{2,4}']+[D_4;\phi_{4,7}']+(D_4;\phi_{1,12}'+D_4;\phi_{8,3}'+D_4;\phi_{9,6}')
\\ D_4;\sigma_2'&\mapsto [D_4;\phi_{4,8}]+[D_4;\phi_{8,9}']+(D_4;\phi_{16,5}+D_4;\phi_{2,16}'+D_4;\phi_{9,10}+D_4;\phi_{9,6}')
\\ D_4;\ep&\mapsto [D_4;\phi_{2,16}'']+[D_4;\phi_{4,13}]+(D_4;\phi_{1,24}+D_4;\phi_{8,9}''+D_4;\phi_{9,10})
\\ E_6[\theta];1&\mapsto [E_6[\theta];\phi_{2,2}]+[E_6[\theta];\phi_{1,0}]+(E_6[\theta];\phi_{1,3}'+E_6[\theta];\phi_{2,1})
\\ E_6[\theta^2];1&\mapsto [E_6[\theta^2];\phi_{2,2}]+[E_6[\theta^2];\phi_{1,0}]+(E_6[\theta^2];\phi_{1,3}'+E_6[\theta^2];\phi_{2,1})
\end{align*}
This completes the proof of the Brauer trees for the first three blocks.

\newpage
\subsection{$d=14$}

For $E_8(q)$ and $d=14$ there are two unipotent blocks of weight $1$, together with 138 unipotent blocks of defect zero.
\\[1cm]
\noindent(i)\;\;\textbf{Block 1:} Cuspidal pair is $(\Phi_{14}.A_1(q),\phi_{2})$, of degree $1$. There are fourteen unipotent characters in the block, two of which -- $E_7[\pm \I];1$ -- are non-real.

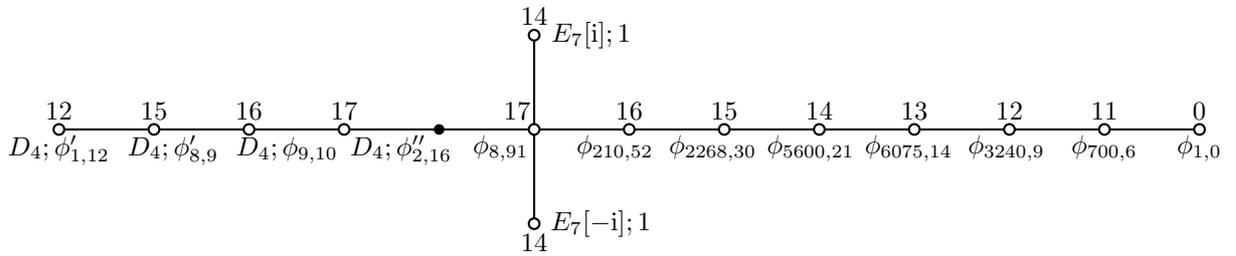
\begin{center}\begin{tabular}{lcccccccc}
\hline Character & $A(-)$ & $\omega_i q^{aA/e}$ & $\kappa=1$ & $\kappa=3$ & $\kappa=5$ & $\kappa=9$ & $\kappa=11$ & $\kappa=13$
\\\hline $\phi_{1,0}$ & $0$ & $q^{0}$ & $0$ & $0$ & $0$ & $0$ & $0$ & $0$
\\ $\phi_{8,91}$ & $119$ & $q^{15}$ & $17$ & $51$ & $85$ & $153$ & $187$ & $221$
\\ $D_4;\phi_{8,9}'$ & $105$ & $-q^{9}$ & $15$ & $45$ & $75$ & $135$ & $165$ & $195$
\\ $D_4;\phi_{9,10}$ & $110$ & $-q^{10}$ & $16$ & $48$ & $78$ & $142$ & $172$ & $204$
\\ $E_7[-\I];1$ & $94$ & $-\I q^{15/2}$ & $14$ & $40$ & $68$ & $120$ & $148$ & $174$
\\ $D_4;\phi_{2,16}''$ & $116$ & $-q^{12}$ & $17$ & $49$ & $83$ & $149$ & $183$ & $215$
\\ $\phi_{700,6}$ & $78$ & $q^{6}$ & $11$ & $33$ & $55$ & $101$ & $123$ & $145$
\\ $\phi_{3240,9}$ & $89$ & $q^{7}$ & $12$ & $38$ & $64$ & $114$ & $140$ & $166$
\\ $\phi_{6075,14}$ & $98$ & $q^{8}$ & $13$ & $43$ & $71$ & $125$ & $153$ & $183$
\\ $\phi_{5600,21}$ & $105$ & $q^{9}$ & $14$ & $44$ & $74$ & $136$ & $166$ & $196$
\\ $\phi_{2268,30}$ & $110$ & $q^{10}$ & $15$ & $47$ & $79$ & $141$ & $173$ & $205$
\\ $E_7[\I];1$ & $94$ & $\I q^{15/2}$ & $14$ & $40$ & $68$ & $120$ & $148$ & $174$
\\ $\phi_{210,52}$ & $116$ & $q^{12}$ & $16$ & $50$ & $82$ & $150$ & $182$ & $216$
\\ $D_4;\phi_{1,12}'$ & $78$ & $-q^{6}$ & $12$ & $34$ & $56$ & $100$ & $122$ & $144$
\\\hline\end{tabular}\end{center}
\begin{center}\begin{tikzpicture}[thick,scale=1.25]
\draw (0,-0.22) node{${D_4;\phi_{1,12}'}$};
\draw (1.2,-0.22) node{$D_4;\phi_{8,9}'$};
\draw (2.4,-0.22) node{$D_4;\phi_{9,10}$};
\draw (3.6,-0.22) node{$D_4;\phi_{2,16}''$};
\draw (4.65,-0.2) node{$\phi_{8,91}$};
\draw (5.85,-0.2) node{$\phi_{210,52}$};
\draw (6.88,-0.2) node{$\phi_{2268,30}$};
\draw (7.91,-0.2) node{$\phi_{5600,21}$};
\draw (8.94,-0.2) node{$\phi_{6075,14}$};
\draw (9.97,-0.2) node{$\phi_{3240,9}$};
\draw (11,-0.2) node{$\phi_{700,6}$};
\draw (12,-0.2) node{$\phi_{1,0}$};

\draw (5.6,1) node {$E_7[\I];1$};
\draw (5.7,-1) node {$E_7[-\I];1$};

\draw (0,0.2) node{$12$};
\draw (1,0.2) node{$15$};
\draw (2,0.2) node{$16$};
\draw (3,0.2) node{$17$};
\draw (4.82,0.2) node{$17$};

\draw (5,-1.2) node{$14$};
\draw (5,1.2) node{$14$};

\draw (6,0.2) node{$16$};
\draw (7,0.2) node{$15$};
\draw (8,0.2) node{$14$};
\draw (9,0.2) node{$13$};
\draw (10,0.2) node{$12$};
\draw (11,0.2) node{$11$};
\draw (12,0.2) node{$0$};

\draw (0,0) -- (12,0);
\draw (5,-1) -- (5,1);

\draw (12,0) node [draw] (l0) {};
\draw (11,0) node [draw] (l0) {};
\draw (10,0) node [draw] (l0) {};
\draw (9,0) node [draw] (l0) {};
\draw (8,0) node [draw] (l0) {};
\draw (7,0) node [draw] (l1) {};
\draw (6,0) node [draw] (l2) {};
\draw (5,0) node [draw] (l2) {};
\draw (5,1) node [draw] (l2) {};
\draw (5,-1) node [draw] (l2) {};
\draw (4,0) node [fill=black!100] (ld) {};
\draw (3,0) node [draw] (l2) {};
\draw (2,0) node [draw] (l2) {};
\draw (1,0) node [draw] (l4) {};
\draw (0,0) node [draw] (l4) {};
\end{tikzpicture}\end{center}

\noindent(ii)\;\; \textbf{Block 2:} Cuspidal pair is $(\Phi_{14}.A_1(q),\phi_{11})$, of degree $q$. There are fourteen unipotent characters in the block, two of which -- $E_7[\pm \I];\ep$ -- are non-real.

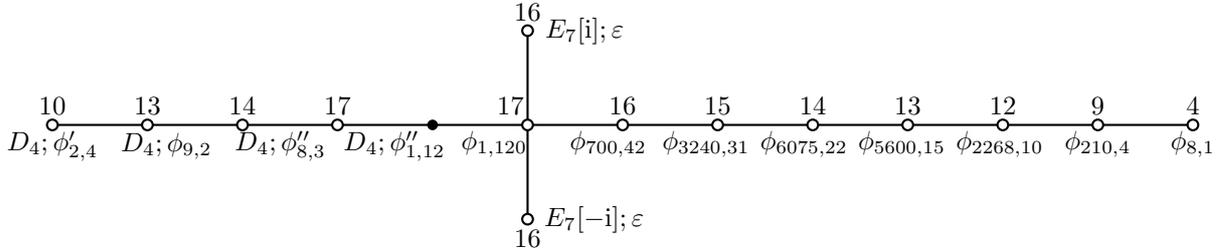
\begin{center}\begin{tabular}{lcccccccc}
\hline Character & $A(-)$ & $\omega_i q^{aA/e}$ & $\kappa=1$ & $\kappa=3$ & $\kappa=5$ & $\kappa=9$ & $\kappa=11$ & $\kappa=13$
\\\hline $\phi_{8,1}$ & $28$ & $q^{2}$ & $4$ & $12$ & $20$ & $36$ & $44$ & $52$
\\ $\phi_{1,120}$ & $119$ & $q^{17}$ & $17$ & $51$ & $85$ & $153$ & $187$ & $221$
\\ $D_4;\phi_{1,12}''$ & $113$ & $-q^{11}$ & $17$ & $49$ & $81$ & $145$ & $177$ & $209$
\\ $\phi_{210,4}$ & $67$ & $q^{5}$ & $9$ & $29$ & $47$ & $87$ & $105$ & $125$
\\ $E_7[-\I];\ep$ & $108$ & $-\I q^{19/2}$ & $16$ & $46$ & $78$ & $138$ & $170$ & $200$
\\ $\phi_{2268,10}$ & $89$ & $q^{7}$ & $12$ & $38$ & $64$ & $114$ & $140$ & $166$
\\ $\phi_{5600,15}$ & $98$ & $q^{8}$ & $13$ & $41$ & $69$ & $127$ & $155$ & $183$
\\ $\phi_{6075,22}$ & $105$ & $q^{9}$ & $14$ & $46$ & $76$ & $134$ & $164$ & $196$
\\ $\phi_{3240,31}$ & $110$ & $q^{10}$ & $15$ & $47$ & $79$ & $141$ & $173$ & $205$
\\ $\phi_{700,42}$ & $113$ & $q^{11}$ & $16$ & $48$ & $80$ & $146$ & $178$ & $210$
\\ $D_4;\phi_{2,4}'$ & $67$ & $-q^{5}$ & $10$ & $28$ & $48$ & $86$ & $106$ & $124$
\\ $E_7[\I];\ep$ & $108$ & $\I q^{19/2}$ & $16$ & $46$ & $78$ & $138$ & $170$ & $200$
\\ $D_4;\phi_{9,2}$ & $89$ & $-q^{7}$ & $13$ & $39$ & $63$ & $115$ & $139$ & $165$
\\ $D_4;\phi_{8,3}''$ & $98$ & $-q^{8}$ & $14$ & $42$ & $70$ & $126$ & $154$ & $182$
\\\hline\end{tabular}\end{center}
\begin{center}\begin{tikzpicture}[thick,scale=1.25]
\draw (0,-0.22) node{${D_4;\phi_{2,4}'}$};
\draw (1.2,-0.22) node{$D_4;\phi_{9,2}$};
\draw (2.4,-0.22) node{$D_4;\phi_{8,3}''$};
\draw (3.6,-0.22) node{$D_4;\phi_{1,12}''$};
\draw (4.65,-0.2) node{$\phi_{1,120}$};
\draw (5.85,-0.2) node{$\phi_{700,42}$};
\draw (6.88,-0.2) node{$\phi_{3240,31}$};
\draw (7.91,-0.2) node{$\phi_{6075,22}$};
\draw (8.94,-0.2) node{$\phi_{5600,15}$};
\draw (9.97,-0.2) node{$\phi_{2268,10}$};
\draw (11,-0.2) node{$\phi_{210,4}$};
\draw (12,-0.2) node{$\phi_{8,1}$};

\draw (5.6,1) node {$E_7[\I];\ep$};
\draw (5.7,-1) node {$E_7[-\I];\ep$};

\draw (0,0.2) node{$10$};
\draw (1,0.2) node{$13$};
\draw (2,0.2) node{$14$};
\draw (3,0.2) node{$17$};
\draw (4.82,0.2) node{$17$};

\draw (5,-1.2) node{$16$};
\draw (5,1.2) node{$16$};

\draw (6,0.2) node{$16$};
\draw (7,0.2) node{$15$};
\draw (8,0.2) node{$14$};
\draw (9,0.2) node{$13$};
\draw (10,0.2) node{$12$};
\draw (11,0.2) node{$9$};
\draw (12,0.2) node{$4$};

\draw (0,0) -- (12,0);
\draw (5,-1) -- (5,1);

\draw (12,0) node [draw] (l0) {};
\draw (11,0) node [draw] (l0) {};
\draw (10,0) node [draw] (l0) {};
\draw (9,0) node [draw] (l0) {};
\draw (8,0) node [draw] (l0) {};
\draw (7,0) node [draw] (l1) {};
\draw (6,0) node [draw] (l2) {};
\draw (5,0) node [draw] (l2) {};
\draw (5,1) node [draw] (l2) {};
\draw (5,-1) node [draw] (l2) {};
\draw (4,0) node [fill=black!100] (ld) {};
\draw (3,0) node [draw] (l2) {};
\draw (2,0) node [draw] (l2) {};
\draw (1,0) node [draw] (l4) {};
\draw (0,0) node [draw] (l4) {};
\end{tikzpicture}\end{center}

\noindent\textbf{Proof of Brauer trees}: We prove that the two unipotent blocks here are Morita equivalent to the principal block of $E_7(q)$ and $d=14$, via a Morita argument.

We take Harish-Chandra induction from the $E_7(q)$-Levi subgroup, showing that the Harish-Chandra induction of each character, cut by each block, is irreducible.

\begin{align*}
\phi_{1,0}&\mapsto [\phi_{1,0}]+[\phi_{8,1}]+(\phi_{35,2}+\phi_{84,4}+\phi_{112,3})
\\\phi_{27,2}&\mapsto [\phi_{700,6}]+[\phi_{210,4}]+(\phi_{35,2}+\phi_{84,4}+\phi_{300,8}+\phi_{567,6}+\phi_{1344,8}+\phi_{112,3}+\phi_{160,7}
\\ &\qquad+\phi_{560,5}+\phi_{1008,9}+\phi_{1400,7})
\\\phi_{105,5}&\mapsto [\phi_{3240,9}]+[\phi_{2268,10}]+(\phi_{300,8}+\phi_{350,14}+\phi_{567,6}+\phi_{1344,8}+\phi_{1575,10}+\phi_{2100,16}+\phi_{4096,12}
\\ &\qquad+\phi_{160,7}+\phi_{560,5}+\phi_{840,13}+\phi_{1008,9}+\phi_{1296,13}+\phi_{1400,7}+\phi_{4096,11})
\\\phi_{189,10}&\mapsto [\phi_{6075,14}]+[\phi_{5600,15}]+(\phi_{350,14}+\phi_{1134,20}+\phi_{1680,22}+\phi_{1575,10}+\phi_{2100,16}+\phi_{5670,18}+\phi_{4096,12}
\\ &\qquad+\phi_{448,25}+\phi_{840,13}+\phi_{1296,13}+\phi_{2400,17}+\phi_{2400,23}+\phi_{5600,19}+\phi_{4096,11})
\\\phi_{189,17}&\mapsto [\phi_{5600,21}]+[\phi_{6075,22}]+ (\phi_{350,38}+\phi_{1134,20}+\phi_{1680,22}+\phi_{1575,34}+\phi_{2100,28}+\phi_{5670,18}+\phi_{4096,26}
\\ &\qquad+\phi_{448,25}+\phi_{840,31}+\phi_{1296,33}+\phi_{2400,17}+\phi_{2400,23}+\phi_{5600,19}+\phi_{4096,27})
\\\phi_{105,26}&\mapsto [\phi_{2268,30}]+[\phi_{3240,31}]+(\phi_{300,8}+\phi_{350,14}+\phi_{567,6}+\phi_{1344,8}+\phi_{1575,10}+\phi_{2100,16}+\phi_{4096,12}
\\ &\qquad+\phi_{160,7}+\phi_{560,5}+\phi_{840,13}+\phi_{1008,9}+\phi_{1296,13}+\phi_{1400,7}+\phi_{4096,11})
\\\phi_{27,37}&\mapsto [\phi_{210,52}]+[\phi_{700,42}]+(\phi_{35,74}+\phi_{84,64}+\phi_{300,44}+\phi_{567,46}+\phi_{1344,38}+\phi_{112,63}+\phi_{160,55}
\\ &\qquad+\phi_{560,47}+\phi_{1008,39}+\phi_{1400,37})
\\\phi_{1,63}&\mapsto [\phi_{8,91}]+[\phi_{1,120}]+(\phi_{35,74}+\phi_{84,64}+\phi_{112,63})
\\ D_4;\ep_1&\mapsto [D_4;\phi_{1,12}']+[D_4;\phi_{2,4}']+(D_4;\phi_{4,7}'+D_4;\phi_{8,3}'+D_4;\phi_{9,6}')
\\ D_4;r\ep_1&\mapsto [D_4;\phi_{8,9}']+[D_4;\phi_{9,2}]+(D_4;\phi_{12,4}+D_4;\phi_{16,5}+D_4;\phi_{4,7}'+D_4;\phi_{6,6}'+D_4;\phi_{8,3}'+D_4;\phi_{9,6}')
\\ D_4;r\ep_2&\mapsto [D_4;\phi_{9,10}]+[D_4;\phi_{8,3}'']+(D_4;\phi_{12,4}+D_4;\phi_{16,5}+D_4;\phi_{4,7}''+D_4;\phi_{6,6}'+D_4;\phi_{8,9}''+D_4;\phi_{9,6}'')
\\ D_4;\ep_2&\mapsto [D_4;\phi_{2,16}'']+[D_4;\phi_{1,12}'']+(D_4;\phi_{4,7}''+D_4;\phi_{8,9}''+D_4;\phi_{9,6}'')
\\ E_7[\I]&\mapsto [E_7[\I];1]+[E_7[\I];\ep]
\\ E_7[-\I]&\mapsto [E_7[-\I];1]+[E_7[-\I];\ep]
\end{align*}
This completes the proof of the Brauer trees, assuming the tree for $E_7(q)$ and $d=14$, which is given in \cite{cdr2012un}.

\newpage

\subsection{$d=15$}

For $E_8(q)$ and $d=15$ there is a single unipotent block of weight $1$, together with 136 unipotent blocks of defect zero.
\\[1cm]
\noindent(i)\;\;\textbf{Block 1:} Cuspidal pair is $(\Phi_{15},1)$, of degree $1$. There are thirty characters in the block, twelve of which are non-real.

\begin{center}\begin{tabular}{lcccccccccc}
\hline Character & $A(-)$ & $\omega_i q^{aA/e}$ & $\kappa=1$ & $\kappa=2$ & $\kappa=4$ & $\kappa=7$ & $\kappa=8$ & $\kappa=11$ & $\kappa=13$ & $\kappa=14$
\\\hline $\phi_{1,0}$ & $0$ & $q^{0}$ & $0$ & $0$ & $0$ & $0$ & $0$ & $0$ & $0$ & $0$
\\ $E_6[\theta^2];\phi_{1,0}$ & $83$ & $-\theta q^{3}$ & $11$ & $23$ & $45$ & $77$ & $89$ & $121$ & $143$ & $155$
\\ $E_8[\zeta^3]$ & $104$ & $\zeta^4 q^{4}$ & $14$ & $28$ & $56$ & $96$ & $112$ & $152$ & $180$ & $194$
\\ $E_6[\theta^2];\phi_{2,2}$ & $104$ & $-\theta q^{4}$ & $14$ & $28$ & $56$ & $96$ & $112$ & $152$ & $180$ & $194$
\\ $\phi_{84,4}$ & $57$ & $q^{2}$ & $7$ & $15$ & $31$ & $55$ & $59$ & $83$ & $99$ & $107$
\\ $E_6[\theta^2];\phi_{1,6}$ & $113$ & $-\theta q^{5}$ & $15$ & $31$ & $61$ & $105$ & $121$ & $165$ & $195$ & $211$
\\ $\phi_{1344,8}$ & $83$ & $q^{3}$ & $10$ & $22$ & $44$ & $78$ & $88$ & $122$ & $144$ & $156$
\\ $\phi_{4096,11}$ & $94$ & $q^{7/2}$ & $11$ & $25$ & $51$ & $87$ & $101$ & $137$ & $163$ & $177$
\\ $\phi_{5670,18}$ & $104$ & $q^{4}$ & $12$ & $28$ & $52$ & $100$ & $108$ & $156$ & $180$ & $196$
\\ $\phi_{4096,27}$ & $109$ & $q^{9/2}$ & $13$ & $29$ & $59$ & $101$ & $117$ & $159$ & $189$ & $205$
\\ $\phi_{1344,38}$ & $113$ & $q^{5}$ & $14$ & $30$ & $60$ & $106$ & $120$ & $166$ & $196$ & $212$
\\ $E_6[\theta];\phi_{1,0}$ & $83$ & $-\theta q^{3}$ & $11$ & $23$ & $45$ & $77$ & $89$ & $121$ & $143$ & $155$
\\ $\phi_{84,64}$ & $117$ & $q^{6}$ & $15$ & $31$ & $63$ & $111$ & $123$ & $171$ & $203$ & $219$
\\ $E_6[\theta];\phi_{2,2}$ & $104$ & $-\theta q^{4}$ & $14$ & $28$ & $56$ & $96$ & $112$ & $152$ & $180$ & $194$
\\ $E_8[\zeta^2]$ & $104$ & $\zeta q^{4}$ & $14$ & $28$ & $56$ & $96$ & $112$ & $152$ & $180$ & $194$
\\ $E_6[\theta];\phi_{1,6}$ & $113$ & $-\theta q^{5}$ & $15$ & $31$ & $61$ & $105$ & $121$ & $165$ & $195$ & $211$
\\ $\phi_{1,120}$ & $120$ & $q^{8}$ & $16$ & $32$ & $64$ & $112$ & $128$ & $176$ & $208$ & $224$
\\ $\phi_{8,1}$ & $29$ & $-q$ & $4$ & $8$ & $16$ & $28$ & $30$ & $42$ & $50$ & $54$
\\ $E_8[\theta^2]$ & $104$ & $\theta q^{4}$ & $15$ & $29$ & $57$ & $97$ & $111$ & $151$ & $179$ & $193$
\\ $\phi_{112,3}$ & $57$ & $-q^{2}$ & $7$ & $15$ & $29$ & $53$ & $61$ & $85$ & $99$ & $107$
\\ $E_8[\zeta^4]$ & $104$ & $\zeta^2 q^{4}$ & $14$ & $28$ & $56$ & $96$ & $112$ & $152$ & $180$ & $194$
\\ $\phi_{1400,7}$ & $83$ & $-q^{3}$ & $10$ & $22$ & $44$ & $80$ & $86$ & $122$ & $144$ & $156$
\\ $\phi_{4096,12}$ & $94$ & $-q^{7/2}$ & $11$ & $25$ & $51$ & $87$ & $101$ & $137$ & $163$ & $177$
\\ $\phi_{5600,19}$ & $104$ & $-q^{4}$ & $12$ & $26$ & $56$ & $98$ & $110$ & $152$ & $182$ & $196$
\\ $\phi_{4096,26}$ & $109$ & $-q^{9/2}$ & $13$ & $29$ & $59$ & $101$ & $117$ & $159$ & $189$ & $205$
\\ $\phi_{1400,37}$ & $113$ & $-q^{5}$ & $14$ & $30$ & $60$ & $108$ & $118$ & $166$ & $196$ & $212$
\\ $E_8[\zeta]$ & $104$ & $\zeta^3 q^{4}$ & $14$ & $28$ & $56$ & $96$ & $112$ & $152$ & $180$ & $194$
\\ $\phi_{112,63}$ & $117$ & $-q^{6}$ & $15$ & $31$ & $61$ & $109$ & $125$ & $173$ & $203$ & $219$
\\ $E_8[\theta]$ & $104$ & $\theta^2 q^{4}$ & $15$ & $29$ & $57$ & $97$ & $111$ & $151$ & $179$ & $193$
\\ $\phi_{8,91}$ & $119$ & $-q^{7}$ & $16$ & $32$ & $64$ & $112$ & $126$ & $174$ & $206$ & $222$
\\\hline\end{tabular}\end{center}
\textbf{WARNING:} this tree is conjectural.
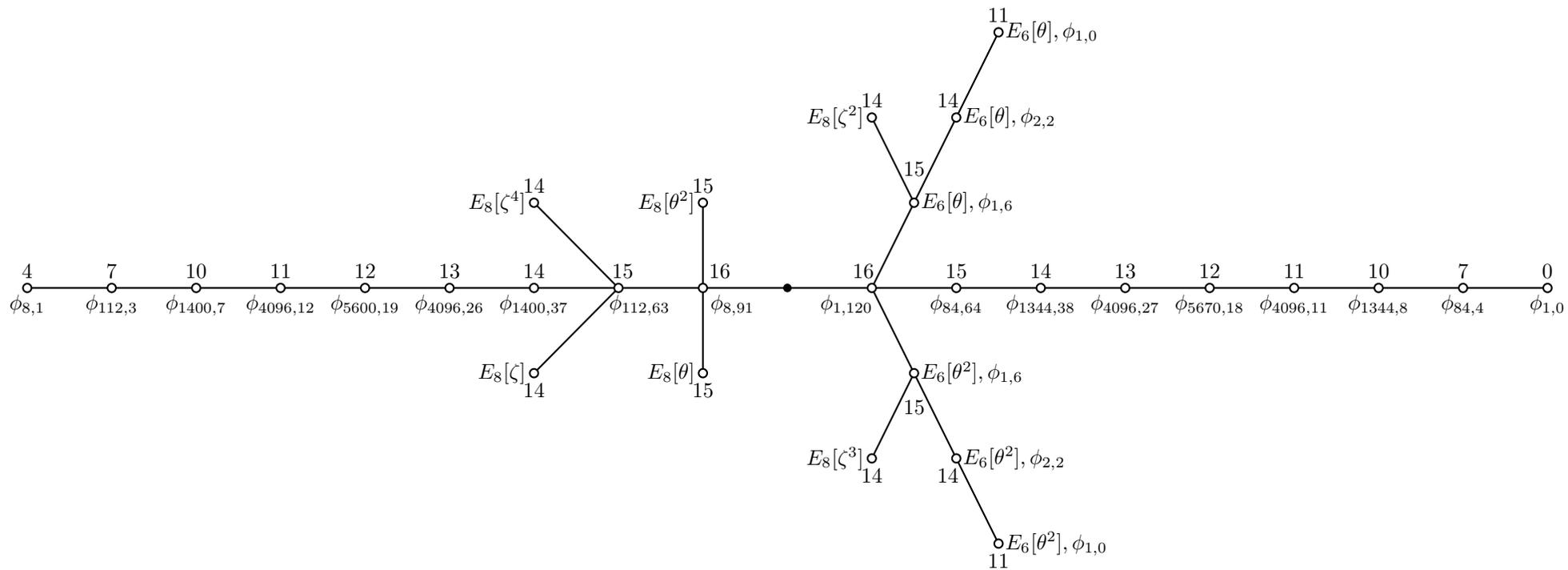
\begin{landscape}
\begin{center}\begin{tikzpicture}[thick,scale=1.4]

\draw (0,0.8) node {$\phi_{1,0}$};
\draw (-1,0.8) node {$\phi_{84,4}$};
\draw (-2,0.8) node {$\phi_{1344,8}$};
\draw (-3,0.8) node {$\phi_{4096,11}$};
\draw (-4,0.8) node {$\phi_{5670,18}$};
\draw (-5,0.8) node {$\phi_{4096,27}$};
\draw (-6,0.8) node {$\phi_{1344,38}$};
\draw (-7,0.8) node {$\phi_{84,64}$};
\draw (-8.3,0.8) node {$\phi_{1,120}$};

\draw (-9.65,0.8) node {$\phi_{8,91}$};
\draw (-10.75,0.8) node {$\phi_{112,63}$};
\draw (-12,0.8) node {$\phi_{1400,37}$};
\draw (-13,0.8) node {$\phi_{4096,26}$};
\draw (-14,0.8) node {$\phi_{5600,19}$};
\draw (-15,0.8) node {$\phi_{4096,12}$};
\draw (-16,0.8) node {$\phi_{1400,7}$};
\draw (-17,0.8) node {$\phi_{112,3}$};
\draw (-18,0.8) node {$\phi_{8,1}$};

\draw (0,1.2) node{$0$};
\draw (-1,1.2) node{$7$};
\draw (-2,1.2) node{$10$};
\draw (-3,1.2) node{$11$};
\draw (-4,1.2) node{$12$};
\draw (-5,1.2) node{$13$};
\draw (-6,1.2) node{$14$};
\draw (-7,1.2) node{$15$};
\draw (-8.1,1.2) node{$16$};

\draw (-9.8,1.2) node{$16$};
\draw (-10.95,1.2) node{$15$};
\draw (-12,1.2) node{$14$};
\draw (-13,1.2) node{$13$};
\draw (-14,1.2) node{$12$};
\draw (-15,1.2) node{$11$};
\draw (-16,1.2) node{$10$};
\draw (-17,1.2) node{$7$};
\draw (-18,1.2) node{$4$};

\draw (-12,-0.2) node{$14$};
\draw (-12,2.2) node{$14$};

\draw (-10,-0.2) node{$15$};
\draw (-10,2.2) node{$15$};

\draw (-8,-1.2) node{$14$};
\draw (-8,3.2) node{$14$};

\draw (-7.5,-0.4) node{$15$};
\draw (-7.5,2.4) node{$15$};
\draw (-7.1,-1.2) node{$14$};
\draw (-7.1,3.2) node{$14$};
\draw (-6.5,-2.2) node{$11$};
\draw (-6.5,4.2) node{$11$};

\draw (0,1) -- (-18,1);
\draw (-10,0) -- (-10,2);

\draw (-11,1) -- (-12,2);
\draw (-11,1) -- (-12,0);

\draw (-8,1) -- (-6.5,4);
\draw (-8,1) -- (-6.5,-2);

\draw (-7.5,2) -- (-8,3);
\draw (-7.5,0) -- (-8,-1);

\draw (-0,1) node [draw] (l0) {};
\draw (-1,1) node [draw] (l1) {};
\draw (-2,1) node [draw] (l1) {};
\draw (-3,1) node [draw] (l1) {};
\draw (-4,1) node [draw] (l1) {};
\draw (-5,1) node [draw] (l1) {};
\draw (-6,1) node [draw] (l1) {};
\draw (-7,1) node [draw] (l1) {};
\draw (-8,1) node [draw] (l1) {};

\draw (-9,1) node [fill=black!100] (ld) {};

\draw (-10,1) node [draw] (l1) {};
\draw (-11,1) node [draw] (l1) {};
\draw (-12,1) node [draw] (l1) {};
\draw (-13,1) node [draw] (l1) {};
\draw (-14,1) node [draw] (l1) {};
\draw (-15,1) node [draw] (l1) {};
\draw (-16,1) node [draw] (l1) {};
\draw (-17,1) node [draw] (l1) {};
\draw (-18,1) node [draw] (l1) {};

\draw (-7.5,2) node [draw,label=right:${E_6[\theta],\phi_{1,6}}$] (l4) {};
\draw (-7,3) node [draw,label=right:${E_6[\theta],\phi_{2,2}}$] (l4) {};
\draw (-6.5,4) node [draw,label=right:${E_6[\theta],\phi_{1,0}}$] (l4) {};

\draw (-7.5,0) node [draw,label=right:${E_6[\theta^2],\phi_{1,6}}$] (l4) {};
\draw (-7,-1) node [draw,label=right:${E_6[\theta^2],\phi_{2,2}}$] (l4) {};
\draw (-6.5,-2) node [draw,label=right:${E_6[\theta^2],\phi_{1,0}}$] (l4) {};

\draw (-12,2) node [draw,label=left:${E_8[\zeta^4]}$] (l4) {};
\draw (-12,0) node [draw,label=left:${E_8[\zeta]}$] (l4) {};
\draw (-10,0) node [draw,label=left:${E_8[\theta]}$] (l4) {};
\draw (-10,2) node [draw,label=left:${E_8[\theta^2]}$] (l4) {};

\draw (-8,-1) node [draw,label=left:${E_8[\zeta^3]}$] (l4) {};
\draw (-8,3) node [draw,label=left:${E_8[\zeta^2]}$] (l4) {};

\end{tikzpicture}\end{center}
\end{landscape}
\newpage
\subsection{$d=18$}

For $E_8(q)$ and $d=18$ there are three unipotent blocks of weight $1$, together with 112 unipotent blocks of defect zero.
\\[1cm]
\noindent(i)\;\;\textbf{Block 1:} Cuspidal pair is $(\Phi_{18}.{}^2\!A_2(q),\phi_3)$, of degree $1$. There are eighteen characters in the block, six of which are non-real.

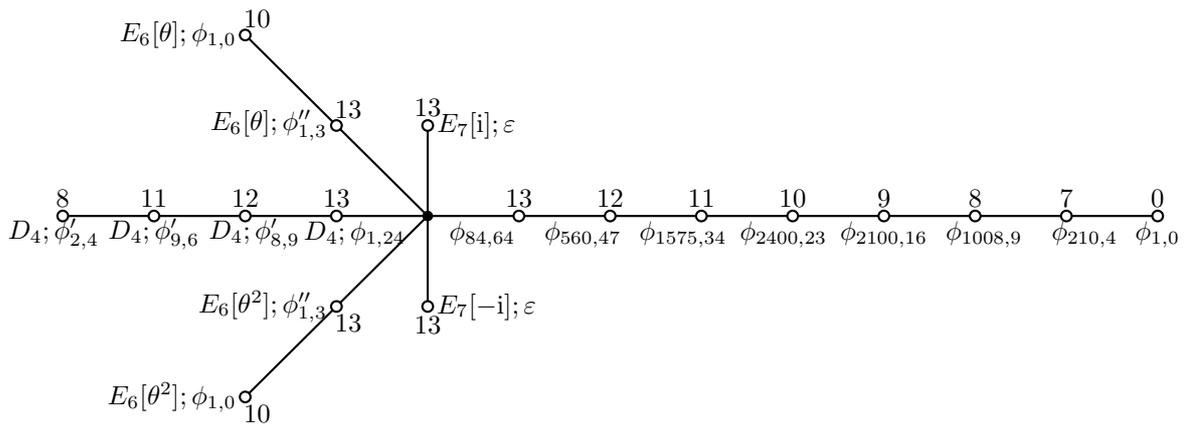
\begin{center}\begin{tabular}{lcccccccc}
\hline Character & $A(-)$ & $\omega_i q^{aA/e}$ & $\kappa=1$ & $\kappa=5$ & $\kappa=7$ & $\kappa=11$ & $\kappa=13$ & $\kappa=17$
\\\hline $\phi_{1,0}$ & $0$ & $q^{0}$ & $0$ & $0$ & $0$ & $0$ & $0$ & $0$
\\ $D_4,\phi_{1,24}$ & $117$ & $-q^{10}$ & $13$ & $65$ & $91$ & $143$ & $169$ & $221$
\\ $E_6[\theta^2];\phi_{1,3}''$ & $112$ & $\theta^2 q^{8}$ & $13$ & $63$ & $87$ & $137$ & $161$ & $211$
\\ $E_7[-\I];\ep$ & $109$ & $-\I q^{15/2}$ & $13$ & $61$ & $85$ & $133$ & $157$ & $205$
\\ $\phi_{210,4}$ & $68$ & $q^{4}$ & $7$ & $37$ & $53$ & $83$ & $99$ & $129$
\\ $\phi_{1008,9}$ & $83$ & $q^{5}$ & $8$ & $46$ & $66$ & $100$ & $120$ & $158$
\\ $\phi_{2100,16}$ & $95$ & $q^{6}$ & $9$ & $53$ & $73$ & $117$ & $137$ & $181$
\\ $\phi_{2400,23}$ & $105$ & $q^{7}$ & $10$ & $60$ & $82$ & $128$ & $150$ & $200$
\\ $\phi_{1575,34}$ & $112$ & $q^{8}$ & $11$ & $61$ & $87$ & $137$ & $163$ & $213$
\\ $\phi_{560,47}$ & $115$ & $q^{9}$ & $12$ & $64$ & $90$ & $140$ & $166$ & $218$
\\ $\phi_{84,64}$ & $117$ & $q^{10}$ & $13$ & $65$ & $91$ & $143$ & $169$ & $221$
\\ $E_6[\theta];\phi_{1,0}$ & $83$ & $\theta q^{5}$ & $10$ & $46$ & $64$ & $102$ & $120$ & $156$
\\ $E_7[\I];\ep$ & $109$ & $\I q^{15/2}$ & $13$ & $61$ & $85$ & $133$ & $157$ & $205$
\\ $D_4,\phi_{2,4}'$ & $68$ & $-q^{4}$ & $8$ & $38$ & $52$ & $84$ & $98$ & $128$
\\ $E_6[\theta];\phi_{1,3}''$ & $112$ & $\theta q^{8}$ & $13$ & $63$ & $87$ & $137$ & $161$ & $211$
\\ $D_4,\phi_{9,6}'$ & $95$ & $-q^{6}$ & $11$ & $51$ & $75$ & $115$ & $139$ & $179$
\\ $D_4,\phi_{8,9}'$ & $105$ & $-q^{7}$ & $12$ & $58$ & $82$ & $128$ & $152$ & $198$
\\ $E_6[\theta^2];\phi_{1,0}$ & $83$ & $\theta^2 q^{5}$ & $10$ & $46$ & $64$ & $102$ & $120$ & $156$
\\\hline\end{tabular}\end{center}

\begin{center}\begin{tikzpicture}[thick,scale=1.2]

\draw (12,0.8) node {$\phi_{1,0}$};
\draw (11.2,0.8) node {$\phi_{210,4}$};
\draw (10.1,0.8) node {$\phi_{1008,9}$};
\draw (9,0.8) node {$\phi_{2100,16}$};
\draw (7.9,0.8) node {$\phi_{2400,23}$};
\draw (6.8,0.8) node {$\phi_{1575,34}$};
\draw (5.7,0.8) node {$\phi_{560,47}$};
\draw (4.6,0.8) node {$\phi_{84,64}$};

\draw (3.2,0.8) node {$D_4;\phi_{1,24}$};
\draw (2.1,0.8) node {$D_4;\phi_{8,9}'$};
\draw (1,0.8) node {$D_4;\phi_{9,6}'$};
\draw (-0.1,0.8) node {$D_4;\phi_{2,4}'$};

\draw (0,1.2) node {$8$};
\draw (1,1.2) node {$11$};
\draw (2,1.2) node {$12$};
\draw (3,1.2) node {$13$};

\draw (5,1.2) node {$13$};
\draw (6,1.2) node {$12$};
\draw (7,1.2) node {$11$};
\draw (8,1.2) node {$10$};
\draw (9,1.2) node {$9$};
\draw (10,1.2) node {$8$};
\draw (11,1.2) node {$7$};
\draw (12,1.2) node {$0$};

\draw (4,2.2) node {$13$};
\draw (4,-0.2) node {$13$};

\draw (3.13,2.18) node {$13$};
\draw (3.13,-0.18) node {$13$};
\draw (2.13,3.18) node {$10$};
\draw (2.13,-1.18) node {$10$};

\draw (0,1) -- (12,1);
\draw (2,-1) -- (4,1);
\draw (2,3) -- (4,1);
\draw (4,0) -- (4,2);

\draw (0,1) node [draw] (l1) {};
\draw (1,1) node [draw] (l0) {};
\draw (2,1) node [draw] (l1) {};
\draw (3,1) node [draw] (l1) {};
\draw (4,1) node [fill=black!100] (ld) {};
\draw (5,1) node [draw] (l1) {};
\draw (6,1) node [draw] (l1) {};
\draw (7,1) node [draw] (l1) {};
\draw (8,1) node [draw] (l1) {};
\draw (9,1) node [draw] (l1) {};
\draw (10,1) node [draw] (l1) {};
\draw (11,1) node [draw] (l1) {};
\draw (12,1) node [draw] (l1) {};

\draw (4,2) node [draw,label=right:${E_7[\I];\ep}$] (l1) {};
\draw (4,0) node [draw,label=right:${E_7[-\I];\ep}$] (l1) {};

\draw (3,2) node [draw,label=left:${E_6[\theta];\phi_{1,3}''}$] (l1) {};
\draw (3,0) node [draw,label=left:${E_6[\theta^2];\phi_{1,3}''}$] (l1) {};
\draw (2,3) node [draw,label=left:${E_6[\theta];\phi_{1,0}}$] (l1) {};
\draw (2,-1) node [draw,label=left:${E_6[\theta^2];\phi_{1,0}}$] (l1) {};
\end{tikzpicture}\end{center}

\noindent(ii)\;\;\textbf{Block 2:} Cuspidal pair is $(\Phi_{18}.{}^2\!A_2(q),\phi_{111})$, of degree $q^3$. There are eighteen characters in the block, four of which are non-real.

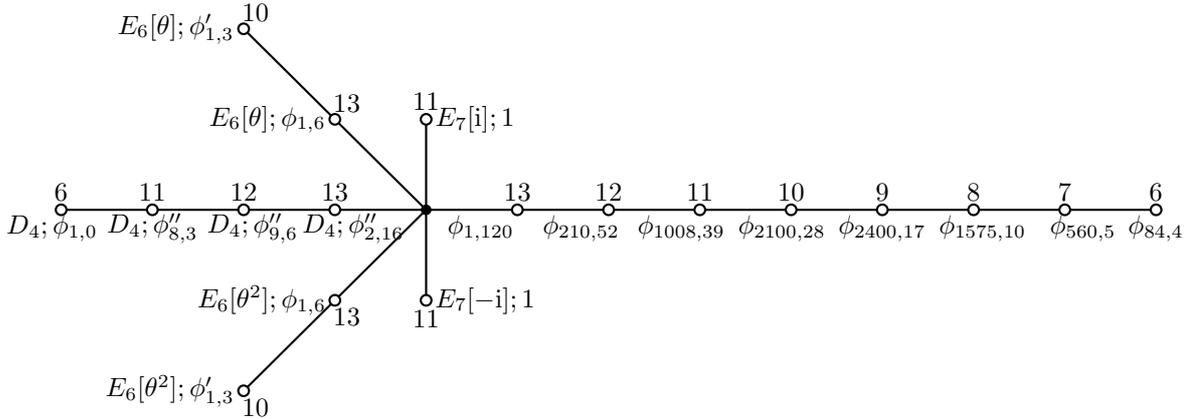
\begin{center}\begin{tabular}{lcccccccc}
\hline Character & $A(-)$ & $\omega_i q^{aA/e}$ & $\kappa=1$ & $\kappa=5$ & $\kappa=7$ & $\kappa=11$ & $\kappa=13$ & $\kappa=17$
\\\hline $\phi_{84,4}$ & $54$ & $q^{3}$ & $6$ & $30$ & $42$ & $66$ & $78$ & $102$
\\ $\phi_{560,5}$ & $70$ & $q^{4}$ & $7$ & $39$ & $55$ & $85$ & $101$ & $133$
\\ $\phi_{1575,10}$ & $85$ & $q^{5}$ & $8$ & $46$ & $66$ & $104$ & $124$ & $162$
\\ $\phi_{2400,17}$ & $96$ & $q^{6}$ & $9$ & $55$ & $75$ & $117$ & $137$ & $183$
\\ $\phi_{2100,28}$ & $104$ & $q^{7}$ & $10$ & $58$ & $80$ & $128$ & $150$ & $198$
\\ $\phi_{1008,39}$ & $110$ & $q^{8}$ & $11$ & $61$ & $87$ & $133$ & $159$ & $209$
\\ $\phi_{210,52}$ & $113$ & $q^{9}$ & $12$ & $62$ & $88$ & $138$ & $164$ & $214$
\\ $E_7[\I];1$ & $91$ & $\I q^{11/2}$ & $11$ & $51$ & $71$ & $111$ & $131$ & $171$
\\ $E_6[\theta];\phi_{1,3}'$ & $85$ & $\theta q^{5}$ & $10$ & $48$ & $66$ & $104$ & $122$ & $160$
\\ $D_4,\phi_{1,0}$ & $54$ & $-q^{3}$ & $6$ & $30$ & $42$ & $66$ & $78$ & $102$
\\ $\phi_{1,120}$ & $120$ & $q^{13}$ & $13$ & $65$ & $91$ & $143$ & $169$ & $221$
\\ $E_6[\theta];\phi_{1,6}$ & $110$ & $\theta q^{8}$ & $13$ & $61$ & $85$ & $135$ & $159$ & $207$
\\ $D_4,\phi_{8,3}''$ & $96$ & $-q^{6}$ & $11$ & $53$ & $75$ & $117$ & $139$ & $181$
\\ $D_4,\phi_{9,6}''$ & $104$ & $-q^{7}$ & $12$ & $56$ & $82$ & $126$ & $152$ & $196$
\\ $E_6[\theta^2];\phi_{1,3}'$ & $85$ & $\theta^2 q^{5}$ & $10$ & $48$ & $66$ & $104$ & $122$ & $160$
\\ $D_4,\phi_{2,16}''$ & $113$ & $-q^{9}$ & $13$ & $63$ & $87$ & $139$ & $163$ & $213$
\\ $E_7[-\I];1$ & $91$ & $-\I q^{11/2}$ & $11$ & $51$ & $71$ & $111$ & $131$ & $171$
\\ $E_6[\theta^2];\phi_{1,6}$ & $110$ & $\theta^2 q^{8}$ & $13$ & $61$ & $85$ & $135$ & $159$ & $207$
\\\hline\end{tabular}\end{center}
\begin{center}\begin{tikzpicture}[thick,scale=1.2]

\draw (12,0.8) node {$\phi_{84,4}$};
\draw (11.2,0.8) node {$\phi_{560,5}$};
\draw (10.1,0.8) node {$\phi_{1575,10}$};
\draw (9,0.8) node {$\phi_{2400,17}$};
\draw (7.9,0.8) node {$\phi_{2100,28}$};
\draw (6.8,0.8) node {$\phi_{1008,39}$};
\draw (5.7,0.8) node {$\phi_{210,52}$};
\draw (4.6,0.8) node {$\phi_{1,120}$};

\draw (3.2,0.8) node {$D_4;\phi_{2,16}''$};
\draw (2.1,0.8) node {$D_4;\phi_{9,6}''$};
\draw (1,0.8) node {$D_4;\phi_{8,3}''$};
\draw (-0.1,0.8) node {$D_4;\phi_{1,0}$};

\draw (0,1.2) node {$6$};
\draw (1,1.2) node {$11$};
\draw (2,1.2) node {$12$};
\draw (3,1.2) node {$13$};

\draw (5,1.2) node {$13$};
\draw (6,1.2) node {$12$};
\draw (7,1.2) node {$11$};
\draw (8,1.2) node {$10$};
\draw (9,1.2) node {$9$};
\draw (10,1.2) node {$8$};
\draw (11,1.2) node {$7$};
\draw (12,1.2) node {$6$};

\draw (4,2.2) node {$11$};
\draw (4,-0.2) node {$11$};

\draw (3.13,2.18) node {$13$};
\draw (3.13,-0.18) node {$13$};
\draw (2.13,3.18) node {$10$};
\draw (2.13,-1.18) node {$10$};

\draw (0,1) -- (12,1);
\draw (2,-1) -- (4,1);
\draw (2,3) -- (4,1);
\draw (4,0) -- (4,2);

\draw (0,1) node [draw] (l1) {};
\draw (1,1) node [draw] (l0) {};
\draw (2,1) node [draw] (l1) {};
\draw (3,1) node [draw] (l1) {};
\draw (4,1) node [fill=black!100] (ld) {};
\draw (5,1) node [draw] (l1) {};
\draw (6,1) node [draw] (l1) {};
\draw (7,1) node [draw] (l1) {};
\draw (8,1) node [draw] (l1) {};
\draw (9,1) node [draw] (l1) {};
\draw (10,1) node [draw] (l1) {};
\draw (11,1) node [draw] (l1) {};
\draw (12,1) node [draw] (l1) {};

\draw (4,2) node [draw,label=right:${E_7[\I];1}$] (l1) {};
\draw (4,0) node [draw,label=right:${E_7[-\I];1}$] (l1) {};

\draw (3,2) node [draw,label=left:${E_6[\theta];\phi_{1,6}}$] (l1) {};
\draw (3,0) node [draw,label=left:${E_6[\theta^2];\phi_{1,6}}$] (l1) {};
\draw (2,3) node [draw,label=left:${E_6[\theta];\phi_{1,3}'}$] (l1) {};
\draw (2,-1) node [draw,label=left:${E_6[\theta^2];\phi_{1,3}'}$] (l1) {};
\end{tikzpicture}\end{center}

\noindent(iii)\;\;\textbf{Block 3:} Cuspidal pair is $(\Phi_{18}.{}^2\!A_2(q),\phi_{21})$, of degree $q\Phi_1$. There are eighteen characters in the block, four of which are non-real.

\begin{center}\begin{tabular}{lcccccccc}
\hline Character & $A(-)$ & $\omega_i q^{aA/e}$ & $\kappa=1$ & $\kappa=5$ & $\kappa=7$ & $\kappa=11$ & $\kappa=13$ & $\kappa=17$
\\\hline $\phi_{8,1}$ & $27$ & $q^{3/2}$ & $3$ & $15$ & $21$ & $33$ & $39$ & $51$
\\ $\phi_{35,2}$ & $44$ & $q^{5/2}$ & $4$ & $24$ & $34$ & $54$ & $64$ & $84$
\\ $E_8[-\theta]$ & $102$ & $-\theta q^{13/2}$ & $11$ & $57$ & $79$ & $125$ & $147$ & $193$
\\ $\phi_{300,8}$ & $76$ & $q^{9/2}$ & $7$ & $43$ & $59$ & $93$ & $109$ & $145$
\\ $\phi_{840,13}$ & $90$ & $q^{11/2}$ & $8$ & $50$ & $70$ & $110$ & $130$ & $172$
\\ $\phi_{1134,20}$ & $102$ & $q^{13/2}$ & $9$ & $55$ & $81$ & $123$ & $149$ & $195$
\\ $\phi_{840,31}$ & $108$ & $q^{15/2}$ & $10$ & $60$ & $84$ & $132$ & $156$ & $206$
\\ $\phi_{300,44}$ & $112$ & $q^{17/2}$ & $11$ & $63$ & $87$ & $137$ & $161$ & $213$
\\ $E_8[-\theta^2]$ & $102$ & $-\theta^2q^{13/2}$ & $11$ & $57$ & $79$ & $125$ & $147$ & $193$
\\ $\phi_{35,74}$ & $116$ & $q^{21/2}$ & $12$ & $64$ & $90$ & $142$ & $168$ & $220$
\\ $\phi_{8,91}$ & $117$ & $q^{23/2}$ & $13$ & $65$ & $91$ & $143$ & $169$ & $221$
\\ $E_8[\theta]$ & $102$ & $\theta q^{13/2}$ & $13$ & $57$ & $79$ & $125$ & $147$ & $191$
\\ $D_4,\phi_{1,12}'$ & $76$ & $-q^{9/2}$ & $9$ & $41$ & $59$ & $93$ & $111$ & $143$
\\ $D_4,\phi_{4,7}'$ & $90$ & $-q^{11/2}$ & $10$ & $50$ & $70$ & $110$ & $130$ & $170$
\\ $D_4,\phi_{6,6}'$ & $102$ & $-q^{13/2}$ & $11$ & $57$ & $79$ & $125$ & $147$ & $193$
\\ $D_4,\phi_{4,7}''$ & $108$ & $-q^{15/2}$ & $12$ & $60$ & $84$ & $132$ & $156$ & $204$
\\ $D_4,\phi_{1,12}''$ & $112$ & $-q^{17/2}$ & $13$ & $61$ & $87$ & $137$ & $163$ & $211$
\\ $E_8[\theta^2]$ & $102$ & $\theta^2 q^{13/2}$ & $13$ & $57$ & $79$ & $125$ & $147$ & $191$
\\\hline\end{tabular}\end{center}
\textbf{WARNING:} this tree is conjectural.
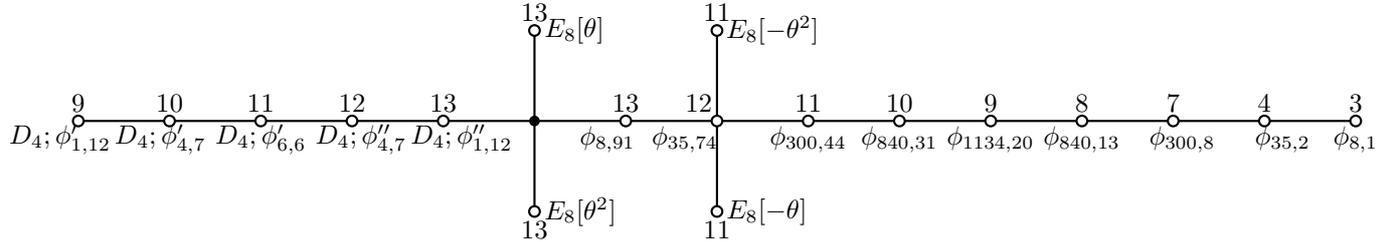
\begin{center}\begin{tikzpicture}[thick,scale=1.2]

\draw (12,0.8) node {$\phi_{8,1}$};
\draw (11.2,0.8) node {$\phi_{35,2}$};
\draw (10.1,0.8) node {$\phi_{300,8}$};
\draw (9,0.8) node {$\phi_{840,13}$};
\draw (8,0.8) node {$\phi_{1134,20}$};
\draw (7,0.8) node {$\phi_{840,31}$};
\draw (6,0.8) node {$\phi_{300,44}$};
\draw (4.65,0.8) node {$\phi_{35,74}$};
\draw (3.8,0.8) node {$\phi_{8,91}$};

\draw (2.2,0.8) node {$D_4;\phi_{1,12}''$};
\draw (1.1,0.8) node {$D_4;\phi_{4,7}''$};
\draw (0,0.8) node {$D_4;\phi_{6,6}'$};
\draw (-1.1,0.8) node {$D_4;\phi_{4,7}'$};
\draw (-2.2,0.8) node {$D_4;\phi_{1,12}'$};

\draw (-2,1.2) node {$9$};
\draw (-1,1.2) node {$10$};
\draw (0,1.2) node {$11$};
\draw (1,1.2) node {$12$};
\draw (2,1.2) node {$13$};

\draw (4,1.2) node {$13$};
\draw (4.8,1.2) node {$12$};
\draw (6,1.2) node {$11$};
\draw (7,1.2) node {$10$};
\draw (8,1.2) node {$9$};
\draw (9,1.2) node {$8$};
\draw (10,1.2) node {$7$};
\draw (11,1.2) node {$4$};
\draw (12,1.2) node {$3$};

\draw (3,2.2) node {$13$};
\draw (3,-0.2) node {$13$};
\draw (5,2.2) node {$11$};
\draw (5,-0.2) node {$11$};

\draw (-2,1) -- (12,1);
\draw (3,0) -- (3,2);
\draw (5,0) -- (5,2);

\draw (-2,1) node [draw] (l1) {};
\draw (-1,1) node [draw] (l0) {};
\draw (0,1) node [draw] (l1) {};
\draw (1,1) node [draw] (l0) {};
\draw (2,1) node [draw] (l1) {};
\draw (4,1) node [draw] (l1) {};
\draw (3,1) node [fill=black!100] (ld) {};
\draw (5,1) node [draw] (l1) {};
\draw (6,1) node [draw] (l1) {};
\draw (7,1) node [draw] (l1) {};
\draw (8,1) node [draw] (l1) {};
\draw (9,1) node [draw] (l1) {};
\draw (10,1) node [draw] (l1) {};
\draw (11,1) node [draw] (l1) {};
\draw (12,1) node [draw] (l1) {};

\draw (3,2) node [draw,label=right:${E_8[\theta]}$] (l1) {};
\draw (3,0) node [draw,label=right:${E_8[\theta^2]}$] (l1) {};

\draw (5,2) node [draw,label=right:${E_8[-\theta^2]}$] (l1) {};
\draw (5,0) node [draw,label=right:${E_8[-\theta]}$] (l1) {};
\end{tikzpicture}\end{center}

\noindent\textbf{Proof of Brauer trees}: We prove that the first two unipotent blocks here are Morita equivalent to the principal block of $E_7(q)$ and $d=18$, via a Morita argument.

We take Harish-Chandra induction from the $E_7(q)$-Levi subgroup, showing that the Harish-Chandra induction of each character, cut by each block, is irreducible. We put square brackets around the contributions to each of the first two blocks, and round brackets around all other characters.

\begin{align*}
\phi_{1,0}&\mapsto [\phi_{1,0}]+[\phi_{84,4}]+(\phi_{35,2}+\phi_{8,1}+\phi_{112,3})
\\ \phi_{7,1}&\mapsto [\phi_{210,4}]+[\phi_{560,5}]+(\phi_{28,8}+\phi_{35,2}+\phi_{567,6}+\phi_{8,1}+\phi_{112,3}+\phi_{160,7})
\\ \phi_{21,6}&\mapsto [\phi_{1008,9}]+[\phi_{1575,10}]+(\phi_{28,8}+\phi_{350,14}+\phi_{567,6}+\phi_{56,19}+\phi_{160,7}+\phi_{1296,13})
\\ \phi_{35,13}&\mapsto [\phi_{2100,16}]+[\phi_{2400,17}]+(\phi_{70,32}+\phi_{350,14}+\phi_{1680,22}+\phi_{56,19}+\phi_{448,25}+\phi_{1296,13})
\\ \phi_{35,22}&\mapsto [\phi_{2400,23}]+[\phi_{2100,28}]+(\phi_{70,32}+\phi_{350,38}+\phi_{1680,22}+\phi_{56,49}+\phi_{448,25}+\phi_{1296,33})
\\ \phi_{21,33}&\mapsto [\phi_{1575,34}]+[\phi_{1008,39}]+(\phi_{28,68}+\phi_{350,38}+\phi_{567,46}+\phi_{56,49}+\phi_{160,55}+\phi_{1296,33})
\\ \phi_{7,46}&\mapsto [\phi_{560,47}]+[\phi_{210,52}]+(\phi_{28,68}+\phi_{35,74}+\phi_{567,46}+\phi_{8,91}+\phi_{112,63}+\phi_{160,55})
\\ \phi_{1,63}&\mapsto [\phi_{84,64}]+[\phi_{1,120}]+(\phi_{8,91}+\phi_{35,74}+\phi_{112,63})
\\ D_4;1&\mapsto [D_4;\phi_{2,4}']+[D_4;\phi_{1,0}]+(D_4;\phi_{4,1}+D_4;\phi_{8,3}'+D_4;\phi_{9,2})
\\ D_4;r&\mapsto [D_4;\phi_{9,6}']+[D_4;\phi_{8,3}'']+(D_4;\phi_{12,4}+D_4;\phi_{16,5}+D_4;\phi_{4,1}+D_4;\phi_{6,6}''+D_4;\phi_{8,3}'+D_4;\phi_{9,2})
\\ D_4;r\ep&\mapsto [D_4;\phi_{8,9}']+[D_4;\phi_{9,6}'']+(D_4;\phi_{12,4}+D_4;\phi_{16,5}+D_4;\phi_{4,13}+D_4;\phi_{6,6}''+D_4;\phi_{8,9}''+D_4;\phi_{9,10})
\\ D_4;\ep&\mapsto [D_4;\phi_{1,24}]+[D_4;\phi_{2,16}'']+(D_4;\phi_{4,13}+D_4;\phi_{8,9}''+D_4;\phi_{9,10})
\\ E_6[\theta];1&\mapsto [E_6[\theta];\phi_{1,0}]+[E_6[\theta];\phi_{1,3}']+(E_6[\theta];\phi_{2,1}+E_6[\theta];\phi_{2,2})
\\ E_6[\theta];\ep&\mapsto [E_6[\theta];\phi_{1,3}'']+[E_6[\theta];\phi_{1,6}]+(E_6[\theta];\phi_{2,1}+E_6[\theta];\phi_{2,2})
\\ E_6[\theta^2];1&\mapsto [E_6[\theta^2];\phi_{1,0}]+[E_6[\theta^2];\phi_{1,3}']+(E_6[\theta^2];\phi_{2,1}+E_6[\theta^2];\phi_{2,2})
\\ E_6[\theta^2];\ep&\mapsto [E_6[\theta^2];\phi_{1,3}'']+[E_6[\theta^2];\phi_{1,6}]+(E_6[\theta^2];\phi_{2,1}+E_6[\theta^2];\phi_{2,2})
\\ E_7[\I]&\mapsto [E_7[\I];\ep]+[E_7[\I];1]
\\ E_7[-\I]&\mapsto [E_7[-\I];\ep]+[E_7[-\I];1]
\end{align*}
This completes the proof of the Brauer trees.

\newpage
\subsection{$d=20$}

For $E_8(q)$ and $d=20$ there is a single unipotent block of weight $1$, together with 146 unipotent blocks of defect zero.
\\[1cm]
\noindent(i)\;\;\textbf{Block 1:} Cuspidal pair is $(\Phi_{20},1)$, of degree $1$. There are twenty characters in the block, six of which -- $E_8[\pm\I]$ and $E_8[\zeta^i]$ -- are non-real.

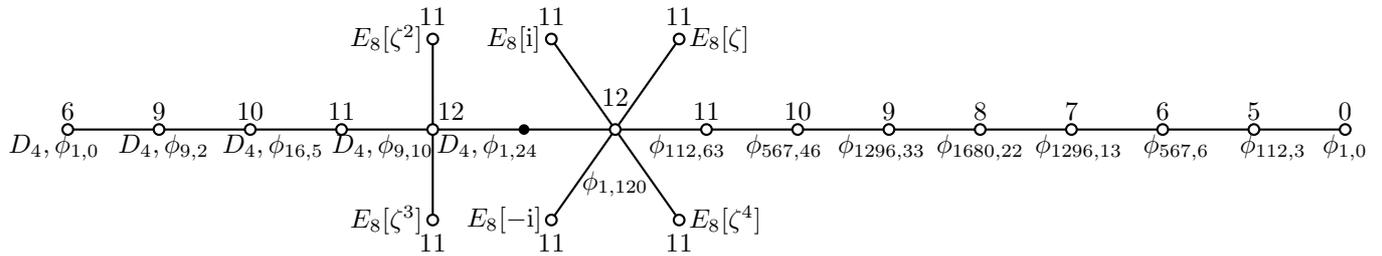
\begin{center}\begin{tabular}{lcccccccccc}
\hline Character & $A(-)$ & $\omega_i q^{aA/e}$ & $\kappa=1$ & $\kappa=3$ & $\kappa=7$ & $\kappa=9$ & $\kappa=11$ & $\kappa=13$ & $\kappa=17$ & $\kappa=19$
\\\hline $\phi_{1,0}$ & $0$ & $q^{0}$ & $0$ & $0$ & $0$ & $0$ & $0$ & $0$ & $0$ & $0$
\\ $E_8[-\I]$ & $104$ & $-\I q^{6}$ & $11$ & $31$ & $73$ & $93$ & $115$ & $135$ & $177$ & $197$
\\ $E_8[\zeta^4]$ & $104$ & $\zeta^4 q^{6}$ & $11$ & $31$ & $73$ & $93$ & $115$ & $135$ & $177$ & $197$
\\ $\phi_{112,3}$ & $57$ & $q^{3}$ & $5$ & $17$ & $39$ & $51$ & $63$ & $75$ & $97$ & $109$
\\ $\phi_{567,6}$ & $74$ & $q^{4}$ & $6$ & $22$ & $52$ & $68$ & $80$ & $96$ & $126$ & $142$
\\ $\phi_{1296,13}$ & $90$ & $q^{5}$ & $7$ & $27$ & $63$ & $81$ & $99$ & $117$ & $153$ & $173$
\\ $\phi_{1680,22}$ & $104$ & $q^{6}$ & $8$ & $32$ & $72$ & $96$ & $112$ & $136$ & $176$ & $200$
\\ $\phi_{1296,33}$ & $110$ & $q^{7}$ & $9$ & $33$ & $77$ & $99$ & $121$ & $143$ & $187$ & $211$
\\ $\phi_{567,46}$ & $114$ & $q^{8}$ & $10$ & $34$ & $80$ & $104$ & $124$ & $148$ & $194$ & $218$
\\ $\phi_{112,63}$ & $117$ & $q^{9}$ & $11$ & $35$ & $81$ & $105$ & $129$ & $153$ & $199$ & $223$
\\ $E_8[\zeta]$ & $104$ & $\zeta q^{6}$ & $11$ & $31$ & $73$ & $93$ & $115$ & $135$ & $177$ & $197$
\\ $E_8[\I]$ & $104$ & $\I q^{6}$ & $11$ & $31$ & $73$ & $93$ & $115$ & $135$ & $177$ & $197$
\\ $\phi_{1,120}$ & $120$ & $q^{12}$ & $12$ & $36$ & $84$ & $108$ & $132$ & $156$ & $204$ & $228$
\\ $D_4;\phi_{1,0}$ & $57$ & $-q^{3}$ & $6$ & $18$ & $40$ & $52$ & $62$ & $74$ & $96$ & $108$
\\ $E_8[\zeta^2]$ & $104$ & $\zeta^2 q^{6}$ & $11$ & $31$ & $73$ & $93$ & $115$ & $135$ & $177$ & $197$
\\ $D_4;\phi_{9,2}$ & $90$ & $-q^{5}$ & $9$ & $27$ & $63$ & $83$ & $97$ & $117$ & $153$ & $171$
\\ $D_4;\phi_{16,5}$ & $104$ & $-q^{6}$ & $10$ & $32$ & $72$ & $94$ & $114$ & $136$ & $176$ & $198$
\\ $D_4;\phi_{9,10}$ & $110$ & $-q^{7}$ & $11$ & $33$ & $77$ & $101$ & $119$ & $143$ & $187$ & $209$
\\ $E_8[\zeta^3]$ & $104$ & $\zeta^3 q^{6}$ & $11$ & $31$ & $73$ & $93$ & $115$ & $135$ & $177$ & $197$
\\ $D_4;\phi_{1,24}$ & $117$ & $-q^{9}$ & $12$ & $36$ & $82$ & $106$ & $128$ & $152$ & $198$ & $222$
\\\hline\end{tabular}\end{center}
\begin{center}\begin{tikzpicture}[thick,scale=1.2]

\draw (0,0.8) node {$\phi_{1,0}$};
\draw (-0.79,0.8) node {$\phi_{112,3}$};
\draw (-1.85,0.8) node {$\phi_{567,6}$};
\draw (-2.92,0.8) node {$\phi_{1296,13}$};
\draw (-4,0.8) node {$\phi_{1680,22}$};
\draw (-5.08,0.8) node {$\phi_{1296,33}$};
\draw (-6.15,0.8) node {$\phi_{567,46}$};
\draw (-7.21,0.8) node {$\phi_{112,63}$};
\draw (-8,0.4) node {$\phi_{1,120}$};

\draw (-9.4,0.8) node {$D_4,\phi_{1,24}$};
\draw (-10.55,0.8) node {$D_4,\phi_{9,10}$};
\draw (-11.75,0.8) node {$D_4,\phi_{16,5}$};
\draw (-12.95,0.8) node {$D_4,\phi_{9,2}$};
\draw (-14.15,0.8) node {$D_4,\phi_{1,0}$};

\draw (0,1) -- (-14,1);
\draw (-10,0) -- (-10,2);

\draw (-8.7,0) -- (-7.3,2);
\draw (-7.3,0) -- (-8.7,2);

\draw (-14,1.2) node {$6$};
\draw (-13,1.2) node {$9$};
\draw (-12,1.2) node {$10$};
\draw (-11,1.2) node {$11$};
\draw (-9.8,1.2) node {$12$};

\draw (0,1.2) node {$0$};
\draw (-1,1.2) node {$5$};
\draw (-2,1.2) node {$6$};
\draw (-3,1.2) node {$7$};
\draw (-4,1.2) node {$8$};
\draw (-5,1.2) node {$9$};
\draw (-6,1.2) node {$10$};
\draw (-7,1.2) node {$11$};
\draw (-8,1.35) node {$12$};

\draw (-10,-0.25) node[rectangle] {$11$};
\draw (-10,2.25) node[rectangle] {$11$};
\draw (-8.7,-0.25) node[rectangle] {$11$};
\draw (-8.7,2.25) node[rectangle] {$11$};
\draw (-7.3,-0.25) node[rectangle] {$11$};
\draw (-7.3,2.25) node[rectangle] {$11$};

\draw (0,1) node [draw] (l0) {};
\draw (-1,1) node [draw] (l1) {};
\draw (-2,1) node [draw] (l1) {};
\draw (-3,1) node [draw] (l1) {};
\draw (-4,1) node [draw] (l1) {};
\draw (-5,1) node [draw] (l1) {};
\draw (-6,1) node [draw] (l1) {};
\draw (-7,1) node [draw] (l1) {};
\draw (-8,1) node [draw] (l1) {};

\draw (-9,1) node [fill=black!100] (ld) {};
\draw (-10,1) node [draw] (l1) {};
\draw (-11,1) node [draw] (l1) {};
\draw (-12,1) node [draw] (l1) {};
\draw (-13,1) node [draw] (l1) {};
\draw (-14,1) node [draw] (l1) {};

\draw (-7.3,0) node [draw,label=right:${E_8[\zeta^4]}$] (l4) {};
\draw (-7.3,2) node [draw,label=right:${E_8[\zeta]}$] (l4) {};
\draw (-8.7,0) node [draw,label=left:${E_8[-\I]}$] (l4) {};
\draw (-8.7,2) node [draw,label=left:${E_8[\I]}$] (l4) {};
\draw (-10,0) node [draw,label=left:${E_8[\zeta^3]}$] (l4) {};
\draw (-10,2) node [draw,label=left:${E_8[\zeta^2]}$] (l4) {};

\end{tikzpicture}\end{center}

\noindent\textbf{Proof of Brauer tree}: This is given in \cite{cdr2012un}.

\newpage
\subsection{$d=24$}

For $E_8(q)$ and $d=24$ there is a single unipotent block of weight $1$, together with 142 unipotent blocks of defect zero.
\\[1cm]
\noindent(i)\;\;\textbf{Block 1:} Cuspidal pair is $(\Phi_{24},1)$, of degree $1$. There are twenty-four characters in the block, ten of which are non-real.

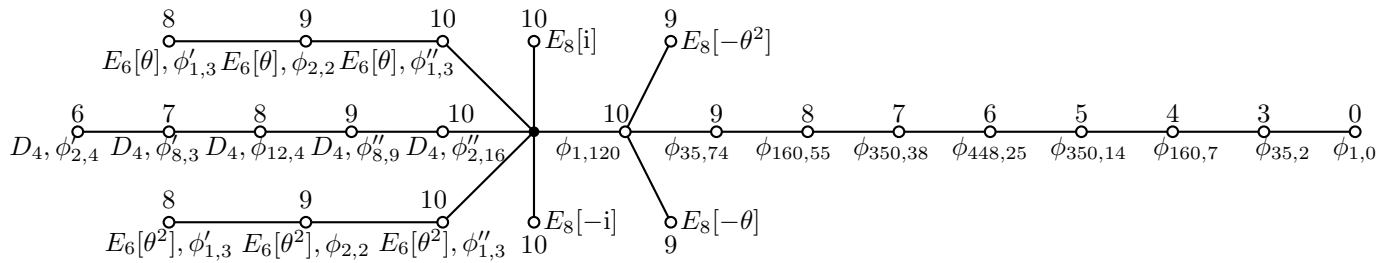
\begin{center}\begin{tabular}{lcccccccccc}
\hline Character & $A(-)$ & $\omega_i q^{aA/e}$ & $\kappa=1$ & $\kappa=5$ & $\kappa=7$ & $\kappa=11$ & $\kappa=13$ & $\kappa=17$ & $\kappa=19$ & $\kappa=23$
\\\hline $\phi_{1,0}$ & $0$ & $q^{0}$ & $0$ & $0$ & $0$ & $0$ & $0$ & $0$ & $0$ & $0$
\\ $E_8[-\theta]$ & $104$ & $-\theta q^{5}$ & $9$ & $43$ & $61$ & $95$ & $113$ & $147$ & $165$ & $199$
\\ $\phi_{35,2}$ & $46$ & $q^{2}$ & $3$ & $19$ & $27$ & $43$ & $49$ & $65$ & $73$ & $89$
\\ $\phi_{160,7}$ & $68$ & $q^{3}$ & $4$ & $28$ & $40$ & $62$ & $74$ & $96$ & $108$ & $132$
\\ $\phi_{350,14}$ & $88$ & $q^{4}$ & $5$ & $37$ & $51$ & $83$ & $93$ & $125$ & $139$ & $171$
\\ $\phi_{448,25}$ & $104$ & $q^{5}$ & $6$ & $44$ & $62$ & $96$ & $112$ & $146$ & $164$ & $202$
\\ $\phi_{350,38}$ & $112$ & $q^{6}$ & $7$ & $47$ & $65$ & $105$ & $119$ & $159$ & $177$ & $217$
\\ $\phi_{160,55}$ & $116$ & $q^{7}$ & $8$ & $48$ & $68$ & $106$ & $126$ & $164$ & $184$ & $224$
\\ $\phi_{35,74}$ & $0$ & $q^{8}$ & $9$ & $49$ & $69$ & $109$ & $127$ & $167$ & $187$ & $227$
\\ $E_8[-\theta^2]$ & $104$ & $-\theta^2 q^{5}$ & $9$ & $43$ & $61$ & $95$ & $113$ & $147$ & $165$ & $199$
\\ $\phi_{1,120}$ & $120$ & $q^{10}$ & $10$ & $50$ & $70$ & $110$ & $130$ & $170$ & $190$ & $230$
\\ $E_8[\I]$ & $104$ & $\I q^{5}$ & $10$ & $44$ & $60$ & $94$ & $114$ & $148$ & $164$ & $198$
\\ $E_6[\theta];\phi_{1,3}'$ & $88$ & $\theta q^{4}$ & $8$ & $36$ & $52$ & $80$ & $96$ & $124$ & $140$ & $168$
\\ $E_6[\theta];\phi_{2,2}$ & $104$ & $\theta q^{5}$ & $9$ & $43$ & $61$ & $95$ & $113$ & $147$ & $165$ & $199$
\\ $E_6[\theta];\phi_{1,3}''$ & $112$ & $\theta q^{6}$ & $10$ & $46$ & $66$ & $102$ & $122$ & $158$ & $178$ & $214$
\\ $D_4,\phi_{2,4}'$ & $68$ & $-q^{3}$ & $6$ & $28$ & $40$ & $64$ & $72$ & $96$ & $108$ & $130$
\\ $D_4,\phi_{8,3}'$ & $88$ & $-q^{4}$ & $7$ & $37$ & $51$ & $81$ & $95$ & $125$ & $139$ & $169$
\\ $D_4,\phi_{12,4}$ & $104$ & $-q^{5}$ & $8$ & $42$ & $60$ & $98$ & $110$ & $148$ & $166$ & $200$
\\ $D_4,\phi_{8,9}''$ & $112$ & $-q^{6}$ & $9$ & $47$ & $65$ & $103$ & $121$ & $159$ & $177$ & $215$
\\ $D_4,\phi_{2,16}''$ & $116$ & $-q^{7}$ & $10$ & $48$ & $68$ & $108$ & $124$ & $164$ & $184$ & $222$
\\ $E_6[\theta^2];\phi_{1,3}'$ & $88$ & $\theta^2 q^{4}$ & $8$ & $36$ & $52$ & $80$ & $96$ & $124$ & $140$ & $168$
\\ $E_6[\theta^2];\phi_{2,2}$ & $104$ & $\theta^2 q^{5}$ & $9$ & $43$ & $61$ & $95$ & $113$ & $147$ & $165$ & $199$
\\ $E_6[\theta^2];\phi_{1,3}''$ & $112$ & $\theta^2 q^{6}$ & $10$ & $46$ & $66$ & $102$ & $122$ & $158$ & $178$ & $214$
\\ $E_8[-\I]$ & $104$ & $-\I q^{5}$ & $10$ & $44$ & $60$ & $94$ & $114$ & $148$ & $164$ & $198$
\\\hline\end{tabular}\end{center}
\begin{center}\begin{tikzpicture}[thick,scale=1.2]

\draw (0,0.8) node {$\phi_{1,0}$};
\draw (-0.79,0.8) node {$\phi_{35,2}$};
\draw (-1.85,0.8) node {$\phi_{160,7}$};
\draw (-2.92,0.8) node {$\phi_{350,14}$};
\draw (-4,0.8) node {$\phi_{448,25}$};
\draw (-5.08,0.8) node {$\phi_{350,38}$};
\draw (-6.15,0.8) node {$\phi_{160,55}$};
\draw (-7.21,0.8) node {$\phi_{35,74}$};
\draw (-8.4,0.8) node {$\phi_{1,120}$};

\draw (-9.85,0.8) node {$D_4,\phi_{2,16}''$};
\draw (-10.95,0.8) node {$D_4,\phi_{8,9}''$};
\draw (-12.05,0.8) node {$D_4,\phi_{12,4}$};
\draw (-13.15,0.8) node {$D_4,\phi_{8,3}'$};
\draw (-14.25,0.8) node {$D_4,\phi_{2,4}'$};

\draw (-10.5,1.75) node {${E_6[\theta],\phi_{1,3}''}$};
\draw (-11.8,1.75) node {${E_6[\theta],\phi_{2,2}}$};
\draw (-13.1,1.75) node {${E_6[\theta],\phi_{1,3}'}$};

\draw (-10,-0.25) node {${E_6[\theta^2],\phi_{1,3}''}$};
\draw (-11.5,-0.25) node {${E_6[\theta^2],\phi_{2,2}}$};
\draw (-13,-0.25) node {${E_6[\theta^2],\phi_{1,3}'}$};

\draw (-14,1.2) node {$6$};
\draw (-13,1.2) node {$7$};
\draw (-12,1.2) node {$8$};
\draw (-11,1.2) node {$9$};
\draw (-9.8,1.2) node {$10$};

\draw (0,1.2) node {$0$};
\draw (-1,1.2) node {$3$};
\draw (-2,1.2) node {$4$};
\draw (-3,1.2) node {$5$};
\draw (-4,1.2) node {$6$};
\draw (-5,1.2) node {$7$};
\draw (-6,1.2) node {$8$};
\draw (-7,1.2) node {$9$};
\draw (-8.1,1.2) node {$10$};

\draw (-9,-0.25) node[rectangle] {$10$};
\draw (-9,2.25) node[rectangle] {$10$};

\draw (-7.5,-0.25) node[rectangle] {$9$};
\draw (-7.5,2.25) node[rectangle] {$9$};

\draw (-10.1,0.25) node[rectangle] {$10$};
\draw (-11.5,0.25) node[rectangle] {$9$};
\draw (-13,0.25) node[rectangle] {$8$};

\draw (-10,2.25) node[rectangle] {$10$};
\draw (-11.5,2.25) node[rectangle] {$9$};
\draw (-13,2.25) node[rectangle] {$8$};

\draw (0,1) -- (-14,1);
\draw (-8,1) -- (-7.5,2);
\draw (-8,1) -- (-7.5,0);
\draw (-9,0) -- (-9,2);
\draw (-9,1) -- (-10,2);
\draw (-10,2) -- (-13,2);
\draw (-9,1) -- (-10,0);
\draw (-10,0) -- (-13,0);

\draw (0,1) node [draw] (l0) {};
\draw (-1,1) node [draw] (l1) {};
\draw (-2,1) node [draw] (l1) {};
\draw (-3,1) node [draw] (l1) {};
\draw (-4,1) node [draw] (l1) {};
\draw (-5,1) node [draw] (l1) {};
\draw (-6,1) node [draw] (l1) {};
\draw (-7,1) node [draw] (l1) {};
\draw (-8,1) node [draw] (l1) {};

\draw (-9,1) node [fill=black!100] (ld) {};
\draw (-10,1) node [draw] (l1) {};
\draw (-11,1) node [draw] (l1) {};
\draw (-12,1) node [draw] (l1) {};
\draw (-13,1) node [draw] (l1) {};
\draw (-14,1) node [draw] (l1) {};

\draw (-7.5,0) node [draw,label=right:${E_8[-\theta]}$] (l4) {};
\draw (-7.5,2) node [draw,label=right:${E_8[-\theta^2]}$] (l4) {};

\draw (-9,0) node [draw,label=right:${E_8[-\I]}$] (l4) {};
\draw (-9,2) node [draw,label=right:${E_8[\I]}$] (l4) {};

\draw (-10,2) node [draw] (l1) {};
\draw (-11.5,2) node [draw] (l1) {};
\draw (-13,2) node [draw] (l1) {};

\draw (-10,0) node [draw] (l1) {};
\draw (-11.5,0) node [draw] (l1) {};
\draw (-13,0) node [draw] (l1) {};

\end{tikzpicture}\end{center}

\noindent\textbf{Proof of Brauer tree}: This is given in \cite{cdr2012un}.

\newpage
\subsection{$d=30$}

For $E_8(q)$ and $d=30$ there is a single unipotent block of weight $1$, together with 136 unipotent blocks of defect zero.
\\[1cm]
\noindent(i)\;\;\textbf{Block 1:} Cuspidal pair is $(\Phi_{30},1)$, of degree $1$. There are thirty characters in the block, sixteen of which are non-real.

\begin{center}\begin{tabular}{lcccccccccc}
\hline Character & $A(-)$ & $\omega_i q^{aA/e}$ & $\kappa=1$ & $\kappa=7$ & $\kappa=11$ & $\kappa=13$ & $\kappa=17$ & $\kappa=19$ & $\kappa=23$ & $\kappa=29$
\\\hline $\phi_{1,0}$ & $0$ & $q^{0}$ & $0$ & $0$ & $0$ & $0$ & $0$ & $0$ & $0$ & $0$
\\ $\phi_{8,1}$ & $29$ & $q$ & $1$ & $13$ & $21$ & $25$ & $33$ & $37$ & $45$ & $57$
\\ $\phi_{28,8}$ & $57$ & $q^{2}$ & $2$ & $26$ & $42$ & $50$ & $64$ & $72$ & $88$ & $112$
\\ $\phi_{56,19}$ & $83$ & $q^{3}$ & $3$ & $39$ & $61$ & $73$ & $93$ & $105$ & $127$ & $163$
\\ $\phi_{70,32}$ & $104$ & $q^{4}$ & $4$ & $52$ & $76$ & $92$ & $116$ & $132$ & $156$ & $204$
\\ $\phi_{56,49}$ & $113$ & $q^{5}$ & $5$ & $53$ & $83$ & $99$ & $127$ & $143$ & $173$ & $221$
\\ $\phi_{28,68}$ & $117$ & $q^{6}$ & $6$ & $54$ & $86$ & $102$ & $132$ & $148$ & $180$ & $228$
\\ $\phi_{8,91}$ & $119$ & $q^{7}$ & $7$ & $55$ & $87$ & $103$ & $135$ & $151$ & $183$ & $231$
\\ $\phi_{1,120}$ & $120$ & $q^{8}$ & $8$ & $56$ & $88$ & $104$ & $136$ & $152$ & $184$ & $232$
\\ $E_8[-\theta^2]$ & $104$ & $-\theta^2 q^{4}$ & $8$ & $48$ & $76$ & $90$ & $118$ & $132$ & $160$ & $200$
\\ $E_8[\zeta]$ & $104$ & $\zeta q^{4}$ & $8$ & $48$ & $76$ & $90$ & $118$ & $132$ & $160$ & $200$
\\ $E_7[\I];1$ & $94$ & $\I q^{7/2}$ & $7$ & $43$ & $69$ & $83$ & $105$ & $119$ & $145$ & $181$
\\ $E_7[\I];\ep$ & $109$ & $\I q^{9/2}$ & $8$ & $50$ & $80$ & $96$ & $122$ & $138$ & $168$ & $210$
\\ $E_6[\theta];\phi_{1,0}$ & $83$ & $\theta q^{3}$ & $6$ & $38$ & $60$ & $72$ & $94$ & $106$ & $128$ & $160$
\\ $E_6[\theta];\phi_{2,1}$ & $104$ & $\theta q^{4}$ & $7$ & $47$ & $75$ & $89$ & $119$ & $133$ & $161$ & $201$
\\ $E_6[\theta];\phi_{1,6}$ & $113$ & $\theta q^{5}$ & $8$ & $52$ & $82$ & $98$ & $128$ & $144$ & $174$ & $218$
\\ $E_8[\zeta^2]$ & $104$ & $\zeta^2 q^{4}$ & $8$ & $48$ & $76$ & $90$ & $118$ & $132$ & $160$ & $200$
\\ $D_4;\phi_{1,0}$ & $57$ & $-q^{2}$ & $4$ & $28$ & $42$ & $50$ & $64$ & $72$ & $86$ & $110$
\\ $D_4;\phi_{4,1}$ & $83$ & $-q^{3}$ & $5$ & $39$ & $61$ & $73$ & $93$ & $105$ & $127$ & $161$
\\ $D_4;\phi_{6,6}''$ & $104$ & $-q^{4}$ & $6$ & $48$ & $78$ & $92$ & $116$ & $130$ & $160$ & $202$
\\ $D_4;\phi_{4,13}$ & $113$ & $-q^{5}$ & $7$ & $53$ & $83$ & $99$ & $127$ & $143$ & $173$ & $219$
\\ $D_4;\phi_{1,24}$ & $117$ & $-q^{6}$ & $8$ & $56$ & $86$ & $102$ & $132$ & $148$ & $178$ & $226$
\\ $E_8[\zeta^3]$ & $104$ & $\zeta^3 q^{4}$ & $8$ & $48$ & $76$ & $90$ & $118$ & $132$ & $160$ & $200$
\\ $E_6[\theta^2];\phi_{1,0}$ & $83$ & $\theta^2 q^{3}$ & $6$ & $38$ & $60$ & $72$ & $94$ & $106$ & $128$ & $160$
\\ $E_6[\theta^2];\phi_{2,1}$ & $104$ & $\theta^2 q^{4}$ & $7$ & $47$ & $75$ & $89$ & $119$ & $133$ & $161$ & $201$
\\ $E_6[\theta^2];\phi_{1,6}$ & $113$ & $\theta^2 q^{5}$ & $8$ & $52$ & $82$ & $98$ & $128$ & $144$ & $174$ & $218$
\\ $E_7[-\I];1$ & $94$ & $-\I q^{7/2}$ & $7$ & $43$ & $69$ & $83$ & $105$ & $119$ & $145$ & $181$
\\ $E_7[-\I];\ep$ & $109$ & $-\I q^{9/2}$ & $8$ & $50$ & $80$ & $96$ & $122$ & $138$ & $168$ & $210$
\\ $E_8[\zeta^4]$ & $104$ & $\zeta^4 q^{4}$ & $8$ & $48$ & $76$ & $90$ & $118$ & $132$ & $160$ & $200$
\\ $E_8[-\theta]$ & $104$ & $-\theta q^{4}$ & $8$ & $48$ & $76$ & $90$ & $118$ & $132$ & $160$ & $200$
\\\hline\end{tabular}\end{center}

\begin{landscape}
\begin{center}\begin{tikzpicture}[thick,scale=1.6]

\draw (14,0.82) node {$\phi_{1,0}$};
\draw (13,0.82) node {$\phi_{8,1}$};
\draw (12,0.82) node {$\phi_{28,8}$};
\draw (11,0.82) node {$\phi_{56,19}$};
\draw (10,0.82) node {$\phi_{70,32}$};
\draw (9,0.82) node {$\phi_{56,49}$};
\draw (8,0.82) node {$\phi_{28,68}$};
\draw (7,0.82) node {$\phi_{8,91}$};
\draw (6,0.82) node {$\phi_{1,120}$};

\draw (4,0.82) node {$D_4;\phi_{1,24}$};
\draw (3,0.82) node {$D_4;\phi_{4,13}$};
\draw (2,0.82) node {$D_4;\phi_{6,6}''$};
\draw (1,0.82) node {$D_4;\phi_{4,1}$};
\draw (0,0.82) node {$D_4;\phi_{1,0}$};

\draw (7.97,2) node {$E_8[-\theta^2]$};
\draw (7.94,0) node {$E_8[-\theta]$};

\draw (6.58,2) node {$E_8[\zeta]$};
\draw (6.63,0) node {$E_8[\zeta^4]$};
\draw (2.62,2) node {$E_8[\zeta^2]$};
\draw (2.62,0) node {$E_8[\zeta^3]$};

\draw (5.4,2) node {$E_7[\I];\ep$};
\draw (5.405,3) node {$E_7[\I];1$};
\draw (5.475,0) node {$E_7[-\I];\ep$};
\draw (5.48,-1) node {$E_7[-\I];1$};

\draw (0,1) -- (14,1);
\draw (5,-1) -- (5,3);
\draw (5,1) -- (6.25,2);
\draw (5,1) -- (7.5,2);
\draw (5,1) -- (6.25,0);
\draw (5,1) -- (7.5,0);
\draw (5,1) -- (2,4);
\draw (5,1) -- (3,2);
\draw (5,1) -- (2,-2);
\draw (5,1) -- (3,0);

\draw (14,1) node [draw] (l0) {};
\draw (13,1) node [draw] (l1) {};
\draw (12,1) node [draw] (l1) {};
\draw (11,1) node [draw] (l1) {};
\draw (10,1) node [draw] (l1) {};
\draw (9,1) node [draw] (l1) {};
\draw (8,1) node [draw] (l1) {};
\draw (7,1) node [draw] (l1) {};
\draw (6,1) node [draw] (l1) {};
\draw (5,1) node [fill=black!100] (ld) {};
\draw (4,1) node [draw] (l0) {};
\draw (3,1) node [draw] (l1) {};
\draw (2,1) node [draw] (l1) {};
\draw (1,1) node [draw] (l1) {};
\draw (0,1) node [draw] (l1) {};

\draw (6.25,2) node [draw] (l1) {};
\draw (6.25,0) node [draw] (l1) {};
\draw (3,2) node [draw] (l1) {};
\draw (3,0) node [draw] (l1) {};

\draw (7.5,2) node [draw] (l1) {};
\draw (7.5,0) node [draw] (l1) {};

\draw (5,-1) node [draw] (l4) {};
\draw (5,0) node [draw] (l4) {};
\draw (5,2) node [draw] (l4) {};
\draw (5,3) node [draw] (l4) {};

\draw (4,2) node [draw] (l1) {};
\draw (3,3) node [draw] (l1) {};
\draw (2,4) node [draw] (l1) {};
\draw (4,0) node [draw] (l1) {};
\draw (3,-1) node [draw] (l1) {};
\draw (2,-2) node [draw] (l1) {};

\draw (14,1.2) node{$0$};
\draw (13,1.2) node{$1$};
\draw (12,1.2) node{$2$};
\draw (11,1.2) node{$3$};
\draw (10,1.2) node{$4$};
\draw (9,1.2) node{$5$};
\draw (8,1.2) node{$6$};
\draw (7,1.2) node{$7$};
\draw (6,1.2) node{$8$};

\draw (4.9,-1.15) node{$7$};
\draw (4.9,-0.15) node{$8$};
\draw (4.9,2.15) node{$8$};
\draw (4.9,3.15) node{$7$};

\draw (4,1.2) node{$8$};
\draw (3,1.2) node{$7$};
\draw (2,1.2) node{$6$};
\draw (1,1.2) node{$5$};
\draw (0,1.2) node{$4$};

\draw (6.25,2.2) node{$8$};
\draw (7.5,2.2) node{$8$};

\draw (6.25,-0.2) node{$8$};
\draw (7.5,-0.2) node{$8$};

\draw (4,2.2) node{$8$};
\draw (3,3.2) node{$7$};
\draw (2,4.2) node{$6$};

\draw (3,2.2) node{$8$};

\draw (4,-0.2) node{$8$};
\draw (3,-1.2) node{$7$};
\draw (2,-2.2) node{$6$};

\draw (3,-0.2) node{$8$};
\end{tikzpicture}\end{center}
\end{landscape}

\bibliography{references}

\end{document}